\documentclass[11pt,leqno]{amsart}
\usepackage[letterpaper,textwidth=6.1in,textheight=8.8in,centering]{geometry}
\usepackage{amsmath,amssymb,mathtools,mathrsfs,bm,enumitem,microtype}
\usepackage{needspace,booktabs,array}
\usepackage{etoolbox}
\usepackage[T1]{fontenc}
\usepackage{lmodern}
\usepackage{xcolor}
\usepackage[colorlinks=true,allcolors=blue,linktoc=page]{hyperref}
\usepackage[nameinlink,noabbrev]{cleveref}
\hypersetup{
  bookmarksdepth=3,
  pdftitle={The q-Miura Transformation and Fundamental q-Characters},
  pdfauthor={Anton M. Zeitlin},
  pdfsubject={The discrete Miura transformation and fundamental q-characters}
}
\numberwithin{equation}{section}
\makeatletter
\renewcommand{\subsection}{\@startsection{subsection}{2}%
  {\z@}{.5\linespacing\@plus.7\linespacing}{.25\linespacing}%
  {\normalfont\bfseries}}
\renewcommand\paragraph{\@startsection{paragraph}{4}{\z@}%
  {.8\baselineskip \@plus .2\baselineskip}%
  {-.5em}%
  {\normalfont\normalsize\bfseries}}
\makeatother
\theoremstyle{plain}
\newtheorem{theorem}{Theorem}[section]
\newtheorem{proposition}[theorem]{Proposition}
\newtheorem{lemma}[theorem]{Lemma}
\newtheorem{corollary}[theorem]{Corollary}

\theoremstyle{definition}
\newtheorem{definition}[theorem]{Definition}
\AtBeginEnvironment{definition}{\pushQED{\qed}}
\AtEndEnvironment{definition}{\popQED}
\theoremstyle{plain}

\theoremstyle{definition}

\theoremstyle{remark}
\newtheorem{remark}[theorem]{Remark}
\newcommand{\sectionoverview}[1]{\noindent#1\par}
\newcommand{\qaff}{\hbar}
\newcommand{\miuraq}{\mathsf{q}}
\newcommand{\C}{\mathbb C}
\newcommand{\Z}{\mathbb Z}

\newcommand{\Ad}{\operatorname{Ad}}

\newcommand{\Frac}{\operatorname{Frac}}

\newcommand{\SL}{\mathrm{SL}}

\title{The $q$-Miura Transformation and Fundamental $q$-Characters}
\author[A. M. Zeitlin]{Anton M. Zeitlin}
\address{School of Mathematics,\newline
  Georgia Institute of Technology,\newline
  686 Cherry Street,\newline
  Atlanta, GA 30332-0160, USA}

\email{\href{mailto:zeitlin@gatech.edu}{zeitlin@gatech.edu}}
\date{}
\begin{document}
\begin{abstract}
Frenkel and Reshetikhin conjectured that under the $q$-Miura transformation
the fundamental $q$-characters express the Steinberg coordinates as Laurent
polynomials in shifted Miura variables. This transformation forgets the
additional flag preserved by the connection of a Miura $q$-oper. We give a
direct proof for every Coxeter element in simply laced type and extend it to
all non-simply-laced types using Miura $(G,\eta,\sigma)$-opers, whose
connections combine the spectral shift $u\mapsto\eta u$ with a diagram
automorphism $\sigma$ of a simply laced group $G$, with the pair
$(G,\sigma)$ prescribed by affine Langlands duality. For ordinary Miura $q$-opers on non-simply-laced groups,
restriction from simply laced groups gives the folded fundamental
$q$-characters.

With the connection and the oper flag fixed, we use generalized B\"acklund
transformations to vary the preserved flag that defines the Miura structure.
On the space of universal Wronskians constructed here, the positive simple-root
subgroups in each $\sigma$-orbit act with a common parameter. We embed the
target module of the corresponding quantum screening operator in the field of
rational functions on this space. Under this embedding, the operator acts on
the Miura variables by differentiation at zero with respect to the common parameter.
The Steinberg coordinates
therefore pull back to elements of the common kernel, where their unique
dominant terms identify them with the fundamental $q$-characters.

These flag variations arise from the action of $G$ on the initial value of
a universal Wronskian. This action gives a
geometric realization of the Frenkel--Mukhin characterization of the
$q$-character ring, with the pullbacks and their shifts generating the common
screening kernel in the Laurent ring of shifted Miura variables.
The same $G$-action lifts the Frenkel--Hernandez
Weyl symmetry to these universal Wronskians, where the pinned representatives of the simple
reflections satisfy the braid relations. On torus-invariant functions they
induce a Weyl group action.
Expanding the transformed Miura variables in each Weyl chamber gives a
uniform geometric proof of the reflection formulas. For folded
$q$-characters, each reflection is induced by the product of the commuting
simple reflections in one orbit of the folding automorphism.

In every type, a meromorphic connection satisfies the oper flag condition
precisely when it admits a unique local Steinberg representative for the
chosen torus trivialization and Coxeter lift. Under the stated normalization
and nonvanishing hypotheses, holomorphic solutions of the $QQ$-system give
holomorphic Steinberg coordinates away from the prescribed singularities. In the $q$DE/ODE/IM
correspondence, the $q$-character formulas express the Steinberg coordinates of twisted monodromy in terms
of principal minors of the central connection matrix relating solutions
normalized at zero and infinity. In XXZ-type models with a formal Cartan twist,
evaluating a fundamental $q$-character at normalized ratios of Baxter operators
gives the corresponding transfer matrix.
\end{abstract}
\maketitle
\begingroup
\linespread{0.97}\selectfont
\tableofcontents
\removelastskip
\endgroup
\clearpage

\section{Introduction}

In 1968 Miura found a transformation taking solutions of the
modified Korteweg--de Vries equation to solutions of the
Korteweg--de Vries equation \cite{Miura}. Drinfeld and Sokolov
placed it in a family of integrable hierarchies associated with
affine Lie algebras \cite{DS}. For $\mathfrak{sl}_N$, the Miura
transformation forms the ordered product of the first-order
operators $\partial_x+v_j(x)$, $1\leq j\leq N$, with
$\sum_j v_j(x)=0$. The sum condition makes the coefficient of
$\partial_x^{N-1}$ vanish, and the remaining coefficients generate
the classical $W$-algebra of $\mathfrak{sl}_N$ as a differential
algebra. In their study of the center of a quantum affine algebra
at the critical level, Frenkel and Reshetikhin introduced a
$q$-deformation of the Miura transformation. In type $A$ they related the
coefficients of the resulting scalar $q$-difference operators to eigenvalues of
transfer matrices, a commuting family of operators in a quantum
integrable model \cite{FRcenter}.

These operators act on functions of the spectral parameter $u$ through
the multiplicative shift $u\mapsto qu$, where $q\in\mathbb C^\times$.
If $f$ solves the associated scalar equation of order $N$, its values at
$u,qu,\ldots,q^{N-1}u$ determine $f(q^Nu)$ at a generic $u$.
Grouping these $N$ values into a column therefore gives a first-order
system that maps the column at $u$ to the column at $qu$. The corresponding
construction for a simple simply connected complex algebraic group $G$
starts with a principal $G$-bundle over $\mathbb P^1$, where the shift
fixes $0$ and $\infty$ \cite[\S2.2]{FKSZ}. A meromorphic $q$-connection
on this bundle assigns to each $u$ where it is defined a $G$-equivariant
isomorphism from the fiber at $u$ to the fiber at $qu$.
In a local trivialization the map sends $h\in G$ to $A(u)h$.
Changing fiber coordinates by left multiplication by $g(u)$ thus
replaces $A(u)$ by $g(qu)A(u)g(u)^{-1}$. Frenkel, Reshetikhin, and
Semenov-Tian-Shansky used unipotent changes of frame in the difference
Drinfeld--Sokolov reduction to construct the deformed classical Virasoro
algebra \cite{FRSTS}. Semenov-Tian-Shansky and Sevostyanov extended
this reduction to semisimple groups \cite{STS}.

To describe the semisimple reduction, we fix opposite Borel subgroups
$B_\pm=HN_\pm$ with common maximal torus $H$ and unipotent radicals
$N_\pm$. A Coxeter element $c$ of the Weyl group is a product
containing each simple reflection once. A pinning specifies an
additive parametrization $x_i(z)$ of each positive simple-root
subgroup. The associated homomorphism $\varphi_i:\mathrm{SL}_2\to G$
takes an upper unitriangular matrix with entry $z$ above the diagonal
to $x_i(z)$. We use
$\dot s_i=\varphi_i\bigl(\begin{smallmatrix}0&-1\\1&0\end{smallmatrix}\bigr)$
as the representative of the simple reflection $s_i$. The product
of these representatives in the chosen order gives the lift
$\dot c$ of $c$. It determines the Coxeter cell $M_c$ and the
Steinberg section $\mathcal A_c$ inside it,
\[
 M_c=N_+\dot cN_+,\qquad
 \mathcal A_c=(N_+\cap\dot cN_-\dot c^{-1})\dot c .
\]
For ordinary conjugation, the section meets each regular conjugacy
class in one point, where regularity means that the centralizer has
dimension $\dim H$ \cite{Steinberg}. For a $q$-connection over the field
$\mathbb C((u))$ of formal Laurent series with $A(u)$ in $M_c$, there is a
unique $N_+$-valued change of frame $n(u)$ for which
$n(qu)A(u)n(u)^{-1}$ lies in $\mathcal A_c$ \cite[Theorem~1.1]{STS}.
We call this representative the \emph{Steinberg form} of the connection.

The same Coxeter element specifies the relative position of flags in
the geometric definition of a $(G,q)$-oper. This definition was developed for
$\mathrm{SL}_N$ in \cite[Definition~4.1 and \S4.2]{KSZ} and for
general simple $G$ in \cite[Definitions~2.1--2.2]{FKSZ}. In the Borel convention
of \cite[Definitions~3.1 and~3.3]{KZ}, the $q$-connection carries a flag
$F_+(u)\in G/B_+$, a full flag in the standard
representation when $G=\mathrm{SL}_N$. Both $F_+(qu)$ and the image
$A(u)F_+(u)$ of the flag at $u$ lie in the fiber at $qu$. For generic $u$, their relative
position is required to be $c$. In a frame in which
$F_+$ is the standard coset $B_+$, this condition reads $A(u)\in B_+\dot cB_+$
\cite[Definition~2.1]{FKSZ} (the same Borel convention is used, for example, in
\cite{QDE}).

A Miura $(G,q)$-oper carries a second flag $F_-(u)\in G/B_-$ which
the connection preserves, so $F_-(qu)=A(u)F_-(u)$. On a dense
open set the two flags are opposite, which means that a single
frame identifies them with the standard cosets $B_+$ and $B_-$
\cite[Theorem~2.3]{FKSZ}. In this frame $A(u)$ takes values in
$B_-$, and its projection to $H=B_-/N_-$ gives the torus
coordinates called the Miura variables. When the chosen frame
puts $A(u)$ in the double coset $M_c$, the two conditions together
give $A(u)\in M_c^-:=M_c\cap B_-$ \cite[Proposition~3.6]{KZ}.
The $q$-Miura transformation forgets the preserved flag. In
coordinates it expresses the Steinberg representative of the
underlying $(G,q)$-oper in terms of the Miura variables and their shifts.

Frenkel and Reshetikhin conjectured that the functions of the
Miura variables obtained in this way are characters of quantum
affine modules after a substitution of spectral variables. The
characters in question are the $q$-characters introduced in their
work on deformed $W$-algebras. They record a refinement of the
ordinary weights of a representation
\cite{FRdeformed,FR}. For a finite-dimensional complex simple Lie
algebra $\mathfrak a$, its loop algebra is
$L\mathfrak a=\mathfrak a\otimes_{\mathbb C}\mathbb C[w,w^{-1}]$. The quantum loop algebra
$U_{\qaff}(L\mathfrak a)$ is obtained from the untwisted quantum
affine algebra by setting the central generator equal to one and
omitting the degree operator. We write $\qaff\in\mathbb C^\times$ for
its multiplicative deformation parameter, assumed to have infinite order,
and reserve $q$ for the base shift $u\mapsto qu$ of an ordinary $q$-connection.

The weight spaces of a finite-dimensional $U_{\qaff}(L\mathfrak a)$-module
$V$ of type one are the simultaneous eigenspaces of the Cartan generators
$k_i$. The type-one condition requires each $k_i$ to act
diagonalizably with eigenvalues in $\qaff_i^{\mathbb Z}$. Here
$\qaff_i=\qaff^{d_i^{\mathrm{sym}}}$, where
$d_i^{\mathrm{sym}}=(\alpha_i,\alpha_i)/2$ for the $i$th simple root
of $\mathfrak a$, with the shortest roots normalized to have squared
length two. In Drinfeld's current presentation of the quantum loop
algebra (see the \hyperref[app:quantum-presentation]{Appendix}), the Cartan currents
$\phi_i^+(x)$ and $\phi_i^-(x)$ have constant terms $k_i$ and $k_i^{-1}$,
respectively. All their coefficients commute. The joint generalized eigenspaces of these coefficients
therefore refine the ordinary weight spaces. The eigenvalues on
such a summand form its \emph{$\ell$-weight}
\cite[\S\S1.2--1.3 and 2.2]{FM}.

To record an $\ell$-weight one uses a Laurent monomial
$m=\prod_{i,b}Y_{i,b}^{n_{i,b}}$, where $i$ indexes a simple root
and $b\in\mathbb C^\times$ is a spectral parameter. The eigenvalues
of $\phi_i^+(x)$ and $\phi_i^-(x)$ are the expansions at zero and
infinity of the same rational function
\[
 \Psi_{i,m}(x)=
 \prod_b\left(\qaff_i\,\frac{1-b\qaff_i^{-1}x}
                            {1-b\qaff_i x}\right)^{n_{i,b}} .
\]
Thus the monomial records eigenvalues that depend on $\qaff$
\cite[Propositions~2.3--2.4]{FM}. If $V_m$ denotes the corresponding
summand of $V$, its dimension is the coefficient of $m$ in the
usual $q$-character. With our notation for the deformation parameter
the formula reads
\[
 \chi_{\qaff}(V)=\sum_{m:\,V_m\ne 0}\dim(V_m)\,m .
\]
Additivity on short exact sequences and multiplicativity on tensor
products make $\chi_{\qaff}$ a homomorphism on the Grothendieck ring
of finite-dimensional type-one modules. Its image is the
$q$-character ring.

For a simple module, the monomial of a generating $\ell$-weight
vector annihilated by the positive Drinfeld generators is called
its \emph{highest monomial}. The fundamental module $L(Y_{i,b})$
is the simple finite-dimensional type-one $U_{\qaff}(L\mathfrak a)$-module
with highest monomial $Y_{i,b}$. The classes of the fundamental modules freely
generate the Grothendieck ring as a polynomial ring \cite[Theorem~3 and Corollary~2]{FR}. For
$\mathrm{SL}_N$, Frenkel and Reshetikhin identified the fundamental
$q$-characters, after a spectral substitution, among the coefficients
of an ordered product of first-order difference operators with shift
$u\mapsto\qaff^2u$. They conjectured that difference
Drinfeld--Sokolov reduction realizes the $q$-character homomorphism
in every simply laced type \cite[\S8.7]{FR}.

To recognize $q$-characters among Laurent polynomials, Frenkel and
Reshetikhin associated to each simple root $i$ a derivation $S_i$
from the Laurent ring to a module generated by symbols $s_{i,b}$,
one for each spectral parameter $b$. This \emph{screening operator}
sends $Y_{i,b}$ to $Y_{i,b}s_{i,b}$ and vanishes on $Y_{j,b}$ for
$j\ne i$. The Cartan matrix determines the Laurent monomial
relating $s_{i,b}$ to $s_{i,b\qaff_i^2}$ in the target module.

Every screening $S_i$ annihilates $\chi_{\qaff}(V)$ for every
finite-dimensional type-one $U_{\qaff}(L\mathfrak a)$-module $V$,
including each fundamental module $L(Y_{j,b})$. The conjecture states
that the common kernel of the screenings in the integral Laurent polynomial
ring is the $\mathbb Z$-linear span of the $q$-characters of all
finite-dimensional type-one modules \cite[\S7, Conjecture~2]{FR}.
Frenkel and Mukhin proved this characterization for all untwisted
quantum affine algebras \cite[Theorem~5.1]{FM}. Within the common
kernel, an element is determined by the coefficients of its
\emph{dominant monomials}, those in which every variable has nonnegative
exponent. The fundamental $q$-character of $L(Y_{i,b})$ is the unique element of the common
kernel whose only dominant monomial is $Y_{i,b}$, with coefficient
one \cite[Corollary~4.5 and Lemma~6.17]{FM}.

In the analytic Bethe Ansatz, this screening criterion applies to the
scalar difference operators constructed in
\cite[Proposition~2.3, Theorem~2.5, and \S5]{KOSY} from the shift
$f(x)\mapsto f(x+1)$ and Laurent monomials in the character variables.
In type $C_n$ the operator has order $2n+2$. Every screening annihilates
its coefficients, among which all fundamental $q$-characters occur up to
signs and spectral shifts. In types $B_n$ and $D_n$ the construction
gives formal series in the same shift.
The generalized Baxter relations give another interpretation of the
same characters in an integrable model. For a transfer matrix
labeled by an auxiliary representation $V$, one obtains its
eigenvalues by evaluating $\chi_{\qaff}(V)$. Each $Y_{i,b}$ is
replaced by a prescribed factor times a ratio of shifted values
of a Baxter polynomial. The roots of the Baxter polynomials are the Bethe roots
\cite[Theorem~5.11 and \S5.6]{FH}.

A similar ratio occurs in the torus coordinates of Miura $(G,q)$-opers. In the Borel convention of \cite[\S4.1]{FKSZ}, the
preserved Borel subgroup acts by upper triangular matrices on the span
of a highest-weight vector in the $i$th fundamental representation
and its image under the simple lowering operator at $i$. In this basis
the connection restricts to a matrix $A_i(u)$. For the nondegenerate $Z$-twisted Miura--Pl\"ucker $(G,q)$-opers of
\cite[Theorem~6.1]{FKSZ}, the rational change of frame
$v_i(u)$ satisfies $A_i(u)v_i(u)=v_i(qu)Z_i$, where $Z_i$ is the
restriction of the torus element $Z$. The first row of $v_i(u)$ is
$(Q_i^+(u),-Q_i^-(u))$. If $\zeta_i$ is the eigenvalue of $Z$ on
the highest-weight line, the upper left entry of the matrix identity
gives the torus coordinate $\zeta_iQ_i^+(qu)/Q_i^+(u)$.

The upper right entries give the $QQ$-system \cite[(6.2)--(6.3)]{FKSZ}.
The prescribed singularities determine the polynomials
$\Lambda_i$, and one seeks polynomial solutions $Q_i^\pm$ with $Q_i^+$
monic. With the simple roots numbered in the order defining $c$ and
$A^G_{ji}=\langle\alpha_i,\alpha_j^\vee\rangle$, the relation at vertex $i$ is
\begin{equation}\label{eq:intro-oper-qq}
 \widetilde\gamma_i Q_i^-(u)Q_i^+(qu)-\gamma_i Q_i^-(qu)Q_i^+(u)
 =\Lambda_i(u)\prod_{j>i}Q_j^+(qu)^{-A^G_{ji}}
               \prod_{j<i}Q_j^+(u)^{-A^G_{ji}} .
\end{equation}
The twist fixes
$\widetilde\gamma_i=\zeta_i\prod_{j>i}\zeta_j^{A^G_{ji}}$ and
$\gamma_i=\zeta_i^{-1}\prod_{j<i}\zeta_j^{-A^G_{ji}}$.
At a root $w$ of $Q_i^+$, evaluating \eqref{eq:intro-oper-qq} at $w$
and $q^{-1}w$ and eliminating the common value $Q_i^-(w)$ gives the
Bethe equation for $w$ under the nondegeneracy assumptions of
\cite[\S6.3]{FKSZ}.

For $G=\mathrm{SL}_N$, the polynomials $Q_i^+$ and $Q_i^-$ in
\eqref{eq:intro-oper-qq} can be expressed as normalized minors of the
$q$-Wronskian matrix formed from a section of the line in the defining
flag of the $(\mathrm{SL}_N,q)$-oper and its successive images under the $q$-connection \cite[\S6]{KZToroidal}, \cite[(2.13)]{KZZoo}.
For a general simple group $G$, the \emph{$(G,q)$-Wronskian} of
\cite[Definition~4.9]{KZ} represents the flag $F_+$ of a $Z$-twisted $(G,q)$-oper by a
meromorphic section in the frame where the connection sends $(u,h)$ to $(qu,Zh)$. Applied to a highest-weight
vector of the $i$th fundamental representation, the section gives a vector
with components $Q_i^+(u)$ along the highest-weight line and $Q_i^-(u)$
along its translate by $s_i$, under the hypotheses of
\cite[Proposition~4.17 and (4.29)]{KZ}. The Wronskian shift identities
on highest-weight lines turn a quadratic relation among these matrix
coefficients, called generalized minors, into the $i$th $QQ$-relation.\footnote{For differential $G$-opers, generalized minors of a $G$-Wronskian give the differential $qq$-system of the twisted Gaudin model \cite{ZWronskians}.}
These minor identities also give exchange relations on double Bruhat
cells, which suggested a cluster interpretation of the $QQ$-systems
\cite[Abstract and \S4.3]{KZ}.

The \emph{universal $(G,q)$-Wronskian} relates a frame in which the
connection is constant to one in which it takes values in the Steinberg
section \cite[\S4.4]{KZ}. An algebraic analogue on a discrete spectral
lattice was developed in \cite[Definition~10.10 and Remark~10.14]{GHL}.
For the ordinary shift in simply laced type, the discrete universal
Wronskians are precisely the sequences that Francone and Leclerc call
$(G,c)$-bands \cite[Definition~1.1]{FLbands}. Their convention reverses
the lattice index and takes the inverse Wronskian with the constant
twist removed, as in Section~\ref{sec:wronskian-descent}. They used these sequences to prove the
discrete Frenkel--Reshetikhin conjecture for every Coxeter element in
types $A$, $D$, and $E$. Their proof compares cluster expansions of
invariant functions with Laurent expansions obtained from the
generalized Baxter relations \cite[\S\S5--6, Theorem~30]{FL}.

We retain the passage to the integer lattice, where the connection in
Steinberg form and one initial value determine the Wronskian. Our
shorter proof of the character identity uses the screenings directly.
Varying the preserved flag shows that every screening annihilates
the pullback of each Steinberg coordinate. Its unique dominant monomial
then identifies the corresponding fundamental $q$-character.
This argument also gives the ordinary fundamental $q$-characters in
types $B_n$, $C_n$, $F_4$, and $G_2$.

Earlier constructions in non-simply-laced type approached the
character formula through the deformed $W$-algebras
$W_{\qaff,t}(\mathfrak a)$ with two deformation parameters. Formulas for
their fundamental fields were
obtained in rank two \cite[Appendix~A]{BP}. Assakaf constructed
fields whose specialization at $t=1$ gives the fundamental
$q$-characters at every vertex of $B_n$, $C_n$, and $G_2$ and at
vertices $1$ and $4$ of $F_4$ in his numbering
\cite[Theorem~4]{Assakaf}.

Our realization of the ordinary fundamental $q$-characters in these
four types uses the lattice analogue of a Miura $(G,\eta,\sigma)$-oper.
In the $q$DE/ODE/IM correspondence, the corresponding meromorphic
connection arises from a fractional
rotation of the differential variable that shifts the spectral parameter
from $u$ to $\eta u$ and acts on the fibers through a diagram automorphism
$\sigma$ of a simply laced group $G$. This construction is described in
\cite[\S7]{QDE} and developed with full proofs in \cite{QDEII}.
In a local frame the connection sends $(u,h)$ to $(\eta u,A(u)\sigma(h))$.
The relative position of the flag at $\eta u$ and the image of the flag
at $u$ is prescribed by an ordered product of simple reflections, one
for each $\sigma$-orbit of simple roots.
In types $A$, $D$, and $E$ we recover the ordinary $(G,q)$-oper by
taking $\sigma=1$ and $\eta=q$.

For both $(G,q)$-opers and $(G,\eta,\sigma)$-opers, changing the
preserved flag in simple-root directions realizes the screening operators
in the character formula. The construction also applies to ordinary
$(G_\nu,q)$-opers for the identity component $G_\nu$ of the subgroup fixed by a diagram
automorphism $\nu$ of a simply laced group. Restriction to this subgroup
gives the folded $q$-characters of \cite{FHR}.

These variations of the preserved flag also relate the screening
operators to the Weyl symmetry of $q$-characters. Frenkel and
Hernandez constructed a Weyl group action on completions of the
Laurent ring of spectral variables and identified the invariant
Laurent polynomials with the $q$-character ring. They obtained each
screening operator as the first-order term in a deformation of the
corresponding simple reflection
\cite[Theorems~3.4 and~5.1, \S5.2]{FHWeyl}.

In our construction the rank-one group at each simple root acts on
the initial value of a discrete universal Wronskian. Its positive root
subgroup realizes the one-parameter deformation. The same rank-one
group contains the pinned representative of the simple reflection.
After expansion in the Weyl chambers, its action on the Miura
variables gives the Frenkel--Hernandez formula for that reflection. In this
way the relations among the representatives give a geometric proof
of the Weyl symmetry on the spectral lattice in every simple Lie type.
Restriction to the fixed subgroup gives the analogous symmetry for
folded $q$-characters.

The following subsections introduce the discrete Miura transformation
for ordinary $(G,q)$-opers and its extension to $(G,\eta,\sigma)$-opers,
then state the main theorems and outline their proofs. We describe
the role of the Wronskian in the geometric proof of the screening-kernel
theorem and in the construction of the Weyl action. We then explain how
the $QQ$-relations give local cancellation of apparent poles in the
Steinberg coordinates. The applications concern two evaluations of
the fundamental $q$-characters. In the $q$DE/ODE/IM correspondence,
the Miura variables obtained from the central connection matrix give
the Steinberg coordinates of twisted monodromy. In quantum integrable
models of XXZ type, normalized ratios of Baxter operators give the
fundamental transfer matrices.

\subsection{From \texorpdfstring{$q$}{q}-opers to the discrete Miura transformation}
\label{sec:intro-lattice}

\sectionoverview{The conjecture compares a Laurent polynomial in independent
character variables with an expression in shifted values of the Miura
variables. In this subsection we explain why we can first compute that
expression with independent variables at each integer and then
evaluate it on a Miura $(G,q)$-oper.}

The reason for introducing the integer lattice is already visible in
the change-of-frame formula $A(u)\mapsto g(qu)A(u)g(u)^{-1}$.
Each connection map uses the frames at two consecutive points of a
multiplicative orbit. When $q$ has infinite order, the orbit through
$u_*\in\mathbb C^\times$ consists of the distinct points
$q^ku_*$, $k\in\mathbb Z$. On an orbit where the connection is
defined, $a_k=A(q^ku_*)$ represents its map from the fiber at
$q^ku_*$ to the fiber at $q^{k+1}u_*$. A change of frame $g$ with
values $g_k=g(q^ku_*)$ replaces $a_k$ by $g_{k+1}a_kg_k^{-1}$.
This formula also makes sense for arbitrary
sequences $(a_k)$ and $(g_k)$ in $G$. Working on $\mathbb Z$ thus allows the
Miura variables at different indices to vary independently while
retaining the change-of-frame formula. An identity obtained in
these variables can then be evaluated on any sequence in the chosen
Miura cell, including a periodic sequence. At $u_*=0$ the orbit consists of a single
point. If the connection and change of frame are regular there, $A(0)$
acts on that one fiber and changes by conjugation by $g(0)$.

The two flags of a Miura $(G,q)$-oper determine which sequences enter the
Miura transformation. On an orbit where the flags are opposite, a
frame in which both are standard puts $a_k$ in $B_-$, since the
connection preserves the second flag. The relative-position condition
on the first flag puts $a_k$ in $B_+\dot c B_+$. With the lift $\dot c$ fixed, a choice of torus
frame at one index determines successive torus changes of frame for
which $A_k^{\mathrm M}\in M_c^-:=B_-\cap N_+\dot cN_+$. On an
infinite orbit each index labels a distinct point. On an orbit of
length $L$, these choices give a frame on the orbit precisely when
the frame after $L$ shifts agrees with the initial one. If a
meromorphic function $A^{\mathrm M}(u)$ already takes values in
$M_c^-$, its values $A^{\mathrm M}(q^ku_*)$ give the sequence
directly, with period $L$ on an orbit of that length.

For simply laced $G$ and any Coxeter element $c$, each sequence
$(A_k^{\mathrm M})$ in $M_c^-$ has a unique change of frame
$(u_k)$ in $N_+$ such that
\[
 A_k=u_{k+1}A_k^{\mathrm M}u_k^{-1}\in\mathcal A_c
 \qquad(k\in\mathbb Z).
\]
The assignment $(A_k^{\mathrm M})\mapsto(A_k)$ is the ordinary
discrete Miura transformation $\mu_c$. To express this transformation in coordinates, we use the same
ordering of the simple roots on both sides.
On the Miura side, an element $A_k^{\mathrm M}$ is determined by
its projection $[A_k^{\mathrm M}]_0$ to $H=B_-/N_-$.
On the Steinberg side, the ordering $c=s_{i_1}\cdots s_{i_s}$
gives the unique factorization
\[
 A_k=\bigl(x_{i_1}(T_{i_1,k})\dot s_{i_1}\bigr)\cdots
     \bigl(x_{i_s}(T_{i_s,k})\dot s_{i_s}\bigr).
\]
Thus $T_{i,k}$ is the parameter of the simple-root subgroup at $i$
in this product.

A matrix coefficient allows us to compare these two sets of
coordinates. The simple coroot $\alpha_j^\vee$ sends $t\in\mathbb C^\times$
to $\varphi_j(\operatorname{diag}(t,t^{-1}))\in H$. The $i$th fundamental
weight $\omega_i$ takes $\alpha_j^\vee(t)$ to $t$ if $j=i$ and to $1$
otherwise. In the irreducible
representation $V(\omega_i)$ of highest weight $\omega_i$, a choice of
highest-weight vector $v_i$ determines a unique linear functional
$v_i^*$ which takes value one on $v_i$ and vanishes on the other
weight spaces. The corresponding
\emph{principal minor} $\Delta_i(g)=v_i^*(gv_i)$ satisfies
\[
 \miuraq_{i,k}:=\omega_i([A_k^{\mathrm M}]_0)=\Delta_i(A_k^{\mathrm M}),
 \qquad T_{i,k}=\Delta_i(A_k).
\]
We use the upright letter $\miuraq$ for the Miura coordinates and $q$
for the multiplier of the shift $u\mapsto qu$. Expressing $T_{i,k}$
in the Miura coordinates $\miuraq_{j,\ell}$ gives the pullback
$\mu_c^*(T_{i,k})$ that we will identify with a fundamental $q$-character.

In type $A$ these coordinates have a familiar interpretation.
For $G=\mathrm{SL}_3$ with its standard pinning and $c=s_1s_2$, the
displayed factorization gives its matrix $A_k$ in the standard
representation. If $f_k$ is the first component of a solution of
$y_{k+1}=A_k y_k$, the equation becomes
\[
 \begin{pmatrix} f_{k+1} \\ f_k \\ f_{k-1} \end{pmatrix}
 =
 \begin{pmatrix}
     T_{1,k} & -T_{2,k} & 1 \\
     1       &  0       & 0 \\
     0       &  1       & 0
 \end{pmatrix}
 \begin{pmatrix} f_k \\ f_{k-1} \\ f_{k-2} \end{pmatrix}.
\]
The first row gives
$f_{k+1}=T_{1,k}f_k-T_{2,k}f_{k-1}+f_{k-2}$. The other rows retain
the two consecutive values needed at the next index. This matrix acting on consecutive values
is called the \emph{companion matrix} of the scalar
recurrence. For $\mathrm{SL}_N$ with the simple reflections in their
linear order, the same construction gives a scalar recurrence of
order $N$.

In every simply laced type, each pullback $\mu_c^*(T_{i,k})$ is a
Laurent polynomial in finitely many independent variables
$\miuraq_{j,\ell}$. A meromorphic connection $A^{\mathrm M}(u)$ in the
cell $M_c^-$ gives functions $\miuraq_j(u)=\Delta_j(A^{\mathrm M}(u))$.
Substituting $\miuraq_{j,\ell}=\miuraq_j(q^\ell u)$ in the Laurent polynomial
gives a meromorphic function wherever the indicated shifts are
defined, or a meromorphic germ if the inputs are germs at zero.
Indeed, $u_k(w)$ denotes the change of frame at index $k$ on the orbit through $w$. The sequences based at $w$ and $qw$
differ by one index. Uniqueness of the change of frame gives
$u_{k+1}(w)=u_k(qw)$ and hence $u_k(w)=u_0(q^kw)$. Thus the
substitution computes the Steinberg coordinate at $q^ku$.

The universal $(G,q)$-Wronskian also has an analogue
on the integer lattice. In the meromorphic setting it compares a
frame in which the connection of the $(G,q)$-oper is represented by
$A(u)\in\mathcal A_c$ with
one in which the connection is the constant $Z_0\in H$. A function
$W(u)$ giving this change of frame satisfies $W(qu)A(u)=Z_0W(u)$.
On a multiplicative orbit its values satisfy
\[
 W_{k+1}A_k=Z_0W_k.
\]
For any sequence $(A_k)$ in the Steinberg section, this recurrence
defines its \emph{discrete universal Wronskian} once $W_0\in G$ is
chosen. Each $A_k$ is invertible, so the recurrence determines $W_k$
in both directions on $\mathbb Z$. For a sequence obtained from an
orbit of length $L$, its discrete Wronskian gives a function on that
orbit precisely when $W_L=W_0$. Over a field of meromorphic functions, existence
of a single $W(u)$ giving the change of frame to the prescribed
$Z_0$ is precisely the $Z_0$-twistedness condition
\cite[Definition~3.1]{FKSZ}. It can fail even for a connection in
Miura form holomorphic at zero. The formula for its Steinberg
coordinates still applies, since it only evaluates finite Laurent
polynomials.

\subsection{Fundamental and folded \texorpdfstring{$q$}{q}-characters}

\sectionoverview{In this subsection we state the theorem identifying the
pullbacks of Steinberg coordinates under the discrete Miura transformation
with ordinary fundamental $q$-characters in every simple Lie type. We then explain
why restricting the ordinary Miura transformation to the subgroup
fixed by a diagram automorphism gives the folded $q$-characters.}

For a simply laced Lie algebra $\mathfrak a$, the fundamental
$q$-characters of $U_{\qaff}(L\mathfrak a)$ are realized by the lattice
version of the ordinary Miura transformation for the simply connected
group $G$ with Lie algebra $\mathfrak a$. To extend this realization to types
$B_n$, $C_n$, $F_4$, and $G_2$, we use the affine Langlands dual of the
untwisted affine algebra of $\mathfrak a$. This dual affine algebra has
a twisted loop realization in $\operatorname{Lie}G$ for a simply
connected, simply laced group $G$. Its loops are equivariant under rotation
of the loop parameter through $1/r$ of a turn and a diagram
automorphism $\sigma$ of order $r$. The same dual affine types occur
in the differential equations associated with these quantum affine
algebras \cite{MRV2,QDEII}. On the difference side, the corresponding
generalization is the \emph{$(G,\eta,\sigma)$-oper}, whose connection
combines each spectral shift with the action of $\sigma$ on the fiber.

On a domain invariant under multiplication by $\eta\in\mathbb C^\times$,
the connection in a local frame sends $h$ in the fiber at $u$ to
$A(u)\sigma(h)$ in the fiber at $\eta u$, where $A$ is a meromorphic
$G$-valued function. The pinning specifies the lift of $\sigma$ to $G$
by $\sigma(x_i(z))=x_{\sigma(i)}(z)$. Thus right multiplication by
$x_i(z)$ in the source fiber corresponds to right multiplication by
$x_{\sigma(i)}(z)$ in the target fiber. After $r$ successive maps the
composite is $G$-equivariant because $\sigma^r=1$, and hence is an
ordinary $q$-connection with $q=\eta^r$. The defining flag $F_+$ of
the $(G,\eta,\sigma)$-oper satisfies a condition for a single map from
$u$ to $\eta u$. At generic $u$ the relative position of $F_+(\eta u)$
and the image $A(u)\sigma(F_+(u))$ of the flag at $u$ must be $c=s_{i_1}\cdots s_{i_s}$,
where $i_1,\ldots,i_s$ are ordered representatives of the
$\sigma$-orbits of simple roots. When $\sigma=1$ and $\eta=q$, each
orbit consists of a single vertex and this gives the ordinary
$(G,q)$-oper condition \cite{QDEII}.

A Miura $(G,\eta,\sigma)$-oper carries a second flag $F_-$ of type
$B_-$, preserved in the sense that $F_-(\eta u)=A(u)\sigma(F_-(u))$.
Forgetting this flag gives the Miura transformation. On a multiplicative
orbit $u_k=\eta^ku_0$ where the two flags are opposite, we can choose
a frame at each $u_k$ in which both are standard. With $\dot c$ fixed,
changes of frame by elements of $H$ then put each connection map in the
form $h\mapsto A_k^{\mathrm M}\sigma(h)$ with $A_k^{\mathrm M}\in M_c^-$.
The lattice formulation allows an arbitrary sequence of maps of this
form. Changing fiber coordinates by $g_k$ at the source and by
$g_{k+1}$ at the target replaces $A_k^{\mathrm M}$ by
$g_{k+1}A_k^{\mathrm M}\sigma(g_k)^{-1}$.

The affine Langlands dual types and the corresponding choices of
$(G,\sigma)$ for these lattice sequences are listed below. The last
column specifies $c$ by giving the simple reflections in their
chosen order, one for each $\sigma$-orbit.
\[
\begin{array}{c|c|c|c|l}
 \text{Quantum affine type}&\text{Dual affine type}&\text{Type of }G
       &\operatorname{ord}(\sigma)&c\\\hline
 B_n^{(1)}&A_{2n-1}^{(2)}&A_{2n-1}&2&s_1\cdots s_n\\
 C_n^{(1)}&D_{n+1}^{(2)}&D_{n+1}&2&s_1\cdots s_n\\
 F_4^{(1)}&E_6^{(2)}&E_6&2&s_1s_2s_3s_4\\
 G_2^{(1)}&D_4^{(3)}&D_4&3&s_1s_2
\end{array}
\]
In the first two rows $n\ge2$, with $D_3=A_3$ when $n=2$.
The involution reverses the $A_{2n-1}$ chain or interchanges
vertices $n$ and $n+1$ of $D_{n+1}$. In $E_6$ it reverses the chain
$1,3,4,5,6$ and fixes vertex $2$, attached to $4$. In the last row $\mathfrak a$ has type $G_2$ and $G$ has type
$D_4$. The order-three symmetry of $D_4$ cycles its three outer
vertices $1,3,4$ and fixes the central vertex $2$. Thus the two orbits are $\{1,3,4\}$ and $\{2\}$. Choosing $1$
from the first and $2$ from the second, in that order, gives
$c=s_1s_2$ in the Weyl group of $D_4$.
For each chosen representative $i$, the size $o_i$ of its
$\sigma$-orbit is the least positive integer for which
$\sigma^{o_i}(i)=i$. In the four types in the table, the corresponding
simple root of $\mathfrak a$ is represented by the sum of the simple
roots of $G$ in that orbit. These roots are pairwise orthogonal and
have squared length two, so their sum has squared length $2o_i$.
Hence $d_i^{\mathrm{sym}}=o_i$ and
$\qaff_i=\qaff^{d_i^{\mathrm{sym}}}=\qaff^{o_i}$. In simply laced
type $\sigma=\mathrm{id}_G$, so $o_i=d_i^{\mathrm{sym}}=1$.

The comparison with $q$-characters requires a spectral label for each
Miura variable $\miuraq_{i,k}$. Increasing the lattice index $k$ by one
multiplies this label by $\qaff^2$. For a $(G,\eta,\sigma)$-connection
the same change of index takes the base point $u$ to $\eta u$.
The two multipliers $\qaff^2$ and $\eta$
may be chosen independently. The labels at different vertices are
fixed relative to one another by integers $\xi_i$. In simply
laced type, the integers satisfy $\xi_j=\xi_i-1$ whenever $i$ and $j$ are
adjacent and $i$ precedes $j$ in the ordering for $c$. This determines
the integers up to a common translation. The corresponding integers for the four
orderings in the table are
\[
\begin{array}{c|l}
 B_n&\xi_i=-2(i-1)\ (i<n),\quad \xi_n=3-2n\\
 C_n&\xi_i=1-i\\
 F_4&(\xi_1,\xi_2,\xi_3,\xi_4)=(0,-2,-2,-3)\\
 G_2&(\xi_1,\xi_2)=(0,-1).
\end{array}
\]
Thus a choice of $a\in\mathbb C^\times$ identifies $\miuraq_{i,k}$ with
$Y_{i,a\qaff^{2k+\xi_i}}$. The pairs $(i,a\qaff^{2k+\xi_i})$ form
the \emph{spectral lattice at $a$}. A monomial is supported on this
lattice if every variable occurring with nonzero exponent has its
label there. The Laurent polynomial identity on the lattice can be
evaluated at invertible elements of any commutative complex
algebra $D$. For an affine variety $X$ such as
$M_c^-$ or $\mathcal A_c$, $X(D)$ denotes the set of algebra
homomorphisms $\mathbb C[X]\to D$, each of which assigns a value in
$D$ to every regular function on $X$. We use complex coefficients
for the character rings by extending scalars from $\mathbb Z$ to
$\mathbb C$.

\Needspace{12\baselineskip}
\begin{theorem}\label{thm:intro}
For each finite-dimensional complex simple Lie algebra $\mathfrak a$,
the ordinary fundamental $q$-characters of $U_{\qaff}(L\mathfrak a)$
are obtained from Miura sequences in a simply connected, simply
laced group $G$ as follows. In simply laced type,
$\operatorname{Lie}G=\mathfrak a$, $\sigma=\mathrm{id}_G$, and $c$
is any Coxeter element. In non-simply-laced type, $G$ and its
automorphism $\sigma$ are specified by the table of affine types.
The $\sigma$-twisted affine algebra of $\operatorname{Lie} G$ is
the Langlands dual of the untwisted affine algebra of $\mathfrak a$, and $c$ is the
product of simple reflections displayed for that type.

For any commutative complex algebra $D$, a sequence
$(A_k^{\mathrm M})_{k\in\mathbb Z}$ in $M_c^-(D)$, where
$M_c^-=B_-\cap N_+\dot cN_+\subset G$, represents the maps
$h\mapsto A_k^{\mathrm M}\sigma(h)$ from index $k$ to index $k+1$.
\begin{enumerate}[label=\textup{(\roman*)},leftmargin=*,nosep]
\item A unique sequence of changes of frame $(u_k)_{k\in\mathbb Z}$ in
$N_+(D)$ takes these maps to the Steinberg section, with
\begin{equation}\label{eq:intro-Miura}
 A_k=u_{k+1}A_k^{\mathrm M}\sigma(u_k)^{-1}\in\mathcal A_c(D),
 \qquad k\in\mathbb Z.
\end{equation}
The discrete Miura transformation $\mu_{c,\sigma}$ sends
$(A_k^{\mathrm M})$ to $(A_k)$.
\item For the Steinberg coordinate
$T_{i,k}=\Delta_i(A_k)$ and the independent Miura coordinates
$\miuraq_{j,d}=\Delta_j(A_d^{\mathrm M})$, the identification
$\miuraq_{j,d}=Y_{j,a\qaff^{2d+\xi_j}}$ gives the Laurent polynomial identity
\begin{equation}\label{eq:intro-character}
 \mu_{c,\sigma}^*(T_{i,k})
      =\chi_{\qaff}\!\left(L(Y_{i,a\qaff^{2k+\xi_i}})\right).
\end{equation}
\item As $i$ ranges over the chosen representatives of the $\sigma$-orbits and $k$ over
$\mathbb Z$, these pullbacks
freely generate the complexified $q$-character
ring of the finite-dimensional type-one $U_{\qaff}(L\mathfrak a)$-modules
whose simple composition factors have highest monomials supported on this lattice.
\end{enumerate}
\end{theorem}

The left-hand side of \eqref{eq:intro-character} expresses the
principal minor $\Delta_i(A_k)$ in the Miura variables. After the
spectral identification, the coefficient of each monomial is the
dimension of the corresponding $\ell$-weight space of the fundamental
module on the right. Replacing $a$ by $a\qaff^{2m}$ for any
$m\in\mathbb Z$ translates $k$ by $m$ and gives the same spectral
lattice. Choosing one $a$ from each orbit under multiplication by
$\qaff^2$ partitions the variables $Y_{i,b}$ among these lattices.
The formula therefore recovers every fundamental
$q$-character of $U_{\qaff}(L\mathfrak a)$.

For folded characters the diagram automorphism enters in a different
place. We start with a simply connected, simply laced group
$\widetilde G$ and a diagram automorphism $\nu$ preserving its pinning.
Folding gives the group $G_\nu=(\widetilde G^\nu)^\circ$, the identity
component of the subgroup fixed by $\nu$. The ordinary connection on
$G_\nu$ sends $h$ to $A_k^{\mathrm M}h$ at index $k$, so its
Steinberg representative is $A_k=u_{k+1}A_k^{\mathrm M}u_k^{-1}$.
For example, reversing the $A_{2n-1}$ diagram of $\mathrm{SL}_{2n}$
gives the fixed group $\mathrm{Sp}_{2n}$ of type $C_n$. In
Theorem~\ref{thm:intro} the same symmetry acts in the connection maps
used to realize the ordinary $q$-characters of type $B_n$.

For a vertex $v$ of the Dynkin diagram of $\widetilde G$, its $\nu$-orbit is
$O=\{\nu^t(v):t\in\mathbb Z\}$. We assume that distinct vertices
of $O$ are nonadjacent. If $\widetilde H$ is the maximal torus in the
pinning of $\widetilde G$, the simple roots at the vertices of $O$
have the same restriction to $H_\nu=(\widetilde H^\nu)^\circ$. This
restriction is a simple root of $G_\nu$, and every simple root of
$G_\nu$ arises from exactly one such $O$. An ordering of these
folded simple roots gives a Coxeter element $c_\nu$ of $G_\nu$.
Replacing the reflection for each folded simple root by the commuting simple
reflections at the vertices of its corresponding $\nu$-orbit gives a Coxeter element of
$\widetilde G$. Its Steinberg section is stable under $\nu$, and
the integers $\xi_i$ have a common value $\xi_O$ for $i\in O$.

Restricting Miura data from $\widetilde G$ to $G_\nu$ identifies the
variables at vertices of each $\nu$-orbit. On the character side
this gives the homomorphism $\pi$ of Laurent polynomial rings
defined by $\pi(\widetilde Y_{i,b})=Y_{O,b}$ for $i\in O$, where
$\widetilde Y_{i,b}$ is a spectral variable for the quantum loop
algebra of $\operatorname{Lie}\widetilde G$. The image of a
$q$-character under $\pi$ is its folded $q$-character
\cite[Theorem~4.3]{FHR}. Each vertex of the simply laced algebra has
$\qaff_i=\qaff$. Since $\pi$ preserves the spectral parameter $b$,
the shifts $b\mapsto b\qaff^{\pm1}$ are the same for every
$\nu$-orbit $O$. The folded spectral lattice consists of the pairs
$(O,b_{O,k})$ with $b_{O,k}=a\qaff^{2k+\xi_O}$. On this lattice
the folded character ring is the common kernel of the folded screening
operators \cite[Proposition~3.2 and Theorem~4.3]{FHR}.

\Needspace{12\baselineskip}
\begin{theorem}\label{thm:intro-folded}
For any ordering of the simple roots of $G_\nu$, every sequence in Miura form for
$G_\nu$ has a unique change of frame with values in its positive
unipotent subgroup taking it to the Steinberg section.
Under the identification $\miuraq_{O',n}=Y_{O',b_{O',n}}$ at every
folded vertex $O'$ and index $n$, the ordinary discrete Miura
transformation $\mu_{c_\nu}$ pulls back the Steinberg coordinate at
$O$ and $k\in\mathbb Z$ to the folded character
\begin{equation}\label{eq:intro-folded-character}
 \mu_{c_\nu}^*(T_{O,k})
   =\pi\!\left(\chi_{\qaff}\!\left(L(\widetilde Y_{i,b_{O,k}})\right)\right)
 \qquad(i\in O).
\end{equation}
These pullbacks freely generate the complexified folded
$q$-character ring on this lattice, the image under $\pi$ of the
complexified simply laced $q$-character ring on the pairs
$(i,b_{O,k})$ with $i\in O$.
\end{theorem}

The folded formula specializes to an ordinary Miura $(G_\nu,q)$-oper
whose connection is represented by a meromorphic function
$A^{\mathrm M}(u)$ with values in the Miura cell $M_{c_\nu}^-$ of $G_\nu$.
If $\miuraq_O(u)$ denotes its
Miura coordinate at $O$, the substitution
$\miuraq_{O,k}=\miuraq_O(q^ku)$ gives the Steinberg coordinates of this
connection. The same shift $u\mapsto qu$ occurs at every vertex of
\eqref{eq:intro-oper-qq}. In the simply laced system, identifying the
singularity polynomials $\Lambda_i$, the twist parameters $\zeta_i$, and the $Q$-functions
within each $\nu$-orbit leaves one pair $Q_O^+,Q_O^-$ for each simple
root of $G_\nu$. For a nondegenerate polynomial solution
in the sense of \cite[\S6.3]{FKSZ}, the roots of the $Q_O^+$ satisfy
the Bethe equations of \cite[(6.16)]{FKSZ}. The derivation of these equations by
folding and the corresponding integrable model are given in
\cite[Proposition~4.4 and \S5]{FHR}.

\subsection{Sketch of the proofs of Theorems~\ref*{thm:intro} and~\ref*{thm:intro-folded}}
\label{sec:intro-minors}

\sectionoverview{Here we explain how reduction to the Steinberg section
expresses each Steinberg coordinate as a Laurent polynomial in the Miura
variables. Varying the preserved flag shows that every screening annihilates
this polynomial. Its unique dominant monomial identifies the ordinary
fundamental character, and restriction to the fixed subgroup gives the folded
character.}

For a $(G,\eta,\sigma)$-oper the prescribed relative position is
$c=s_{i_1}\cdots s_{i_s}$. When $\sigma\ne1$, the list
$i_1,\ldots,i_s$ consists of one chosen vertex from each orbit of
$\sigma$ on the simple roots of $G$. For an ordinary $(G,q)$-oper,
$\sigma=\mathrm{id}_G$ and each orbit is a single vertex, so the list
contains all the simple roots. The factorization in $M_c^-$ uses one
factor for each reflection in this product. At index $k$ it gives
\[
 A_k^{\mathrm M}=P_{i_1}(\miuraq_{i_1,k})\cdots P_{i_s}(\miuraq_{i_s,k}),
 \qquad
 P_j(t)=\varphi_j\!\begin{pmatrix}t&0\\1&t^{-1}\end{pmatrix}.
\]
Factoring products of consecutive connection maps determines the unique
changes of frame $u_k\in N_+$ for which
$A_k=u_{k+1}A_k^{\mathrm M}\sigma(u_k)^{-1}$ lies in the Steinberg
section at every index. Each $u_k$ depends regularly on finitely many
Miura variables. The pullback $\mu_{c,\sigma}^*T_{i,k}$ of the principal
minor $T_{i,k}=\Delta_i(A_k)$ is therefore a regular function on a
finite product of tori, hence a Laurent polynomial in the Miura variables.

\subsubsection{Varying the preserved flag}
The value of this Laurent polynomial is constant as the preserved flag
varies over a fixed sequence $(A_k)$ in the Steinberg section. In the
Steinberg frame the flag at index $k$ is $u_kB_-$. To vary this flag
we use changes of frame $v_k$ in positive simple-root
subgroups. We require $v_k$ to vary to first order from the identity
and the transformed sequence to remain in Miura form, so that
\[
 v_{k+1}A_k^{\mathrm M}=\widehat A_k^{\mathrm M}\sigma(v_k),
 \qquad \widehat A_k^{\mathrm M}\in M_c^- .
\]
Indeed, $u_kv_k^{-1}$ then takes the varied sequence back to $(A_k)$,
and the new preserved flag in the Steinberg frame is $u_kv_k^{-1}B_-$.
Thus the derivative of each Steinberg coordinate along this variation
is zero. For the differential Miura transformation, differentiating
along changes of Miura structure gives the screening operators
\cite[\S9.9]{Young}.

The B\"acklund transformations for the additive shift $u\mapsto u+h$
\cite[Theorem~3.1]{MV} and the multiplicative shift $u\mapsto qu$
\cite[Proposition~7.1]{FKSZ} both use a change of frame in a simple-root
subgroup that preserves Miura form. To obtain the screening operator,
we linearize this condition at the identity by factoring
$v_{k+1}A_k^{\mathrm M}$ from left to right. For
$x_i(z)=\varphi_i\bigl(\begin{smallmatrix}1&z\\0&1\end{smallmatrix}\bigr)$
the two identities we need are
\[
 \begin{aligned}
 x_i(z)P_i(t)&=P_i(t+z)\,x_i\!\left(\frac{z}{t(t+z)}\right),\\
 x_i(z)P_j(t)&=P_j(t)\,x_i\!\left(t^{-\langle\alpha_i,\alpha_j^\vee\rangle}z\right)
                &&(j\ne i).
 \end{aligned}
\]
The first follows by multiplication in $\mathrm{SL}_2$ for $t$ and
$t+z$ invertible. In the second, $x_i(z)$ commutes with the negative
simple-root subgroup at $j$, and conjugation by the inverse torus
factor $\alpha_j^\vee(t)^{-1}$ rescales $z$ by the indicated power.
Thus the Miura coordinate changes only when the two roots coincide.
In either case the positive factor on the right belongs to the same
root subgroup as the one on the left.

The first-order change of frame has the form
$v_k=x_{h_k}(\varepsilon r_k)$ with $\varepsilon^2=0$, where the vertex $h_k$ specifies the positive simple-root subgroup containing $v_k$. As $v_{k+1}$ crosses the factors of $A_k^{\mathrm M}$,
the subgroup stays the same by the two identities. The element left
on the right must be
$\sigma(v_k)=x_{\sigma(h_k)}(\varepsilon r_k)$, so for a nonzero
variation $h_{k+1}=\sigma(h_k)$. The indices $h_k$ therefore run through
a single $\sigma$-orbit. This orbit has $o_i$ vertices, exactly one of which occurs in the list $i_1,\ldots,i_s$. We denote that vertex by $i$. There is a unique $0\le\rho<o_i$ for which $h_{\rho+1}=i$.
It follows that the root subgroup of $v_{k+1}$ is the image of
$x_{\sigma^{k-\rho}(i)}$. At index $k$, the factorization changes a
Miura coordinate only if the subgroup of $v_{k+1}$ matches a factor
in $A_k^{\mathrm M}$. In this orbit the only such factor is
$P_i(\miuraq_{i,k})$, and the match occurs exactly when $k=\rho+n o_i$ for
some $n\in\mathbb Z$.

For instance, the nontrivial diagram automorphism in type $A_3$ interchanges $1$ and $3$ and fixes $2$. The ordered representatives $1,2$ give
$A_k^{\mathrm M}=P_1(\miuraq_{1,k})P_2(\miuraq_{2,k})$. If $h_1=1$, the subgroup
of $v_{k+1}$ is at root $1$ for even $k$ and at root $3$ for odd $k$.
At odd $k$ the second identity applies to both factors and only
rescales the parameter of $v_{k+1}$. Thus this variation changes
$\miuraq_{1,k}$ at even $k$. For a fixed vertex of $\sigma$, the orbit has
size one and the variation changes its Miura coordinate at every index.

\subsubsection{From variations of the flag to quantum screening operators}
At each index $k=\rho+n o_i$ where the Miura coordinate changes,
the positive root element immediately before $P_i(\miuraq_{i,k})$ has
the form $x_i(\varepsilon\miuraq_{i,k}d_{i,k})$ for a unique $d_{i,k}$,
since $\miuraq_{i,k}$ is invertible. The identity for $x_i(z)P_i(t)$
then gives $\miuraq_{i,k}(1+\varepsilon d_{i,k})$ as the varied
coordinate, so $d_{i,k}$ measures its first-order change relative to
$\miuraq_{i,k}$. At index $\rho$, each factor preceding
$P_i(\miuraq_{i,\rho})$ multiplies the root subgroup parameter by an invertible
monomial. We choose the parameter of $v_{\rho+1}$ so that the element
immediately before $P_i(\miuraq_{i,\rho})$ is $x_i(\varepsilon \miuraq_{i,\rho})$.
This sets $d_{i,\rho}=1$.

The first factorization identity gives the outgoing parameter
$\varepsilon d_{i,k+o_i}/\miuraq_{i,k+o_i}$ just after $P_i(\miuraq_{i,k+o_i})$.
The positive element then crosses the remaining factors at $k+o_i$,
all factors at the intermediate indices, and the factors preceding
$P_i(\miuraq_{i,k})$ at $k$. The factor leaving $A_r^{\mathrm M}$ is
$\sigma(v_r)$, while the factor entering $A_{r-1}^{\mathrm M}$ is $v_r$.
Thus the transition from $r$ to $r-1$ changes the root label by
$\sigma^{-1}$ and preserves the scalar parameter.
We denote by $M_{i,k}$ the product of the rescaling factors
$t^{-\langle\alpha_\ell,\alpha_j^\vee\rangle}$ from the second identity,
where $\ell$ is the current root label and $t$ is the Miura coordinate
of the factor at $j$ being crossed. The resulting parameter at
$P_i(\miuraq_{i,k})$ is $\varepsilon \miuraq_{i,k}d_{i,k}$, so
\[
 \miuraq_{i,k}d_{i,k}=M_{i,k}\frac{d_{i,k+o_i}}{\miuraq_{i,k+o_i}},\qquad
 d_{i,k+o_i}=\mathsf R_{i,k}d_{i,k},\qquad
 \mathsf R_{i,k}=\frac{\miuraq_{i,k}\miuraq_{i,k+o_i}}{M_{i,k}} .
\]
With $d_{i,\rho}=1$, the recurrence gives
$d_{i,\rho+o_i}=\mathsf R_{i,\rho}$ and
$d_{i,\rho-o_i}=\mathsf R_{i,\rho-o_i}^{-1}$. Iteration determines
every $d_{i,\rho+n o_i}$ since each multiplier is invertible.
The coefficient of $\varepsilon$ in $F$ evaluated at the varied
coordinates defines the vector field
\[
 \mathscr D_{i,\rho}F
   =\sum_{n\in\mathbb Z}
       d_{i,\rho+n o_i}\,\miuraq_{i,\rho+n o_i}
       \frac{\partial F}{\partial \miuraq_{i,\rho+n o_i}} .
\]
Only the finitely many coordinates occurring in the Laurent polynomial
$F$ contribute to this sum.

Under $\miuraq_{j,\ell}=Y_{j,a\qaff^{2\ell+\xi_j}}$, the multiplier $\mathsf R_{i,k}$ is
the simple-root monomial $\mathsf R_i(b_k\qaff_i)$, where
$b_k=a\qaff^{2k+\xi_i}$ and $\qaff_i=\qaff^{o_i}$. Indeed,
$\mathsf R_i(b)$ has numerator $Y_{i,b\qaff_i^{-1}}Y_{i,b\qaff_i}$. For
each neighbor $j$ with $C_{ji}=-m$ in the Cartan matrix of
$\mathfrak a$, its denominator contains
$\prod_{r=0}^{m-1}Y_{j,b\qaff^{m-1-2r}}$. Thus $\mathsf R_i(b)$ is
the monomial denoted by $A_{i,b}$ in \cite[(2.15)]{FM}.
At $b=b_k\qaff_i$ these factors are the variables encountered at
the neighboring roots.

This recurrence identifies the vector field with a component of the
quantum screening $S_i$. The values of $S_i$ lie in a module over
the Laurent ring in the spectral variables. The module has generators
$s_{i,b}$ subject to
$s_{i,b\qaff_i^2}=\mathsf R_i(b\qaff_i)s_{i,b}$ \cite[\S7.1]{FR}.
The action on a generator and on a Laurent polynomial is
\[
 S_i(Y_{j,b})=
 \begin{cases}
   Y_{i,b}s_{i,b},& j=i,\\
   0,& j\ne i,
 \end{cases}
 \qquad
 S_i(F)=\sum_{b\in\mathbb C^\times}s_{i,b}Y_{i,b}
                 \frac{\partial F}{\partial Y_{i,b}} .
\]
Since $b_{k+o_i}=b_k\qaff_i^2$, the symbols and logarithmic derivatives
obey the same recurrence. Starting at $\rho=0,\ldots,o_i-1$ gives
$s_{i,b_k}=d_{i,k}s_{i,b_\rho}$ for $k=\rho+n o_i$. In particular
$S_i(\miuraq_{i,k})=\miuraq_{i,k}d_{i,k}s_{i,b_\rho}$ at these indices.
Grouping the terms of $S_i(F)$ by their starting index now gives
\[
 S_i(F)=\sum_{\rho=0}^{o_i-1}s_{i,b_\rho}\,\mathscr D_{i,\rho}F .
\]
Since $\qaff$ has infinite order, the parameters $b_0,\ldots,b_{o_i-1}$
lie in distinct orbits under multiplication by $\qaff_i^2$.
The relations determine every symbol on an orbit from its initial
symbol and give no relations between different orbits. Thus
$S_i(F)=0$ precisely when $\mathscr D_{i,\rho}F=0$ for every $\rho$.
Each vector field differentiates the choice of preserved flag while the
Steinberg sequence stays fixed. Hence $S_i(\mu_{c,\sigma}^*T_{j,\ell})=0$
for all $i,j,\ell$.

\subsubsection{Identification of the fundamental \texorpdfstring{$q$}{q}-character}
The calculation above puts $\mu_{c,\sigma}^*T_{i,k}$ in the common
screening kernel. To recognize this one polynomial, we use the
dominant-monomial criterion of \cite[Lemmas~5.6 and~6.17]{FM}.
Two elements of the kernel agree if the coefficients of every dominant
monomial agree, since every nonzero element contains a dominant monomial.
The fundamental character $\chi_{\qaff}(L(Y_{i,b_k}))$ belongs to
the kernel and has $Y_{i,b_k}$ as its only dominant monomial, with
coefficient one \cite[Corollaries~4.5 and~5.5]{FM}.

In $\mu_{c,\sigma}^*T_{i,k}$ the monomial $\miuraq_{i,k}$ occurs with
coefficient one and is the only dominant monomial. For the chain ordering in type $A$ this follows
from the companion matrix formula. For bipartite orderings in $D,E$ and
the four twisted cases, we compute $\Delta_i(A_k)=v_i^*(A_kv_i)$.
In the Miura product, $P_i(\miuraq_{i,k})v_i$ has highest-weight component
$\miuraq_{i,k}v_i$. The other factors and the unipotent changes of frame fix
$v_i$, so this component contributes $\miuraq_{i,k}$ to the minor.
The remaining terms involve only Miura variables from preceding factors
at index $k$ and from the later indices on which $u_{k+1}$ depends.
No nonnegative integral combination of the torus weights of these variables
equals the weight of the minor, so every remaining monomial has a negative
exponent. Changes of simply laced Coxeter ordering
preserve dominance. Under $\miuraq_{i,k}=Y_{i,b_k}$ the pullback and
$\chi_{\qaff}(L(Y_{i,b_k}))$ thus have the same dominant term. The
criterion gives \eqref{eq:intro-character}.

\subsubsection{Folded \texorpdfstring{$q$}{q}-characters by restriction to the fixed subgroup}
The simply laced identity \eqref{eq:intro-character} applies in
particular to a Miura sequence in $G_\nu=(\widetilde G^\nu)^\circ$.
The chosen ordering of the folded simple roots gives a $\nu$-stable
Steinberg section in $\widetilde G$. Since the sequence is fixed by
$\nu$, both $(u_k)$ and $(\nu(u_k))$ take it to this section.
Uniqueness gives $\nu(u_k)=u_k$ for every $k$, so the Steinberg
representative is fixed by $\nu$. By uniqueness of its ordered
factorization, its coordinates agree at all vertices of each
$\nu$-orbit $O$. Restriction of the Miura variables to $G_\nu$ induces
the homomorphism $\pi$ on spectral variables. Applying it to
\eqref{eq:intro-character} thus identifies the pullback of the
coordinate at $O$ with the folded fundamental character
in~\eqref{eq:intro-folded-character}, as asserted in
Theorem~\ref{thm:intro-folded}.

\subsection{Discrete universal Wronskians and the common screening kernel}
\label{sec:intro-applications}

\sectionoverview{Here we explain how varying a discrete universal
Wronskian shows that every Laurent polynomial annihilated by all the
screenings is a polynomial in the Steinberg coordinates.}

On the spectral lattice of Theorem~\ref{thm:intro}, the ordinary
fundamental $q$-characters generate the common screening kernel.
Their identification with the pullbacks of the Steinberg coordinates
therefore gives a geometric realization of the Frenkel--Mukhin theorem
on this lattice. To see directly why every Laurent polynomial
annihilated by the screenings is such a pullback, we vary the initial
value of a discrete universal Wronskian.

For fixed $(A_k)$ and $Z_0\in H$, the equation $W_{k+1}A_k=Z_0\sigma(W_k)$
determines $(W_k)$ in both directions from any $W_0\in G$. The map
$h\mapsto Z_0\sigma(h)$ preserves $B_-$, so $W_k^{-1}B_-$ is a
preserved flag for $(A_k)$. Varying $W_0$ gives every choice of
the flag at index zero. Where the preserved flag is opposite to the
standard flag $B_+$ at both $k$ and $k+1$, both $W_k$ and $W_{k+1}$ factor in
$B_-N_+$. The torus parts of these factorizations give
\[
 \miuraq_{i,k}=\omega_i(Z_0)
       \frac{\Delta_{\sigma^{-1}(i)}(W_k)}{\Delta_i(W_{k+1})}.
\]
When a meromorphic Wronskian exists, its restriction to a spectral
orbit gives the same Miura coordinate wherever both factorizations exist
\cite{QDE,QDEII}.

Substituting the ratios of minors into a Laurent polynomial $F$ gives
a rational function of $(A_k)$ and $g=W_0^{-1}$. Its dependence on
$g$ is through the initial flag $gB_-$. Varying $g$ to $gx_j(z)$
in a positive simple-root subgroup and differentiating at $z=0$ gives,
where the factorizations exist, an invertible multiple of the
corresponding normalized derivative $\mathscr D_{i,\rho}$. If all
screenings annihilate the polynomial, each of these derivatives is zero.
The same calculation starting at $gx_j(z)$ shows that the rational
function of $z$ has zero derivative and is therefore constant.
It is also invariant under right multiplication by $B_-$, which fixes
$gB_-$. Since $B_-$ and the positive simple-root subgroups generate
$G$, the function depends only on $(A_k)$.

This independence of $W_0$ also gives regularity. For any fixed
sequence $(A_k)$, only finitely many minors $\Delta_j(W_k)$ occur
in the denominator obtained by substituting the ratios of minors
into $F$. Each is nonzero on a nonempty open set of initial values,
since $W_0\mapsto W_k$ is an isomorphism of $G$ and
$\Delta_j(1)=1$. These open sets intersect because $G$ is irreducible.
With $W_0$ chosen in their intersection, the expression in minors
represents the descended function on a neighborhood of $(A_k)$ where
its denominator is nonzero. The function is therefore regular at
every sequence and hence polynomial in the finitely many Steinberg
coordinates on which it depends.

\subsection{Weyl group symmetry of the Miura transformation}
\label{sec:intro-weyl-symmetry}

\sectionoverview{Here we explain how the action of $G$ on the initial value
of a discrete universal Wronskian lifts the Frenkel--Hernandez reflections
and gives their formulas in each Weyl chamber. We also show how the
quantum screening operators describe the infinitesimal changes of the
Miura variables as the positive simple-root subgroups vary the preserved
flag.}

For a fixed sequence $(A_k)$ in the Steinberg section, the initial
value $W_0$ determines the discrete universal Wronskian. The flag at
index zero is $W_0^{-1}B_-$. For $h\in G$, right multiplication of
$W_0^{-1}$ by $h$ replaces this flag by $W_0^{-1}hB_-$ and leaves
every $A_k$ fixed. The connection carries the new initial flag to a preserved flag at each
index. These flags give Miura coordinates at $k$ when they are opposite to the standard
flag $B_+$ at both $k$ and $k+1$. For $h=x_j(t)$
in a positive simple-root subgroup, differentiation at $t=0$ acts on
functions of the Miura variables as a nonzero rational multiple of
the corresponding screening component.
The elements $x_j(t)$ for $j$ in the $\sigma$-orbit of $i$ commute.
Differentiating the action of their product at $t=0$ sums the
derivatives over this orbit. To identify the sum with the screening
operator $S_i$, we represent the symbols in its target module by rational
functions of the Steinberg sequence and the initial value $W_0$.
Under this embedding, the derivative of a Laurent polynomial $F$ in the
Miura variables along the product action is the image of $S_iF$.

The rank-one homomorphism $\mathrm{SL}_2\to G$ fixed by the pinning
at $j$ gives both $x_j(t)$ and the representative $\dot s_j$ of the
simple reflection $s_j$. In the twisted
construction the product of these representatives over the
$\sigma$-orbit of $i$ represents a simple reflection of the connected
fixed subgroup $(G^\sigma)^\circ$. The Cartan matrix of this subgroup is the transpose of that of $\mathfrak a$, so the
two Weyl groups have the same braid relations, including the relations
of lengths four and six at double and triple edges. The pinned
representatives satisfy these relations as group elements. Thus the
two products in each relation give the same Wronskian and the same
total unipotent change of frame. The square of each representative
is the corresponding simple coroot evaluated at $-1$. Consequently
right translation by the representative acts as an involution on
torus-invariant functions.

To compare these reflections with those of Frenkel and Hernandez, we
compute the action of the pinned representative at $i$ on the Miura
variables. It fixes the variables at the other vertices. If
$\miuraq'_{i,k}$ denotes its value on $\miuraq_{i,k}$, Gauss factorization
in the corresponding rank-one groups gives a rational function
$\Sigma_{i,k}$ of the Wronskian for which
\[
 \frac{\miuraq'_{i,k}}{\miuraq_{i,k}}=1-\frac{1}{\Sigma_{i,k}},
 \qquad
 \Sigma_{i,k+o_i}=1+\mathsf R_{i,k}^{-1}\Sigma_{i,k}.
\]
Here $o_i$ is the size of the $\sigma$-orbit of $i$, and $\mathsf R_{i,k}$
is the root monomial in the recurrence for the $i$th screening.
Frenkel and Hernandez use the same recurrence to define a reflection
on series in the spectral variables. The sign of the simple root in a
Weyl chamber determines which direction of iteration gives a series
with weights bounded above in that chamber. In
Theorem~\ref{fh:intertwining} we construct these expansions
simultaneously for all simple roots. Evaluation in the two adjacent
chambers $w$ and $ws_i$ then identifies the action of the pinned
representative with the Frenkel--Hernandez reflection between their
completed Laurent rings. The square and braid relations on
torus-invariant functions of Wronskians thus give the Weyl group
relations on these completions. Theorem~\ref{fh:fixed-functions}
shows that a Laurent polynomial is fixed by the reflection at $i$
precisely when $S_i$ annihilates it. The common screening kernel
described in the preceding subsection is therefore the ring of
invariant Laurent polynomials.

For the folded characters of Theorem~\ref{thm:intro-folded}, the
connection takes values in the fixed subgroup $G_\nu$ and the shift
has no automorphism on the fibers. The simple reflection at a
$\nu$-orbit is induced by the product of the commuting simply laced
reflections at its vertices. After the spectral variables at these
vertices are identified, this product gives the reflection between
the corresponding completions of the folded Laurent ring. The Laurent polynomials fixed by
all these reflections form the folded $q$-character ring on the
chosen lattice (Proposition~\ref{fh:folded-weyl}).

\subsection{Cancellation of apparent poles by the \texorpdfstring{$QQ$}{QQ}-relations}
\label{sec:intro-local-poles}

\sectionoverview{Here we specialize the independent Miura variables
to meromorphic functions of a common parameter $u$. We show how a
single screening equation and the corresponding $QQ$-relation give
holomorphic Steinberg coordinates at a shifted zero of a $Q$-function
under local nonvanishing conditions.}

In this application the singularity factors $\Lambda_i$ may be
meromorphic. A holomorphic solution of the $QQ$-system, together with
prescribed scalar functions $s_i$, gives the specialization
\[
 \miuraq_{i,t}(u)=s_i(\eta^t u)
      \frac{Q_i^+(\eta^t u)}{Q_i^+(\eta^{t+o_i}u)} .
\]
If $Q_i^+(\eta^{t+o_i}u_0)=0$, an individual monomial in the Laurent
expression for a Steinberg coordinate $T_{j,k}(u)$ can have a pole at
$u_0$. The question is whether $T_{j,k}(u)$ is holomorphic there after
the monomials are summed.
In the folded formula the
connection uses the ordinary shift, so $\eta=q$ and $o_i=1$. At a vertex indexed by a folding
orbit $O$, the possible pole therefore occurs when
$Q_O^+(q^{t+1}u_0)=0$, whatever the size of $O$.

To examine the terms with this denominator, we use the fact that the
character lies in $\ker S_i$. It is therefore a polynomial in the sums
\[
 E_{i,t}=\miuraq_{i,t}+\frac{M_{i,t}}{\miuraq_{i,t+o_i}},
\]
with Laurent coefficients in the variables at the other vertices
\cite{FM,FHR}.
Here $M_{i,t}$ is the product of Miura variables at neighboring vertices
that occurs in the root monomial $\mathsf R_{i,t}$.
After substitution, the first summand has $Q_i^+(x)$ in its denominator
at $x=\eta^{t+o_i}u$. The reciprocal in the second summand puts the
same $Q_i^+(x)$ in its denominator.

To isolate a possible pole at vertex $i$, we assume the $Q$-functions
are holomorphic at the shifted arguments used and the required scalar
factors and shifted $Q_\ell^+$ at vertices $\ell\ne i$ are holomorphic
and nonzero there. We also assume that the right-hand side of the
$i$th $QQ$-relation is nonzero at each zero of $Q_i^+(x)$ in these
denominators. With the relative coefficient of the two summands fixed
by the scalar normalization, we combine the $i$th $QQ$-relation at
$x$ with the same relation at $x/\eta^{o_i}$. Eliminating $Q_i^-(x)$
expresses the numerator of $E_{i,t}(u)$ as $Q_i^+(x)$
times a holomorphic function. Thus $E_{i,t}(u)$ is holomorphic at
$u_0$ even for a multiple zero. The Laurent coefficients in the
character are holomorphic there by the assumptions at the other
vertices. It follows that $T_{j,k}(u)$ and its quotient by
$s_j(\eta^k u)$ are holomorphic at $u_0$. The divisibility holds in
the local ring of holomorphic germs, so the argument applies equally to
polynomial and to holomorphic $Q$-functions. If a neighboring shifted
zero also makes the right-hand side of the $QQ$-relation vanish, the
remaining denominator requires a separate analysis.

When these hypotheses hold for all coordinates of the Steinberg
product, they give a holomorphic representative near $u_0$. For an
ordinary $(G,q)$-oper with simply laced $G$, this leads to a local comparison with
the construction of Elliott and Pestun. With a fixed torus
trivialization and the same order and lifts of the simple
reflections, the flag description gives the Steinberg product of their
multiplicative Hitchin section \cite[Definitions~9.1--9.2]{EP}. They
relate the apparent poles of its triangular presentation to the Bethe
equations \cite[Remark~9.10]{EP}. Their global $(G,q)$-oper locus in the
complex structure for $q$-connections is the framed Hitchin section,
with prescribed singularities and framing at infinity
\cite[Remarks~2.19 and~5.7]{EP}.

\subsection{Application to twisted monodromy}
\label{sec:intro-monodromy}

\sectionoverview{In this subsection we apply Theorem~\ref{thm:intro} to the
twisted monodromy of a family of differential equations in the
$q$DE/ODE/IM correspondence. When the
monodromy is in Steinberg form, the central connection matrix comparing
solutions normalized at zero and at infinity satisfies the Wronskian
equation. Its principal minors give the Miura variables, and evaluating
the fundamental $q$-characters at these variables gives the Steinberg
coordinates of twisted monodromy.}

Affine opers provide families of differential equations for the
$q$DE/ODE/IM correspondence. In the loop realization, an affine
oper gives a differential equation in $z$ for each value of a
spectral parameter $u$. Feigin and Frenkel conjectured that affine
opers for the Langlands dual affine algebra encode the spectra of
the quantum KdV Hamiltonians \cite[\S\S4--5]{FFaffine}.
For affine opers with singularities only at zero and infinity,
Masoero, Raimondo, and Valeri derived the Bethe equations by
studying relations among solutions of the differential equations
in simply laced and non-simply-laced types \cite{MRV1,MRV2}.
Frenkel and Hernandez proved the corresponding functional relations
on the representation-theoretic side and gave explicit forms of the
affine opers, including those for twisted affine algebras
\cite[\S\S3--4 and 8]{FHaffine}. Masoero and Raimondo extended the
construction to affine opers with additional singularities
of trivial local monodromy, obtaining solutions of the $QQ$-system
from entries of the central connection matrix \cite{MRaffine}.

The comparison uses two normalized $G$-frames of solutions of the
same differential equation. The local behavior at $z=0$ determines
the Frobenius frame, and the asymptotic behavior in a sector at
$z=\infty$ determines the sectorial frame. We assume that both
continue to a common domain in $z$ and depend meromorphically on $u$.
For generic untwisted simply laced affine opers the analytic hypotheses
are verified in \cite{QDE}. On that common domain,
the change from the sectorial to the Frobenius frame is independent
of $z$ and gives a meromorphic function $C(u)$ with values in $G$.
In a linear representation of $G$, the
matrix of this change of frame is the central connection matrix.

These two frames describe the monodromy associated with a rotation of $z$.
If $r$ is the order of $\sigma$, a rotation through $1/r$ of a turn
induces a pullback from solutions at $\eta u$ to solutions at $u$.
Its inverse is the twisted monodromy from $u$ to $\eta u$.
Composing $r$ such maps gives the monodromy for a full turn, from
$u$ to $qu$ with $q=\eta^r$. The map for a single rotation
acts in the sectorial frame by $h\mapsto\mathscr M(u)\sigma(h)$, and in the
Frobenius frame it acts by $h\mapsto Z_0\sigma(h)$ for a fixed $Z_0\in H$
\cite{QDE,QDEII}.
Changing frames at both ends gives
\[
 C(\eta u)\mathscr M(u)=Z_0\sigma(C(u)).
\]
This identity is the meromorphic change of frame to the constant
action $h\mapsto Z_0\sigma(h)$. When $\mathscr M(u)$ lies in the
Steinberg section, $C(u)$ is the universal Wronskian. Along the
orbit through $u_*$, the sequence $W_k=C(\eta^ku_*)$ is therefore
the discrete universal Wronskian with initial value $W_0=C(u_*)$.

The Miura variables can now be read from $C$ along any orbit on which
it is defined and all its principal minors are nonzero. Indeed, in
Frobenius coordinates the flag defined by the sectorial frame is
$C(u)B_+$ and the preserved flag is $B_-$. They are opposite precisely where these minors are nonzero. If $\mathscr M(u)$ belongs to $M_c=N_+\dot cN_+$ on this
orbit, the Miura variable at vertex $i$ is given by
$\omega_i(Z_0)\Delta_{\sigma^{-1}(i)}(C(u))/\Delta_i(C(\eta u))$. More generally, when $\mathscr M(u)$ lies in
the larger cell $B_+\dot cB_+$, torus changes of frame $s_k\in H$
put the integer-indexed connection maps in $N_+\dot cN_+$. The ratio
at index $k$ acquires the factor
$\omega_i(s_{k+1})/\omega_{\sigma^{-1}(i)}(s_k)$. Evaluating the
fundamental $q$-characters at these ratios gives the Steinberg
coordinates in the chosen frames. On an orbit of length $L$, the
torus frames descend to its points precisely when $s_L=s_0$.
A meromorphic torus normalization on the spectral domain requires a
single meromorphic function $s$ satisfying the normalization equation
there.

The fundamental characters used in this calculation also describe the
transfer matrices of XXZ-type models. In such a model the auxiliary
module $V$ specifies the Laurent polynomial $\chi_{\qaff}(V)$, while
the quantum space $W$ and a twist specify the commuting operators
substituted for its variables \cite[Theorems~4.8 and~5.3]{FH}. For a
tensor product of finite-dimensional simple quantum loop modules and
a formal twist, this substitution gives the fundamental transfer
matrices as Steinberg coordinates. On a joint eigenvector in $W$, its
arguments become ratios of shifted Baxter polynomials with prescribed
scalar factors \cite[Theorem~5.11]{FH}. Comparing a differential
family with a particular model therefore requires an identification
of these ratios with the Miura variables obtained from $C$. When they
agree, the two evaluations of the fundamental character coincide.
In this model, dividing the transfer matrix by the product of the
scalar factors of its $R$-matrices leaves a twisted trace of rational
matrices. At every nonzero spectral parameter where these matrices are
regular, all the normalized transfer eigenvalues are regular, even
when shifted Baxter roots coincide.

\subsection{Structure of the paper}
\label{sec:intro-structure}

Section~\ref{sec:group-preliminaries} expresses the Miura variables
and the Steinberg coordinates in terms of principal minors. It defines
ordinary $(G,q)$-opers and the more general $(G,\eta,\sigma)$-opers, then
relates their Miura reductions to sequences along a multiplicative orbit.
To recognize the Laurent polynomials that arise in these variables,
Section~\ref{sec:quantum-preliminaries} recalls how the screening
operators characterize the $q$-character ring and how the coefficients
of the dominant monomials determine an element of that ring.

In simply laced type, Section~\ref{ade:normalization-section}
constructs the unique $N_+$-valued change of frame that puts a sequence
in the Coxeter cell into Steinberg form, for every ordering of the
simple reflections. With the Steinberg representative fixed,
Section~\ref{ade:screening-section} varies the preserved flag and
identifies the resulting vector fields on the Miura variables with the
components of the screening operators. Using these screenings,
Section~\ref{ade:character-section} shows that the pullback of each
Steinberg coordinate is the $q$-character of the corresponding
fundamental representation. Restriction to the identity component
of the fixed subgroup then gives the folded characters.
For the ordinary characters in types $B_n$, $C_n$, $F_4$, and $G_2$,
Section~\ref{tw:section} incorporates a diagram automorphism into the
connection maps. Section~\ref{tw:b2-section} works out this construction
in $\mathrm{SL}_4$ and gives both fundamental $q$-characters of type $B_2$.

Section~\ref{sec:wronskian-descent} uses discrete universal Wronskians
to prove that every Laurent polynomial on the spectral lattice
annihilated by all screenings is a polynomial in the Steinberg
coordinates pulled back by the Miura transformation.
Section~\ref{sec:fh-weyl} uses this description of the screening
kernel to identify the invariant Laurent polynomials for the
Frenkel--Hernandez reflections. It constructs the lifts of the
reflections to Wronskians and treats folded characters by restriction.
Section~\ref{sec:difference-opers} passes from independent lattice
variables to functions of a spectral parameter. It compares the flag
and Steinberg descriptions of meromorphic $(G,q)$-opers and
$(G,\eta,\sigma)$-opers. It derives local cancellation of apparent poles
from the $QQ$-relations and computes the Steinberg
coordinates of twisted monodromy from the central connection matrix.
The section concludes by evaluating the same characters at normalized ratios of
Baxter operators to obtain the fundamental transfer matrices of
XXZ-type models. The \hyperref[app:quantum-presentation]{Appendix} records the
Drinfeld current presentation of the quantum loop algebra and the
coproduct used for tensor products.

\section*{Acknowledgments}

The author is partially supported by NSF grant DMS-2526435
(formerly DMS-2203823) and by the Simons Foundation through its
Travel Support for Mathematicians program.

\paragraph{Use of AI tools.}
The author formulated the research questions and developed the method of
proof. ChatGPT assisted with constructing proofs, checking mathematical
arguments, and refining the exposition. The author takes full responsibility
for the manuscript, including the correctness of all results and proofs
and the accuracy of the references.

\section{Root subgroups and the Steinberg section}
\label{sec:group-preliminaries}

\sectionoverview{In this section, we construct coordinates on the Steinberg section and on
the intersection of the Coxeter cell with the opposite Borel subgroup.
We express both sets of coordinates as principal minors. The flag
descriptions of Miura $(G,q)$-opers and Miura $(G,\eta,\sigma)$-opers
then relate these varieties to connection maps along a multiplicative
orbit. We describe the regular functions on the resulting spaces of sequences.}

\subsection{Root data and pinned root subgroups}
\label{sec:pinned-root-subgroups}

\sectionoverview{In this subsection, we fix the root data of a simple simply connected
complex algebraic group $G$ and describe the automorphisms that preserve
its pinning. We construct the two rank-one factors associated with each
simple root and prove the root subgroup factorization used to describe
the Coxeter cell.}

Throughout this section, $G$ is a simple simply connected complex algebraic group with opposite Borel subgroups $B_+$ and $B_-$. Their intersection is a maximal
torus $H$. If $N_\pm$ are their unipotent radicals, then $B_\pm=HN_\pm$.
We work over a commutative unital $\mathbb C$-algebra $D$.
A point $g\in G(D)$ is a
$\mathbb C$-algebra homomorphism $\mathbb C[G]\to D$ that assigns
a value in $D$ to each regular function on $G$.
In a faithful representation $G\subset\mathrm{GL}_m$, it is an
invertible matrix over $D$ satisfying the defining equations of
$G$. For example, $\mathrm{SL}_N(D)$ consists of the matrices over
$D$ of determinant one. We write $\mathbb G_a(D)=D$ with addition
and $\mathbb G_m(D)=D^\times$, the group of invertible elements of $D$,
with multiplication.

Taking $D$ to be a Laurent polynomial ring allows the torus
coordinates to vary independently. To compute first-order changes of
regular functions, we use the algebra of dual numbers
$D[\varepsilon]/(\varepsilon^2)$. In $a+\varepsilon b$, the
coefficient $b$ specifies a first-order variation of $a$. The
product rule follows from
\[
 (a+\varepsilon b)(c+\varepsilon d)
 =ac+\varepsilon(ad+bc).
\]

The roots $\Phi$ are the nonzero weights for the action of $H$
on $\mathfrak g=\operatorname{Lie}G$ obtained by differentiating
conjugation. We denote the weight space of a root $\beta$ by
$\mathfrak g_\beta$.
The weights of $\mathfrak n_\pm=\operatorname{Lie}(N_\pm)$ give the
decomposition $\Phi=\Delta_+\sqcup\Delta_-$ into positive and negative
roots, with $\Delta_-=-\Delta_+$. The simple roots $\alpha_i$ and their
coroots $\alpha_i^\vee$ are indexed by $I$. The characters
$\lambda:H\to\mathbb G_m$ and cocharacters $\nu:\mathbb G_m\to H$
form the lattices $X^*(H)$ and $X_*(H)$, respectively. Their
composition defines the integer pairing
$\lambda(\nu(\zeta))=\zeta^{\langle\lambda,\nu\rangle}$.
Since $G$ is simply connected, the fundamental weights $\omega_i$,
defined by $\langle\omega_i,\alpha_j^\vee\rangle=\delta_{ij}$,
form a basis of $X^*(H)$. Here $\delta_{ij}$ denotes the Kronecker
delta. Our convention for the Cartan matrix of $G$ is
$A^{G}_{ij}=\langle\alpha_j,\alpha_i^\vee\rangle$.

The Weyl group $W=N_G(H)/H$ acts on characters by
$(w\lambda)(h)=\lambda(n^{-1}hn)$ for any representative $n$ of $w$
in $N_G(H)$. The simple reflection at $i$ sends
$\lambda$ to $s_i\lambda=\lambda-\langle\lambda,\alpha_i^\vee\rangle\alpha_i$.
An expression $w=s_{i_1}\cdots s_{i_s}$ is reduced if it uses the
least possible number of simple reflections. This number is the length $\ell(w)$.
The longest element $w_0$ sends $\Delta_+$ to $\Delta_-$.
A product containing each simple reflection exactly once is a Coxeter
element. Its order is the Coxeter number $h$.

For a root $\alpha$, the root subgroup $U_\alpha$ is the
one-dimensional unipotent subgroup with Lie algebra the root space
of weight $\alpha$. A pinning consists of $H$, $B_+$, and a
parametrization $x_i:\mathbb G_a\xrightarrow{\sim}U_{\alpha_i}$
for every simple root. The derivative in the additive coordinate is
$e_i=\left.\tfrac{d}{dz}x_i(z)\right|_{z=0}$. It belongs to $T_{1_G}G=\operatorname{Lie}G$
because $x_i(0)=1_G$. The pinning determines the rank-one
homomorphism $\varphi_i:\mathrm{SL}_2\to G$ whose differential sends
the raising matrix $E_{12}$ to $e_i$ and for which
$\varphi_i(\operatorname{diag}(\zeta,\zeta^{-1}))=\alpha_i^\vee(\zeta)$. For $z\in D$
the elements in the positive and negative simple-root subgroups
are
\[
 x_i(z)=\varphi_i\!\begin{pmatrix}1&z\\0&1\end{pmatrix},
 \qquad y_i(z)=\varphi_i\!\begin{pmatrix}1&0\\z&1\end{pmatrix}.
\]
Both are parametrized by addition in $D$. The negative simple-root
vector $f_i=\left.\tfrac{d}{dz}y_i(z)\right|_{z=0}$ is defined by
the same differentiation, using $y_i(0)=1_G$. Our representative
of $s_i$ is
\[
 \dot s_i=\varphi_i\begin{pmatrix}0&-1\\1&0\end{pmatrix}
          =x_i(-1)y_i(1)x_i(-1).
\]
For a product $c=s_{i_1}\cdots s_{i_s}$ of distinct simple
reflections we use the representative $\dot c=\dot s_{i_1}\cdots\dot s_{i_s}$.
The two rank-one factors
for $\zeta\in D^\times$ and $z\in D$ are
\begin{equation}\label{eq:rank-one-factors}
 P_i(\zeta)=y_i(\zeta^{-1})\alpha_i^\vee(\zeta)
       =\varphi_i\begin{pmatrix}\zeta&0\\1&\zeta^{-1}\end{pmatrix},\qquad
 X_i(z)=x_i(z)\dot s_i
       =\varphi_i\begin{pmatrix}z&-1\\1&0\end{pmatrix}.
\end{equation}

For products involving more than one reflection, the braid relations
ensure that the product $\dot w$ of the $\dot s_i$ along a reduced expression depends
only on $w$ \cite[\S1.4, (1.7)--(1.9)]{FZ}. In particular
$\dot w\dot w'=\dot{ww'}$ whenever $\ell(ww')=\ell(w)+\ell(w')$. In rank one
$\dot s_i^2=\alpha_i^\vee(-1)$. Conjugation by $\dot w$ sends
$U_\alpha$ to $U_{w\alpha}$, and for $h\in H(D)$ one has
$h x_i(z)h^{-1}=x_i(\alpha_i(h)z)$.

A diagram automorphism is a permutation $\sigma$ of $I$ that preserves
the Cartan matrix $A^G$. It lifts uniquely to an automorphism of $G$
that preserves $H$ and $B_+$ and satisfies
$\sigma(x_i(z))=x_{\sigma(i)}(z)$ \cite[\S32]{Humphreys}.
We use the same symbol for this lift, which preserves the pinning.
The rank-one homomorphisms satisfy
$\sigma\circ\varphi_i=\varphi_{\sigma(i)}$, so the lift also permutes
the negative simple-root subgroups and preserves $B_-$.

When $c$ is a Coxeter element, its Coxeter cell is the double coset
$N_+\dot cN_+$. To factor elements in this cell, we distinguish the
positive root subgroups which $\dot c^{-1}$ conjugates into $N_-$.
Their roots form $\Delta_+\cap c\Delta_-$. The same set
$\Delta_+\cap w\Delta_-$ for any $w\in W$ can be read from a
reduced expression by the Weyl group length criterion. For any $w\in W$
and simple reflection $s_i$, one has $\ell(ws_i)=\ell(w)+1$ when
$w\alpha_i$ is positive and $\ell(ws_i)=\ell(w)-1$ when it is negative
\cite[\S1.6]{HumphreysCoxeter}.

\begin{lemma}\label{lem:inversion-roots}
For $w=s_{i_1}\cdots s_{i_s}$ reduced, the roots
$s_{i_1}\cdots s_{i_{p-1}}\alpha_{i_p}$, $1\le p\le s$, are distinct
and form $\Delta_+\cap w\Delta_-$. In particular
$\ell(w)=|\Delta_+\cap w\Delta_-|$ and $\ell(w_0)=|\Delta_+|$.
An expression in pairwise distinct simple reflections is reduced.
\end{lemma}

\begin{proof}
When $\ell(ws_i)=\ell(w)+1$, the length criterion gives
$w\alpha_i\in\Delta_+$. The simple reflection $s_i$ permutes
$\Delta_-\setminus\{-\alpha_i\}$ and sends $-\alpha_i$ to $\alpha_i$, so
\[
 \Delta_+\cap ws_i\Delta_-
  =(\Delta_+\cap w\Delta_-)\sqcup\{w\alpha_i\}.
\]
Iteration along $w=s_{i_1}\cdots s_{i_s}$ gives the asserted roots
and their number. For $w_0$ the intersection is all of $\Delta_+$.

If $i_1,\ldots,i_s$ are distinct, each reflection $s_{i_a}$ with
$a<p$ changes only the coefficient of $\alpha_{i_a}$ in the
simple-root expansion. Thus the coefficient of $\alpha_{i_p}$ in
$s_{i_1}\cdots s_{i_{p-1}}\alpha_{i_p}$ is one. All simple-root coefficients of a negative root are nonpositive,
so this root is positive. The length criterion now shows that
adjoining $s_{i_p}$ increases the length by one for every $p$.
\end{proof}

We use additive parametrizations of the positive root subgroups $U_\alpha$, with $x_i$ as the parametrization at $\alpha_i$. For $S\subset\Delta_+$, the subgroups $U_\alpha$ with $\alpha\in S$ generate a subgroup denoted by $U_S$. We call
$S$ closed if $\alpha,\beta\in S$ and $\alpha+\beta\in\Phi$
imply $\alpha+\beta\in S$.
For such a set $S$, the ordered multiplication theorem for root subgroups
says that in any ordering the map $\prod_{\alpha\in S}U_\alpha\to G$
is an isomorphism onto the smooth closed subgroup $U_S$. Its image is
independent of the ordering \cite[Proposition~5.1.16(1)]{Conrad}.

\begin{lemma}\label{lem:root-factorization}
If a closed set $S$ is the disjoint union $S_1\sqcup S_2$ of closed subsets, multiplication
identifies $U_{S_1}\times U_{S_2}$ with $U_S$. For closed subsets $S_1,S_2$
and any $w\in W$,
\[
 U_{S_1}\cap U_{S_2}=U_{S_1\cap S_2},\qquad
 N_+\cap\dot w N_+\dot w^{-1}=U_{\Delta_+\cap w\Delta_+}.
\]
\end{lemma}

\begin{proof}
Applying the ordered multiplication theorem with the roots of $S_1$
first and those of $S_2$ second gives the asserted factorization of $U_S$.
In the coordinates for any fixed order of $\Delta_+$,
an element belongs to $U_{S_a}$ precisely when its coordinates outside
$S_a$ vanish. Belonging to both subgroups therefore means that all
coordinates outside $S_1\cap S_2$ vanish.

For the second intersection, each $g\in N_+(D)$ has a unique factorization $g=pv$, with $p\in U_{\Delta_+\cap w\Delta_-}(D)$ and
$v\in U_{\Delta_+\cap w\Delta_+}(D)$. Conjugation by $\dot w^{-1}$
takes $p$ into $N_-(D)$ and $v$ into $N_+(D)$. If $\dot w^{-1}g\dot w$
lies in $N_+(D)$, then $\dot w^{-1}p\dot w$ belongs to
$N_-(D)\cap N_+(D)=\{1\}$. Thus $g=v$, and every such $v$ lies
in the intersection. This argument holds for every coefficient algebra $D$.
\end{proof}

For $U_i=U_{\Delta_+\setminus\{\alpha_i\}}$, the factorization gives
$N_+=x_i(\mathbb G_a)U_i$. Moreover, $\varphi_i(\mathrm{SL}_2)$
normalizes $U_i$. Indeed, a positive root different from $\alpha_i$
has a positive coefficient at some $j\ne i$, and this coefficient is
unchanged by adding or subtracting $\alpha_i$. Hence taking the Lie
bracket with $e_i$ or $f_i$ preserves $\operatorname{Lie}U_i$.
Exponentiating shows that the two simple-root subgroups normalize
$U_i$. Their complex points generate $\varphi_i(\mathrm{SL}_2(\mathbb C))$.
For the conjugation map $\mathrm{SL}_2\times U_i\to G$, every
defining equation of $U_i$ therefore vanishes on all complex points
of $\mathrm{SL}_2\times U_i$. Since this variety is reduced, the
pullbacks of these equations are zero as regular functions. Thus
$\varphi_i(\mathrm{SL}_2(D))$ normalizes $U_i(D)$ for every $D$.
In particular $U_i$ is normal in $N_+$, and the quotient
$N_+/U_i\simeq\mathbb G_a$ gives the simple-root coordinate at $i$.

\subsection{Generalized minors and Steinberg coordinates}

\sectionoverview{In this subsection, we define generalized minors as matrix
coefficients of the fundamental representations. We use the principal
minors to characterize the Gauss cell and to give coordinates on the
Steinberg section and the intersection of the Coxeter cell with $B_-$.}

An integral weight $\lambda$ is dominant if
$\langle\lambda,\alpha_i^\vee\rangle\geq0$ for every simple coroot.
For such a weight, $V(\lambda)$ denotes the irreducible rational
$G$-module of highest weight $\lambda$. Its
highest-weight space is an $N_+$-fixed line, and all its weights
lie in $\lambda-\sum_i\mathbb Z_{\geq 0}\alpha_i$
\cite[Theorems~22.18 and~22.20]{Milne}. A representative $\dot w$
of a Weyl group element identifies the weight spaces of weights
$\mu$ and $w\mu$. In particular, each extremal weight space
$V(\lambda)_{w\lambda}$ is a line.

For each $i$ we use a highest-weight vector $v_i\in V(\omega_i)$ and
the covector $v_i^*$ taking value one on $v_i$ and zero on the other
weight spaces. In the rank-one subalgebra at $j$, the vector $v_i$ has highest
weight $m=\langle\omega_i,\alpha_j^\vee\rangle=\delta_{ij}$.
Since there are finitely many weights, there is a least $p>0$
with $f_j^p v_i=0$. The $\mathfrak{sl}_2$ relations give
\[
 0=e_j f_j^p v_i=p(m-p+1)f_j^{p-1}v_i,\qquad f_j^{p-1}v_i\ne0.
\]
Thus $p=m+1$, and $v_i,f_jv_i,\ldots,f_j^m v_i$ span the
$\mathfrak{sl}_2$-submodule generated by $v_i$.
For $j\ne i$ this submodule is trivial, so $\varphi_j(\mathrm{SL}_2)$
fixes $v_i$. For $j=i$ the
matrices in \eqref{eq:rank-one-factors} act in the basis $v_i,f_iv_i$,
giving
\begin{equation}\label{eq:rank-one-action-on-highest-weight}
 P_i(\zeta)v_i=\zeta v_i+f_iv_i,\quad X_i(z)v_i=z v_i+f_iv_i,\quad
 P_j(\zeta)v_i=X_j(z)v_i=v_i\quad(j\ne i).
\end{equation}
Applying $v_i^*$ to these identities extracts the highest-weight
coefficient, which is $\zeta$ for $P_i(\zeta)$ and $z$ for $X_i(z)$.
This same matrix coefficient and its Weyl translates define regular
functions on $G$.

\begin{definition}
\label{def:fundamental-minors}
For $g\in G$, the principal minor at $i$ is
$\Delta_i(g)=v_i^*(gv_i)$. For $w,w'\in W$, the corresponding
generalized minor in the chosen pinning is \cite[\S1.4]{FZ}
\begin{equation}\label{eq:generalized-minor}
 \Delta_{w\omega_i,w'\omega_i}(g)=v_i^*(\dot w^{-1}g\dot w'v_i).
\end{equation}
\end{definition}

These functions extract the $\dot w v_i$ coefficient
of $g\dot w'v_i$. The parabolic decomposition in the Weyl group shows why the two weights determine the minor. The stabilizer
of a dominant weight $\lambda$ is generated by the $s_j$ with
$\langle\lambda,\alpha_j^\vee\rangle=0$. If $W_J$ is the subgroup
generated by the reflections with labels in $J$, each coset $wW_J$
has a unique shortest element $u$, and
$\ell(uv)=\ell(u)+\ell(v)$ for every $v\in W_J$
\cite[\S\S1.10 and 1.12]{HumphreysCoxeter}.

For $\lambda=\omega_i$ the stabilizer is thus generated by the
$s_j$ with $j\ne i$. For these $j$, the representative
$\dot s_j$ fixes $v_i$ by \eqref{eq:rank-one-action-on-highest-weight}.
It fixes $v_i^*$ as well, since $s_j$ fixes $\omega_i$ and permutes
the other weights. In the decomposition $w=uv$ for this stabilizer,
the weight $w\omega_i$ determines $u$. Length additivity gives $\dot w=\dot u\dot v$
and $v_i^*\dot w^{-1}=v_i^*\dot u^{-1}$. The same decomposition of
$w'$ shows that $\dot w'v_i$ depends only on $w'\omega_i$.

The signs in these matrix coefficients depend on the chosen lifts.
The extremal vectors in \cite{QDE} are defined using
$\widehat s_i=x_i(1)y_i(-1)x_i(1)=\dot s_i^{-1}$ and the corresponding
products $\widehat w$ along reduced expressions. The lifts denoted
there by $\bar s_i$ and $\bar c$ agree with our $\dot s_i$ and $\dot c$.
If $\dot w=\widehat w h_w$ with $h_w\in H$, the minor
\eqref{eq:generalized-minor} differs from the one in that convention
by the factor $\omega_i(h_{w'})/\omega_i(h_w)$. Thus the two conventions
give the same minor when the two extremal weights coincide.

A Gauss factorization of $g$ is an expression $g=n_-hn_+$ with
$n_\pm\in N_\pm$ and $h\in H$. For $b=n_-h\in B_-(D)$, we write
$[b]_0=h$ for its torus factor. We use two facts for a split reductive
group with opposite Borels. The Bruhat decomposition partitions $G(\mathbb C)$
into the double cosets $B_-(\mathbb C)\dot w B_+(\mathbb C)$,
$w\in W$. Multiplication $N_-\times H\times N_+\to G$ is an open
immersion with image $B_-B_+$ \cite[Theorems~21.73 and~21.84]{Milne}.
The principal minors describe this open set and recover its torus
coordinate.

\begin{proposition}\label{prop:gauss-cell}
The Gauss cell $G_0=B_-B_+=N_-HN_+$ is the open subscheme of $G$
on which every $\Delta_i$ is invertible. In particular, for a
commutative complex algebra $D$, an element $g\in G(D)$ has a
factorization
\[
 g=n_- h n_+,\qquad n_\pm\in N_\pm(D),\qquad h\in H(D),
\]
if and only if every $\Delta_i(g)$ belongs to $D^\times$. The factors
are unique and satisfy $\Delta_i(g)=\omega_i(h)$.

For $g=ub$ in the inverse cell $G_0^{-1}=N_+B_-$, with
$u\in N_+(D)$ and $b\in B_-(D)$, the torus factor of $b$ satisfies
$\omega_i([b]_0)=\Delta_i(g^{-1})^{-1}$.
\end{proposition}

\begin{proof}
The positive unipotent subgroup fixes $v_i$. The negative root
subgroups add only lower-weight components, so they preserve the
coefficient of $v_i$, that is, $v_i^*n_-=v_i^*$. It follows that on
$G_0=N_-HN_+$,
\[
 \Delta_i(n_-hn_+)=\omega_i(h).
\]
For the converse over $\mathbb C$, the Bruhat decomposition gives
$g=b_-\dot w b_+$ with $b_-=n_-h_-$ and $b_+=h_+n_+$. Hence
\[
 \Delta_i(g)=\omega_i(h_-h_+)\,v_i^*(\dot w v_i),
 \qquad \dot w v_i\in V(\omega_i)_{w\omega_i}.
\]
If $w\ne1$, then $w\omega_i\ne\omega_i$ for some $i$, since the
fundamental weights form a basis. Since $v_i^*$ vanishes on every weight space other than
$V(\omega_i)_{\omega_i}$, the corresponding minor is zero.

Thus $G_0$ and the principal open defined by $f=\prod_i\Delta_i$
have the same complex points. Their difference in either direction
is locally closed in a complex variety of finite type, hence contains
a complex point if nonempty. The opens therefore coincide. For a
commutative complex algebra $D$, a point $g:\mathbb C[G]\to D$
factors through this open exactly when it extends to
$\mathbb C[G][f^{-1}]$. Equivalently,
\[
 g\in G_0(D)\quad\Longleftrightarrow\quad f(g)\in D^\times
 \quad\Longleftrightarrow\quad \Delta_i(g)\in D^\times\text{ for all }i.
\]
Multiplication identifies $N_-\times H\times N_+$ with this open, so
its inverse gives the unique factors over $D$. If $g=ub$ with
$b=n_-h$, then $h^{-1}n_-^{-1}h\in N_-(D)$, and the factorization of
$g^{-1}$ gives
\[
 g^{-1}=(h^{-1}n_-^{-1}h)\,h^{-1}u^{-1},
 \qquad \Delta_i(g^{-1})=\omega_i(h^{-1})=\omega_i([b]_0)^{-1}.
\]
\end{proof}

\begin{definition}
\label{def:coxeter-cells}
For a product $c=s_{i_1}\cdots s_{i_s}$ of distinct simple reflections, its
support is $J=\{i_1,\ldots,i_s\}$ and its inversion set is
$R_c=\Delta_+\cap c\Delta_-$. The associated subgroup and cells are
\begin{equation}\label{eq:universal-cells}
 \begin{gathered}
 U_c=N_+\cap\dot c N_-\dot c^{-1},\qquad M_c=N_+\dot c N_+,\\
 \mathcal A_c=U_c\dot c,\qquad M_c^-=M_c\cap B_- .
 \end{gathered}
\end{equation}
When $J=I$, $M_c$ is the Coxeter cell and $\mathcal A_c$ is the
Steinberg section for $c$.
\end{definition}

Since $w_0\Delta_+=\Delta_-$ and
$N_-=\dot w_0N_+\dot w_0^{-1}$, the intersection formula of
Lemma~\ref{lem:root-factorization} gives $U_c=U_{R_c}$.
The intersection defining this subgroup is independent of the lift of
$c$ because $H$ normalizes $N_-$.

For $J=I$, Steinberg's theorem for a simply connected semisimple
group over $\mathbb C$ says that the ordered product
$\prod_{i\in I}x_i(T_i)\dot s_i$ meets each regular conjugacy class
in exactly one point. Here an element is regular if its centralizer
has dimension $\dim H$. The regular functions on $G$ invariant
under conjugation form the ring $\mathbb C[G]^G$. The quotient map
$G\to\operatorname{Spec}\mathbb C[G]^G$ restricts to an isomorphism
on this section \cite[Theorem~1.4 and Corollary~7.16]{Steinberg}.
For any support $J$, the principal minors give coordinates on both
$\mathcal A_c$ and $M_c^-$. On $M_c^-$ the $i$th Miura coordinate is
$\miuraq_i=\Delta_i|_{M_c^-}$ for $i\in J$. On a sequence we write
$\miuraq_{i,k}=\Delta_i(A_k^{\mathrm M})$ for its value at index $k$.

\begin{proposition}
\label{prop:cell-parametrizations}
For $T=(T_i)_{i\in J}$ and $\boldsymbol{\miuraq}=(\miuraq_i)_{i\in J}$, the ordered products
\[
 a(T)=X_{i_1}(T_{i_1})\cdots X_{i_s}(T_{i_s}),\qquad
 A^{\mathrm M}(\boldsymbol{\miuraq})=P_{i_1}(\miuraq_{i_1})\cdots P_{i_s}(\miuraq_{i_s})
\]
define isomorphisms
\[
 \mathbb A^J\times N_+\longrightarrow M_c,\quad (T,w)\longmapsto a(T)w,
 \qquad
 (\mathbb G_m)^J\longrightarrow M_c^-,\quad \boldsymbol{\miuraq}\longmapsto A^{\mathrm M}(\boldsymbol{\miuraq}).
\]
The first identifies $\mathbb A^J\times\{1\}$
with $\mathcal A_c$. For $i\in J$,
$\Delta_i(a(T))=T_i$ and $\Delta_i(A^{\mathrm M}(\boldsymbol{\miuraq}))=\omega_i([A^{\mathrm M}(\boldsymbol{\miuraq})]_0)=\miuraq_i$.
For $j\notin J$ both minors equal one.
\end{proposition}

\begin{proof}
The sign of $c^{-1}\alpha$ partitions $\Delta_+$ into the two
closed subsets $R_c$ and $\Delta_+\cap c\Delta_+$. Thus, with
$V=N_+\cap\dot cN_+\dot c^{-1}$, Lemma~\ref{lem:root-factorization}
gives $N_+=U_cV$. Since $\dot c^{-1}V\dot c\subset N_+$,
this gives $M_c=U_c\dot c N_+$. For $g=p\dot c w$ in this product,
$\dot c^{-1}g=(\dot c^{-1}p\dot c)w$ is its Gauss factorization with
torus factor $1$. Thus, in the Gauss coordinates
$\dot c^{-1}g=n_-hn_+$ on $\dot cG_0$, one has
\[
 g\in M_c\quad\Longleftrightarrow\quad
 h=1,\quad n_-\in\dot c^{-1}U_c\dot c .
\]
These conditions define a closed subset of $\dot cG_0$, and the inverse of
$(p,w)\mapsto p\dot c w$ is the regular map
$g\mapsto(\dot c n_-\dot c^{-1},n_+)$ over every $D$.

For $w_{p-1}=s_{i_1}\cdots s_{i_{p-1}}$, we have
\[
 a(T)\dot c^{-1}=\prod_{p=1}^s
   \dot w_{p-1}x_{i_p}(T_{i_p})\dot w_{p-1}^{-1}.
\]
Each conjugated parametrization is an isomorphism onto
$U_{w_{p-1}\alpha_{i_p}}$. By Lemma~\ref{lem:inversion-roots} the roots
$w_{p-1}\alpha_{i_p}$ enumerate $R_c$, so ordered multiplication
identifies $\mathbb A^J$ with $U_c$. Translation by $\dot c$ gives
the section $\mathcal A_c$, and right multiplication by $N_+$ gives
the first isomorphism.

For $i=i_p$, every $X_j(T_j)$ with $j\ne i$ fixes $v_i$. Applying
\eqref{eq:rank-one-action-on-highest-weight} at the $p$th factor gives
\[
 a(T)v_i=T_i v_i+
     \left(\prod_{a<p}X_{i_a}(T_{i_a})\right)f_iv_i.
\]
The second term has weights in
$\omega_i-\alpha_i+\sum_{a<p}\mathbb Z\alpha_{i_a}$, since
$X_j(T_j)$ shifts weights only by multiples of $\alpha_j$.
This set omits $\omega_i$ because the $i_a$ are distinct. Thus
$\Delta_i(a(T))=T_i$. For $i\notin J$, all factors fix $v_i$.

If $g=a(T)w\in M_c^-(D)$, right $N_+$-invariance and
Proposition~\ref{prop:gauss-cell} give
\[
 T_i=\Delta_i(a(T))=\Delta_i(g)=\omega_i([g]_0)\in D^\times
 \qquad(i\in J).
\]
The rank-one matrices in \eqref{eq:rank-one-factors} yield
$X_i(\zeta)=P_i(\zeta)x_i(-\zeta^{-1})$. For $i\ne j$, $x_i$
commutes with $y_j$ by $[e_i,f_j]=0$, and conjugation by
$\alpha_j^\vee(\zeta)$ gives
\begin{equation}\label{ade:rank-one-distinct}
 x_i(s)P_j(\zeta)
 =P_j(\zeta)x_i\bigl(\zeta^{-\langle\alpha_i,\alpha_j^\vee\rangle}s\bigr),
 \qquad s\in D,\quad \zeta\in D^\times .
\end{equation}
For any tuple of units $\boldsymbol{\miuraq}$, these identities move each
$x_{i_p}$ to the right of all the later $P_{i_a}$, since the labels are
distinct. The resulting product is
\[
 a(\boldsymbol{\miuraq})=A^{\mathrm M}(\boldsymbol{\miuraq})\prod_{p=1}^s
 x_{i_p}\!\left(-\miuraq_{i_p}^{-1}
  \prod_{a=p+1}^s \miuraq_{i_a}^{-\langle\alpha_{i_p},\alpha_{i_a}^\vee\rangle}\right).
\]
Thus $A^{\mathrm M}(\boldsymbol{\miuraq})$ belongs to $M_c(D)$,
and it lies in $B_-(D)$ because each $P_j$ does. For the tuple $T$
associated with $g$ above, both $g$ and $A^{\mathrm M}(T)$ lie in
$a(T)N_+(D)$ and in $B_-(D)$. Their quotient belongs to
$N_+(D)\cap B_-(D)=\{1\}$, so $g=A^{\mathrm M}(T)$.
Conversely, projection $B_-\to H$ gives
\[
 [A^{\mathrm M}(\boldsymbol{\miuraq})]_0
  =\prod_{j\in J}\alpha_j^\vee(\miuraq_j),\qquad
 \Delta_i(A^{\mathrm M}(\boldsymbol{\miuraq}))
  =\begin{cases}\miuraq_i,&i\in J,\\1,&i\notin J.\end{cases}
\]
Thus $g\mapsto(\Delta_i(g))_{i\in J}$ is the inverse over every $D$.
\end{proof}

For $J=I$ the description of $M_c^-$ over a field recovers the
factorization of \cite[Theorem~2.5 and \S2.5]{FKSZ}, after interchanging
the two Borels. For a proper subset $J$ the corresponding field
factorization takes place in the semisimple subgroup generated by the
rank-one groups at its vertices.

\subsection{From Miura \texorpdfstring{$(G,q)$}{(G,q)}-opers and Miura
\texorpdfstring{$(G,\eta,\sigma)$}{(G,η,σ)}-opers to sequences in \texorpdfstring{$M_c^-$}{M\string_c\string^-}}
\label{sec:miura-opers-orbit}

\sectionoverview{In this subsection we define ordinary $(G,q)$-opers
and explain how a diagram automorphism enters the definition of a
$(G,\eta,\sigma)$-oper. In both cases a preserved flag gives the
Miura reduction. We relate its restriction to a multiplicative orbit
to the sequences in $M_c^-$ and describe the coordinate rings of
sequences in Miura and Steinberg form.}

Both definitions use the relative position of two flags in the same
fiber of a principal right $G$-bundle $P\to U$. A meromorphic reduction to $B_\pm$ is a
meromorphic section of $P/B_\pm$, so its value at $u$ is a flag in
the fiber $P_u$. If two flags of type $B_+$ are represented in a
common frame by $xB_+$ and $yB_+$, their relative position is the
Weyl group element $w$ characterized by $x^{-1}y\in B_+\dot wB_+$.
For an ordinary $q$-connection, $q\in\C^\times$ and the base
$U\subset\C$ is a domain satisfying $qU=U$.

\begin{definition}
\label{def:ordinary-flag-oper}
A meromorphic $q$-connection on $P$ is a meromorphic family of
$G$-equivariant isomorphisms $\mathcal A_u:P_u\to P_{qu}$. For a
Coxeter element $c$ of the Weyl group of $G$, a $(G,q)$-oper consists
of $P$, such a connection, and a meromorphic $B_+$-reduction $F_+$
with
\[
 \operatorname{pos}\bigl(F_+(qu),\mathcal A_u(F_+(u))\bigr)=c
\]
on a dense open subset of $U$. A Miura $(G,q)$-oper also has a
meromorphic $B_-$-reduction $F_-$ satisfying
$\mathcal A_u(F_-(u))=F_-(qu)$ wherever both sides are defined.
\end{definition}

Our convention for the two reductions interchanges $B_+$ and $B_-$
in \cite[Definitions~2.1--2.2]{FKSZ}. To extend the definition, we
allow a pinned diagram automorphism $\sigma$ of $G$ to act on the
fibers as the base point moves from $u$ to $\eta u$. Here
$\eta\in\C^\times$ and $\eta U=U$. Equivariance for the right
action of $G$ on $P$ then takes the following form.

\begin{definition}
\label{def:eta-sigma-connection}
A meromorphic $(\eta,\sigma)$-connection on $P$ is a meromorphic
family of isomorphisms $\mathcal A_u:P_u\to P_{\eta u}$ satisfying
\[
 \mathcal A_u(pg)=\mathcal A_u(p)\sigma(g)
 \qquad (p\in P_u,\ g\in G).
\]
In a meromorphic frame $e$, the identity
$\mathcal A_u(e(u))=e(\eta u)A(u)$ determines the $G$-valued function
$A$. In this frame the connection sends $(u,h)$ to
$(\eta u,A(u)\sigma(h))$.
\end{definition}

Since $\sigma(B_\pm)=B_\pm$, the connection sends a flag $xB_\pm$
at $u$ to $A(u)\sigma(x)B_\pm$ at $\eta u$. For a $(G,\eta,\sigma)$-oper,
we compare $F_+(\eta u)$ with the image of $F_+(u)$ under this single
map. Their relative position is specified by an ordered list $J=(i_1,\ldots,i_s)$ containing one
simple root from each $\sigma$-orbit of simple roots. We use the
product and its pinned lift
\[
 c=s_{i_1}\cdots s_{i_s},\qquad
 \dot c=\dot s_{i_1}\cdots\dot s_{i_s}.
\]

With this choice of $c$, the two flag conditions are those of
\cite[\S7.2]{QDE} and \cite{QDEII}.

\begin{definition}
\label{def:eta-sigma-oper}
For this choice of $c$, a meromorphic $(G,\eta,\sigma)$-oper consists
of $P$, a meromorphic $(\eta,\sigma)$-connection $\mathcal A$ on $P$,
and a meromorphic $B_+$-reduction $F_+$ satisfying
\[
 \operatorname{pos}\bigl(F_+(\eta u),\mathcal A_u(F_+(u))\bigr)=c
\]
on a dense open subset of $U$. A Miura reduction is a meromorphic
$B_-$-reduction $F_-$ such that $\mathcal A_u(F_-(u))=F_-(\eta u)$
wherever both sides are defined. The data $(P,\mathcal A,F_+,F_-)$
then form a Miura $(G,\eta,\sigma)$-oper.
\end{definition}

When $\sigma=1$, each orbit consists of one simple root, so $c$ is a
Coxeter element of the Weyl group of $G$. With $\eta=q$ this gives
the ordinary $(G,q)$-oper of Definition~\ref{def:ordinary-flag-oper}. More
generally, if $\sigma$ has order $r$, the composite of $r$ consecutive
connection maps is $G$-equivariant and covers $u\mapsto\eta^r u$.
It is therefore an ordinary $q$-connection for $q=\eta^r$.

In a frame adapted to $F_+$, the image of $F_+(u)$ under the connection is
$A(u)B_+$, so the relative position condition reads
$A(u)\in B_+\dot cB_+$. Where the two flags are opposite, a frame
adapted to both also gives $A(u)\in B_-$ by preservation of $F_-$.
For an orbit along which these conditions hold, write
$A(\eta^ku_*)\in N_+h_k\dot cN_+$ with $h_k\in H$, and set
$\gamma(s)=\dot c\,\sigma(s)\dot c^{-1}$ for $s\in H$. Torus changes of
frame $(s_k)$ replace $h_k$ by $s_{k+1}h_k\gamma(s_k)^{-1}$. Thus the
recurrence $s_{k+1}=\gamma(s_k)h_k^{-1}$ puts the transformed
connection maps in $M_c^-$ and determines all $s_k$ from $s_0$ in
both directions on $\mathbb Z$. For an orbit of length $L$, the
sequence $(h_k)$ is periodic. Uniqueness in the recurrence then shows
that the chosen torus frames descend to that orbit if and only if
$s_L=s_0$. On the integer lattice the discrete Miura transformation
is defined for arbitrary sequences in $M_c^-$ at independent indices.

For an affine scheme $X=\operatorname{Spec}A$, the product of a copy
of $X$ at each $k\in\mathbb Z$ has coordinate ring
$\bigotimes_{k\in\mathbb Z}A$, the direct limit of the tensor products
over finite subsets of $\mathbb Z$. Thus a $D$-point of $X^{\mathbb Z}$
is a sequence in $X(D)$, and every regular function depends on only
finitely many entries. For $X=\mathcal A_c$, the element at index $k$ is $A_k=a((T_{i,k})_{i\in J})$. For $X=M_c^-$, it is the ordered product
$A_k^{\mathrm M}=\prod_{p=1}^s P_{i_p}(\miuraq_{i_p,k})$. Applying
Proposition~\ref{prop:cell-parametrizations} at each index gives
\[
 \begin{aligned}
  \mathbb C[\mathcal A_c^{\mathbb Z}]&=\mathbb C[T_{i,k}:i\in J,k\in\mathbb Z],
   &T_{i,k}&=\Delta_i(A_k),\\
  \mathbb C[(M_c^-)^{\mathbb Z}]&=\mathbb C[\miuraq_{i,k}^{\pm1}:i\in J,k\in\mathbb Z],
   &\miuraq_{i,k}&=\omega_i([A_k^{\mathrm M}]_0).
 \end{aligned}
\]
Translation acts on a sequence $x$ by $(\tau x)_k=x_{k+1}$.
Consequently the pullbacks of the coordinate functions are
$\tau^*T_{i,k}=T_{i,k+1}$ and $\tau^*\miuraq_{i,k}=\miuraq_{i,k+1}$.
On the multiplicative orbit $u_*\eta^{\mathbb Z}$, a function $f$
gives the sequence $x_k=f(\eta^k u_*)$. Its translate is
$(\tau x)_k=f(\eta^{k+1}u_*)$, with $\eta=q$ in the ordinary case.

\section{Quantum loop algebras and screenings}
\label{sec:quantum-preliminaries}

\sectionoverview{In this section, we state the representation-theoretic criteria used to
recognize the Steinberg coordinates as fundamental $q$-characters.
The screening operators characterize the image of the character map,
and the coefficients of the dominant monomials determine each element of that image.
We also determine when these criteria restrict to a prescribed set of
spectral variables.}

\subsection{Fundamental modules and their characters}

\sectionoverview{In this subsection, we define the $q$-character from the joint
generalized eigenspaces of the Cartan currents. We explain how
Drinfeld polynomials label the simple modules and recall why the
classes of the fundamental modules generate the Grothendieck ring.}

The simple roots $\alpha_i^{\rm q}$ and fundamental weights $\varpi_i$ of a simple complex Lie algebra $\mathfrak a$ are indexed by $K$. The Cartan matrix entry $C_{ij}$ is the pairing of $\alpha_j^{\rm q}$ with the
$i$th simple coroot. These are the roots denoted by $\beta_i$ in
\cite{QDEII}. The relatively prime positive integers
$d_i^{\mathrm{sym}}$ are determined by the symmetry of
$B_{ij}=d_i^{\mathrm{sym}}C_{ij}$. With the shortest roots of
$\mathfrak a$ normalized to have squared length two, these integers
are $d_i^{\mathrm{sym}}=(\alpha_i^{\rm q},\alpha_i^{\rm q})/2$.
We fix a multiplicative
deformation parameter $\qaff\in\mathbb C^\times$ of infinite order.
The simple root at $i$ determines the scalar
$\qaff_i=\qaff^{d_i^{\mathrm{sym}}}$.

The notation $L\mathfrak a=\mathfrak a\otimes_{\mathbb C}\mathbb C[w,w^{-1}]$
denotes the Lie algebra of Laurent polynomial loops in $\mathfrak a$.
The quantum loop algebra $U_{\qaff}(L\mathfrak a)$ is the quotient of
the untwisted quantum affine algebra, with the degree generator omitted,
by the relation $\mathbf c^{1/2}=1$ on its central generator. We use the
Drinfeld--Jimbo coproduct for tensor products and write
$U_{\qaff}(\widehat{\mathfrak a})$ when referring to this algebra
at level zero. In Drinfeld's current presentation the generators are
$k_i^{\pm1}$, $h_{i,m}$ for $m\in\mathbb Z\setminus\{0\}$, and
$x_{i,r}^{\pm}$ for $r\in\mathbb Z$ \cite[Theorem 1.2]{FM}.
The Cartan currents in the formal variable $u$ are
\begin{equation}\label{eq:cartan-currents}
 \phi_i^{\pm}(u)=k_i^{\pm1}
 \exp\left(\pm(\qaff-\qaff^{-1})\sum_{m>0}h_{i,\pm m}u^{\pm m}\right).
\end{equation}
The positive current is a series in $u$ with constant term $k_i$,
and the negative current is a series in $u^{-1}$ with constant term
$k_i^{-1}$. At level zero the elements $k_i^{\pm1}$ and $h_{j,m}$
all commute. The coefficients of the currents therefore admit
simultaneous generalized eigenspaces in every finite-dimensional module.

The category $\mathcal C$ consists of finite-dimensional $U_{\qaff}(L\mathfrak a)$-modules
of type one. This means that each $k_i$ acts diagonally with eigenvalues in
$\qaff_i^{\mathbb Z}$.
The Grothendieck ring $K_0(\mathcal C)$ is generated additively by classes
$[V]$, with $[V]=[V']+[V'']$ for every short exact sequence
$0\to V'\to V\to V''\to0$. Multiplication is induced by tensor product.
By the Jordan--H\"older theorem, each class is the sum of the simple
classes with their composition multiplicities. The simple classes
therefore form an additive basis of $K_0(\mathcal C)$.

An $\ell$-weight is a simultaneous generalized eigenvalue of the
coefficients of \eqref{eq:cartan-currents}. To describe these
eigenvalues we use the Laurent rings
\[
 \mathcal Y_{\mathbb Z}=\mathbb Z[Y_{i,b}^{\pm1}:i\in K,b\in\mathbb C^\times],
 \qquad \mathcal Y=\mathbb C\otimes_{\mathbb Z}\mathcal Y_{\mathbb Z}.
\]
A Laurent monomial in these independent variables,
$M=\prod_{i,b}Y_{i,b}^{n_{i,b}}$, has only finitely many nonzero
integer exponents. The set of pairs $(i,b)$ for which $n_{i,b}\ne0$
is its support, denoted by $\operatorname{supp}(M)$. At vertex $i$
the monomial determines the rational function
\begin{equation}\label{eq:ellweight-monomial}
 \Psi_i(u)=\qaff_i^{\sum_b n_{i,b}}
 \prod_b\left(\frac{1-b\qaff_i^{-1}u}{1-b\qaff_i u}\right)^{n_{i,b}}.
\end{equation}
Every $\ell$-weight of a module in $\mathcal C$ arises in this way. Its eigenvalues for $\phi_i^+$ and $\phi_i^-$ are the
expansions of $\Psi_i$ at zero and infinity, respectively
\cite[Propositions 2.3 and 2.4]{FM}. To see that they determine $M$,
it suffices by multiplicativity of \eqref{eq:ellweight-monomial}
to consider the case $\Psi_i=1$ for all $i$. The order of vanishing
at $u=b^{-1}$ is then zero. In the product for $\Psi_i$, the numerator
and denominator at this point give
\[
 0=\operatorname{ord}_{u=b^{-1}}\Psi_i
   =n_{i,b\qaff_i}-n_{i,b\qaff_i^{-1}}
 \qquad(i\in K,\ b\in\mathbb C^\times).
\]
Thus $n_{i,b}=n_{i,b\qaff_i^{2n}}$ for every $n\in\mathbb Z$.
These spectral parameters are distinct, and only finitely many
exponents are nonzero. Hence $n_{i,b}=0$ for every $(i,b)$ and $M=1$.

For $V\in\mathcal C$, we denote by $V_M$ the common generalized eigenspace in $V$ for the two expansions of \eqref{eq:ellweight-monomial}
at every vertex. The commuting current coefficients generate a
finite-dimensional commutative subalgebra of $\operatorname{End}(V)$,
so $V=\bigoplus_M V_M$, with finitely many nonzero summands.

\begin{definition}
\label{def:q-character}
The $q$-character of a finite-dimensional type-one
$U_{\qaff}(L\mathfrak a)$-module $V$ is the Laurent polynomial
\begin{equation}\label{eq:qcharacter-definition}
 \chi_{\qaff}(V)=\sum_M\dim(V_M)M .
\end{equation}
Here $M$ runs over the $\ell$-weights of $V$ and $V_M$ is its joint
generalized eigenspace of $\ell$-weight $M$.
\end{definition}

In this sum $\qaff$ enters through the eigenvalues encoded by $M$.
For example, the coefficient of $Y_{i,b}$ is the dimension of the
joint generalized eigenspace with
$\Psi_i(u)=\qaff_i(1-b\qaff_i^{-1}u)/(1-b\qaff_i u)$ and
$\Psi_j(u)=1$ for $j\ne i$.

The ordinary weight can be recovered by retaining just the constant
coefficient of each positive current. From \eqref{eq:ellweight-monomial},
$k_i$ acts on $V_M$ by $\qaff_i^{\sum_b n_{i,b}}$. Thus the weight
is $\operatorname{wt}(M)=\sum_{i,b}n_{i,b}\varpi_i$. We order
ordinary weights by $\mu<\lambda$ when $\lambda-\mu$ is a nonzero
nonnegative integral combination of the $\alpha_i^{\rm q}$.

\begin{definition}
\label{def:dominant-monomial}
The monomial
$M=\prod_{i,b}Y_{i,b}^{n_{i,b}}$ is \emph{dominant} if every
$n_{i,b}$, the power of $Y_{i,b}$ after cancellation, is
nonnegative \cite[\S2.2, after (2.13)]{FM}.
\end{definition}

For example,
$Y_{i,a}/Y_{i,b}$ with $a\ne b$ has ordinary weight zero, a
dominant weight, although its exponent at $(i,b)$ is $-1$.

\begin{definition}
\label{def:drinfeld-polynomials}
A highest $\ell$-weight vector is a nonzero simultaneous eigenvector of
the Cartan currents annihilated by all $x_{i,r}^+$. For a simple module generated by such a vector, the monomial
encoding its $\ell$-weight is called the \emph{highest monomial} of
the module. If this monomial is $M=\prod_{i,b}Y_{i,b}^{n_{i,b}}$ with
$n_{i,b}\geq 0$, its \emph{Drinfeld polynomials} are
\[
 \pi_i(u)=\prod_{b\in\mathbb C^\times}(1-bu)^{n_{i,b}},
 \qquad i\in K.
\]
Each polynomial has constant term one, and the multiplicity of
$b^{-1}$ as its root recovers $n_{i,b}$.
\end{definition}

Substitution in
\eqref{eq:ellweight-monomial} expresses the current eigenvalues as
\[
 \Psi_i(u)=\qaff_i^{\deg\pi_i}
           \frac{\pi_i(u\qaff_i^{-1})}{\pi_i(u\qaff_i)}
 \qquad(i\in K)
\]
in the normalization of \cite[Theorem~1.3]{FM}. In particular, for
$M=Y_{i,b}$ the only nonconstant polynomial is $\pi_i(u)=1-bu$.
The simple module with this highest $\ell$-weight is called the
fundamental module at vertex $i$ and parameter $b$. The following
theorem combines the classification by Chari and Pressley
\cite{CPguide,CPaffine} with the character theorem of Frenkel and
Reshetikhin \cite[Theorem~3 and Corollary~2]{FR}.

\begin{theorem}
\label{thm:qcharacter-facts}
The simple objects of $\mathcal C$ are indexed by dominant monomials,
or equivalently by their tuples of Drinfeld polynomials.
The simple module $L(M)$ is generated by a highest $\ell$-weight vector
of monomial $M$ and its highest $\ell$-weight space is one-dimensional.
The $q$-character defines an injective ring homomorphism
$\chi_{\qaff}:K_0(\mathcal C)\to\mathcal Y_{\mathbb Z}$.
The classes of the fundamental modules $L(Y_{i,b})$ freely generate
$K_0(\mathcal C)$ as a polynomial ring over $\mathbb Z$.
\end{theorem}

The dependence of these generators on $b$ comes from rescaling the
loop parameter. For $z\in\mathbb C^\times$, the spectral automorphism
$\tau_z$ multiplies a generator of loop degree $r$ by $z^r$. The degrees of $k_i$,
$h_{i,m}$, and $x_{i,r}^\pm$ are $0$, $m$, and $r$, respectively \cite[(1.4)]{FM}. On a module pulled back by
$\tau_z$, the current $\phi_i^\pm(u)$ acts as $\phi_i^\pm(zu)$.
Replacing $u$ by $zu$ in \eqref{eq:ellweight-monomial} therefore gives
\begin{equation}\label{eq:spectral-translation}
 \chi_{\qaff}(\tau_z^* V)=
 \left.\chi_{\qaff}(V)\right|_{Y_{i,b}\mapsto Y_{i,zb}},
 \qquad \tau_z^*L(Y_{i,b})\simeq L(Y_{i,zb}).
\end{equation}
The second identity follows because pullback preserves simplicity
and sends the highest monomial to $Y_{i,zb}$.

\subsection{Screening operators}

\sectionoverview{In this subsection, we define the screening operators and use the
dominant-monomial criterion of Frenkel and Mukhin to recognize
an individual Laurent polynomial in their common kernel. We also describe when a restriction
on the spectral parameters gives a ring generated by fundamental $q$-characters.}

The screening operator at vertex $i$ differentiates the variables
$Y_{i,b}$ with values in a module that relates the derivatives at
$b$ and $b\qaff_i^2$. The relation is expressed in terms of the root
monomial associated with column $i$ of the Cartan matrix. At spectral
parameter $b$, this monomial is denoted by $A_{i,b}$ in
\cite[(2.15)]{FM}. In our notation it is
\begin{align}\label{eq:quantum-root}
 \mathsf R_i(b)={}&Y_{i,b\qaff_i}Y_{i,b\qaff_i^{-1}}\prod_{j:C_{ji}=-1}Y_{j,b}^{-1}
 \notag\\
 &{}\cdot\prod_{j:C_{ji}=-2}(Y_{j,b\qaff}Y_{j,b\qaff^{-1}})^{-1}
 \prod_{j:C_{ji}=-3}(Y_{j,b\qaff^2}Y_{j,b}Y_{j,b\qaff^{-2}})^{-1}.
\end{align}
Its total exponent at vertex $j$ is $C_{ji}$, so its weight is
$\alpha_i^{\rm q}$.

\begin{definition}
\label{def:quantum-screening}
The target $\mathcal S_i$ of the $i$th screening
is the $\mathcal Y$-module generated by $s_{i,b}$ with relations
\begin{equation}\label{eq:quantum-screening-module}
 s_{i,b\qaff_i^2}=\mathsf R_i(b\qaff_i)s_{i,b}.
\end{equation}
The action of $S_i$ on
the Laurent polynomial ring is specified by the Leibniz rule and
\begin{equation}\label{eq:quantum-screening}
 S_i:\mathcal Y\longrightarrow\mathcal S_i,\qquad
 S_i(Y_{j,b})=\delta_{ij}Y_{i,b}s_{i,b}.
\end{equation}
The same generators and relations over $\mathcal Y_{\mathbb Z}$ give
the integral screening $S_i^{\mathbb Z}$.
\end{definition}

For a monomial the action
is
\[
 S_i(M)=M\sum_b n_{i,b}s_{i,b}
 \qquad\left(M=\prod_{j,b}Y_{j,b}^{n_{j,b}}\right).
\]

Starting from a parameter $b_0$, the recurrence
\eqref{eq:quantum-screening-module} relates the symbols at
$b_0\qaff_i^{2n}$ for $n\in\mathbb Z$. Every multiplier is
invertible, so the recurrence determines the symbols in both
directions from $s_{i,b_0}$. Explicitly, for $n>0$ we have
\[
 \begin{split}
 s_{i,b_0\qaff_i^{2n}}&=
  \left(\prod_{a=0}^{n-1}\mathsf R_i(b_0\qaff_i^{\,2a+1})\right)s_{i,b_0},\\
 s_{i,b_0\qaff_i^{-2n}}&=
  \left(\prod_{a=1}^{n}\mathsf R_i(b_0\qaff_i^{-2a+1})^{-1}\right)s_{i,b_0}.
 \end{split}
\]
These formulas express every generator in terms of $s_{i,b_0}$ and
satisfy all the relations on its orbit. Distinct orbits have no
relations in common, so the symbols at one representative per orbit
form a $\mathcal Y$-basis of $\mathcal S_i$.

In $\mathsf R_i(b\qaff_i)$ the factors with vertex $i$ are $Y_{i,b}$ and
$Y_{i,b\qaff_i^2}$. It follows from \eqref{eq:quantum-screening-module}
that
\[
 S_i\bigl(Y_{i,b}+Y_{i,b}\mathsf R_i(b\qaff_i)^{-1}\bigr)
 =Y_{i,b}s_{i,b}-Y_{i,b}\mathsf R_i(b\qaff_i)^{-1}s_{i,b\qaff_i^2}=0.
\]
Frenkel and Mukhin determine the full kernel of each $S_i$ and their
intersection. They also describe the possible monomials of simple
and fundamental characters
\cite[Theorems~4.1 and~5.1, Proposition~5.2, Corollary~4.5, and Lemmas~6.1 and~6.5]{FM}.

\begin{theorem}\label{thm:FM}
\begin{enumerate}[label=\textup{(\roman*)},leftmargin=*,nosep]
\item For a fixed vertex $i$, $\ker S_i^{\mathbb Z}$ is the algebra over
$\mathbb Z[Y_{j,b}^{\pm1}:j\ne i,\ b\in\mathbb C^\times]$ generated by
$Y_{i,b}(1+\mathsf R_i(b\qaff_i)^{-1})$, $b\in\mathbb C^\times$.
\item For each dominant monomial $M$, its coefficient in
$\chi_{\qaff}(L(M))$ is one. Every other monomial is $M$ times a
nonempty product of inverse root monomials $\mathsf R_j(b)^{-1}$.
\item In the fundamental character $\chi_{\qaff}(L(Y_{i,b}))$,
$Y_{i,b}$ is the unique dominant monomial, and all spectral parameters
lie in $b\qaff^{\mathbb Z}$. If $m\ne Y_{i,b}$ is a monomial of
this character and $s$ is the greatest integer for which $m$ contains
a variable with parameter $b\qaff^s$, every variable of $m$ at that
parameter has negative exponent.
\item The integral and complexified character homomorphisms satisfy
\[
 \begin{aligned}
 \chi_{\qaff}(K_0(\mathcal C))&=\bigcap_{i\in K}\ker S_i^{\mathbb Z}
   &&\text{inside }\mathcal Y_{\mathbb Z},\\
 \chi_{\qaff}(\mathbb C\otimes_{\mathbb Z}K_0(\mathcal C))&=
 \bigcap_{i\in K}\ker S_i&&\text{inside }\mathcal Y.
 \end{aligned}
\]
\end{enumerate}
\end{theorem}

Tensoring with $\mathbb C$ preserves exactness, so it gives the same
generators for each complex kernel. Applied to the kernel sequence
of $F\mapsto(S_i^{\mathbb Z}F)_{i\in K}$ into the direct sum of the
screening modules, it gives the assertion for their intersection.
The monomial
description also orders the terms of each simple character by ordinary
weight. Each factor $\mathsf R_j(b)^{-1}$ subtracts the simple root
$\alpha_j^{\rm q}$, so every term other than the highest monomial has
strictly smaller weight. The following lemma gives the corresponding
uniqueness criterion \cite[Lemma~6.17]{FM}.

\begin{lemma}
\label{lem:screened-uniqueness}
A Laurent polynomial in $\bigcap_i\ker S_i$ is determined by the
coefficients of its dominant monomials. In particular, if its only
dominant monomial is $Y_{i,b}$ with coefficient one, the polynomial
equals $\chi_{\qaff}(L(Y_{i,b}))$.
\end{lemma}

\begin{proof}
A nonzero integral Laurent polynomial annihilated by all screenings
has a dominant monomial \cite[Lemma~5.6]{FM}. To see why the same
holds over $\mathbb C$, we express a Laurent polynomial as
$F=\sum_{a=1}^r c_a F_a$, where the $c_a\in\mathbb C$ are linearly
independent over $\mathbb Q$ and the $F_a$ have integer coefficients.
Such an expression follows by expanding the finitely many coefficients
in a rational basis for their span and clearing denominators.
In the basis of screening symbols, each $S_iF_a$ has integer Laurent
coefficients. Hence $S_i F=0$ implies $S_iF_a=0$ for every $a$.
If all dominant coefficients of $F$ vanish, those of every $F_a$
vanish by the same independence. The integral result gives $F_a=0$
for all $a$. Applied to the difference of two elements of the kernel,
this proves the first assertion.
The fundamental character belongs to the kernel by
\cite[Corollary~5.5]{FM} and has $Y_{i,b}$ with coefficient one as its
only dominant term by \cite[Corollary~4.5]{FM}. This gives the second assertion.
\end{proof}

These facts suffice for the character
identification. The equality between the whole kernel and the image of
the character map is the stronger assertion of \cite[Theorem~5.1]{FM}.

For a prescribed set $\Lambda\subset K\times\mathbb C^\times$ of spectral variables, the corresponding Laurent ring is $\mathcal Y_\Lambda=\mathbb C[Y_{i,b}^{\pm1}:(i,b)\in\Lambda]$. The category $\mathcal C_\Lambda$ consists of the modules in $\mathcal C$
for which the highest monomial of every simple composition factor
involves only variables indexed by $\Lambda$. In a short exact
sequence $0\to U\to V\to W\to0$, the composition factors of $V$
are those of $U$ and $W$, counted with multiplicity. Thus $V$
belongs to $\mathcal C_\Lambda$ precisely when both $U$ and $W$ do.

\begin{lemma}\label{lem:spectral-subring}
If $\chi_{\qaff}(L(Y_{i,b}))\in\mathcal Y_\Lambda$
for every $(i,b)\in\Lambda$, then $\mathcal C_\Lambda$
is closed under tensor product and
\[
 \mathcal Y_\Lambda\cap\bigcap_{i\in K}\ker S_i
 =\chi_{\qaff}\bigl(\mathbb C\otimes_{\mathbb Z}K_0(\mathcal C_\Lambda)\bigr).
\]
The right-hand side is a polynomial algebra freely generated by the
fundamental characters indexed by $\Lambda$.
\end{lemma}

\begin{proof}
For a dominant monomial $M=\prod_{(i,b)\in\Lambda}Y_{i,b}^{n_{i,b}}$,
we form the tensor product $V$ of the corresponding fundamental modules
in any fixed order. Multiplicativity of the character gives
\[
 V=\bigotimes_{(i,b)\in\Lambda}L(Y_{i,b})^{\otimes n_{i,b}},\qquad
 \chi_{\qaff}(V)=\prod_{(i,b)\in\Lambda}
       \chi_{\qaff}(L(Y_{i,b}))^{n_{i,b}}\in\mathcal Y_\Lambda .
\]
By Theorem~\ref{thm:FM}\textup{(ii)}, this product has term $M$ with
coefficient one and all other terms have strictly lower weight.
Every constituent contributes its character with nonnegative
coefficients. The constituent containing $M$ must therefore have
highest monomial $M$ and occur once. Every other constituent $L(M')$
has $\operatorname{wt}(M')<\operatorname{wt}(M)$. Furthermore,
every monomial of each constituent occurs in the product, so
$\chi_{\qaff}(L(M))\in\mathcal Y_\Lambda$ for every dominant $M$
supported on $\Lambda$. By additivity the same support assertion holds
for every object in $\mathcal C_\Lambda$. Multiplicativity then gives,
for $V_1,V_2\in\mathcal C_\Lambda$,
\[
 \chi_{\qaff}(V_1\otimes V_2)
   =\chi_{\qaff}(V_1)\chi_{\qaff}(V_2)\in\mathcal Y_\Lambda .
\]
The highest monomial of each constituent occurs in this product,
so $V_1\otimes V_2\in\mathcal C_\Lambda$.

In $K_0(\mathcal C_\Lambda)$ the composition multiplicities give
\[
 [L(M)]=\prod_{(i,b)\in\Lambda}[L(Y_{i,b})]^{n_{i,b}}
   -\sum_{\operatorname{wt}(M')<\operatorname{wt}(M)}
        [V:L(M')]\,[L(M')],
\]
where $[V:L(M')]$ denotes the multiplicity of $L(M')$ in $V$.
To induct on this identity, we use the integer
$d(M)=\langle\operatorname{wt}(M),2\rho^\vee\rangle$, where
$2\rho^\vee$ is the sum of the positive coroots.
Each $s_i$ negates its simple coroot and permutes the other positive coroots,
so $\langle\alpha_i^{\rm q},2\rho^\vee\rangle=2$.
Hence $d(M)\in\mathbb Z_{\ge0}$ for dominant $M$, and $d(M')<d(M)$
in the sum. Induction gives generation by the fundamental classes.
Their algebraic independence follows from
Theorem~\ref{thm:qcharacter-facts}, since the simple classes of
$\mathcal C_\Lambda$ embed in the basis of $K_0(\mathcal C)$.

If $0\ne F\in\mathcal Y_\Lambda\cap\bigcap_i\ker S_i$, its dominant
monomials form a nonempty finite set by Lemma~\ref{lem:screened-uniqueness}.
For the largest height $d$ among these monomials, we eliminate
the dominant terms of height $d$ by subtracting simple characters. The remainder is
\[
 F'=F-\sum_{\substack{M\ \mathrm{dominant}\\d(M)=d}}
      [M]F\,\chi_{\qaff}(L(M)),
\]
where $[M]F$ denotes the coefficient of $M$ and the sum is over the
monomials of $F$. These simple characters belong to
$\mathcal Y_\Lambda\cap\bigcap_i\ker S_i$ by the support calculation
above and \cite[Corollary~5.5]{FM}. Theorem~\ref{thm:FM}\textup{(ii)}
gives $[M']F'=0$ for every dominant $M'$ of height at least $d$.
Thus either $F'=0$ or its maximal dominant height is smaller than $d$.
Since these heights are nonnegative integers, iteration expresses $F$ as
a linear combination of the characters of $L(M)$ with $M$ supported
on $\Lambda$. Each of these characters belongs to the intersection,
which gives the reverse inclusion.
\end{proof}

\section{Reduction to the Steinberg section in simply laced type}
\label{ade:normalization-section}

\sectionoverview{In this section, we prove that for every sequence $(g_k)$ in a Coxeter
cell of type $A$, $D$, or $E$ there is a unique sequence $(u_k)$ in $N_+$
such that $u_{k+1}g_ku_k^{-1}$ lies in the Steinberg section for all $k$.
We construct $u_k$ from finitely many consecutive terms of $(g_k)$ for one
ordering in each type, then pass to arbitrary orderings of the simple reflections.}

\subsection{Reduction by root subgroups}
\label{sec:root-subgroup-reduction}

\sectionoverview{In this subsection we show how factorization in root
subgroups determines the change of frame at $k$ from finitely many
consecutive elements of the Coxeter cell. The factorization lemma allows
a diagram automorphism to act on the source of each connection map.}

For $G$ of type $A$, $D$ or $E$, an ordering
$\mathbf i=(i_1,\ldots,i_r)$ of the simple roots determines
the Coxeter element $c$, its pinned representative and its inversion set by
\[
 c=s_{i_1}\cdots s_{i_r},\qquad \dot c=\dot s_{i_1}\cdots\dot s_{i_r},
 \qquad R_c=\Delta_+\cap c\Delta_-.
\]
The cell $M_c=N_+\dot cN_+$ contains both $M_c^-=M_c\cap B_-$ and
the Steinberg section $\mathcal A_c=U_c\dot c$. Over a commutative
complex algebra $D$, Proposition~\ref{prop:cell-parametrizations}
parametrizes sequences in these two subvarieties by
\[
 A_k^{\mathrm M}=\prod_{p=1}^rP_{i_p}(\miuraq_{i_p,k}),\qquad
 A_k=\prod_{p=1}^rX_{i_p}(T_{i_p,k}),
\]
where $\miuraq_{i,k}\in D^\times$ and $T_{i,k}\in D$. For a sequence
$(g_k)$ in $M_c(D)$, left multiplication by $g_k$ sends the copy of
$G(D)$ at $k$ to the copy at $k+1$. A change of coordinates by
$u_k\in N_+(D)$ at each $k$ therefore transforms $g_k$ into
$u_{k+1}g_ku_k^{-1}$. For $g_k=A_k^{\mathrm M}$ we seek
$(u_k)\in N_+(D)^{\mathbb Z}$ satisfying
\begin{equation}\label{ade:gauge-equation}
 u_{k+1}A_k^{\mathrm M}u_k^{-1}=A_k\in\mathcal A_c(D)\qquad(k\in\mathbb Z).
\end{equation}

The factorization lemma also applies when the connection map from $k$ to
$k+1$ first applies a fixed diagram automorphism $\sigma$ in the source
fiber and then multiplies by $g_k$. Here $\sigma$ is the lift preserving
the pinning from Section~\ref{sec:pinned-root-subgroups}. The map in the
original frames is $h\mapsto g_k\sigma(h)$. If $h$ is a coordinate in
the new frame at $k$, its coordinate in the original frame is $u_k^{-1}h$.
Its image in the new frame at $k+1$ is
\[
 u_{k+1}g_k\sigma(u_k^{-1}h)
  =\bigl(u_{k+1}g_k\sigma(u_k)^{-1}\bigr)\sigma(h).
\]
For $\sigma=\mathrm{id}_G$ this is the ordinary change of frame.

\begin{definition}
\label{def:discrete-miura-map}
For a sequence $(A_k^{\mathrm M})$ in $M_c^-(D)$ admitting a
unique sequence $(u_k)$ in $N_+(D)$ such that
\[
 A_k=u_{k+1}A_k^{\mathrm M}\sigma(u_k)^{-1}\in\mathcal A_c(D)
 \qquad(k\in\mathbb Z),
\]
the discrete Miura transformation $\mu_{c,\sigma}$ sends
$(A_k^{\mathrm M})$ to $(A_k)$. For $\sigma=1$ we write $\mu_c$.
\end{definition}

The condition for the transformed group element to lie in the Steinberg
section is $u_{k+1}g_k=A_k\sigma(u_k)$ with $A_k\in\mathcal A_c(D)$.
For a prescribed $u_{k+1}$, the factorization in $\mathcal A_cN_+$
determines both factors on the right. Writing $A_k=a_k\dot c$ with
$a_k\in U_c(D)$ gives $A_k\sigma(g)=a_k\Gamma(g)\dot c$, where
$\Gamma(g)=\dot c\,\sigma(g)\dot c^{-1}$. On the roots this is the action
of $\gamma=c\sigma$, and $R_c=\Delta_+\cap\gamma\Delta_-$ because
$\sigma$ preserves $\Delta_-$. The following lemma gives a condition on
the positive roots under which iterating the factorization from right
to left makes $u_k$ independent of the choice at the right endpoint.

\begin{lemma}\label{ade:contraction}
Let $c$ be a product of distinct simple reflections, let $\sigma$ be an
automorphism of $G$ preserving its pinning, and write $\gamma=c\sigma$. For the
ordinary connection $\sigma=\mathrm{id}_G$ and $\gamma=c$.
For some $m\ge1$, let $S_0,\ldots,S_m\subset\Delta_+$ be closed
sets of roots, with $S_0=\Delta_+$ and $S_m=\varnothing$. Suppose
that for $0\le j<m$ one has
\begin{equation}\label{ade:filtration-hypothesis}
 E_j:=R_c\cup S_j=R_c\sqcup\gamma S_{j+1},
 \qquad E_j\text{ is closed}.
\end{equation}
For each $g=(g_k)\in M_c(D)^{\mathbb Z}$ there is a unique
$u=(u_k)\in N_+(D)^{\mathbb Z}$ satisfying, for every $k$,
\[
 A_k=u_{k+1}g_k\sigma(u_k)^{-1}\in\mathcal A_c(D).
\]
In the affine coordinates of Proposition~\ref{prop:cell-parametrizations},
$u_k$ is polynomial in $g_k,\ldots,g_{k+m-1}$, and $A_k$ is
polynomial in $g_k,\ldots,g_{k+m}$. The map
\[
 (u,A)\longmapsto\bigl(u_{k+1}^{-1}A_k\sigma(u_k)\bigr)_k
\]
identifies $N_+^{\mathbb Z}\times\mathcal A_c^{\mathbb Z}$ with $M_c^{\mathbb Z}$.
\end{lemma}

\begin{proof}
We first factor $u_{k+1}g_k$ for an arbitrary choice of $u_{k+1}$.
The subgroup $V=N_+\cap\Gamma^{-1}(N_+)$ consists of the positive
unipotent elements that remain in $N_+$ under
$\Gamma(g)=\dot c\,\sigma(g)\dot c^{-1}$. The equality
$R_c=\Delta_+\cap\gamma\Delta_-$ and
\eqref{ade:filtration-hypothesis} for $j=0$ give
\[
 \gamma S_1=\Delta_+\setminus R_c=\Delta_+\cap\gamma\Delta_+,
 \qquad S_1=\gamma^{-1}\Delta_+\cap\Delta_+,\qquad V=U_{S_1}.
\]
For every $j<m$, $\gamma S_{j+1}\subseteq\Delta_+\cap\gamma\Delta_+$,
so $U_{S_{j+1}}\subseteq V$. The disjoint union in
\eqref{ade:filtration-hypothesis} and Lemma~\ref{lem:root-factorization}
give multiplication isomorphisms
\begin{equation}\label{ade:step-factorizations}
 U_c\times\Gamma(V)\simeq N_+,\qquad
 U_c\times\Gamma(U_{S_{j+1}})\simeq U_{E_j}.
\end{equation}
For $p\in N_+(D)$, Proposition~\ref{prop:cell-parametrizations} and the
first isomorphism of \eqref{ade:step-factorizations} give the unique factors
\[
 g_k=a_k\dot c\,r_k,\qquad pa_k=a'\Gamma(v),\qquad
 a_k,a'\in U_c(D),\quad r_k\in N_+(D),\quad v\in V(D).
\]
Since $\Gamma(v)\dot c=\dot c\,\sigma(v)$,
\[
 pg_k=(a'\dot c)\,\sigma\bigl(\mathcal R_k(p)\bigr),\qquad
 \mathcal R_k(p)=v\sigma^{-1}(r_k).
\]
For $u_{k+1}=p$ this determines $A_k=a'\dot c$ and $u_k=\mathcal R_k(p)$,
both polynomial in the affine coordinates of $(p,g_k)$.

Suppose $u_{k+1}$ and $\hat u_{k+1}$ give the respective factors
$a'_k\dot c$ and $a''_k\dot c$ in the section. Comparing the two
factorizations of $g_k$ and writing $d_s=\hat u_s u_s^{-1}$ gives
\begin{equation}\label{ade:relative-gauge}
 d_{k+1} a'_k=a''_k\Gamma(d_k).
\end{equation}
In particular $\Gamma(d_k)\in N_+(D)$, so $d_k\in V(D)$.
The second factorization in \eqref{ade:step-factorizations} restricts the
first to $U_{E_j}$. By uniqueness of its factors,
\[
 d_{k+1}\in U_{S_j}(D)
 \ \Longrightarrow\ d_{k+1}a'_k\in U_{E_j}(D)
 \ \Longrightarrow\ d_k\in U_{S_{j+1}}(D).
\]
For any two choices of the unipotent element at $k+m$, their quotient
$d_{k+m}$ lies in $U_{S_0}(D)=N_+(D)$. Iteration gives
$d_{k+m-j}\in U_{S_j}(D)$ for $0\le j\le m$. Since $S_m=\varnothing$,
we have $d_k=1$. Thus for fixed $g_k,\ldots,g_{k+m-1}$ the composite
$\mathcal R_k\circ\cdots\circ\mathcal R_{k+m-1}$ takes the same
value $u_k$ at every choice of terminal element. In particular, we can compute $u_{k+1}$ with terminal value $1$
at $k+m+1$. It then satisfies
\[
 \mathcal R_k(u_{k+1})
 =(\mathcal R_k\circ\cdots\circ\mathcal R_{k+m-1})(\mathcal R_{k+m}(1))
 =u_k .
\]
These equalities give a solution at every index. For any other solution
$\hat u$, the recurrence gives
$\hat u_k=(\mathcal R_k\circ\cdots\circ\mathcal R_{k+m-1})(\hat u_{k+m})=u_k$, which proves uniqueness.
The expressions for $u_k$ and $A_k=u_{k+1}g_k\sigma(u_k)^{-1}$ are
polynomial in $g_k,\ldots,g_{k+m-1}$ and $g_k,\ldots,g_{k+m}$,
respectively. They define the inverse morphism to
$(u,A)\mapsto(u_{k+1}^{-1}A_k\sigma(u_k))_k$.
\end{proof}

The dependence on $g_{k+m-1}$ in Lemma~\ref{ade:contraction} is
necessary whenever there is a root $\beta$ for which
$\gamma^a\beta\in\Delta_+$ for $0\le a<m$. Indeed, a nonidentity element
$x\in U_\beta(\mathbb C)$ determines a sequence with
$g_s=\dot c$ for $s\ne k+m-1$ and
\[
 g_{k+m-1}=\dot c\,\sigma\bigl(\Gamma^{m-1}(x)\bigr).
\]
Since $g_s=\dot c$ for $s\ge k+m$, we have $u_{k+m}=1$. The factorization
at $k+m-1$ gives $u_{k+m-1}=\Gamma^{m-1}(x)$. For $1\le a\le m-1$,
\[
 \Gamma^a(x)\dot c=\dot c\,\sigma(\Gamma^{a-1}(x)),\qquad
 \Gamma^{a-1}(x),\Gamma^a(x)\in N_+(\mathbb C).
\]
Successive factorizations thus give $u_{k+a}=\Gamma^a(x)$ for
$0\le a<m$. In particular $u_k=x$, while the constant sequence gives $u_k=1$.

\subsection{Companion matrices for the linear ordering in type \texorpdfstring{$A$}{A}}

\sectionoverview{In this subsection, we construct the gauge transformation to companion
matrices for the linear ordering in $\mathrm{SL}_N$. We show that
$g_k,\ldots,g_{k+N-2}$ determine $u_k$. For sequences of lower triangular
matrices, we also give Laurent polynomial formulas for the scalar recurrence.}

We now apply Lemma~\ref{ade:contraction} with $\sigma=\mathrm{id}_G$.
For $G=\mathrm{SL}_N$, $N\ge2$, the ordering $s_1,\ldots,s_{N-1}$
gives the Coxeter element $c=s_1\cdots s_{N-1}$.
In the standard representation, the elements of the Steinberg section
are the companion matrices
\[
 A_k=
 \begin{pmatrix}
 T_{1,k} & -T_{2,k} & \cdots & (-1)^{N-2}T_{N-1,k} & (-1)^{N-1} \\
 1 & 0 & \cdots & 0 & 0 \\
 0 & 1 & \ddots & \vdots & \vdots \\
 \vdots & \ddots & \ddots & 0 & 0 \\
 0 & \cdots & 0 & 1 & 0
 \end{pmatrix}.
\]
Indeed, in the first row of $X_1(T_1)\cdots X_{N-1}(T_{N-1})$,
right multiplication by $X_j(T_j)$ sends the pair $((-1)^{j-1},0)$
in columns $j,j+1$ to $((-1)^{j-1}T_j,(-1)^j)$. For $j\ge2$,
the $j$th row is the product of the standard row $\mathbf e_j^T$ with
the same factors. Those preceding $X_{j-1}$ fix this row,
$X_{j-1}$ sends it to $\mathbf e_{j-1}^T$, and the subsequent factors
fix $\mathbf e_{j-1}^T$. Since $M_c(D)=\mathcal A_c(D)N_+(D)$,
right multiplication by an upper unitriangular matrix shows that every
matrix in the cell has ones on the subdiagonal and zeros below it.

For $v_k=(f_k,f_{k-1},\ldots,f_{k-N+1})^T$, the lower rows of
$v_{k+1}=A_kv_k$ copy the $N-1$ entries common to both columns.
The first row is the scalar recurrence
\[
 f_{k+1}=\sum_{j=1}^{N-1}(-1)^{j-1} T_{j,k} f_{k+1-j}
       +(-1)^{N-1}f_{k-N+1}.
\]

The action of $c$ on roots determines how many terms of $(g_k)$ are
needed for the change to companion form. It permutes the coordinate weights of the diagonal torus
by the cycle $(1\ 2\ \cdots\ N)$. Writing these weights as
$\epsilon_a$, the positive root $\epsilon_a-\epsilon_b$ for $a<b$
first becomes negative when $b$ reaches $1$ after $N-b+1$
applications of $c$. The positive roots that become negative under
$c^{-1}$ are the $\epsilon_1-\epsilon_b$. Hence
\[
 \begin{aligned}
 R_c&=\{\epsilon_1-\epsilon_b:2\le b\le N\},\\
 S_j&:=\bigcap_{d=0}^{j}c^{-d}\Delta_+
 =\{\epsilon_a-\epsilon_b:1\le a<b\le N-j\},\qquad 0\le j\le N-1.
 \end{aligned}
\]
In particular $S_{N-1}=\varnothing$. For $j<N-1$, $cS_{j+1}$ consists of the roots
$\epsilon_a-\epsilon_b\in S_j$ with $a\ge2$.
Thus $R_c\cup S_j=R_c\sqcup cS_{j+1}$. Apart from the order
of the summands, a sum of positive roots is a root precisely when
it has the form $(\epsilon_a-\epsilon_b)+(\epsilon_b-\epsilon_d)=\epsilon_a-\epsilon_d$, with $a<b<d$.
For two such summands in $R_c\cup S_j$, either $a=1$ and
$\epsilon_a-\epsilon_d\in R_c$, or $a>1$ and both summands lie in $S_j$, hence
$d\le N-j$ and $\epsilon_a-\epsilon_d\in S_j$. Thus $R_c\cup S_j$ is closed.
Lemma~\ref{ade:contraction} with $m=N-1$ therefore gives $u_k$ as a
polynomial in $g_k,\ldots,g_{k+N-2}$. This bound is sharp on $M_c$, since $c^a(\epsilon_1-\epsilon_2)=\epsilon_{1+a}-\epsilon_{2+a}$ is positive for
$0\le a\le N-2$.

With $\miuraq_{0,k}=\miuraq_{N,k}=1$ and $a_{j,k}=\miuraq_{j,k}/\miuraq_{j-1,k}$, the product
of the factors $P_j(\miuraq_{j,k})$ on the lower triangular locus is
\begin{equation}\label{ade:bidiagonal}
 A_k^{\mathrm M}=\begin{pmatrix}
 a_{1,k}&0&\cdots&0\\
 1&a_{2,k}&\ddots&\vdots\\
 0&\ddots&\ddots&0\\
 0&\cdots&1&a_{N,k}
 \end{pmatrix}.
\end{equation}
\begin{proposition}\label{ade:companion}
For $(g_k)\in M_c(D)^{\mathbb Z}$, the matrices $u_k\in N_+(D)$ for which
$u_{k+1}g_k u_k^{-1}$ is in companion form have rows
\begin{equation}\label{ade:companion-gauge}
 \mathbf e_j^Tu_k=\mathbf e_N^Tg_{k+N-j-1}\cdots g_k\qquad(1\le j\le N),
\end{equation}
with the empty product for $j=N$. For $g_k=A_k^{\mathrm M}$ as in
\eqref{ade:bidiagonal}, let $A_k=u_{k+1}A_k^{\mathrm M}u_k^{-1}$.
Its $i$th Steinberg coordinate $T_{i,k}=\Delta_i(A_k)$ is
\begin{equation}\label{ade:companion-coefficients}
 T_{i,k}=
 \sum_{1\le j_1<\cdots<j_i\le N}
 \prod_{p=1}^{i}
       \frac{\miuraq_{j_p,k+j_p-p}}{\miuraq_{j_p-1,k+j_p-p}}
       \qquad(1\le i<N).
\end{equation}
As a Laurent polynomial in the independent variables $\miuraq_{j,s}$, this
coordinate depends only on indices $k\le s\le k+N-i$, and both
endpoints occur.
\end{proposition}

\begin{proof}
Every $g\in M_c(D)$ is a companion matrix times an upper unitriangular
matrix, so $g_{ab}=0$ for $a>b+1$ and $g_{b+1,b}=1$.
If a row $v$ satisfies $v_1=\cdots=v_{s-1}=0$, $v_s=1$, with $s>1$,
then
\[
 (vg)_b=\sum_{a=s}^{N}v_a g_{ab}=\delta_{b,s-1}
 \qquad(1\le b\le s-1).
\]
Each multiplication therefore moves the first nonzero entry of the
row one place to the left and leaves it equal to $1$. Starting with
$\mathbf e_N^T$ and multiplying $N-j$ times gives a row whose first
nonzero entry is $1$ in column $j$. Thus the rows in
\eqref{ade:companion-gauge} form an upper unitriangular matrix $u_k$
and satisfy
\[
 \mathbf e_j^Tu_{k+1}g_k=\mathbf e_N^Tg_{k+N-j}\cdots g_k=\mathbf e_{j-1}^Tu_k
 \qquad(2\le j\le N).
\]
Thus $A_k=u_{k+1}g_k u_k^{-1}$ has rows $\mathbf e_{j-1}^T$ for $2\le j\le N$
and $(A_k)_{1N}=(-1)^{N-1}\det A_k=(-1)^{N-1}$, as required. Conversely, the same row equations and
$\mathbf e_N^Tu_k=\mathbf e_N^T$ force \eqref{ade:companion-gauge}.

It remains to compute the first row when $g_k=A_k^{\mathrm M}$.
For a solution of $y_{k+1}=A_k^{\mathrm M}y_k$, we write $y_{j,k}$
for the $j$th component of $y_k$. On scalar sequences $f$ the shift
$D_+$ and multiplication by $a_j$ act by $(D_+f)_k=f_{k+1}$ and
$(a_j f)_k=a_{j,k}f_k$. The rows of
\eqref{ade:bidiagonal} give
\[
 (D_+-a_j)y_j=y_{j-1}\quad(2\le j\le N),\qquad
 (D_+-a_1)y_1=0,
\]
where $y_j=(y_{j,k})_{k\in\mathbb Z}$. Consequently $y_N$ is annihilated
by $(D_+-a_1)\cdots(D_+-a_N)$. For a multiplication operator $b$,
$(D_+^s(bf))_k=b_{k+s}f_{k+s}$. If the multiplication factors occupy
positions $j_1<\cdots<j_i$, then $j_p-p$ shifts precede $a_{j_p}$. Hence
\[
 \left[\left(\prod_{j=1}^N(D_+-a_j)\right)f\right]_k
 =\sum_{i=0}^N(-1)^i
   \left(\sum_{1\le j_1<\cdots<j_i\le N}
         \prod_{p=1}^i a_{j_p,k+j_p-p}\right)f_{k+N-i}.
\]
To identify these coefficients with the Steinberg coordinates, we use
the same scalar solution $y_N$ in the companion form.
Equation~\eqref{ade:companion-gauge} gives
$(u_ky_k)_j=y_{N,k+N-j}$. The first row of
$u_{k+1}y_{k+1}=A_ku_ky_k$ is therefore
\[
 y_{N,k+N}+\sum_{i=1}^{N-1}(-1)^i T_{i,k}y_{N,k+N-i}
             +(-1)^N y_{N,k}=0 .
\]
The operator expansion for $f=y_N$ has the same coefficients at
$y_{N,k+N}$ and $y_{N,k}$, namely $1$ and
$(-1)^N\prod_j a_{j,k}=(-1)^N$. All $A_s^{\mathrm M}$ are invertible,
so $y_k$ can be prescribed arbitrarily in $D^N$. Consequently
$u_ky_k=(y_{N,k+N-1},\ldots,y_{N,k})^T$ also ranges over $D^N$.
The coefficients at its entries must agree in the two recurrences,
which proves
\eqref{ade:companion-coefficients}. The inequalities
$p\le j_p\le N-i+p$ bound its indices by
$k\le k+j_p-p\le k+N-i$. The subsets $\{1,\ldots,i\}$ and $\{N-i+1,\ldots,N\}$ contribute
$\miuraq_{i,k}$ and $\miuraq_{N-i,k+N-i}^{-1}$, respectively. Every summand has coefficient one, so both endpoints occur.
\end{proof}

For $N=3$ the subsets in \eqref{ade:companion-coefficients}
give three terms in each Steinberg coordinate. In this case the
companion matrix and its scalar recurrence are
\[
 A_k=\begin{pmatrix}
       T_{1,k}&-T_{2,k}&1\\
       1&0&0\\
       0&1&0
     \end{pmatrix},\qquad
 f_{k+1}=T_{1,k}f_k-T_{2,k}f_{k-1}+f_{k-2},
\]
where
\[
 T_{1,k}=\miuraq_{1,k}+\frac{\miuraq_{2,k+1}}{\miuraq_{1,k+1}}
                         +\frac1{\miuraq_{2,k+2}},\qquad
 T_{2,k}=\miuraq_{2,k}+\frac{\miuraq_{1,k}}{\miuraq_{2,k+1}}
                         +\frac1{\miuraq_{1,k+1}}.
\]
In each polynomial only the first term has nonnegative exponents in
every variable. Thus $T_{1,k}$ and $T_{2,k}$ have the unique dominant
monomials $\miuraq_{1,k}$ and $\miuraq_{2,k}$, respectively.

\subsection{Bipartite Coxeter elements in types \texorpdfstring{$D$ and $E$}{D and E}}

\sectionoverview{In this subsection, we prove that for a bipartite Coxeter element of
order $h$ in type $D$ or $E$, the $h/2$ terms of the sequence starting with $g_k$
determine $u_k$. We identify the sets of positive roots that give this
bound and show that changing the last of those terms can change $u_k$.}

A bipartition $I=I_+\sqcup I_-$ places the two endpoints of each
edge in opposite parts. The reflections with labels in the same part
therefore commute, so their order in the bipartite Coxeter element
\[
 c=\left(\prod_{i\in I_+}s_i\right)
    \left(\prod_{i\in I_-}s_i\right)
\]
is immaterial. The root calculation below uses Steinberg's expression
for the longest element $w_0$ in this ordering \cite[Section~2.2]{BKW}.

\begin{theorem}\label{ade:bipartite-longest}
If a bipartite Coxeter element $c$ of a finite irreducible Coxeter group
has even order $h$, then $c^{h/2}=w_0$. If $c=s_{i_1}\cdots s_{i_r}$ in the bipartite ordering, then the
product $(s_{i_1}\cdots s_{i_r})^{h/2}$ is a reduced expression for $w_0$.
\end{theorem}

In types $D$ and $E$ the Coxeter numbers are even
\cite[Table~1]{MRV2}. The values of $h$ and $m$ in $c^m=w_0$ are
\[
\begin{array}{c|c|c}
 G&h&m=h/2\\\hline
 D_n\ (n\ge4)&2n-2&n-1\\
 E_6&12&6\\
 E_7&18&9\\
 E_8&30&15
\end{array}.
\]

\begin{lemma}\label{ade:blocks}
For a bipartite Coxeter element $c$ of type $D$ or $E$, with $m=h/2$,
\[
 \Delta_+=\bigsqcup_{a=0}^{m-1}c^aR_c,\qquad w_0R_c=-R_c,
\]
and
\begin{equation}\label{ade:block-filtration}
 S_j:=\bigcap_{a=0}^{j}c^{-a}\Delta_+
   =\bigsqcup_{a=0}^{m-1-j}c^aR_c\qquad(0\le j\le m).
\end{equation}
In particular $S_m=\varnothing$ and $S_j=R_c\sqcup cS_{j+1}$ for
$j<m$. For $(g_s)_s\in M_c(D)^{\mathbb Z}$, the solution $(u_s)_s$
of $u_{s+1}g_su_s^{-1}\in\mathcal A_c(D)$ has $u_k$ polynomial in
$g_k,\ldots,g_{k+m-1}$.
\end{lemma}

\begin{proof}
The roots $\beta_p=s_{i_1}\cdots s_{i_{p-1}}\alpha_{i_p}$,
$1\le p\le r$, satisfy $R_c=\{\beta_1,\ldots,\beta_r\}$ by
Lemma~\ref{lem:inversion-roots}. Theorem~\ref{ade:bipartite-longest}
gives a reduced expression for $w_0=c^m$ consisting of $m$ copies of
$s_{i_1}\cdots s_{i_r}$. At position $ar+p$, the first $a$ copies act
on $\beta_p$, so the corresponding positive root is $c^a\beta_p$.
Since $w_0\Delta_-=\Delta_+$, Lemma~\ref{lem:inversion-roots} gives
\[
 \Delta_+=\bigsqcup_{a=0}^{m-1}c^a R_c,\qquad
 w_0 R_c=w_0(\Delta_+\cap c\Delta_-)=\Delta_-\cap c\Delta_+=-R_c .
\]
The formula for $w_0R_c$ uses $cw_0=w_0c$, which follows from $w_0=c^m$.

For $\alpha=c^a\beta$ with $0\le a<m$ and $\beta\in R_c$,
\[
 c^b\alpha\in c^{a+b}R_c\subset\Delta_+\quad(0\le b<m-a),
 \qquad c^{m-a}\alpha=w_0\beta\in\Delta_- .
\]
Thus $\alpha\in S_j$ if and only if $a\le m-1-j$, proving
\eqref{ade:block-filtration}. If $\alpha,\alpha'\in S_j$ and
$\alpha+\alpha'$ is a root, then
$c^b(\alpha+\alpha')=c^b\alpha+c^b\alpha'\in\Delta_+$ for
$0\le b\le j$. Thus $S_j$ is closed. For $j<m$, \eqref{ade:block-filtration} gives
$R_c\cup S_j=S_j=R_c\sqcup cS_{j+1}$, so
Lemma~\ref{ade:contraction} applies with this $m$.
\end{proof}

For any $\beta\in R_c$ the roots $\beta,c\beta,\ldots,c^{m-1}\beta$
are all positive, so the dependence on the $m$ entries beginning at $g_k$ is sharp
on $M_c$. Since $c^{2m}=1$ and $w_0$ sends the positive root
$c^{m-a}\beta$ to a negative root for $1\le a\le m$, the backward
iterates satisfy
\begin{equation}\label{ade:negative-past-roots}
 c^{-a}\beta=w_0c^{m-a}\beta\in\Delta_-\qquad
 (\beta\in R_c,\ 1\le a\le m).
\end{equation}
In particular \eqref{ade:negative-past-roots} holds in types
$D_{2n+1}$ and $E_6$ because the argument uses only $w_0\Delta_+=\Delta_-$,
which allows $w_0$ to permute the simple-root labels up to sign.

\subsection{Changing the Coxeter element}

\sectionoverview{In this subsection we prove that reduction to the
Steinberg section works for every Coxeter element in types $A$, $D$,
and $E$. We first determine how moving the first simple reflection to the
end changes the Miura and Steinberg coordinates. The orientations of the
Dynkin diagram then allow us to pass from this single change to an
arbitrary ordering of the simple reflections.}

\begin{lemma}\label{ade:rotation}
For $c=s_i s_{i_2}\cdots s_{i_r}$ and
$c'=s_{i_2}\cdots s_{i_r}s_i=s_i c s_i$, there is a regular
bijection between solutions of \eqref{ade:gauge-equation} for $c$
and $c'$ over every $D$.
The transformed coordinates are
\[
 \miuraq'_{i,k}=\miuraq_{i,k-1},\quad T'_{i,k}=T_{i,k-1},\qquad
 \miuraq'_{j,k}=\miuraq_{j,k},\quad T'_{j,k}=T_{j,k}\quad(j\ne i),
\]
and the corresponding elements $u'_k\in N_+(D)$ are
\begin{equation}\label{ade:rotation-gauge}
 u'_k=X_i(T_{i,k-1})^{-1}u_kP_i(\miuraq_{i,k-1}).
\end{equation}
\end{lemma}

\begin{proof}
To determine how the first factor changes under rotation, we factor
$u_{k+1}$ in the simple-root subgroup at $i$ and the subgroup $U_i$
generated by the other positive root subgroups. By Lemma~\ref{lem:root-factorization},
$u_{k+1}=x_i(s_{k+1})n_{k+1}$ with $n_{k+1}\in U_i(D)$.
For a highest-weight vector $v_i\in V(\omega_i)$ and
$E_\beta\in\mathfrak g_\beta$ with $\beta>0$, $\beta\ne\alpha_i$,
\[
 E_\beta v_i=0,\qquad E_\beta f_i v_i=[E_\beta,f_i]v_i=0,
\]
since $\beta-\alpha_i$, if a root, is positive. Thus $U_i$ fixes
$v_i,f_iv_i$. The identity $\langle\omega_i,\alpha_j^\vee\rangle=\delta_{ij}$
shows that the rank-one subgroup at $i$ acts on their span by the
standard representation, while each such subgroup at $j\ne i$ fixes
$v_i$. Hence $u_{k+1}A_k^{\mathrm M}v_i=A_ku_kv_i$ reads
\[
 (\miuraq_{i,k}+s_{k+1})v_i+f_iv_i=T_{i,k}v_i+f_iv_i,
 \qquad T_{i,k}=\miuraq_{i,k}+s_{k+1}.
\]
In the rank-one subgroup the identity
\[
 \begin{pmatrix}1&s\\0&1\end{pmatrix}
 \begin{pmatrix}\zeta&0\\1&\zeta^{-1}\end{pmatrix}
 =\begin{pmatrix}\zeta+s&-1\\1&0\end{pmatrix}
  \begin{pmatrix}1&\zeta^{-1}\\0&1\end{pmatrix}
\]
reads $x_i(s)P_i(\zeta)=X_i(\zeta+s)x_i(\zeta^{-1})$.
Since the $\alpha_i$-string through each positive root other than
$\alpha_i$ consists of positive roots other than $\alpha_i$,
the rank-one subgroup normalizes $U_i$. Hence
\[
 u'_{k+1}=X_i(T_{i,k})^{-1}u_{k+1}P_i(\miuraq_{i,k})
 =x_i(\miuraq_{i,k}^{-1})\,
   \operatorname{Ad}(P_i(\miuraq_{i,k})^{-1})(n_{k+1})\in N_+(D).
\]
This expression is regular because $\miuraq_{i,k}$ is a unit.
In the rotated ordering the $i$th factor occurs last, with its
parameter from index $k-1$. Denoting the reordered Miura and Steinberg
products by $\widehat A_k^{\mathrm M}$ and $A'_k$, respectively, we have
\[
 \widehat A_k^{\mathrm M}
   =P_i(\miuraq_{i,k})^{-1}A_k^{\mathrm M}P_i(\miuraq_{i,k-1}),\qquad
 A'_k=X_i(T_{i,k})^{-1}A_kX_i(T_{i,k-1}).
\]
Consequently
\[
 u'_{k+1}\widehat A_k^{\mathrm M}
 =X_i(T_{i,k})^{-1}u_{k+1}A_k^{\mathrm M}P_i(\miuraq_{i,k-1})
 =X_i(T_{i,k})^{-1}A_ku_kP_i(\miuraq_{i,k-1})
 =A'_ku'_k.
\]
After $r$ rotations each of the $r$ labels has moved once, so
all coordinates at index $k$ have their original values at $k-1$, and
\eqref{ade:rotation-gauge} gives
\[
 u^{(r)}_k
 =\Bigl(\prod_{a=1}^r X_{i_a}(T_{i_a,k-1})\Bigr)^{-1}
   u_k\prod_{a=1}^r P_{i_a}(\miuraq_{i_a,k-1})
 =A_{k-1}^{-1}u_kA_{k-1}^{\mathrm M}=u_{k-1},
\]
where $i_1=i$. Thus the remaining $r-1$ rotations followed by the
shift $k\mapsto k+1$ invert the first rotation in both orders.
\end{proof}

An ordering of the simple reflections orients each edge of the Dynkin
diagram from the earlier vertex to the later one. Moving the first
reflection to the end reverses the arrows at its vertex. The Dynkin
diagrams of types $A$, $D$, and $E$ are trees, so the possible changes
of orientation are described by the following theorem of Eriksson
and Eriksson \cite[Proposition~2.3, tree case]{ErikssonEriksson}.
A sink of an oriented tree is a vertex at which all incident arrows
point inward. A reversal at a sink reverses each of these arrows.

\begin{theorem}
\label{thm:eriksson-tree}
Any orientation of a finite tree can be obtained from any other by a
sequence of reversals at sinks.
\end{theorem}

\begin{corollary}
\label{ade:coxeter-coordinate-shifts}
For any two ordered products defining Coxeter elements $c$ and $c'$
in types $A$, $D$, and $E$, the first product can be changed into the
second by interchanging consecutive commuting reflections and moving
the first reflection to the end. Each such sequence induces a regular
bijection between the solutions of \eqref{ade:gauge-equation} for
$c$ and $c'$. If the sequence moves $s_i$ from first to last $n_i$
times, then
\[
 \miuraq'_{i,k}=\miuraq_{i,k-n_i},\qquad T'_{i,k}=T_{i,k-n_i}
 \qquad(i\in I,\ k\in\mathbb Z).
\]
\end{corollary}

\begin{proof}
At a source $i$ every incident arrow points outwards. Every preceding
vertex is therefore nonadjacent to $i$, so its reflection commutes
with $s_i$. Once $s_i$ has been moved to the front, moving it to the
end reverses precisely the arrows at $i$. Theorem~\ref{thm:eriksson-tree}
applied to the opposite orientations joins any two orientations by
such reversals at sources. Once the orientations agree, the first
vertex in the target ordering is a source in both. Its reflection
can be commuted to the front, and deletion of this vertex reduces
the comparison to the remaining orderings. Induction gives the
target ordering.

At nonadjacent vertices the rank-one subgroups commute, so
interchanging the factors preserves both products in \eqref{ade:gauge-equation}.
By Lemma~\ref{ade:rotation}, each rotation gives a regular bijection
and decreases the index of the coordinates at the moved vertex by one.
After $s_i$ has moved $n_i$ times, the coordinates at $i$ and index $k$
therefore equal their original values at $k-n_i$.
\end{proof}

Corollary~\ref{ade:coxeter-coordinate-shifts} uses the ordered
factors of a sequence in Miura form. For a general sequence
$(g_k)$ in the Coxeter cell, each $g_k$ also admits a factorization
into elements of the cells for the individual simple reflections.
This extends reduction to the Steinberg section to the whole cell.

\begin{theorem}
\label{thm:ade-full-cell}
For every Coxeter element of type $A$, $D$, or $E$, the map
\[
 N_+^{\Z}\times\mathcal A_c^{\Z}\longrightarrow M_c^{\Z},\qquad
 (u,A)\longmapsto (u_{k+1}^{-1}A_k u_k)_{k\in\Z},
\]
is an isomorphism of affine schemes.
\end{theorem}

\begin{proof}
For the linear ordering in type $A$, Proposition~\ref{ade:companion}
gives the unique $u_k$ as a polynomial in finitely many $g_b$,
and hence also gives $A_k=u_{k+1}g_ku_k^{-1}$. For the bipartite
orderings in types $D,E$, Lemma~\ref{ade:blocks} verifies the
hypotheses of Lemma~\ref{ade:contraction}, which gives the same
existence, uniqueness, and polynomial dependence. These formulas
define the inverse morphism in the stated orderings.

Suppose the assertion holds for $c=sw$, where $s=s_i$ and
$w=s_{i_2}\cdots s_{i_r}$. The rank-one case is already covered, so $r\ge2$ and the rotated
product is $c'=ws$. For each of
$d=s,w,c,c'$, Proposition~\ref{prop:cell-parametrizations} gives
regular maps $a_d:M_d\to\mathcal A_d$ and $n_d:M_d\to N_+$ with
$x=a_d(x)n_d(x)$. The ordered products defining the sections give
$\mathcal A_s\mathcal A_w=\mathcal A_c$ and
$\mathcal A_w\mathcal A_s=\mathcal A_{c'}$. Thus for $p\in N_+(D)$,
$f\in M_s(D)$ and $h\in M_w(D)$,
\[
 pfh=\underbrace{a_s(pf)a_w(n_s(pf)h)}_{\in\mathcal A_c(D)}
       n_w(n_s(pf)h).
\]
Uniqueness of the factors, first for $c=sw$ and then for $c'=ws$, gives
\begin{equation}\label{ade:elementary-normalization}
 \begin{aligned}
 n_c(pfh)&=n_w\bigl(n_s(pf)h\bigr),&
 n_{c'}(phf)&=n_s\bigl(n_w(ph)f\bigr).
 \end{aligned}
\end{equation}

For a prescribed sequence $(g'_k)$ in $M_{c'}(D)$, the coordinates
$T'_{j,k}$ of $a_{c'}(g'_k)$ give regular factors
\[
 g'_k=h_k f_{k-1},\qquad
 h_k=\prod_{b=2}^{r}X_{i_b}(T'_{i_b,k}),\qquad
 f_{k-1}=X_i(T'_{i,k})n_{c'}(g'_k).
\]
Then $g_k=f_k h_k\in M_c(D)$ by the preceding factorization.
The equation $u_{k+1}g_k=A_ku_k$ says that $u_k$ is exactly the
right $N_+$ factor of $u_{k+1}g_k$. The equations for the two
sequences are therefore
\[
 u_k=n_c(u_{k+1}g_k),\qquad u'_k=n_{c'}(u'_{k+1}g'_k).
\]
Equation~\eqref{ade:elementary-normalization} identifies their solutions by
\[
 u'_{k+1}=n_s(u_{k+1}f_k),\qquad
 u_k=n_w(u'_{k+1}h_k).
\]
Indeed, substitution of either formula into the other gives
\[
 n_w\bigl(n_s(u_{k+1}f_k)h_k\bigr)=n_c(u_{k+1}g_k)=u_k,\qquad
 n_s\bigl(n_w(u'_{k+1}h_k)f_{k-1}\bigr)=n_{c'}(u'_{k+1}g'_k)=u'_k .
\]
Thus the unique solution for $(g_k)$ gives the unique solution for
$(g'_k)$. The formulas for $h_b$ and $f_b$ express
$g_b=f_bh_b$ as a regular function of $g'_b$ and $g'_{b+1}$. If $u_{k+1}$
depends regularly on $g_{k+1},\ldots,g_{k+L}$, then
$u'_{k+1}=n_s(u_{k+1}f_k)$ depends regularly on
$g'_{k+1},\ldots,g'_{k+L+1}$. This finite dependence holds in the
initial orderings. The calculation thus proves it by induction for
each subsequent rotation. The other component of the inverse is
$A'_k=a_{c'}(u'_{k+1}g'_k)$, which is regular as well.
Interchanging commuting reflections preserves $\dot c$ and the
ordered product defining $\mathcal A_c$. By
Corollary~\ref{ade:coxeter-coordinate-shifts}, these interchanges
and rotations reach every Coxeter ordering.
\end{proof}

Each time the first reflection moves to the end, the range of indices
in the formula for $u_k$ grows by at most one. After $d$ such moves
and any number of interchanges of commuting reflections, $u_k$ depends only on
$g_k,\ldots,g_{k+L_0+d-1}$, where $L_0=N-1$ in $\mathrm{SL}_N$
and $L_0=h/2$ in types $D$ and $E$.

Under the isomorphism of Theorem~\ref{thm:ade-full-cell}, the map
$\pi_c:M_c^{\mathbb Z}\to\mathcal A_c^{\mathbb Z}$ taking a sequence
to its representative in $\mathcal A_c^{\mathbb Z}$ is projection to the second factor. The gauge transformation
determined by $v=(v_k)\in N_+(D)^{\mathbb Z}$ acts in product coordinates
by $(u,A)\mapsto(uv^{-1},A)$, where $(uv^{-1})_k=u_kv_k^{-1}$. Thus the action is
free and each orbit meets $\mathcal A_c^{\mathbb Z}$ exactly once.
Its invariant functions are
\[
 \mathbb C[M_c^{\mathbb Z}]^{N_+^{\mathbb Z}}
    =\pi_c^*\mathbb C[\mathcal A_c^{\mathbb Z}].
\]
Indeed, evaluating the invariance identity at $v=u$ gives
$F(u,A)=F(1,A)$ over every coefficient algebra.

\begin{corollary}\label{ade:normalization}
For every Coxeter ordering in types $A$, $D$ and $E$, each sequence
$A^{\mathrm M}\in M_c^-(D)^{\mathbb Z}$ determines a unique
$(u_k)\in N_+(D)^{\mathbb Z}$ satisfying \eqref{ade:gauge-equation}. The coordinates of $u_k$ and $A_k$
are regular functions of finitely many $\miuraq_{j,s}$. The
assignment $A^{\mathrm M}\mapsto A$ defines a morphism
\[
 \mu_c:(M_c^-)^{\mathbb Z}\longrightarrow\mathcal A_c^{\mathbb Z}.
\]
It commutes with the lattice translation $k\mapsto k+1$.
\end{corollary}

\begin{proof}
The inverse isomorphism of Theorem~\ref{thm:ade-full-cell} assigns
to a sequence $(g_s)$ in $M_c$ its unique pair $(u,A)$. Restricting
to $g_s=A_s^{\mathrm M}$ gives $\mu_c=\pi_c|_{(M_c^-)^{\mathbb Z}}$.
At each $k$, the coordinates of $u_k$ and $A_k$ are regular functions
of finitely many $g_s$. By Proposition~\ref{prop:cell-parametrizations},
each $A_s^{\mathrm M}$ depends regularly on the invertible coordinates
$\miuraq_{j,s}$, so these functions become Laurent polynomials in the
Miura variables at those indices. Translating \eqref{ade:gauge-equation}
by one index takes the pair $(u,A)$ for $A^{\mathrm M}$ to
$(\tau u,\tau A)$ for $\tau A^{\mathrm M}$. Uniqueness then gives
$\mu_c(\tau A^{\mathrm M})=\tau\mu_c(A^{\mathrm M})$.
\end{proof}

\begin{remark}
\label{rem:coxeter-fksz}
For the Miura transformation of Corollary~\ref{ade:normalization},
a rotation shifts the index at $i$ by $k\mapsto k-1$
(Lemma~\ref{ade:rotation}). Along the orbit $u=q^k u_*$ this is the
change of spectral argument $u\mapsto q^{-1}u$ in
\cite[Proposition~4.10]{FKSZ}.
For a $Z$-twisted Miura $(G,q)$-oper with regular singularities,
the ordered product in that convention has the form
$A(u)=F_i(u)R(u)$, where $F_i$
is the first simple-root factor and $R$ is the product of the remaining
factors. Moving $s_i$ to the end is realized by the rational gauge
transformation
\[
 g(u)=F_i(q^{-1}u)^{-1},\qquad
 g(qu)A(u)g(u)^{-1}=R(u)F_i(q^{-1}u).
\]
This identifies the corresponding spaces of $Z$-twisted Miura
$(G,q)$-opers and preserves nondegeneracy. Only the
singularity polynomial at $i$ changes, from $\Lambda_i(u)$ to
$\Lambda_i(q^{-1}u)$.

In types $A$, $D$, and $E$, Corollary~\ref{ade:coxeter-coordinate-shifts}
connects any two orderings by rotations and interchanges of commuting
reflections. The interchanges preserve the spectral arguments, so the
composite isomorphism takes the space for $c$ with singularity
polynomials $\Lambda_j(u)$ to the space for $c'$ with polynomials
$\Lambda_j(q^{-n_j}u)$, where $n_j$ counts the rotations that move
$s_j$ to the end. At the level of $QQ$-systems,
\cite[(7.33) and Theorem~7.12]{FKSZ} gives the equivalence for any two
Coxeter elements at fixed twist parameters by shifting the argument
at vertex $j$ by a factor $q^{d_j}$ and rescaling the coefficients.
The singularity polynomials thus become $\Lambda_j(q^{d_j}u)$ for
suitable integers $d_j$.
\end{remark}

\section{Screenings of the Miura coordinates in simply laced type}
\label{ade:screening-section}

\sectionoverview{In this section, we construct infinitesimal variations
of a sequence in Miura form that leave its Steinberg representative fixed.
We identify the resulting vector fields with the quantum screenings
on a spectral lattice.}

In this section we take $\mathfrak a=\operatorname{Lie}G$. The indices,
Cartan matrix and fundamental weights of Section~\ref{sec:quantum-preliminaries}
then agree with $I$, $A^G$ and $\omega_i$, respectively. For the ordinary
connection on the simply laced group $G$, Corollary~\ref{ade:normalization}
expresses each Steinberg coordinate $\mu_c^*T_{i,k}$ as a Laurent
polynomial in
\[
 \mathcal Y_c=\mathbb C[\miuraq_{j,s}^{\pm1}:j\in I,\ s\in\mathbb Z].
\]
For a fixed simple root $i$, a change of frame $v_k=x_i(\beta_k)$
preserves Miura form precisely when
\[
 x_i(\beta_{k+1}) A_k^{\mathrm M}
    =\widehat A_k^{\mathrm M}x_i(\beta_k),
 \qquad \widehat A_k^{\mathrm M}\in M_c^- .
\]

To express the infinitesimal change in Miura coordinates, we write
$j\sim i$ for adjacent vertices and $j\prec i$ when $s_j$ precedes
$s_i$ in the chosen ordering. In the product $A_{k+1}^{\mathrm M}A_k^{\mathrm M}$,
the factors between $P_i(\miuraq_{i,k+1})$ and $P_i(\miuraq_{i,k})$ come from
vertices following $i$ at index $k+1$ and vertices preceding $i$ at
index $k$. Of these, only the neighbors of $i$ affect the parameter
of the root subgroup at $i$. The Laurent monomial relating the two
variations is
\begin{equation}\label{ade:geometric-root}
 \mathsf R_{i,k}=
 \frac{\miuraq_{i,k}\miuraq_{i,k+1}}
 {\displaystyle
  \prod_{j\sim i,\ j\prec i}\miuraq_{j,k}
  \prod_{j\sim i,\ i\prec j}\miuraq_{j,k+1}} .
\end{equation}
With the convention that an empty product equals one, each
$\mathsf R_{i,k}$ is a unit. We use it to define the $\mathcal Y_c$-module
$\mathcal M_i$ on generators $s_{i,k}$, $k\in\mathbb Z$, with relations
\begin{equation}\label{ade:screening-recurrence}
 s_{i,k+1}=\mathsf R_{i,k}s_{i,k}\qquad(k\in\mathbb Z).
\end{equation}
These relations identify $\mathcal M_i$ with the free module on $s_{i,0}$ by
$s_{i,k}=d_{i,k}s_{i,0}$, where
\[
 d_{i,0}=1,\qquad
 d_{i,k}=\begin{cases}
 \prod_{s=0}^{k-1}\mathsf R_{i,s}, & k>0,\\
 \prod_{s=k}^{-1}\mathsf R_{i,s}^{-1}, & k<0 .
 \end{cases}
\]
An algebraic vector field on the affine scheme of Miura sequences is
a $\mathbb C$-linear derivation of its coordinate ring $\mathcal Y_c$
with values in $\mathcal Y_c$. The coefficients $d_{i,k}$ define such
a field by
\[
 \mathscr D_i F=\sum_{k\in\mathbb Z} d_{i,k}\,\miuraq_{i,k}
                        \frac{\partial F}{\partial \miuraq_{i,k}},
 \qquad F\in\mathcal Y_c .
\]
The sum is finite for each $F$. In terms of the generators
$s_{i,k}$, the same field gives the module-valued derivation $\mathscr S_i:\mathcal Y_c\to\mathcal M_i$ by
\begin{equation}\label{ade:geometric-screening}
 \mathscr S_i(F)=(\mathscr D_i F)s_{i,0},\qquad
 \mathscr S_i(\miuraq_{j,k})=\delta_{ij}\miuraq_{i,k}s_{i,k}.
\end{equation}
Choosing $s_{i,p}$ as generator instead multiplies $\mathscr D_i$
by the unit $d_{i,p}^{-1}$ and gives the same kernel.

\begin{proposition}\label{ade:screening-invariance}
For each $i\in I$, the vector field $\mathscr D_i$ annihilates every
regular function pulled back by $\mu_c$. In particular, for $j\in I$ and
$k\in\mathbb Z$,
\[
 \mathscr D_i(\mu_c^*T_{j,k})=0,\qquad
 \mathscr S_i(\mu_c^*T_{j,k})=0 .
\]
\end{proposition}

\begin{proof}
We compute the change in the Miura variables induced by
$v_k=x_i(\beta_k)$ and identify its derivative with $\mathscr D_i$.
At the factor $P_i(\zeta)$, for $\zeta$ and $\zeta+z$ invertible, the
required factorization follows from the identity in $\mathrm{SL}_2$
\[
 \begin{pmatrix}1&z\\0&1\end{pmatrix}
 \begin{pmatrix}\zeta&0\\1&\zeta^{-1}\end{pmatrix}
 =
 \begin{pmatrix}\zeta+z&0\\1&(\zeta+z)^{-1}\end{pmatrix}
 \begin{pmatrix}1&z/[\zeta(\zeta+z)]\\0&1\end{pmatrix}.
\]
Applying $\varphi_i$ yields
\begin{equation}
 x_i(z)P_i(\zeta)=P_i(\zeta+z)x_i\!\left(\frac{z}{\zeta(\zeta+z)}\right).
                                                        \label{ade:rank-one-same}
\end{equation}

In the ordered product $x_i(\beta_{k+1})A_k^{\mathrm M}$, moving
$x_i(z)$ past $P_j(\miuraq_{j,k})$ for $j\ne i$ multiplies $z$ by
$\miuraq_{j,k}$ if $j\sim i$ and leaves $z$ unchanged otherwise, by
\eqref{ade:rank-one-distinct}. The parameter reaching the factor
$P_i(\miuraq_{i,k})$ is therefore
$\beta_{k+1}\prod_{j\sim i,\,j\prec i}\miuraq_{j,k}$.
Equation~\eqref{ade:rank-one-same} adds this parameter to
$\miuraq_{i,k}$. For the variation to preserve Miura form there must be
a Miura product $\widehat A_k^{\mathrm M}$ such that
\[
 x_i(\beta_{k+1})A_k^{\mathrm M}
   =\widehat A_k^{\mathrm M}x_i(\beta_k).
\]
Thus the upper factor remaining after the last multiplication must
have parameter $\beta_k$. The change at $P_i$ and this final parameter
are given by
\[
 \widehat{\miuraq}_{i,k}=\miuraq_{i,k}+\beta_{k+1}\prod_{j\sim i,\ j\prec i}\miuraq_{j,k},
 \qquad
 \beta_k=\frac{\beta_{k+1}\prod_{j\sim i}\miuraq_{j,k}}
                    {\miuraq_{i,k}\widehat{\miuraq}_{i,k}},
\]
with $\widehat{\miuraq}_{j,k}=\miuraq_{j,k}$ for $j\ne i$, wherever
$\widehat{\miuraq}_{i,k}$ is a unit.
Over $\mathcal Y_c[\varepsilon]/(\varepsilon^2)$, substitution of
$\beta_k=\varepsilon r_k$ and
$\widehat{\miuraq}_{i,k}=\miuraq_{i,k}+\varepsilon\dot{\miuraq}_{i,k}$
gives
\[
 \begin{aligned}
 \dot{\miuraq}_{i,k}&=r_{k+1}\prod_{j\sim i,\ j\prec i}\miuraq_{j,k},\\
 r_k&=\frac{\dot{\miuraq}_{i,k}}{\miuraq_{i,k}^2}\prod_{j\sim i,\ i\prec j}\miuraq_{j,k}
      =\frac{\prod_{j\sim i}\miuraq_{j,k}}{\miuraq_{i,k}^2}\,r_{k+1}.
 \end{aligned}
\]
The normalization $\dot{\miuraq}_{i,0}=\miuraq_{i,0}$ fixes
$r_1=\miuraq_{i,0}\prod_{j\sim i,\ j\prec i}\miuraq_{j,0}^{-1}$.
The recurrence for $r_k$ has invertible multipliers and thus
determines $r_k$ from $r_1$ in both directions.
Eliminating $r_{k+1}$ between indices $k$ and $k+1$ gives
\[
 \frac{\dot{\miuraq}_{i,k+1}}{\miuraq_{i,k+1}}
 =\frac{\miuraq_{i,k}\miuraq_{i,k+1}}
 {\prod_{j\sim i,\ j\prec i}\miuraq_{j,k}
  \prod_{j\sim i,\ i\prec j}\miuraq_{j,k+1}}\,
  \frac{\dot{\miuraq}_{i,k}}{\miuraq_{i,k}}
 =\mathsf R_{i,k}\frac{\dot{\miuraq}_{i,k}}{\miuraq_{i,k}} .
\]
Thus $\dot{\miuraq}_{i,k}/\miuraq_{i,k}$ and $d_{i,k}$ satisfy
the same recurrence with value $1$ at $k=0$. Since the other Miura
coordinates remain fixed, we have
$\widehat{\miuraq}_{j,k}=\miuraq_{j,k}(1+\varepsilon\delta_{ij}d_{i,k})$
for every $j,k$. Since $(1+\varepsilon d)^n=1+n\varepsilon d$
also for $n<0$, the corresponding substitution on $\mathcal Y_c$ is
$F\mapsto F+\varepsilon\mathscr D_iF$.

If $u_{k+1}A_k^{\mathrm M}u_k^{-1}=A_k$ is the Steinberg representative,
then $v_k=x_i(\varepsilon r_k)$ satisfies
\[
 (u_{k+1}v_{k+1}^{-1})\,\widehat A_k^{\mathrm M}\,(u_kv_k^{-1})^{-1}=A_k .
\]
Corollary~\ref{ade:normalization} gives
$\mu_c(\widehat A^{\mathrm M})=\mu_c(A^{\mathrm M})$. Hence for every
regular function $H$ on $\mathcal A_c^{\mathbb Z}$,
\[
 0=H(\mu_c(\widehat A^{\mathrm M}))-H(\mu_c(A^{\mathrm M}))
   =\varepsilon\mathscr D_i(\mu_c^*H).
\]
The coefficient of $\varepsilon$ vanishes. For $H=T_{j,k}$,
\eqref{ade:geometric-screening} gives the second assertion.
\end{proof}

In $\mathrm{SL}_2$ the cancellation can be read directly from
$\mu_c^*T_{1,k}=\miuraq_{1,k}+\miuraq_{1,k+1}^{-1}$
(Proposition~\ref{ade:companion}). Since
$d_{1,k+1}=\miuraq_{1,k}\miuraq_{1,k+1}d_{1,k}$, we have
\[
 \mathscr D_1(\mu_c^*T_{1,k})
   =\miuraq_{1,k}d_{1,k}-\miuraq_{1,k+1}^{-1}d_{1,k+1}=0 .
\]

The finite change of frame in the proof of
Proposition~\ref{ade:screening-invariance} also describes changes of
Miura form for functions of a spectral parameter. If $A_k^{\mathrm M}$
and $x_i(\beta_k)$ are the values of $A$ and $g_i$ along an orbit of
a shift $\tau$, their defining equation becomes
$g_i(\tau(u))A(u)=\widehat A(u)g_i(u)$. For the additive shift
$\tau(u)=u+h$, the polynomial solutions of the reproduction equation in
\cite[Theorem~3.1]{MV} give rational parameters $b(u)$ for
$g_i(u)=\exp(b(u)e_i)$. For the multiplicative shift $\tau(u)=qu$,
the B\"acklund transformation of \cite[Proposition~7.1]{FKSZ} has the
form $g_i(u)=\exp(\mu_i(u)f_i)$ with rational $\mu_i$. The negative
generator $f_i$ occurs because $B_-$ stabilizes the distinguished flag
of the Miura--Pl\"ucker $(G,q)$-oper in that convention. For a differential
Miura--Pl\"ucker $G$-oper, a change of frame in the same simple-root subgroup replaces the polynomial
$q_+^i$ by $q_-^i$ and reflects the Cartan twist by $s_i$
\cite[Proposition~6.1]{BSZ}. With our stabilizer $B_+$ the corresponding
subgroup is $x_i(\beta_k)=\exp(\beta_k e_i)$.

To identify the infinitesimal action of these changes of frame with the
quantum screening, we assign a spectral parameter to each $\miuraq_{i,k}$. A height function $\xi:I\to\mathbb Z$ satisfies
$\xi_j=\xi_i-1$ on every arrow $i\to j$ in the orientation determined
by $c$. Since the diagram is a tree, the height at one vertex determines
all the others. For $a\in\mathbb C^*$ and the quantum parameter $\qaff$
of Section~\ref{sec:quantum-preliminaries}, the homomorphism assigning
the spectral variables is
\begin{equation}\label{ade:spectral-dictionary}
 \jmath_{a,\xi}:\mathcal Y_c\longrightarrow\mathcal Y,\qquad
 \miuraq_{i,k}\longmapsto Y_{i,a\qaff^{2k+\xi_i}}.
\end{equation}
Increasing $k$ by one multiplies the spectral parameter by $\qaff^2$.
Since $\qaff$ has infinite order, distinct indices give distinct
parameters at fixed $i$. Variables at different vertices are
independent. Hence $\jmath_{a,\xi}$ identifies $\mathcal Y_c$ with
the Laurent ring $\mathcal Y_\Lambda$ on
$\Lambda=\{(i,a\qaff^{2k+\xi_i}):i\in I,\ k\in\mathbb Z\}$.

\begin{lemma}\label{ade:screening-comparison}
Under \eqref{ade:spectral-dictionary}, the derivation $\mathscr S_i$
agrees with the quantum screening $S_i$ of
\eqref{eq:quantum-screening} restricted to $\mathcal Y_\Lambda$.
The target $\mathcal M_i$ embeds in the full screening module, and
\[
 \jmath_{a,\xi}(\mathsf R_{i,k})=\mathsf R_i(a\qaff^{2k+\xi_i+1}).
\]
In particular $\mathscr S_i(F)=0$ if and only if
$S_i(\jmath_{a,\xi}(F))=0$, for every $F\in\mathcal Y_c$.
\end{lemma}

\begin{proof}
At vertex $i$ the spectral parameter is $b_k=a\qaff^{2k+\xi_i}$.
For adjacent vertices $i,j$, the height decreases by one from the
earlier to the later vertex in the ordering for $c$. Thus
\[
 \begin{aligned}
 a\qaff^{2k+\xi_j}
   &=b_k\qaff &&(j\prec i,\ \xi_j=\xi_i+1),\\
 a\qaff^{2(k+1)+\xi_j}
   &=b_k\qaff &&(i\prec j,\ \xi_j=\xi_i-1).
 \end{aligned}
\]
Since $b_{k+1}=b_k\qaff^2$, substitution in
\eqref{ade:geometric-root} gives
\[
 \jmath_{a,\xi}(\mathsf R_{i,k})
 =Y_{i,b_k}Y_{i,b_k\qaff^2}\prod_{j\sim i}Y_{j,b_k\qaff}^{-1}
 =\mathsf R_i(b_k\qaff) .
\]
The relation $s_{i,b_{k+1}}=\mathsf R_i(b_k\qaff)s_{i,b_k}$
therefore defines a map $\mathcal M_i\to\mathcal S_i$ sending
$F s_{i,k}$ to $\jmath_{a,\xi}(F)s_{i,b_k}$. The module $\mathcal M_i$
is free on $s_{i,0}$, and $s_{i,b_0}$ is a basis element of
$\mathcal S_i$ for the orbit $\{b_0\qaff^{2n}:n\in\mathbb Z\}$.
Since $\jmath_{a,\xi}$ is injective, so is its extension to
the modules. On each Laurent generator $\miuraq_{j,k}$, the image
of $\mathscr S_i(\miuraq_{j,k})$ under this extension is
$\delta_{ij}Y_{i,b_k}s_{i,b_k}=S_i(\jmath_{a,\xi}(\miuraq_{j,k}))$.
The Leibniz rule gives the equality on all of $\mathcal Y_c$, and the
injectivity of the module map gives the stated equivalence of vanishing.
\end{proof}

\section{Fundamental \texorpdfstring{$q$}{q}-characters in simply laced type}
\label{ade:character-section}

\sectionoverview{In this section we identify the pullbacks of the Steinberg
coordinates with fundamental $q$-characters in simply laced type and
with folded characters for the fixed subgroup. We also express the
simply laced formula in lowest-weight minors.}

By Proposition~\ref{ade:screening-invariance} and
Lemma~\ref{ade:screening-comparison}, every quantum screening
annihilates $\jmath_{a,\xi}(\mu_c^*T_{i,k})$. Among Laurent
polynomials with this property, the coefficients of the dominant
monomials determine the polynomial by Lemma~\ref{lem:screened-uniqueness}. Since the substitution
\eqref{ade:spectral-dictionary} preserves exponents, it suffices
to prove that the only dominant monomial in $\mu_c^*T_{i,k}$ is
$\miuraq_{i,k}$, with coefficient one.

\begin{lemma}\label{ade:A-dominance}
For $c=s_1\cdots s_{N-1}$ in type $A_{N-1}$, the only dominant
monomial of $\mu_c^*T_{i,k}$ is $\miuraq_{i,k}$, with coefficient $1$.
\end{lemma}

\begin{proof}
Each term of \eqref{ade:companion-coefficients} corresponds to a subset
$\{j_1,\ldots,j_i\}\subseteq\{1,\ldots,N\}$ with $j_1<\cdots<j_i$.
If $j_{p_0}=b,\ldots,j_{p_1}=e$ is a maximal run of consecutive
integers in this subset, then $j_p=b+p-p_0$ for $p_0\le p\le p_1$.
The corresponding factors telescope to
\[
 \prod_{p=p_0}^{p_1}
 \frac{\miuraq_{j_p,k+j_p-p}}{\miuraq_{j_p-1,k+j_p-p}}
 =\prod_{d=b}^e\frac{\miuraq_{d,k+b-p_0}}{\miuraq_{d-1,k+b-p_0}}
 =\frac{\miuraq_{e,k+b-p_0}}{\miuraq_{b-1,k+b-p_0}}.
\]
At a gap $j_{p+1}\ge j_p+2$, one has
$k+j_{p+1}-p-1>k+j_p-p$. The lattice index is constant on each interval and increases at every
gap. The denominator of an interval beginning at $b>1$ therefore
cannot cancel with a numerator from another interval. Its exponent in
the full monomial is $-1$, also when $e=N$ and $\miuraq_{N,k+b-p_0}=1$.
Thus a dominant term can only come from the single interval
$\{1,\ldots,i\}$. Its product is
$\miuraq_{i,k}/\miuraq_{0,k}=\miuraq_{i,k}$, and it occurs once in
\eqref{ade:companion-coefficients}.
\end{proof}

The type $A$ calculation uses the explicit coefficients of the scalar
operator. In types $D$ and $E$, the torus action on the Miura variables
rules out dominant monomials other than $\miuraq_{i,k}$. For these ordinary
connections the map from $k$ to $k+1$ is $h\mapsto A_k^{\mathrm M}h$,
so $\sigma=\mathrm{id}_G$ in the following lemma. As in
Lemma~\ref{ade:contraction}, we include maps
$h\mapsto A_k^{\mathrm M}\sigma(h)$ in the statement so the argument
also applies to the non-simply-laced character formulas.

\begin{lemma}\label{lem:weight-criterion}
Let $G$ be simple and simply connected, let
$c=s_{i_1}\cdots s_{i_s}$ have distinct simple reflections, and let
$\sigma$ preserve the pinning. Suppose that each sequence
$(A_k^{\mathrm M})$ in $M_c^-$ determines a unique sequence
$(u_k)$ in $N_+$ for which
\[
 A_k=u_{k+1}A_k^{\mathrm M}\sigma(u_k)^{-1}\in\mathcal A_c
 \qquad (k\in\mathbb Z).
\]
This defines $\mu_{c,\sigma}(A^{\mathrm M})=A$. Its coordinate at $(i,k)$
is $(\mu_{c,\sigma}^*T_{i,k})(A^{\mathrm M})=\Delta_i(A_k)$.
A monomial in the $\miuraq_{j,d}$ is dominant if all its exponents are
nonnegative.
Suppose that for some integer $m\geq 0$ each $u_{k+1}$ depends regularly
only on the Miura variables at $k+1,\ldots,k+m$.
If $\gamma=c\sigma$ and $\beta_j=(1-c)\omega_j$ satisfy
\[
 \gamma^{-b}\beta_j\in\Delta_-
 \qquad(1\le b\le m,\ j\in\{i_1,\ldots,i_s\}),
\]
then the only dominant monomial of $\mu_{c,\sigma}^*T_{i,k}$ is
$\miuraq_{i,k}$, with coefficient one, for every $i\in\{i_1,\ldots,i_s\}$
and $k\in\mathbb Z$.
\end{lemma}

\begin{proof}
We first isolate the term $\miuraq_{i,k}$ in the principal minor,
and then use the torus action to determine the weights of the
remaining terms. For $i=i_p$, by \eqref{eq:rank-one-action-on-highest-weight},
\[
 P_j(\zeta)v_i=v_i\quad(j\ne i),\qquad
 P_i(\zeta)v_i=\zeta v_i+f_iv_i .
\]
The factors following $P_i$ fix $v_i$, as do the elements $u_{k+1}$
and $\sigma(u_k)^{-1}$ of $N_+$. The principal-minor formula of
Proposition~\ref{prop:cell-parametrizations} therefore gives
\begin{equation}\label{ade:leading-term}
 \begin{aligned}
 \mu_{c,\sigma}^*T_{i,k}
 &=v_i^*(u_{k+1}A_k^{\mathrm M}v_i)=\miuraq_{i,k}+F_{i,k},\\
 F_{i,k}
 &=v_i^*\!\left(u_{k+1}
       P_{i_1}(\miuraq_{i_1,k})\cdots P_{i_{p-1}}(\miuraq_{i_{p-1},k}) f_iv_i\right).
 \end{aligned}
\end{equation}
By the hypothesis on $u_{k+1}$,
\[
 F_{i,k}\in
 \mathbb C\bigl[\miuraq_{i_a,k}^{\pm1}\ (a<p),\;
 \miuraq_{i_a,k+b}^{\pm1}\ (1\le a\le s,\;1\le b\le m)\bigr].
\]
The variable $\miuraq_{i,k}$ does not occur in this ring, so its
coefficient in \eqref{ade:leading-term} is one.

To compute the weights, we let $h\in H$ act at index $d$ by
$g_d\mapsto h_{d+1}g_d\sigma(h_d)^{-1}$, where
$h_d=\gamma^{d-1}(h)$ and $\gamma(h)=\dot c\sigma(h)\dot c^{-1}$.
The identity $h_{d+1}\dot c=\dot c\sigma(h_d)$ gives
\[
 h_{d+1}(n\dot c n')\sigma(h_d)^{-1}
  =(h_{d+1}nh_{d+1}^{-1})\dot c\,
    (\sigma(h_d)n'\sigma(h_d)^{-1}).
\]
Since $H$ normalizes $N_+$, $U_c$ and $B_-$, this action preserves
$M_c$, $\mathcal A_c$ and $M_c^-$. The sequence $(h_du_dh_d^{-1})$ takes the transformed Miura
sequence to $(h_{d+1}A_d\sigma(h_d)^{-1})$ in the Steinberg section.
Uniqueness of the change of frame thus makes $\mu_{c,\sigma}$
$H$-equivariant. To compute the weight of $\Delta_j$, we use
$\sigma(h_d)=\dot c^{-1}h_{d+1}\dot c$ and the action
$(\gamma\lambda)(h)=\lambda(\gamma^{-1}(h))$ on characters. This gives
\[
 \Delta_j\!\left(h_{d+1}g_d\sigma(h_d)^{-1}\right)
 =\frac{\omega_j(h_{d+1})}{(c\omega_j)(h_{d+1})}\Delta_j(g_d)
 =(\gamma^{-d}\beta_j)(h)\Delta_j(g_d).
\]
The minor $\Delta_j$ restricts to $\miuraq_{j,d}$ on $M_c^-$ and
to $T_{j,d}$ on $\mathcal A_c$. Thus $\miuraq_{i,k}$ and
$\mu_{c,\sigma}^*T_{i,k}$ have the same weight $\gamma^{-k}\beta_i$.
The torus acts diagonally on Laurent monomials, so each monomial of
their difference $F_{i,k}$ has this weight as well.

For a root lattice element $\lambda$, let $[\alpha_i]\lambda$ denote
the coefficient of $\alpha_i$ in its expansion in simple roots.
In $\beta_{i_a}=s_{i_1}\cdots s_{i_{a-1}}\alpha_{i_a}$, each preceding
reflection changes only the coefficient of its own simple root.
Since the reflections are distinct, $[\alpha_i]\beta_i=1$ and
$[\alpha_i]\beta_{i_a}=0$ for $a<p$.
If a monomial of $F_{i,k}$ had exponents $n_{j,d}\ge0$, applying
$\gamma^k$ to its weight and comparing these coefficients would give
\[
 1=[\alpha_i]\beta_i
 =\sum_{a<p} n_{i_a,k}[\alpha_i]\beta_{i_a}
  +\sum_{b=1}^m\sum_{a=1}^s
    n_{i_a,k+b}[\alpha_i](\gamma^{-b}\beta_{i_a})\le0 ,
\]
because each $\gamma^{-b}\beta_{i_a}$ is a negative root.
Thus $F_{i,k}$ has no dominant monomial. By \eqref{ade:leading-term},
the only dominant monomial of the principal minor is $\miuraq_{i,k}$
with coefficient one.
\end{proof}

\begin{corollary}\label{ade:DE-dominance}
For a bipartite Coxeter element of type $D$ or $E$, the unique dominant monomial of $\mu_c^*T_{i,k}$ is
$\miuraq_{i,k}$, with coefficient one.
\end{corollary}
\begin{proof}
For $m=h/2$, Lemma~\ref{ade:blocks} gives both the dependence of
$u_{k+1}$ on $A_{k+1}^{\mathrm M},\ldots,A_{k+m}^{\mathrm M}$ and
the decomposition $\Delta_+=\bigsqcup_{a=0}^{m-1}c^aR_c$.
Since $\beta_j\in R_c$ and
$c^m=w_0=w_0^{-1}$, for $1\le b\le m$ we have
\[
 c^{m-b}\beta_j\in c^{m-b}R_c\subset\Delta_+,
 \qquad
 c^{-b}\beta_j=w_0c^{m-b}\beta_j\in w_0\Delta_+=\Delta_-.
\]
Thus the roots $c^{-b}\beta_j$ are negative over the full range
of indices on which $u_{k+1}$ depends. Lemma~\ref{lem:weight-criterion}
with $\sigma=1$ shows that $\miuraq_{i,k}$ is the only dominant term
and has coefficient one.
\end{proof}

\begin{proposition}\label{ade:all-dominance}
For every Coxeter ordering in types $A$, $D$, and $E$, $\mu_c^*T_{i,k}$ has exactly one
dominant monomial, namely $\miuraq_{i,k}$, with coefficient $1$.
\end{proposition}

\begin{proof}
For $c'$ in the linear ordering in type $A$ or a bipartite ordering
in type $D$ or $E$, Corollary~\ref{ade:coxeter-coordinate-shifts}
gives a bijection $\rho$ from Miura sequences for $c$ to those for
$c'$ satisfying
\[
 \rho^*(\miuraq'_{i,k})=\miuraq_{i,k-n_i},\qquad
 \rho^*(\mu_{c'}^*T'_{i,k})=\mu_c^*T_{i,k-n_i},
\]
where $n_i$ counts the moves of $s_i$ to the end.
The map $\rho^*$ only shifts the lattice indices of the Miura
variables. It therefore preserves monomial coefficients and all
exponents. By Lemma~\ref{ade:A-dominance} and
Corollary~\ref{ade:DE-dominance}, the unique dominant term for $c'$
is $\miuraq'_{i,k}$ with coefficient one. Its image is
$\miuraq_{i,k-n_i}$, and $k-n_i$ ranges over all integers.
\end{proof}

\begin{theorem}\label{ade:main}
Let $G$ be a simple simply connected group of simply laced type,
$c$ a Coxeter element, and $\jmath_{a,\xi}$ the spectral identification
\eqref{ade:spectral-dictionary}. For each $i\in I$ and $k\in\mathbb Z$,
\begin{equation}\label{ade:raw-character}
 \jmath_{a,\xi}\!\left(\mu_c^*(T_{i,k})\right)=
 \chi_{\qaff}\bigl(L(Y_{i,a\qaff^{2k+\xi_i}})\bigr).
\end{equation}
For $\mathfrak a=\operatorname{Lie}G$ and the category
$\mathcal C_\Lambda$ of Section~\ref{sec:quantum-preliminaries},
the pullback $\mu_c^*$ is injective with image given by
\begin{equation}\label{ade:character-image}
 \jmath_{a,\xi}\!\left(\mu_c^*\mathbb C[\mathcal A_c^{\mathbb Z}]\right)
 =\bigcap_{i\in I}\ker( S_i|_{\mathcal Y_\Lambda})
 =\chi_{\qaff}\bigl(\mathbb C\otimes_{\mathbb Z}K_0(\mathcal C_\Lambda)\bigr).
\end{equation}
\end{theorem}

\begin{proof}
We compare the pullback with the fundamental character by taking
their difference at $b=a\qaff^{2k+\xi_i}$,
\[
 F=\jmath_{a,\xi}(\mu_c^*T_{i,k})-\chi_{\qaff}(L(Y_{i,b})).
\]
Every $S_j$ annihilates the pullback by
Proposition~\ref{ade:screening-invariance} and
Lemma~\ref{ade:screening-comparison}. The fundamental character also
belongs to every screening kernel \cite[Corollary~5.5]{FM}. The pullback has the unique
dominant term $Y_{i,b}$ with coefficient one by
Proposition~\ref{ade:all-dominance}, and the character has the same
dominant term by Theorem~\ref{thm:FM}\textup{(iii)}. Thus $F\in\bigcap_j\ker S_j$ has no dominant
monomial, and $F=0$ by Lemma~\ref{lem:screened-uniqueness}.

In particular each $\chi_{\qaff}(L(Y_{i,b}))$ with $(i,b)\in\Lambda$
belongs to $\mathcal Y_\Lambda$. By Lemma~\ref{lem:spectral-subring},
the common kernel in that ring equals the image of the complexified
Grothendieck ring and is freely generated by these characters.
The source is the polynomial ring in the Steinberg coordinates
$T_{i,k}$. Equation~\eqref{ade:raw-character} identifies its generators
bijectively with these fundamental characters and thus gives the isomorphism
\[
 \mathbb C[T_{i,k}:i\in I,\ k\in\mathbb Z]
 \ \xrightarrow{\ \jmath_{a,\xi}\circ\mu_c^*\ }\
 \mathbb C[\chi_{\qaff}(L(Y_{i,b})):(i,b)\in\Lambda],
\]
which proves \eqref{ade:character-image}. Under the inverse
identification $[L(Y_{i,b})]$ corresponds to $T_{i,k}$.
\end{proof}

\begin{remark}\label{ade:height-rotation}
If $s_i$ is moved from first to last in the product defining $c$, the transformed
heights are $\xi'_i=\xi_i-2$ and $\xi'_j=\xi_j$ for $j\ne i$.
Indeed every neighbor originally has height $\xi_i-1$ and must
precede $i$ in the new ordering. This agrees with
$\miuraq'_{i,k}=\miuraq_{i,k-1}$, since $2k+\xi'_i=2(k-1)+\xi_i$. All
other variables retain their spectral labels.
\end{remark}

\subsection{Spectral parameters for lowest-weight minors}
\label{ade:fl-comparison}

\sectionoverview{In this subsection, we express the Steinberg coordinates of the
Miura transformation in terms of the lowest-weight minors of its input.
We derive the spectral parameter of each minor from the duality of
fundamental representations.}

Let $A^{\mathrm M}=(A_s^{\mathrm M})_{s\in\mathbb Z}$ be a sequence in $M_c^-$, and let
$w_0$ be the longest Weyl group element. The involution $j\mapsto j^*$
of the Dynkin diagram is determined by $w_0\omega_j=-\omega_{j^*}$.
In the factorization $A_s^{\mathrm M}=n_{-,s}[A_s^{\mathrm M}]_0$ the lower unipotent factor
fixes a lowest-weight vector in $V(\omega_j)$. The torus factor acts
on this vector by its weight, so
\begin{equation}\label{ade:lowest-minor}
 \Delta_{w_0\omega_j,w_0\omega_j}(A_s^{\mathrm M})
 =\omega_{j^*}([A_s^{\mathrm M}]_0)^{-1}=\miuraq_{j^*,s}^{-1}.
\end{equation}
We derive the corresponding change of spectral parameter for any
finite-dimensional complex simple Lie algebra $\mathfrak a$. Its opposition involution $\iota$ is characterized by
$w_0^{\mathfrak a}\varpi_j=-\varpi_{\iota(j)}$, where
$w_0^{\mathfrak a}$ is its longest Weyl group element. We use $H_{\mathfrak a}=r^\vee h^\vee$, where $r^\vee$ is the ratio of the
largest to the smallest squared root length and $h^\vee$ is the dual
Coxeter number. For a finite-dimensional module $V$, its contragredient
dual $V^*$ is the linear dual with action $(x\varphi)(v)=\varphi(S(x)v)$,
where $S$ is the antipode for the coproduct in
the \hyperref[app:quantum-presentation]{Appendix}. In this convention \cite[Proposition~6.18]{FM} gives the
following duality formula.

\begin{proposition}
\label{prop:fm-fundamental-dual}
For $b\in\mathbb C^\times$, the contragredient dual of $L(Y_{i,b})$
is isomorphic to $L(Y_{\iota(i),b\qaff^{-H_{\mathfrak a}}})$. Its
$q$-character is
\[
 \chi_{\qaff}\bigl(L(Y_{i,b})^*\bigr)
  =\left.\chi_{\qaff}\bigl(L(Y_{i,b})\bigr)\right|_{
    Y_{j,b\qaff^n}\,=\,Y_{j,b\qaff^{-n}}^{-1}}
 \qquad(j\in K,\ n\in\mathbb Z).
\]
\end{proposition}

For fixed $a\in\mathbb C^\times$, the involution
$\mathcal D$ of the Laurent ring on $a\qaff^{\mathbb Z}$ given by
$\mathcal D(Y_{j,a\qaff^n})=Y_{\iota(j),a\qaff^{H_{\mathfrak a}-n}}^{-1}$
satisfies
\begin{equation}\label{eq:fundamental-spectral-duality}
 \mathcal D\!\left(\chi_{\qaff}(L(Y_{i,a\qaff^p}))\right)
   =\chi_{\qaff}(L(Y_{i,a\qaff^{-p}})),\qquad p\in\mathbb Z .
\end{equation}
Indeed, at $b=a\qaff^p$ the dual has highest monomial
$Y_{\iota(i),a\qaff^{p-H_{\mathfrak a}}}$ by
Proposition~\ref{prop:fm-fundamental-dual}. Relabeling the vertices
by $\iota$ and multiplying all spectral parameters by
$\qaff^{H_{\mathfrak a}-2p}$ gives highest monomial $Y_{i,a\qaff^{-p}}$.
The composite substitution sends
$Y_{j,a\qaff^{p+n}}$ to $Y_{\iota(j),a\qaff^{H_{\mathfrak a}-p-n}}^{-1}$,
as in \eqref{eq:fundamental-spectral-duality}. In particular, substituting
$Y_{\iota(j),a\qaff^{H_{\mathfrak a}-2s-\xi_j}}=\miuraq_{j,s}^{-1}$
in $\mathcal D(Y_{j,a\qaff^{2s+\xi_j}})$ gives $\miuraq_{j,s}$.

In simply laced type $\mathfrak a=\operatorname{Lie}G$, so
$H_{\mathfrak a}$ is the Coxeter number $h$ and $\iota(j)=j^*$. By
\eqref{ade:lowest-minor}, $\miuraq_{j,s}^{-1}$ is the lowest-weight
minor of $A_s^{\mathrm M}$ in $V(\omega_{j^*})$. Substituting these
minors in \eqref{eq:fundamental-spectral-duality} gives the following
formula.

\begin{corollary}\label{ade:fl-exact}
For $a\in\mathbb C^\times$ and the heights $\xi$ of
\eqref{ade:spectral-dictionary}, the $i$th coordinate of
$\mu_c(A^{\mathrm M})$ at index $k$ is
\begin{equation}\label{eq:ade-forward-fundamental}
 T_{i,k}\bigl(\mu_c(A^{\mathrm M})\bigr)=
 \left.\chi_{\qaff}\!\left(L(Y_{i,a\qaff^{-2k-\xi_i}})\right)
 \right|_{Y_{j,a\qaff^{h-2s-\xi_{j^*}}}
           =\Delta_{w_0\omega_j,w_0\omega_j}(A_s^{\mathrm M})} .
\end{equation}
The substitution runs over $j\in I$ and $s\in\mathbb Z$.
\end{corollary}

\begin{proof}
For $p=2k+\xi_i$, the involution $\mathcal D$ carries the character
in \eqref{ade:raw-character} to the character at $a\qaff^{-p}$
by \eqref{eq:fundamental-spectral-duality}. On its variables,
$\mathcal D(Y_{j,a\qaff^{2s+\xi_j}})
=Y_{j^*,a\qaff^{h-2s-\xi_j}}^{-1}$. To recover the original value
$\miuraq_{j,s}$ after applying $\mathcal D$, the variable
$Y_{j^*,a\qaff^{h-2s-\xi_j}}$ must therefore take the value
$\miuraq_{j,s}^{-1}$. By \eqref{ade:lowest-minor} this is the
lowest-weight minor of $A_s^{\mathrm M}$ in $V(\omega_{j^*})$.
Relabeling $j^*$ as $j$ gives \eqref{eq:ade-forward-fundamental}.
\end{proof}

For $\mathfrak{sl}_2$ the Coxeter number is $h=2$, the opposition
involution fixes the single vertex, and we may choose $\xi_1=0$.
The fundamental character at $b$ is $Y_{1,b}+Y_{1,b\qaff^2}^{-1}$.
In \eqref{eq:ade-forward-fundamental}, $b=a\qaff^{-2k}$ and
$Y_{1,a\qaff^{2-2s}}=\miuraq_{1,s}^{-1}$. Its two terms therefore give
\[
 Y_{1,a\qaff^{-2k}}+Y_{1,a\qaff^{2-2k}}^{-1}
 \longmapsto \miuraq_{1,k+1}^{-1}+\miuraq_{1,k}.
\]
This is the same coordinate as in \eqref{ade:raw-character}, where the
first term gives $\miuraq_{1,k}$ and the second gives
$\miuraq_{1,k+1}^{-1}$.

The shift $h-\xi_{j^*}$ in \eqref{eq:ade-forward-fundamental} can
also be read from the action of $c$ on the roots. Along each arrow
$i\to j$ the heights in \eqref{ade:spectral-dictionary} satisfy
$\xi_j=\xi_i-1$, so $-\xi_j$ increases by one in the direction of the
arrow. This fixes the sign of the height in the root coordinates
of Kirillov and Thind
\cite[Theorem~5.4, Proposition~7.3(3)--(4), and Lemma~8.4]{KT}.

\begin{proposition}
\label{prop:kt-root-coordinates}
For a Coxeter element $c$ in simply laced type, $1-c$ identifies the
weight lattice with the root lattice. With $h$ the Coxeter number,
the roots are in bijection with the pairs $(j,n)$ for
$j\in I$ and $n\in\mathbb Z/(2h)\mathbb Z$ satisfying
$n\equiv-\xi_j\pmod 2$. The origin of the second coordinate can be chosen
so that $\beta_j=(1-c)\omega_j$ corresponds to $(j,-\xi_j)$ for every
$j$. If $\alpha$ corresponds to $(j,n)$, then $c\alpha$ corresponds
to $(j,n+2)$ and $-\alpha$ to $(j^*,n+h)$.
\end{proposition}

The least positive integer $m_j$ such that
$c^{m_j}\omega_j=w_0\omega_j$ therefore satisfies
\begin{equation}\label{ade:fl-height}
 2m_j=h+\xi_j-\xi_{j^*}.
\end{equation}
Indeed, $w_0\omega_j=-\omega_{j^*}$. Applying $1-c$ to the two
weights and comparing their root coordinates in
Proposition~\ref{prop:kt-root-coordinates} gives
\[
 c^m\omega_j=w_0\omega_j
 \quad\Longleftrightarrow\quad c^m\beta_j=-\beta_{j^*}
 \quad\Longleftrightarrow\quad
 -\xi_j+2m\equiv h-\xi_{j^*}\pmod{2h}.
\]
The parity condition at $j$ makes $h+\xi_j-\xi_{j^*}$ even.
Moreover $|\xi_j-\xi_{j^*}|\leq\operatorname{dist}(j,j^*)\leq |I|-1<h$,
where $\operatorname{dist}$ is distance in the Dynkin diagram.
Thus the value of $m_j$ in \eqref{ade:fl-height} lies strictly
between $0$ and $h$ and is the least positive solution. In particular, the exponent
of the input variable in \eqref{eq:ade-forward-fundamental} is also
$2(m_j-s)-\xi_j$. For the bipartite ordering of \cite{QDE}, the convention
$a=\qaff$ and $\xi_i=0$ on the sources gives its spectral parameters.
There $\xi_j$ is denoted by $\operatorname{ht}_c(j)$,
$m_j$ by $\mathfrak m_j$, and $T_{i,k}$ by $\Theta_{i,1}^{(k)}$.

For comparison with the spectral convention of
\cite[\S6.2, Theorem~30]{FL}, we set $a=\qaff$ and reverse the lattice
index. In our coordinates, the fundamental module for
$T_{i,-s}$ then has highest monomial $Y_{i,\qaff^{2s+1-\xi_i}}$.
For every $r\in\mathbb Z$ the variable $Y_{j,\qaff^{2(r+m_j)+1-\xi_j}}$
is evaluated at the lowest-weight minor
$\Delta_{w_0\omega_j,w_0\omega_j}(A_{-r}^{\mathrm M})$.
The resulting exponents and minors agree with that convention when
the same height function is used. For an arbitrary ordering, the
differences $\xi_j-\xi_i=-1$ along the arrows determine all heights
from the height at one source. The other sources can have different
heights. For example, in
$1\to2\to3\leftarrow4$ the normalization $\xi_1=0$ gives
$(0,-1,-2,-1)$, so the other source, vertex $4$, has height $-1$.
Normalizing at a single source works for every ordering. Adding a
common integer $d$ to the heights multiplies every spectral parameter
in this substitution by $\qaff^{-d}$.

\subsection{Folded \texorpdfstring{$q$}{q}-characters and ordinary Miura
\texorpdfstring{$q$}{q}-opers}\label{sec:folded-characters}

\sectionoverview{In this subsection, we restrict the ordinary Miura
transformation for a simply laced group $\widetilde G$ to the
identity component $G_\nu$ of the subgroup fixed by a diagram automorphism
$\nu$. This gives the folded character formula, the folded screening
operators, and the $QQ$-relations for Miura $(G_\nu,q)$-opers.}

Let $\widetilde G$ be a simply connected, simply laced simple group,
and let $\nu$ be an automorphism preserving its pinning. For a vertex
$v$ of the Dynkin diagram of $\widetilde G$, we denote its orbit by
$O=\{\nu^t(v):t\in\mathbb Z\}$. We assume that distinct vertices
in the same $\nu$-orbit are nonadjacent.
For the connected fixed subgroup $G_\nu=(\widetilde G^\nu)^\circ$,
the four nontrivial foldings give the following pairs of types
$(\widetilde G,G_\nu)$
\[
 (A_{2n-1},C_n),\qquad (D_{n+1},B_n),\qquad
 (E_6,F_4),\qquad (D_4,G_2).
\]
Here $n\geq 2$, with $D_3=A_3$. For the pair $(D_4,G_2)$, the
automorphism $\nu$ cyclically permutes the three outer vertices of $D_4$.
In the folding construction each $\nu$-orbit $O$ labels one simple
root $\alpha_O$ of $G_\nu$. If $\widetilde H$ is the maximal torus of the
pinning and $H_\nu=(\widetilde H^\nu)^\circ$, then
$\alpha_O=\widetilde\alpha_i|_{H_\nu}$ for every $i\in O$. These
restrictions agree because $\nu$ fixes $H_\nu$ pointwise. The pinned
rank-one subgroups at vertices of $O$ commute, and the corresponding
homomorphism from $\mathrm{SL}_2$ to $G_\nu$ is
$g\mapsto\prod_{i\in O}\widetilde\varphi_i(g)$
\cite[\S2.3, Lemma~2.2]{HongShen}. Its coroot is the sum of the
simple coroots at the vertices of $O$. Indeed, the cocharacter
lattice of $H_\nu$ is
$X_*(\widetilde H)^\nu$. Since $\nu$ permutes the simple coroot
basis of $X_*(\widetilde H)$, the orbit sums form a basis of
$X_*(H_\nu)$. Thus $G_\nu$ is simply connected. Evaluating the rank-one
homomorphism on the two matrices in \eqref{eq:rank-one-factors} gives
\begin{equation}\label{eq:folded-pinning}
 P_O(\zeta)=\prod_{i\in O}\widetilde P_i(\zeta),\qquad
 X_O(T)=\prod_{i\in O}\widetilde X_i(T).
\end{equation}
The positive unipotent subgroup of $G_\nu$ is $\widetilde N_+^\nu$,
also on points with values in a commutative complex algebra. Indeed,
for any complex unipotent group $U$ the exponential is a natural
isomorphism of varieties $\operatorname{Lie}U\to U$
\cite[Proposition~14.32]{Milne}. It identifies $\widetilde N_+^\nu$
with the vector space $(\operatorname{Lie}\widetilde N_+)^\nu$, so
this fixed subgroup is connected and belongs to $G_\nu$. Its Lie
algebra is the sum of the positive root spaces of $G_\nu$ for the
pinning just described. This proves the equality as algebraic groups.

An ordering of the orbits determines a Coxeter element $c_\nu$ of $G_\nu$.
Replacing the simple reflection for each orbit $O$ by the product of
the commuting simple reflections at its vertices gives a Coxeter
element $\widetilde c$ of $\widetilde G$. The resulting orientation
of the simply laced diagram is $\nu$-invariant.
Since adjacent heights differ by one in the direction of the
orientation, $\xi\circ\nu-\xi$ is constant on the connected diagram.
Iterating $\nu$ shows that this constant is zero. Thus the height
function $\xi$ takes a single value $\xi_O$ on each orbit.

We denote the spectral variables of the quantum loop algebra of
$\operatorname{Lie}\widetilde G$ by $\widetilde Y_{i,b}$.

\begin{definition}
\label{def:folded-character}
Folding is the homomorphism of Laurent polynomial rings
\[
 \pi(\widetilde Y_{i,b})=Y_{O,b}\quad(i\in O).
\]
For a finite-dimensional type-one module $V$ over this quantum loop
algebra, its folded $q$-character is $\pi(\chi_{\qaff}(V))$
\cite[Theorem~4.3]{FHR}.
\end{definition}

For a fixed $a\in\mathbb C^\times$,
\eqref{ade:spectral-dictionary} identifies the Miura variable
$\miuraq_{O,k}$ of $G_\nu$ with $Y_{O,b_{O,k}}$, where
$b_{O,k}=a\qaff^{2k+\xi_O}$.

To obtain the folded screening at $O$, we first compute the image
of the simply laced root monomial $\widetilde{\mathsf R}_i(b)$ for $i\in O$.
For a second $\nu$-orbit $O'\ne O$, the integer $m_{O'O}$ counts the neighbors in $O'$ of a vertex $i\in O$. This number is independent of
the choice of $i$, since $\nu$ preserves adjacency and permutes $O'$
and $O$. The Cartan entry of $G_\nu$ is
$A^{G_\nu}_{O'O}=\langle\alpha_O,\alpha_{O'}^\vee\rangle$. Pairing
$\widetilde\alpha_i$ with the sum of the simple coroots in $O'$
gives $A^{G_\nu}_{O'O}=-m_{O'O}$.
Consequently
\begin{equation}\label{eq:folded-root-monomial}
 \mathsf R_O(b):=\pi(\widetilde{\mathsf R}_i(b))
 =Y_{O,b\qaff^{-1}}Y_{O,b\qaff}
    \prod_{O'\ne O}Y_{O',b}^{-m_{O'O}} .
\end{equation}
For the order-three automorphism of $D_4$, the central vertex has
three neighbors in the orbit $O'$ of outer vertices. Its folded root
monomial therefore has denominator $Y_{O',b}^{3}$. In \eqref{eq:quantum-root}, the three
neighboring factors at a triple edge occur at $b\qaff^{-2},b,b\qaff^2$.
The common argument $b$ in \eqref{eq:folded-root-monomial} comes from
identifying vertices at a fixed spectral parameter.

\begin{definition}
\label{def:folded-screening}
For each $\nu$-orbit $O$, the folded screening $S_O$ is a derivation
on the Laurent ring in the variables $Y_{O',b}$, where $O'$ runs
over all $\nu$-orbits. Its values lie in the module generated by
the symbols $s_{O,b}$ over that ring. The derivation on generators
and the relations in its target are
\begin{equation}\label{eq:folded-screening}
 S_O(Y_{O',b})=\delta_{O,O'}Y_{O',b}s_{O,b},\qquad
 s_{O,b\qaff^2}=\mathsf R_O(b\qaff)s_{O,b}.
\end{equation}
\end{definition}

On a monomial $M=\prod_{O',b}Y_{O',b}^{u_{O',b}}$ the Leibniz rule gives
$S_O(M)=M\sum_b u_{O,b}s_{O,b}$. The root monomials are units and
$\qaff$ has infinite order. For any $b_0\in\mathbb C^\times$, iteration
of the relation expresses all symbols $s_{O,b_0\qaff^{2n}}$,
$n\in\mathbb Z$, as multiples of $s_{O,b_0}$. Choosing one $b_0$
from each orbit under multiplication by $\qaff^2$ gives a basis of
the target module consisting of these symbols $s_{O,b_0}$. Each relation between
successive generators gives an element of the kernel,
\[
 E_{O,b}=Y_{O,b}+Y_{O,b}\mathsf R_O(b\qaff)^{-1}.
\]
By \eqref{eq:folded-screening}, the two contributions to $S_O(E_{O,b})$
cancel exactly as in the calculation following \eqref{eq:quantum-screening}.
The folded screening theorem of Frenkel, Hernandez and Reshetikhin
gives the converse and identifies the common kernel
\cite[Propositions~3.2 and~3.3(1), Theorem~4.3]{FHR}.
We state it on a single orbit $a\qaff^{\mathbb Z}$.

\begin{theorem}
\label{thm:fhr-folded-screening}
In the Laurent ring on the parameters $b\in a\qaff^{\mathbb Z}$,
$\ker S_O$ is generated over $\mathbb C$ by the $E_{O,b}$ and the
$Y_{O',b}^{\pm1}$ with $O'\ne O$. The intersection of these kernels
is the image under $\pi$ of the simply laced $q$-character ring on
the same spectral parameters. An element of this intersection is
determined by the coefficients of its dominant monomials.
\end{theorem}

The positive simple-root subgroup of $G_\nu$ at $O$ is
$z\mapsto\prod_{i\in O}\widetilde x_i(z)$. Its tangent vector at $z=0$
is the sum of the simple-root vectors at the vertices of $O$.
To express the corresponding identity for screenings, we map the
target module of each simply laced screening $\widetilde S_i$,
$i\in O$, to the target module of $S_O$. The rule
$\pi_i(f\widetilde s_{i,b})=\pi(f)s_{O,b}$ defines a map on its
target module. Indeed the relation for $\widetilde s_{i,b\qaff^2}$
maps to the second relation in \eqref{eq:folded-screening}, by
\eqref{eq:folded-root-monomial}. For every simply laced Laurent
polynomial $f$ one has
\begin{equation}\label{eq:folded-orbit-screening}
 S_O(\pi(f))=\sum_{i\in O}\pi_i(\widetilde S_i(f)).
\end{equation}
Indeed, for $f=\widetilde Y_{j,b}$ both sides equal
$Y_{O,b}s_{O,b}$ if $j\in O$ and vanish otherwise. The Leibniz rule
with coefficients mapped by $\pi$ extends the equality to Laurent
polynomials.

On the lattice $b_{O,k}=a\qaff^{2k+\xi_O}$ the relations in
\eqref{eq:folded-screening} express every $s_{O,b_{O,k}}$ as a
multiple of $s_{O,b_{O,0}}$. Thus
$S_O(f)=(\mathscr D_O f)s_{O,b_{O,0}}$ defines a derivation
$\mathscr D_O$ of the folded Laurent ring. For $i\in O$, the simply laced derivation $\widetilde{\mathscr D}_i$ is defined by $\widetilde S_i(f)=(\widetilde{\mathscr D}_i f)\widetilde s_{i,b_{O,0}}$. Since $\pi_i(\widetilde s_{i,b_{O,0}})=s_{O,b_{O,0}}$, \eqref{eq:folded-orbit-screening} gives
\[
 \mathscr D_O(\pi(f))=
       \pi\!\left(\sum_{i\in O}\widetilde{\mathscr D}_i f\right)
\]
for every Laurent polynomial $f$ on the simply laced lattice. Thus
simultaneous variation at the vertices of $O$ induces the vector field $\mathscr D_O$ on the folded Miura variables. In particular, any polynomial killed
by the simply laced screenings restricts to the common folded kernel.

\begin{corollary}\label{cor:folded-miura}
For any ordering of the $\nu$-orbits on the Dynkin diagram of
$\widetilde G$, let $\mu_{c_\nu}$ be the ordinary discrete Miura
transformation for $G_\nu$. Its pullback of the Steinberg coordinate
$T_{O,k}$ at the simple root $\alpha_O$ and lattice index $k$ is
given, for every $i\in O$, by
\begin{equation}\label{eq:folded-character}
 \mu_{c_\nu}^*(T_{O,k})
 =\pi\!\left(\chi_{\qaff}
             \bigl(L(\widetilde Y_{i,b_{O,k}})\bigr)\right).
\end{equation}
Its only dominant monomial is $Y_{O,b_{O,k}}$, with coefficient one.
These pullbacks freely generate the common kernel of the folded
screening operators on the Laurent ring of the variables $Y_{O,b_{O,k}}$.
\end{corollary}

\begin{proof}
For a Miura sequence of $G_\nu$, the pinning in
\eqref{eq:folded-pinning} identifies the Miura variables at all vertices
of each orbit. In $\widetilde G$ the sequence and its unique Steinberg
representative therefore have the forms
\[
 A_k^{\mathrm M}
   =\prod_O\prod_{i\in O}\widetilde P_i(\miuraq_{O,k}),\qquad
 \widetilde u_{k+1}A_k^{\mathrm M}\widetilde u_k^{-1}
   =\prod_O\prod_{i\in O}\widetilde X_i(\widetilde T_{i,k}).
\]
The automorphism $\nu$ fixes $A_k^{\mathrm M}$ and permutes the
commuting factors within each orbit of the Steinberg product. Thus
$(\nu(\widetilde u_k))$ is another sequence of changes of frame taking
the same Miura sequence into the same section. Uniqueness of the
change of frame, followed by uniqueness of the factorization in that
section, gives
\[
 \nu(\widetilde u_k)=\widetilde u_k\in\widetilde N_+^\nu,\qquad
 \widetilde T_{\nu(i),k}=\widetilde T_{i,k}=:T_{O,k}\quad(i\in O).
\]
Here $\widetilde N_+^\nu$ is the positive unipotent subgroup of $G_\nu$,
and the resulting product is $\prod_O X_O(T_{O,k})$ by
\eqref{eq:folded-pinning}. Uniqueness in $\widetilde G$ implies
uniqueness in $G_\nu$. Since $\xi_i=\xi_O$ for $i\in O$,
Theorem~\ref{ade:main} gives
\[
 \mu_{c_\nu}^*T_{O,k}
   =\pi(\mu_{\widetilde c}^*\widetilde T_{i,k})
   =\pi\!\left(\chi_{\qaff}(L(\widetilde Y_{i,b_{O,k}}))\right).
\]
Equation~\eqref{eq:folded-orbit-screening} places each pullback in
every folded screening kernel.

For a monomial other than the highest in the simply laced fundamental
character, let $n_{j,s}$ be its exponent at
$\widetilde Y_{j,b_{O,k}\qaff^s}$. At the largest $s$ with some
$n_{j,s}\ne0$, Theorem~\ref{thm:FM}\textup{(iii)}
gives $n_{j,s}\le0$ for every $j$, with strict inequality for some
$j_0$. In the folded monomial the exponent at the orbit $O'$ of $j_0$ is
$\sum_{j\in O'}n_{j,s}<0$. Thus only the highest monomial has dominant
image, namely $Y_{O,b_{O,k}}$ with coefficient one.

For algebraic independence, we compare the highest weights in
products of the pullbacks $F_{O,k}=\mu_{c_\nu}^*T_{O,k}$. Each term
other than the highest in a simply laced fundamental character has
weight below the highest weight by a nonzero sum of simple roots.
Since $\widetilde\alpha_j$ restricts to $\alpha_{O'}$ for $j\in O'$,
this weight difference remains a nonzero sum of simple roots after
folding. Hence for each
finitely supported family of exponents $e_{O,k}\ge0$,
\[
 \prod_{O,k}F_{O,k}^{e_{O,k}}=
 \prod_{O,k}Y_{O,b_{O,k}}^{e_{O,k}}
    +\text{terms of strictly lower weight}.
\]
Distinct exponent families give distinct highest monomials. Among the
products occurring with nonzero coefficient in a polynomial relation,
any one of maximal highest weight would therefore contribute a
monomial that no other product could cancel. This proves algebraic
independence.

To prove generation on the stated lattice, we first restrict the
screenings from the full parameter set $a\qaff^{\mathbb Z}$. At
each orbit $O$ this set splits as
\[
 a\qaff^{\mathbb Z}
 =\{b_{O,k}:k\in\mathbb Z\}\sqcup\{\qaff b_{O,k}:k\in\mathbb Z\}.
\]
We denote by $\mathcal Y_a$ and $\mathcal Y_{a\qaff}$ the Laurent
rings in the folded variables on the first and second sets,
respectively, as $O$ varies. At a fixed $O$, the recurrence
\eqref{eq:folded-screening} only relates symbols within the same
set. Its multipliers are units, so the two sets give independent
generators $s_{O,b_{O,0}}$ and $s_{O,\qaff b_{O,0}}$ of the target
module over the full ring. On the first set, every multiplier
$\mathsf R_O(b_{O,k}\qaff)$ belongs to $\mathcal Y_a$ because
neighboring heights differ by one. The target of the restricted
screening thus embeds in the full target by sending its generator
to $s_{O,b_{O,0}}$, and
\[
 \bigcap_O\ker(S_O|_{\mathcal Y_a})
   =\mathcal Y_a\cap\bigcap_O\ker S_O ,
\]
where the kernels on the right are taken on $a\qaff^{\mathbb Z}$.

Theorem~\ref{thm:fhr-folded-screening} identifies the full kernel
with the image of the simply laced character ring under $\pi$.
By Lemma~\ref{lem:spectral-subring}, this image is generated by the
folded fundamental characters on the two lattices. Each of these
characters uses only variables on the lattice of its highest monomial,
by Theorem~\ref{ade:main} and the definition of $\pi$. Consequently,
setting the variables of $\mathcal Y_{a\qaff}$ equal to $1$ sends its
fundamental characters to scalars and fixes the generators $F_{O,k}$
on the first lattice. This specialization fixes every element of the
displayed intersection and sends the full kernel into
$\mathbb C[F_{O,k}:O,\ k\in\mathbb Z]$. Since the $F_{O,k}$ belong
to the intersection, they generate it.
\end{proof}

For example, under $A_3\to C_2$ we label the orbits $\{1,3\}$ and
$\{2\}$ by $1$ and $2$, respectively, and order them as $1<2$.
With $\xi_1=1$ and $\xi_2=0$,
\eqref{eq:folded-character} for the second orbit gives
\[
 \mu_{c_\nu}^*(T_{2,k})=
 \miuraq_{2,k}+\frac{\miuraq_{1,k}^2}{\miuraq_{2,k+1}}
       +2\frac{\miuraq_{1,k}}{\miuraq_{1,k+1}}
       +\frac{\miuraq_{2,k+1}}{\miuraq_{1,k+1}^{\,2}}
       +\frac{1}{\miuraq_{2,k+2}} .
\]
The coefficient two comes from the distinct $A_3$ monomials
$\widetilde Y_{1,b\qaff}/\widetilde Y_{3,b\qaff^3}$ and
$\widetilde Y_{3,b\qaff}/\widetilde Y_{1,b\qaff^3}$, where
$b=a\qaff^{2k}$, because both restrict to $\miuraq_{1,k}/\miuraq_{1,k+1}$.
Setting all variables equal to one gives six, the dimension of the
second fundamental $A_3$ module. The ordinary second fundamental $q$-character of
$C_2^{(1)}$ gives five under the same evaluation
\cite[Example~6.7]{FHR}.

To verify the two screening equations, the polynomials $E_{O,b}$ in
this example are
\[
 E_{1,b\qaff}=\miuraq_{1,k}+\frac{\miuraq_{2,k+1}}{\miuraq_{1,k+1}},\qquad
 E_{2,b}=\miuraq_{2,k}+\frac{\miuraq_{1,k}^{\,2}}{\miuraq_{2,k+1}}.
\]
Grouping the terms of the coordinate in two ways gives
\[
 \begin{aligned}
 \mu_{c_\nu}^*(T_{2,k})
   &=\miuraq_{2,k}+\frac{E_{1,b\qaff}^{\,2}}{\miuraq_{2,k+1}}
             +\miuraq_{2,k+2}^{-1}\\
   &=E_{2,b}+2\frac{\miuraq_{1,k}}{\miuraq_{1,k+1}}
       +\frac{E_{2,b\qaff^2}}{\miuraq_{1,k+1}^{\,2}},
 \end{aligned}
\]
where $E_{2,b\qaff^2}=\miuraq_{2,k+1}+\miuraq_{1,k+1}^{\,2}/\miuraq_{2,k+2}$.
The first expression belongs to $\ker S_1$ because $S_1$ annihilates
$E_{1,b\qaff}$ and all variables at vertex $2$. The second belongs
to $\ker S_2$ because $S_2$ annihilates $E_{2,b}$,
$E_{2,b\qaff^2}$ and all variables at vertex $1$. Once the other
four terms are fixed, the first screening forces the coefficient
of $\miuraq_{1,k}/\miuraq_{1,k+1}$ to be two. Changing it by $\lambda$ changes
the first screening by
\[
 \lambda\frac{\miuraq_{1,k}}{\miuraq_{1,k+1}}
 \left(1-\frac{\miuraq_{1,k}\miuraq_{1,k+1}}{\miuraq_{2,k+1}}\right)s_{1,b\qaff},
\]
which vanishes only for $\lambda=0$ since the Miura variables are
independent.

\paragraph{Ordinary \texorpdfstring{$(G_\nu,q)$}{(G-nu,q)}-opers and their \texorpdfstring{$QQ$}{QQ}-relations.}
The restriction from $\widetilde G$ to $G_\nu$ that gives
\eqref{eq:folded-character} also applies to the Wronskian recurrence.
For a sequence $(A_k)$ in the Steinberg section of $G_\nu$ and a twist
$Z\in H_\nu$, the discrete universal Wronskian is the group-valued
sequence satisfying $W_{k+1}A_k=ZW_k$ with given $W_0$. If $W_0\in G_\nu$,
both recursions
\[
 W_{k+1}=ZW_kA_k^{-1},\qquad W_k=Z^{-1}W_{k+1}A_k
\]
take values in $G_\nu$. Thus prescribing the same sequence $(A_k)$ and
initial value $W_0$ in the two groups gives the same Wronskian.

The geometric input for the folded $QQ$-system is an ordinary Miura
$(G_\nu,q)$-oper as in Definition~\ref{def:ordinary-flag-oper}.
In the differential setting, a full turn for an untwisted vacuum affine
oper of non-simply-laced type gives such a Miura $(G_\nu,q)$-oper.
Its Wronskian minors satisfy the folded $QQ$-relations
\cite[\S7.4]{QDE}. In a local frame the connection sends $h$ in the fiber at $u$ to
$A(u)h$ in the fiber at $qu$, so along an orbit it is represented
by $A_k=A(q^ku)$. We assume that $q$ has infinite order.
In the Borel convention of \cite[\S4.1]{FKSZ} the preserved
reduction is to the positive Borel. In a frame adapted to this
reduction, $A(u)$ acts on the two-dimensional space spanned by
a highest vector $w_O$ and $f_Ow_O$ in the fundamental representation
of $G_\nu$ at the simple root indexed by $O$. These restrictions give
the rank-two connections used in the $QQ$-system.

For the polynomial $QQ$-system we use the rational setting of
\cite[\S6.1]{FKSZ} on $\mathbb P^1$. The singularities are prescribed
by monic polynomials $\Lambda_O$, and the regular semisimple element
$Z=\prod_O\alpha_O^\vee(\zeta_O)$ satisfies
$\alpha_O(Z)\notin q^{\mathbb Z}$ for every simple root. A rational
$Z$-twisted Miura $(G_\nu,q)$-oper has $A(u)=v(qu)Zv(u)^{-1}$ for some
$v\in B_+(\mathbb C(u))$ in a frame adapted to the preserved flag.
Writing $A_O$, $v_O$ and $Z_O$ for the restrictions to the
two-dimensional space at $O$ gives
$A_O(u)v_O(u)=v_O(qu)Z_O$. The $Z$-twisted Miura--Pl\"ucker condition
requires this identity only on these two-dimensional spaces, with all
$v_O$ induced by a single rational Borel-valued map
\cite[Definitions~3.2 and~4.3]{FKSZ}.

For a nondegenerate polynomial solution, the construction of
\cite[Theorem~6.1, (6.7) and (6.11)--(6.14)]{FKSZ} gives rational
frames $v_O$ whose first rows are $(Q_O^+,-Q_O^-)$. Here $Q_O^+$ is
monic. In the factorization of $v_O$ as a diagonal matrix times an
upper unipotent matrix, the latter has off-diagonal entry $-Q_O^-/Q_O^+$. Comparing the upper left entries
in $A_O(u)v_O(u)=v_O(qu)Z_O$ gives the torus coordinate
$\zeta_O Q_O^+(qu)/Q_O^+(u)$ of $A(u)$ at $O$. The upper right
entries give the $QQ$-relation.
In the chosen order $<$ of the orbits, the two constants in the equation at $O$ are
\[
 \widetilde\gamma_O=\zeta_O\prod_{O'>O}\zeta_{O'}^{-m_{O'O}},
 \qquad
 \gamma_O=\zeta_O^{-1}\prod_{O'<O}\zeta_{O'}^{m_{O'O}}.
\]
Any ordering $\prec$ of the vertices of $\widetilde G$ that lists the orbits in the chosen order gives, at each $i\in O$, the simply laced equation with right-hand side
\[
 \widetilde\Lambda_i(u)
 \prod_{j\sim i,\, i\prec j}\widetilde Q_j^+(qu)
 \prod_{j\sim i,\, j\prec i}\widetilde Q_j^+(u),
\]
where the tildes denote the polynomials for $\widetilde G$.
We identify the polynomials $\widetilde Q^\pm$, the singularity
polynomials and the torus parameters within each orbit.
Each neighboring orbit $O'$ then contributes $m_{O'O}$ equal factors,
and the two coefficients on the left become
$\widetilde\gamma_O$ and $\gamma_O$. The resulting equation is
\begin{align}
 &\widetilde\gamma_O Q_O^-(u)Q_O^+(qu)
       -\gamma_O Q_O^-(qu)Q_O^+(u) \notag\\
 &\hspace{2em}=\Lambda_O(u)
       \prod_{O'>O} Q_{O'}^+(qu)^{m_{O'O}}
       \prod_{O'<O} Q_{O'}^+(u)^{m_{O'O}} .
 \label{eq:ordinary-folded-qq}
\end{align}
This is \eqref{eq:intro-oper-qq} for $G_\nu$, with the Cartan exponents
given by the multiplicities $m_{O'O}$.

\begin{remark}
\label{rem:folded-polynomial-qq}
For fixed $Q_O^+$ and fixed right-hand side, two solutions for $Q_O^-$
in a field of meromorphic functions differ by $Q_O^+d$, where
\[
 d(qu)=(\widetilde\gamma_O/\gamma_O)d(u)=\alpha_O(Z)d(u).
\]
Indeed, this equation follows by subtracting their $QQ$-relations
and dividing by $Q_O^+(u)Q_O^+(qu)$. A nonzero rational $d$ of order
$n$ at zero would give $\alpha_O(Z)=q^n$ by comparison of leading
terms. Thus the polynomial partner, when it exists, is the unique
rational partner. Over meromorphic functions on $\mathbb C^\times$
the homogeneous equation can have nonzero solutions. For example,
when $0<|q|<1$ the function $d(u)=\theta_q(u)/\theta_q(\alpha_O(Z)u)$
has this property, where
$\theta_q(u)=\prod_{n\ge0}(1-q^n u)(1-q^{n+1}/u)$ satisfies
$\theta_q(qu)=-u^{-1}\theta_q(u)$. Hence the polynomial entries above
refer to the normalized rational frames.
\end{remark}

For a polynomial solution, the right-hand side $F_O(u)$ of
\eqref{eq:ordinary-folded-qq} can be evaluated at a zero
$w\in\mathbb C^\times$ of $Q_O^+$ and at $q^{-1}w$.
Eliminating $Q_O^-(w)$ from these two identities gives
\[
 \gamma_O F_O(w)Q_O^+(q^{-1}w)
 +\widetilde\gamma_O F_O(q^{-1}w)Q_O^+(qw)=0 .
\]
Where the denominators are nonzero, division gives the folded Bethe
equation \cite[(6.16) and Theorem~6.4]{FKSZ}. The nondegeneracy
condition comparing zeros uses the vertices of $G_\nu$ after
folding. If $O_1$ and $O_2$ are distinct and some vertex $P$ is equal
or adjacent to each, the zeros of $Q_{O_1}^+$, $Q_{O_2}^+$ and
$\Lambda_P$ must lie on disjoint multiplicative $q$-orbits when compared
across any two of the three sets
\cite[Proposition~4.8]{FKSZ}. For example,
the two outer nodes of $A_3$ share a neighbor. Assigning
them the same nonconstant $Q^+$ gives coincident roots in the
simply laced system, although they represent a single polynomial
in the folded system.

The same multiplicative shift evaluates the character formula for
a connection represented in Miura form by
$A^{\mathrm M}(u)=\prod_O P_O(\miuraq_O(u))$, with $\miuraq_O$ nonzero
meromorphic functions. At lattice index $k$ we use the group element
$A^{\mathrm M}(q^k u)$, so \eqref{eq:folded-character} is evaluated
by $Y_{O,a\qaff^{2k+\xi_O}}\mapsto \miuraq_O(q^k u)$.

\section{Twisted Miura maps in types \texorpdfstring{$B,C,F_4,G_2$}{B, C, F4, G2}}
\label{tw:section}

\sectionoverview{In this section we prove that the ordinary fundamental
$q$-characters in types $B_n$, $C_n$, $F_4$ and $G_2$ are the
Steinberg coordinates of a Miura transformation on a simply laced
group. We construct the transformation using a diagram automorphism
in each connection map and recover the $q$-character homomorphism by
varying the spectral lattice.}

\subsection{Diagram automorphisms and changes of frame on the lattice}

\sectionoverview{In this subsection we specify the simply laced group
$G$, the automorphism $\sigma$, and the simple reflections defining
the relative position of the flags of a $(G,\eta,\sigma)$-oper for each
of the four non-simply-laced types. For each
choice we describe the Steinberg section and the action of changes of
frame on the corresponding sequences of connection maps.}

In Section~\ref{sec:folded-characters}, the automorphism $\nu$
determined the fixed subgroup $G_\nu$. The connection on this subgroup
was an ordinary $q$-connection, just as in the simply laced case.
For the ordinary $q$-characters in non-simply-laced type we instead
use the $(G,\eta,\sigma)$-opers of Definition~\ref{def:eta-sigma-oper}.
Their connection maps act on the simply laced group $G$ by
$h\mapsto A(u)\sigma(h)$. The relative position of the flags at
$\eta u$ is prescribed by a product containing one simple reflection
from each $\sigma$-orbit.

For each finite type $\mathfrak a$, the table gives the pair $(G,\sigma)$
and an ordered list $J$ of representatives of the $\sigma$-orbits
of simple roots. The order in $J$ specifies the order of the
simple reflections in $c$.
The last column records
$m=|\Delta_+(G)|/|J|$, the number of positive roots divided by the
number of vertices in $J$.
\begin{equation}\label{tw:data}
\begin{array}{c|c|c|c|c}
 \mathfrak a&\text{type of }G&\sigma&J&m\\\hline
 B_n&A_{2n-1}&j\mapsto2n-j&(1,\ldots,n)&2n-1\\
 C_n&D_{n+1}&(n\ n+1)&(1,\ldots,n)&n+1\\
 F_4&E_6&(1\ 6)(3\ 5)&(1,2,3,4)&9\\
 G_2&D_4&(1\ 3\ 4)&(1,2)&6
\end{array}
\end{equation}
In the two classical families $n\ge2$. The nodes of $A_{2n-1}$
are numbered along the chain. In $D_{n+1}$ the terminal nodes $n$ and
$n+1$ are adjacent to $n-1$, so $1$ is the central
node of $D_3=A_3$. The edges of $E_6$ join the pairs
$(1,3),(3,4),(4,5),(5,6),(2,4)$, and $2$ is the central node of $D_4$.
The twisted affine types of the four pairs $(G,\sigma)$ are
$A_{2n-1}^{(2)},D_{n+1}^{(2)},E_6^{(2)},D_4^{(3)}$, respectively.
Their transposed Cartan matrices, up to simultaneous relabeling of rows
and columns, are those of the untwisted affine types $B_n^{(1)},C_n^{(1)},F_4^{(1)},G_2^{(1)}$ associated with the
finite types in the first column
\cite[Table~1]{MRV2}, \cite[Ch.~4]{Kac}.

We write $A^G$ for the Cartan matrix of the simply laced group in
the second column and $\omega_j$ for its fundamental weights. On
weights, $(\sigma\lambda)(h)=\lambda(\sigma^{-1}h)$ for $h\in H$.
For the ordering in each row, $c$ is a Coxeter element of the
parabolic subgroup generated by $s_i$, $i\in J$. The twisted Coxeter
element $\gamma=c\sigma$ acts on roots, and its lift
$\Gamma=\operatorname{Ad}(\dot c)\circ\sigma$ sends $U_\alpha$ to
$U_{\gamma\alpha}$. By Lemma~\ref{lem:inversion-roots} the displayed
product for $c$ is reduced. Since $\sigma$ preserves the positive roots, the inversion set is
$R_c=\Delta_+\cap c\Delta_-=\Delta_+\cap\gamma\Delta_-$.

The Steinberg section of the cell $M_c=N_+\dot cN_+$ is
$\mathcal A_c=U_c\dot c$, with $U_c=N_+\cap\dot cN_-\dot c^{-1}$.
Proposition~\ref{prop:cell-parametrizations} gives the coordinates
\begin{equation}\label{tw:sections}
 \mathcal A_c=\left\{\prod_{p=1}^{\ell}X_{i_p}(T_{i_p})\right\}
   \simeq\mathbb A^\ell,\qquad
 M_c^-=M_c\cap B_-=\left\{\prod_{p=1}^{\ell}P_{i_p}(\miuraq_{i_p})\right\}
   \simeq(\mathbb G_m)^\ell .
\end{equation}

Along a spectral orbit $u^{(k)}=\eta^ku_*$, the connection map
from the fiber at $u^{(k)}$ to the next fiber is represented by
$h\mapsto g_k\sigma(h)$ with $g_k=A(u^{(k)})$. To pass from
$M_c^-$ to $\mathcal A_c$ we use the action of $N_+$-valued
changes of frame derived in Section~\ref{sec:root-subgroup-reduction}.
For sequences in $M_c(D)$ over a commutative complex algebra $D$ it is
\begin{equation}\label{tw:gauge}
 (u\cdot g)_k=u_{k+1}g_k\sigma(u_k)^{-1},
 \qquad u\in N_+(D)^{\Z},\quad g\in M_c(D)^{\Z}.
\end{equation}
If $0$ is in the domain and both $A$ and a change of frame $b$ are
regular there, the two fibers coincide. The same change-of-frame
formula then replaces $A(0)$ by $b(0)A(0)\sigma(b(0))^{-1}$.
The flag condition in
Definition~\ref{def:eta-sigma-oper} is imposed at generic $u$.

\subsection{Reduction to the Steinberg section in twisted type}

\sectionoverview{In this subsection, we prove that a unique twisted unipotent
gauge transformation sends a given sequence in $M_c$ to a sequence
in $\mathcal A_c$. We determine, from the action of $c\sigma$ on
positive roots, how many terms of the given sequence are needed to compute
the transformation at each index.}

The root calculation uses the following theorem of Oh and Suh
\cite[\S3]{OhSuh}, in the formulation of \cite[Theorem~2.2(1)]{FujitaOh}.

\begin{theorem}
\label{thm:oh-suh}
For a simply laced root system with one of the nontrivial automorphisms
$\sigma$ in \eqref{tw:data}, let $r$ be its order and $h^\vee$ the
dual Coxeter number of the corresponding finite type $\mathfrak a$ in
that table. If $c$ is the product, in any order, of one simple reflection
from each $\sigma$-orbit and $m=rh^\vee/2$, then the operator
$\gamma=c\sigma$ satisfies $\gamma^m=w_0\sigma^m=-1$ on the real
span of the roots. Moreover,
\begin{equation}\label{tw:longest-word}
 w_0=c\,\sigma(c)\cdots\sigma^{m-1}(c)
\end{equation}
is a reduced expression when each factor is expanded in the chosen
simple reflections and their images under $\sigma^a$.
\end{theorem}

\begin{lemma}\label{tw:blocks}
For each row of \eqref{tw:data} and $0\le j\le m$, the sets
$S_j=\bigcap_{a=0}^{j}\gamma^{-a}\Delta_+$ satisfy
\begin{equation}\label{tw:filtration}
 S_j=\bigsqcup_{a=0}^{m-1-j}\gamma^aR_c .
\end{equation}
In particular $S_m=\varnothing$ and
$S_j=R_c\sqcup\gamma S_{j+1}$ for $0\le j<m$.
\end{lemma}

\begin{proof}
In each row of \eqref{tw:data}, $J$ contains one vertex from each
$\sigma$-orbit. The values of $m=rh^\vee/2$ in types $B_n$, $C_n$,
$F_4$, and $G_2$ are $2n-1$, $n+1$, $9$, and $6$, respectively. Theorem~\ref{thm:oh-suh} thus gives
$\gamma^m=-1$ and the reduced expression \eqref{tw:longest-word} for $w_0$.
For $0\le a<m$, the product of the first $a$ factors is
$p_a=c\,\sigma(c)\cdots\sigma^{a-1}(c)$, with $p_0=1$.
The inversion roots contributed by the next factor $\sigma^a(c)$
are therefore $p_a\sigma^aR_c=\gamma^aR_c$ by
Lemma~\ref{lem:inversion-roots}. Since a reduced expression for $w_0$
lists each positive root exactly once among its inversion roots,
\[
 \Delta_+=\bigsqcup_{a=0}^{m-1}\gamma^a R_c .
\]
For $\alpha=\gamma^a\beta$ with $\beta\in R_c$, the first $m-a$
roots $\alpha,\gamma\alpha,\ldots,\gamma^{m-a-1}\alpha$ belong to this
union, whereas $\gamma^{m-a}\alpha=\gamma^m\beta=-\beta$ is negative.
Thus $\alpha\in S_j$ exactly when $a\le m-1-j$. This proves
\eqref{tw:filtration}, including $S_m=\varnothing$. For $j<m$,
the summand at $a=0$ is $R_c$ and the others form $\gamma S_{j+1}$.
\end{proof}

\begin{proposition}\label{tw:normalization}
The morphism $(v,A)\mapsto(v_{k+1}A_k\sigma(v_k)^{-1})_k$
from $N_+^{\Z}\times\mathcal A_c^{\Z}$ to $M_c^{\Z}$ is an
isomorphism of affine schemes. For every
$g\in M_c(D)^{\Z}$ there is a unique $u\in N_+(D)^{\Z}$ such that
$u\cdot g\in\mathcal A_c(D)^{\Z}$. Each $u_k$ depends polynomially on
the root coordinates of $g_k,\ldots,g_{k+m-1}$, where $m$ is given by
\eqref{tw:data}. On $(M_c^-)^{\Z}$ the resulting Miura transformation
$\mu_{c,\sigma}:(M_c^-)^{\Z}\to\mathcal A_c^{\Z}$ is a morphism with
Laurent-polynomial coordinates.
\end{proposition}

\begin{proof}
The sets $\gamma^{-a}\Delta_+$ are closed under addition whenever
the sum is a root, so their intersection $S_j$ has the same property.
By Lemma~\ref{tw:blocks}, $S_0=\Delta_+$, $S_m=\varnothing$, and
$S_j=R_c\sqcup\gamma S_{j+1}$ for $0\le j<m$. In particular,
$R_c\cup S_j=S_j$, as required in \eqref{ade:filtration-hypothesis}.
For each $(g_k)$, Lemma~\ref{ade:contraction} therefore gives a unique
sequence $(u_k)$ for which $u_{k+1}g_k\sigma(u_k)^{-1}$ belongs to
$\mathcal A_c(D)$ at every $k$. The same lemma gives the polynomial
dependence of $u_k$ on $g_k,\ldots,g_{k+m-1}$. Hence the inverse to the
morphism in the proposition is regular and is given by
\[
 g\longmapsto\left((u_k^{-1})_k,\,
       \bigl(u_{k+1}g_k\sigma(u_k)^{-1}\bigr)_k\right).
\]
For $g_k\in M_c^-$, the product
$g_k=\prod_{j\in J}P_j(\miuraq_{j,k})$ in \eqref{tw:sections} has
root coordinates Laurent polynomial in the $\miuraq_{j,k}$.
Substitution in the second component gives the assertion for
$\mu_{c,\sigma}$.
\end{proof}

Proposition~\ref{tw:normalization} gives a Steinberg section for
twisted conjugation on sequences. To compare it with the sections
studied by He and Lusztig \cite[\S3]{HeLusztig}, we combine the
diagram automorphism and the lattice shift into one automorphism. If
$\vartheta(v)_k=\sigma(v_{k-1})$ and $v_k=u_{k+1}$, then
$(v g\vartheta(v)^{-1})_k=u_{k+1}g_k\sigma(u_k)^{-1}$.
Malten's cross-section theorem allows the translation of the lattice
as the automorphism of the base \cite[Definition~1.10 and (1.4)]{Malten}.
For $\gamma^{-1}$, moving the diagram automorphisms to the right
in $(\gamma^{-1})^m$ leaves the reverse of \eqref{tw:longest-word}
with $\sigma^{-m}$ applied to each simple reflection. This is again
a reduced expression for $w_0$, so Malten's main lemma, part~(ii),
gives condition~(1.3) of his theorem. Inversion takes his right
section for $\gamma^{-1}$ to $U_c\dot c$ in our coordinates.

\subsection{The finite Cartan matrix and the spectral lattice}

\sectionoverview{In this subsection, we determine the Cartan matrix and its symmetrizing
integers from the action of $\sigma$ on the simple roots of $G$.
We then give the spectral parameter of each Miura variable, thereby
identifying it with a variable of the $q$-character ring.}

For each $i\in J$, the size $o_i$ of its $\sigma$-orbit is the least positive integer such that $\sigma^{o_i}(i)=i$. The root used to compute
the Cartan matrix $C$ of $\mathfrak a$ is the sum
$\sum_{t=0}^{o_i-1}\alpha_{\sigma^t(i)}$ in the real span of the
roots of $G$. We normalize the Weyl-invariant inner product by
$(\alpha_v,\alpha_v)=2$. The distinct vertices $\sigma^t(i)$ in
this sum are pairwise nonadjacent in every diagram of \eqref{tw:data},
so the squared length of the sum is $2o_i$. Pairing the sum at $i$
with the coroot of the sum at $j$ and using $\sigma$-invariance
of the inner product gives
\[
 C_{ji}=\frac{1}{o_j}
    \sum_{s=0}^{o_j-1}\sum_{t=0}^{o_i-1}
        (\alpha_{\sigma^s(j)},\alpha_{\sigma^t(i)})
       =\sum_{t=0}^{o_i-1}(\alpha_j,\alpha_{\sigma^t(i)}) .
\]
For $i\ne j$, $-C_{ji}$ therefore counts the neighbors of $j$ among
the vertices $\sigma^t(i)$. To compare with folding, suppose $\nu$
acts by the same permutation on the same simply laced diagram. In
the fixed group $G_\nu$ of Section~\ref{sec:folded-characters}, the
simple root corresponding to $i$ is the restriction of
$\widetilde\alpha_i$ to $H_\nu$, and the simple coroot corresponding
to $j$ is the sum of the ambient coroots at the vertices $\nu^s(j)$,
one for each distinct vertex. Thus the negative of the Cartan entry
of $G_\nu$ in row $j$, column $i$
counts the neighbors of $i$ among the vertices $\nu^s(j)$ for $i\ne j$.
The two Cartan matrices are transposes.
For example, the central vertex of $D_4$ has three neighbors in the
outer orbit, and each outer vertex has one central neighbor.

These counts give the Cartan matrices of the Lie algebras in the first column of
\eqref{tw:data}. For $B_n$ and $C_n$ each edge of the chain is simple
except the last, for which the entries are
\[
 \begin{array}{c|cc|c}
 &C_{n-1,n}&C_{n,n-1}&(o_1,\ldots,o_n)\\\hline
 B_n&-1&-2&(2,\ldots,2,1)\\
 C_n&-2&-1&(1,\ldots,1,2).
 \end{array}
\]
For the exceptional types the labeling in \eqref{tw:data} gives
\[
 C_{F_4}=\begin{pmatrix}
 2&0&-1&0\\0&2&0&-1\\-1&0&2&-1\\0&-1&-2&2
 \end{pmatrix},\quad (o_i)=(2,1,2,1),\qquad
 C_{G_2}=\begin{pmatrix}2&-1\\-3&2\end{pmatrix},\quad(o_i)=(3,1).
\]
In each case $\operatorname{diag}(o_i)C$ is symmetric. Each
automorphism has a fixed vertex, so the orbit sizes have greatest
common divisor one. Thus $o_i=d_i^{\mathrm{sym}}$ in the notation
of Section~\ref{sec:quantum-preliminaries}, and $\qaff_i=\qaff^{o_i}$.

For a fixed $a\in\C^\times$, we seek a spectral parameter $a\qaff^{2k+\xi_i}$ for $\miuraq_{i,k}$, where $\qaff$ has infinite multiplicative order. In the root monomial
\eqref{eq:quantum-root}, the spectral exponents contributed by any
one neighboring vertex are symmetric about the midpoint of the two
numerator exponents. This condition determines the differences
$\xi_j-\xi_i$ once we know which variables $\miuraq_{j,k+s}$ occur in the
denominator. We determine their indices by factoring a simple-root
gauge transformation through successive elements in Miura form.

For a commutative complex algebra $D$, a first-order variation of an
element in Miura form is an element in Miura form over
$D[\varepsilon]/(\varepsilon^2)$ that reduces to the original element
modulo $\varepsilon$. The factors $P_j(\zeta)$,
$j\in J$, lie in the simply laced group $G$, so the Cartan entry in
\eqref{ade:rank-one-distinct} for $x_h(z)$ and $P_j(\zeta)$ is
$A^G_{jh}=\langle\alpha_h,\alpha_j^\vee\rangle$. For $j=h$,
$\zeta\in D^\times$, and $z\in\varepsilon D$, the parameter of the
right-hand factor in \eqref{ade:rank-one-same} is
$z/[\zeta(\zeta+z)]=z/\zeta^2$, since $z^2=0$.
In either case the right-hand unipotent factor lies in the root subgroup
at $h$. For $r_{k+1}\in D$, moving $v_{k+1}=x_h(\varepsilon r_{k+1})$
through the ordered product $A_k^{\mathrm M}$ therefore gives
\[
 v_{k+1}A_k^{\mathrm M}=\widehat A_k^{\mathrm M}\sigma(v_k),
 \qquad v_k\in U_{\alpha_{\sigma^{-1}(h)}} ,
\]
where $\widehat A_k^{\mathrm M}$ is again in Miura form. In particular,
the simple-root indices at consecutive lattice positions satisfy
$h_{k+1}=\sigma(h_k)$.

Suppose $h_{k+1}=i$ for a representative $i\in J$. The element
multiplying $A_{k+s}^{\mathrm M}$ on the left belongs to the root
subgroup at $\sigma^s(i)$, so its simple-root index first returns to
$i$ at $s=o_i$. For an adjacent representative $j$ following $i$ in $J$, we need the positions at which its factor $P_j$ changes the subgroup parameter. The spectral exponents of $\miuraq_{i,k}$ and $\miuraq_{i,k+o_i}$ have midpoint
$2k+\xi_i+o_i$. Between the two corresponding factors $P_i$,
the neighboring factors $P_j$ that rescale the root parameter
occur at indices $k+s$ with $j$ adjacent to $\sigma^s(i)$ and
$1\le s\le o_i$, by \eqref{ade:rank-one-distinct}. The endpoint $s=o_i$ is included
and $s=0$ is excluded because $j$ follows $i$ in each product.
The representatives of adjacent orbits in \eqref{tw:data} are themselves
adjacent. If $j$ is fixed, it is adjacent to every $\sigma^s(i)$,
so $s$ takes every value from $1$ to $o_i$. If $j$ is moved and $i$ is fixed,
the only possible value of $s$ is $1=o_i$. When both are moved, the edges
between the orbits are precisely the translates of the edge joining $i$
and $j$, as one sees on the $A_{2n-1}$ chain and on the $E_6$ edges
$(1,3),(6,5)$. Thus $j$ is adjacent to $\sigma^s(i)$ only at $s=o_i$ in the
given range.
In these cases the smallest and largest values of $s$ are
$\min(o_i,o_j)$ and $o_i$, respectively. Thus the midpoint of
the neighboring exponents $2(k+s)+\xi_j$ is
$2k+\xi_j+\min(o_i,o_j)+o_i$. Equating it to $2k+\xi_i+o_i$
gives $\xi_j-\xi_i=-\min(o_i,o_j)$.
The orbit diagram is a tree, so $\xi_1=0$ fixes all heights.
A common shift of the heights only changes $a$. With this normalization,
the identifications of the Miura variables are
\begin{equation}\label{tw:dictionary}
 \miuraq_{i,k}=Y_{i,a\qaff^{2k+\xi_i}},\qquad
 \begin{array}{c|l}
 B_n&\xi_i=-2(i-1)\ (i<n),\quad\xi_n=-2n+3,\\
 C_n&\xi_i=-(i-1),\\
 F_4&(\xi_1,\xi_2,\xi_3,\xi_4)=(0,-2,-2,-3),\\
 G_2&(\xi_1,\xi_2)=(0,-1).
 \end{array}
\end{equation}
Under this identification, increasing $k$ by one multiplies the
spectral parameter by $\qaff^2$. The pairs $(i,a\qaff^{2k+\xi_i})$ form
a set $\Lambda$ with Laurent ring $\mathcal Y_\Lambda=\C[\miuraq_{i,k}^{\pm1}:i\in J,k\in\Z]$.
The category $\mathcal C_\Lambda$ consists of the finite-dimensional
type-one $U_{\qaff}(\widehat{\mathfrak a})$-modules whose simple
constituents have highest monomial supported on $\Lambda$.

\subsection{Twisted gauge transformations and quantum screenings}

\sectionoverview{In this subsection, we differentiate twisted gauge
transformations that preserve the Steinberg coordinates. When the
orbit of $i$ has $o_i$ vertices, each variation changes a Miura
coordinate at $i$ every $o_i$ indices. We show how the $o_i$ choices
of starting index give the components of the quantum screening at $i$.}

For a nonzero first-order change of frame
$v_k=x_{h_k}(\varepsilon r_k)$, with $\varepsilon^2=0$, the
factorization used to obtain \eqref{tw:dictionary} gives
$h_{k+1}=\sigma(h_k)$. The set $J$ contains exactly one vertex $i$
in the $\sigma$-orbit of $h_1$. Since this orbit has $o_i$ vertices,
there is a unique $0\le\rho<o_i$ for which $h_{\rho+1}=i$. At index
$k$, the element $v_{k+1}$ acting on $A_k^{\mathrm M}$ thus belongs
to the positive simple-root subgroup at $\sigma^{k-\rho}(i)$.
Since $P_i$ is the only factor of
$A_k^{\mathrm M}$ indexed by a vertex in this orbit, the only Miura
coordinates that can vary are $\miuraq_{i,k}$ with
$k=\rho+n o_i$, $n\in\mathbb Z$.

At such an index we write the positive root element immediately to
the left of $P_i(\miuraq_{i,k})$ as $x_i(\varepsilon \miuraq_{i,k}d_{i,k})$. By
\eqref{ade:rank-one-same} the new Miura coordinate is
$\miuraq'_{i,k}=\miuraq_{i,k}(1+\varepsilon d_{i,k})$, and the factor to the
right of $P_i(\miuraq'_{i,k})$ is $x_i(\varepsilon d_{i,k}/\miuraq_{i,k})$.

Moving $v_{k+o_i+1}$ through $A_{k+o_i}^{\mathrm M}$ leaves
$\sigma(v_{k+o_i})$ on its right. The element $v_{k+o_i}$ then
multiplies $A_{k+o_i-1}^{\mathrm M}$ on the left. Thus between the
two factors $P_i(\miuraq_{i,k+o_i})$ and $P_i(\miuraq_{i,k})$ the lattice index
decreases, and within each product the factors occur in the order
$J$. At index $k+o_i$ only the factors after $P_i$ intervene, and
at index $k$ only those before $P_i$. At each intermediate index $r$,
the positive root element belongs to the subgroup at $\sigma^{r-k}(i)$
and crosses every factor in $A_r^{\mathrm M}$. By
\eqref{ade:rank-one-distinct}, these factors multiply its parameter by
\[
 M_{i,k}=
 \prod_{j\text{ after }i}\miuraq_{j,k+o_i}^{-A^G_{ji}}\,
 \prod_{r=k+1}^{k+o_i-1}\ \prod_{j\in J}
 \miuraq_{j,r}^{-A^G_{j,\sigma^{r-k}(i)}}\,
 \prod_{j\text{ before }i}\miuraq_{j,k}^{-A^G_{ji}} .
\]
The outgoing parameter at $P_i(\miuraq_{i,k+o_i})$ and the incoming
parameter at $P_i(\miuraq_{i,k})$ therefore satisfy
$\miuraq_{i,k}d_{i,k}=M_{i,k}d_{i,k+o_i}/\miuraq_{i,k+o_i}$.
Solving for $d_{i,k+o_i}$ gives the recurrence
$d_{i,k+o_i}=\mathsf R_{i,k}d_{i,k}$. We introduce a
$\mathcal Y_\Lambda$-module whose generators $s_{i,k}$ satisfy
this same recurrence,
\begin{equation}\label{tw:screen-recurrence}
 s_{i,k+o_i}=\mathsf R_{i,k}s_{i,k},\qquad
 \mathsf R_{i,k}=\frac{\miuraq_{i,k}\miuraq_{i,k+o_i}}{M_{i,k}} .
\end{equation}
In the product defining $M_{i,k}$ the exponent is one for a neighbor
of $\sigma^{r-k}(i)$ and zero otherwise. Thus the four choices of
$(G,\sigma)$ give the following denominators. In types $B_n$ and
$C_n$ the vertices are numbered from $1$ to $n$, so a factor with
first index $0$ is omitted.
\begin{align*}
 B_n:\quad
 M_{i,k}&=\miuraq_{i-1,k}\miuraq_{i+1,k+2} &&(i<n-1),\\
 M_{n-1,k}&=\miuraq_{n-2,k}\miuraq_{n,k+2}\miuraq_{n,k+1},&
 M_{n,k}&=\miuraq_{n-1,k};\\[2pt]
 C_n:\quad
 M_{i,k}&=\miuraq_{i-1,k}\miuraq_{i+1,k+1} &&(i<n),\\
 M_{n,k}&=\miuraq_{n-1,k+1}\miuraq_{n-1,k};\\[2pt]
 F_4:\quad
 M_{1,k}&=\miuraq_{3,k+2},&M_{2,k}&=\miuraq_{4,k+1},\\
 M_{3,k}&=\miuraq_{1,k}\miuraq_{4,k+2}\miuraq_{4,k+1},&
 M_{4,k}&=\miuraq_{2,k}\miuraq_{3,k};\\[2pt]
 G_2:\quad
 M_{1,k}&=\miuraq_{2,k+3}\miuraq_{2,k+2}\miuraq_{2,k+1},&
 M_{2,k}&=\miuraq_{1,k}.
\end{align*}
For example, in $E_6$ the successive root indices for $i=3$ are
$3,5,3$. The factors contributing to $M_{3,k}$ occur at $(k+2,4)$
after $(k+2,3)$, at $(k+1,4)$ adjacent to $5$, and at $(k,1)$
before $(k,3)$. In $D_4$ the indices $1,4,3$ are all adjacent to
$2$, giving one factor from each of $A_{k+3}^{\mathrm M},A_{k+2}^{\mathrm M},A_{k+1}^{\mathrm M}$.

All the multipliers $\mathsf R_{i,k}$ are units. For each starting
index $0\le\rho<o_i$, the recurrence with $d_{i,\rho}=1$ gives
\[
 d_{i,\rho+n o_i}=
 \begin{cases}
   \displaystyle\prod_{s=0}^{n-1}\mathsf R_{i,\rho+s o_i},& n>0,\\[4pt]
   1,& n=0,\\[2pt]
   \displaystyle\prod_{s=n}^{-1}\mathsf R_{i,\rho+s o_i}^{-1},& n<0 .
 \end{cases}
\]
Each $k\in\mathbb Z$ determines its starting index $\rho$ by the
unique expression $k=\rho+n o_i$ with $0\le\rho<o_i$. The associated
vector field on the space of Miura sequences differentiates a Laurent
polynomial $F$ by
\[
 \mathscr D_{i,\rho}F=
   \sum_{n\in\mathbb Z}d_{i,\rho+n o_i}\,\miuraq_{i,\rho+n o_i}
         \frac{\partial F}{\partial \miuraq_{i,\rho+n o_i}} .
\]
Each Laurent polynomial involves only finitely many coordinates, so
this formula defines a derivation of $\mathcal Y_\Lambda$. At every
$k=\rho+n o_i$, its logarithmic derivative
$\mathscr D_{i,\rho}(\miuraq_{i,k})/\miuraq_{i,k}$ is precisely $d_{i,k}$.

\begin{lemma}\label{tw:screen-identification}
The restrictions of the quantum screenings to $\mathcal Y_\Lambda$
are the derivations
\begin{equation}\label{tw:screen-derivation}
 S_i(\miuraq_{j,k})=\delta_{ij}\miuraq_{i,k}s_{i,k},
\end{equation}
where $s_{i,k}$ denotes the screening symbol $s_{i,a\qaff^{2k+\xi_i}}$.
These symbols generate the target over $\mathcal Y_\Lambda$ with
relations \eqref{tw:screen-recurrence}. In the basis
$s_{i,0},\ldots,s_{i,o_i-1}$ the screening is given by
\[
 S_i(F)=\sum_{\rho=0}^{o_i-1}(\mathscr D_{i,\rho}F)s_{i,\rho}.
\]
For each $(i,\rho)$, the first-order substitution
$F\mapsto F+\varepsilon\mathscr D_{i,\rho}F$ is induced by a twisted
gauge transformation \eqref{tw:gauge} whose parameters are Laurent
polynomials in the Miura variables.
\end{lemma}

\begin{proof}
For $e=2k+\xi_i+o_i$, equation \eqref{tw:dictionary} gives
\[
 \begin{aligned}
 \miuraq_{i,k}\miuraq_{i,k+o_i}
   &=Y_{i,a\qaff^{e+o_i}}Y_{i,a\qaff^{e-o_i}},\\
 \miuraq_{j,k+s}
   &=Y_{j,a\qaff^{e+\delta}},\qquad
     \delta=\xi_j-\xi_i-o_i+2s .
 \end{aligned}
\]
The first line gives the numerator of $\mathsf R_i(a\qaff^e)$.
To compare the denominators, we use the second line for every
factor $\miuraq_{j,k+s}$ of $M_{i,k}$. Thus $s$ is the displacement
from $k$ in the Miura product and $\delta$ is the displacement from
$e$ in the exponent of $\qaff$.
On an edge with $o_i=o_j=d$, the pair $(s,\xi_j-\xi_i)$ equals
$(d,-d)$ when $j$ follows $i$ and $(0,d)$ when $j$ precedes $i$.
In either case $\delta=0$.
For a multiple edge, the table lists the pairs $(s,\delta)$ in order
for each choice of its endpoints $i$ and $j$.
\[
\begin{array}{c|c|c|c@{\qquad}c|c|c}
 & (i,j)&s&\delta&(i,j)&s&\delta\\\hline
 B_n&(n-1,n)&1,2&-1,1&(n,n-1)&0&0\\
 C_n&(n,n-1)&0,1&-1,1&(n-1,n)&1&0\\
 F_4&(3,4)&1,2&-1,1&(4,3)&0&0\\
 G_2&(1,2)&1,2,3&-2,0,2&(2,1)&0&0
\end{array}
\]
Thus the values of $\delta$ for $C_{ji}=-1,-2,-3$ are
$\{0\},\{-1,1\},\{-2,0,2\}$, respectively. Comparison with
\eqref{eq:quantum-root} gives
\begin{equation}\label{tw:root-identification}
 \mathsf R_{i,k}=\mathsf R_i(a\qaff^{2k+\xi_i+o_i}) .
\end{equation}
With $b=a\qaff^{e-o_i}$ and
$\qaff_i=\qaff^{o_i}$, relation \eqref{eq:quantum-screening-module}
becomes
\[
 s_{i,b\qaff_i^2}=\mathsf R_i(b\qaff_i)s_{i,b}
 \quad\Longleftrightarrow\quad
 s_{i,k+o_i}=\mathsf R_{i,k}s_{i,k}.
\]
Since $\qaff$ has infinite multiplicative order, for $n\in\mathbb Z$ we have
\[
 a\qaff^{2k+\xi_i}=(a\qaff^{2h+\xi_i})\qaff_i^{2n}
 \quad\Longleftrightarrow\quad k=h+n o_i .
\]

The indices $0\le\rho<o_i$ thus represent distinct orbits of spectral
parameters under multiplication by $\qaff_i^2$. The defining relations
of the screening module relate symbols only within a single orbit, so
$s_{i,0},\ldots,s_{i,o_i-1}$ are independent. For $k=\rho+no_i$
the recurrence has invertible multipliers and gives
$s_{i,k}=d_{i,k}s_{i,\rho}$, since $d_{i,\rho}=1$ and
$d_{i,k+o_i}=\mathsf R_{i,k}d_{i,k}$. These symbols therefore also
span the target on the lattice. For $F\in\mathcal Y_\Lambda$,
\[
 S_iF=\sum_{k\in\mathbb Z}\miuraq_{i,k}
               \frac{\partial F}{\partial\miuraq_{i,k}}s_{i,k}
 =\sum_{\rho=0}^{o_i-1}(\mathscr D_{i,\rho}F)s_{i,\rho}
\]
by \eqref{eq:quantum-screening} and the Leibniz rule.

For the component $\mathscr D_{i,\rho}$, the required gauge
transformation over $\mathcal Y_\Lambda[\varepsilon]/(\varepsilon^2)$
has components $v_k=x_{h_k}(\varepsilon r_k)$, where
\[
 h_k=\sigma^{k-\rho-1}(i),\qquad
 r_{\rho+1}=\miuraq_{i,\rho}\prod_{j\prec i}\miuraq_{j,\rho}^{A^G_{ji}},
 \qquad
 r_{k+1}=r_k\prod_{j\in J}\miuraq_{j,k}^{A^G_{j,h_{k+1}}}.
\]
Here $j\prec i$ refers to the order in $J$. The unit multipliers
determine $r_k$ in both directions. By \eqref{ade:rank-one-distinct}
and \eqref{ade:rank-one-same}, moving $v_{k+1}$ through the Miura product
leaves on the right
\[
 x_{h_{k+1}}\!\left(\varepsilon r_{k+1}\prod_{j\in J}
          \miuraq_{j,k}^{-A^G_{j,h_{k+1}}}\right)
     =x_{\sigma(h_k)}(\varepsilon r_k)=\sigma(v_k).
\]
At $k=\rho$, the parameter just before $P_i$ is
$\varepsilon r_{\rho+1}\prod_{j\prec i}
\miuraq_{j,\rho}^{-A^G_{ji}}=\varepsilon\miuraq_{i,\rho}$,
so \eqref{ade:rank-one-same} replaces $\miuraq_{i,\rho}$ by
$\miuraq_{i,\rho}(1+\varepsilon)$. The same factorization at successive
occurrences of $i$ gives \eqref{tw:screen-recurrence}, hence the relative
variation of $\miuraq_{i,k}$ is $\varepsilon d_{i,k}$ at
$k=\rho+no_i$. At every other index $h_{k+1}$ lies in the same
$\sigma$-orbit but differs from its sole representative $i\in J$.
No factor then changes its Miura coordinate. This gives
$F\mapsto F+\varepsilon\mathscr D_{i,\rho}F$ on the Laurent ring.
\end{proof}

\begin{corollary}\label{tw:screen-invariance}
Every pullback $\mu_{c,\sigma}^*T_{i,k}$ belongs to
$\bigcap_{j\in J}\ker S_j$.
\end{corollary}
\begin{proof}
We work over $\mathcal Y_\Lambda[\varepsilon]/(\varepsilon^2)$ with the
sequence $A^{\mathrm M}$ whose Miura coordinates are the Laurent generators. For each pair
$(j,\rho)$, Lemma~\ref{tw:screen-identification} gives the variation
$\widehat A^{\mathrm M}=v\cdot A^{\mathrm M}$. If the gauge transformation
$u$ sends $A^{\mathrm M}$ to the Steinberg section, then $uv^{-1}$ sends
$\widehat A^{\mathrm M}$ to the same point of the section at every index.
Uniqueness in Proposition~\ref{tw:normalization} therefore gives, for
$F=\mu_{c,\sigma}^*T_{i,k}$,
\[
 \mu_{c,\sigma}(\widehat A^{\mathrm M})=\mu_{c,\sigma}(A^{\mathrm M}),
 \qquad F(\widehat A^{\mathrm M})-F(A^{\mathrm M})
       =\varepsilon(\mathscr D_{j,\rho}F)(A^{\mathrm M})=0 .
\]
Comparison of the coefficients of $\varepsilon$ gives
$\mathscr D_{j,\rho}F=0$ for every $\rho$, and the decomposition
of $S_j$ in Lemma~\ref{tw:screen-identification} gives $S_jF=0$
for every $j\in J$.
\end{proof}

\subsection{The fundamental \texorpdfstring{$q$}{q}-characters}

\sectionoverview{In this subsection, we prove that the twisted Miura pullback sends each
Steinberg coordinate to a fundamental $q$-character and that its image
is the common screening kernel on the spectral lattice. The remaining
point is to determine the dominant monomials of each pullback.}

\begin{lemma}\label{tw:dominance}
The monomial $\miuraq_{i,k}$ occurs in $\mu_{c,\sigma}^*T_{i,k}$ with coefficient one. Every other monomial has a negative exponent.
\end{lemma}
\begin{proof}
For each $j=i_p\in J$ the inversion root
$\beta_j=(1-c)\omega_j=s_{i_1}\cdots s_{i_{p-1}}\alpha_j$
belongs to $R_c$. For $1\le b\le m$, the positive-root partition in
Lemma~\ref{tw:blocks} and the identity $\gamma^m=-1$ give
\[
 \gamma^{m-b}\beta_j\in\Delta_+,\qquad
 \gamma^{-b}\beta_j=\gamma^{-m}\gamma^{m-b}\beta_j
                   =-\gamma^{m-b}\beta_j\in\Delta_- .
\]
Proposition~\ref{tw:normalization} shows that the change of frame $u_{k+1}$
depends regularly only on the Miura variables
at $k+1,\ldots,k+m$. This dependence and the negativity just proved
are the two hypotheses of Lemma~\ref{lem:weight-criterion}. It follows that $\miuraq_{i,k}$ has coefficient
one and is the only dominant monomial of $\mu_{c,\sigma}^*T_{i,k}$.
\end{proof}

\begin{theorem}\label{tw:character-theorem}
For the ordered orbit representatives $J$ in \eqref{tw:data} and the
spectral variables in \eqref{tw:dictionary}, the pullback satisfies
\begin{equation}\label{tw:raw-character}
 \mu_{c,\sigma}^*(T_{i,k})=\chi_{\qaff}\bigl(L(Y_{i,a\qaff^{2k+\xi_i}})\bigr)
 \qquad(i\in J,\ k\in\mathbb Z).
\end{equation}
It is injective, with image in $\mathcal Y_\Lambda$ given by
\[
 \mu_{c,\sigma}^*\C[\mathcal A_c^{\Z}]
   =\bigcap_{i\in J}\ker(S_i|_{\mathcal Y_\Lambda})
   =\chi_{\qaff}\bigl(\C\otimes_{\Z}K_0(\mathcal C_\Lambda)\bigr).
\]
\end{theorem}
\begin{proof}
By Corollary~\ref{tw:screen-invariance}, every screening annihilates
$\mu_{c,\sigma}^*T_{i,k}$. Lemma~\ref{tw:dominance} identifies its
only dominant monomial as $Y_{i,b}$ with coefficient one, where
$b=a\qaff^{2k+\xi_i}$. By Lemma~\ref{lem:screened-uniqueness} these
dominant coefficients determine the element of the common kernel.
It is therefore $\chi_{\qaff}(L(Y_{i,b}))$, proving \eqref{tw:raw-character}.

In particular all fundamental characters indexed by $\Lambda$ lie in
$\mathcal Y_\Lambda$. Lemma~\ref{lem:spectral-subring} now shows that
they freely generate both the restricted screening kernel and the
complexified character ring of $\mathcal C_\Lambda$. The pullback maps the
polynomial generators $T_{i,k}$ of $\mathbb C[\mathcal A_c^{\mathbb Z}]$
to these free generators, so it is an isomorphism onto the kernel.
Proposition~\ref{tw:normalization} proves the reduction assertion of
Theorem~\ref{thm:intro} for these four types.
Corollary~\ref{ade:normalization} and Theorem~\ref{ade:main}
give its simply laced cases.
\end{proof}

\subsection{Steinberg coordinates in lowest-weight minors}

\sectionoverview{In this subsection, we express the coordinates of the twisted Miura
transformation in terms of the lowest-weight principal minors of the
elements in Miura form. Duality for the fundamental modules determines the spectral parameter
attached to each such minor.}

For the input $(A_n^{\mathrm M})$ and its representative $(A_n)=\mu_{c,\sigma}(A^{\mathrm M})$,
the Steinberg coordinate of the map from $n$ to $n+1$
is $T_{i,n}(A)=\Delta_i(A_n)$.
To apply \eqref{eq:fundamental-spectral-duality}, we determine the
shift $H_{\mathfrak a}=r^\vee h^\vee$ for the four Lie algebras $\mathfrak a$.
If $\theta$ is the highest root and $\rho$ is half the
sum of the positive roots of $\mathfrak a$, then
$h^\vee=1+\langle\rho,\theta^\vee\rangle$. In the order of
\eqref{tw:data}, the pairs $(r^\vee,h^\vee)$ are
$(2,2n-1),(2,n+1),(2,9),(3,4)$ \cite[Ch.~4]{Kac}, giving the
products $4n-2,2n+2,18,12$. Thus $H_{\mathfrak a}=2m$.
The opposition involution $j\mapsto j^*$ of $G$ is characterized
by $w_0\omega_{j^*}=-\omega_j$. For $j\in J$ the corresponding
lowest-weight minor of $A_s^{\mathrm M}\in B_-$ is
\[
 y_{j,s}:=\Delta_{w_0\omega_{j^*},w_0\omega_{j^*}}(A_s^{\mathrm M})
        =\omega_j([A_s^{\mathrm M}]_0)^{-1}=\miuraq_{j,s}^{-1}.
\]
Here $j^*$ is a node of $G$ and may lie outside $J$.
The opposition involution $\iota$ in \eqref{eq:fundamental-spectral-duality}
belongs to the quantum Lie algebra $\mathfrak a$. It fixes every node
in \eqref{tw:data} because the longest element acts as $-1$
\cite[Ch.~4]{Kac}. Indeed, the Weyl groups of types $B$ and $C$
contain every coordinate sign change. In $F_4$ the product of
reflections in the four coordinate roots is $-1$, and in $G_2$
the third power of a Coxeter element is rotation by $\pi$.

\begin{corollary}
\label{tw:forward-character}
With $a$ and $\xi$ as in
\eqref{tw:dictionary}, one has
\begin{equation}\label{tw:forward-formula}
 T_{i,n}(A)
  =\left.\chi_{\qaff}\!\left(L(Y_{i,a\qaff^{-2n-\xi_i}})\right)
    \right|_{Y_{j,a\qaff^{H_{\mathfrak a}-2s-\xi_j}}=y_{j,s}
               \ (j\in J,\ s\in\Z)} .
\end{equation}
Increasing $n$ or $s$ by one multiplies the corresponding spectral
parameter by $\qaff^{-2}$.
\end{corollary}

\begin{proof}
For $p=2n+\xi_i$, the involution $\mathcal D$ in
\eqref{eq:fundamental-spectral-duality} sends the fundamental
character in \eqref{tw:raw-character} to the one with spectral
parameter $a\qaff^{-2n-\xi_i}$. Since $\iota(j)=j$, it sends
each Laurent generator $Y_{j,a\qaff^{2s+\xi_j}}$ to
$Y_{j,a\qaff^{H_{\mathfrak a}-2s-\xi_j}}^{-1}$. Under the substitution in \eqref{tw:forward-formula}, the image of
this generator evaluates to $y_{j,s}^{-1}=\miuraq_{j,s}$. Thus the transformed character evaluates to the same
coordinate $T_{i,n}(A)$.
\end{proof}

\subsection{The character homomorphism on all spectral parameters}
\label{sec:all-spectral-parameters}

\sectionoverview{In this subsection, we combine the character formulas of
Theorems~\ref{ade:main} and~\ref{tw:character-theorem} to recover
the $q$-character homomorphism on the full Grothendieck ring.
Each spectral lattice contributes one tensor factor.}

With $K$ indexing the simple roots of $\mathfrak a$, the lattice with base parameter $a\in\mathbb C^\times$ is
\[
 \Lambda_a=\{(i,a\qaff^{2k+\xi_i}):i\in K,\ k\in\mathbb Z\}.
\]
Two such lattices agree if their base parameters differ by a factor
$\qaff^{2n}$ for some $n\in\mathbb Z$, and otherwise they are disjoint.
If $E$ contains one parameter from each orbit under multiplication by
$\qaff^2$, then every pair $(i,b)\in K\times\mathbb C^\times$
belongs to exactly one $\Lambda_a$ with $a\in E$. The integer $k$
giving $b=a\qaff^{2k+\xi_i}$ is unique because $\qaff$ has infinite
multiplicative order. By
Theorems~\ref{ade:main} and~\ref{tw:character-theorem}, the
$q$-character of $L(Y_{i,b})$ belongs to $\mathcal Y_{\Lambda_a}$.
The induction on composition multiplicities in the proof of
Lemma~\ref{lem:spectral-subring} expresses every simple class of
$\mathcal C_{\Lambda_a}$ as an integer polynomial in the
fundamental classes at those parameters. By
Theorem~\ref{thm:qcharacter-facts}, all fundamental classes freely
generate $K_0(\mathcal C)$ over $\mathbb Z$. Partitioning these
generators by their spectral lattices gives the first of the
following multiplication isomorphisms. Partitioning the Laurent
variables in the same way gives the second.
\[
 K_0(\mathcal C)\simeq
   \bigotimes_{a\in E} K_0(\mathcal C_{\Lambda_a}),
 \qquad
 \mathcal Y\simeq\bigotimes_{a\in E}\mathcal Y_{\Lambda_a} .
\]
The tensor products are over $\mathbb Z$ and $\mathbb C$, respectively.
Each is the union of the tensor products over finite subsets of $E$,
with the map to a larger tensor product given by inserting units.
The second
decomposition also holds over $\mathbb Z$ for the Laurent rings with
integer coefficients. Under these integral decompositions, the
$q$-character homomorphism is the tensor product of its restrictions,
since both maps agree on every fundamental class.

\section{The two fundamental \texorpdfstring{$q$}{q}-characters in type \texorpdfstring{$B_2$}{B2}}
\label{tw:b2-section}

\sectionoverview{In this section, we compute the $q$-characters of the two fundamental
$U_{\qaff}(L\mathfrak{so}_5)$-modules by a twisted unipotent gauge
transformation in $\mathrm{SL}_4$.
The two Steinberg coordinates give Laurent polynomials with five and
four terms, respectively. For the four-term polynomial we verify
directly that both screening operators annihilate it.}

For $\mathfrak a=\mathfrak{so}_5$ of type $B_2$, \eqref{tw:data} gives
$G=\mathrm{SL}_4$ with its standard upper and lower triangular
Borel subgroups. The involution of the $A_3$ diagram sends $i$ to
$4-i$. Its lift to the pinned group is
\[
 J_0=\begin{pmatrix}
 0&0&0&1\\0&0&-1&0\\0&1&0&0\\-1&0&0&0
 \end{pmatrix},\qquad
 \sigma(g)=J_0g^{-\mathsf T}J_0^{-1},
\]
where $g^{-\mathsf T}=(g^{-1})^{\mathsf T}$. Here $E_{ab}$ denotes the matrix with its only nonzero entry equal to one in row $a$ and column $b$. Then $J_0^{\mathsf T}=-J_0$ gives $\sigma^2=1$, and
$J_0(-E_{i+1,i})J_0^{-1}=E_{4-i,5-i}$ gives
$\sigma(x_i(z))=x_{4-i}(z)$ for $x_i(z)=I+zE_{i,i+1}$.
The corresponding calculation on the negative simple-root subgroup
gives $\sigma(y_i(z))=y_{4-i}(z)$.

The ordered orbit representatives $J=(1,2)$ determine $c=s_1s_2$.
At index $k$ the element in Miura form is
\[
 A_k^{\mathrm M}=P_1(\miuraq_{1,k})P_2(\miuraq_{2,k})=
 \begin{pmatrix}
 \miuraq_{1,k}&0&0&0\\1&\miuraq_{2,k}/\miuraq_{1,k}&0&0\\0&1&\miuraq_{2,k}^{-1}&0\\0&0&0&1
 \end{pmatrix},
\]
and the matrices in the Steinberg section for the lift $\dot s_1\dot s_2$ are
\[
 a(T,V)=X_1(T)X_2(V)=
 \begin{pmatrix}
 T&-V&1&0\\1&0&0&0\\0&1&0&0\\0&0&0&1
 \end{pmatrix}.
\]

\begin{proposition}\label{tw:b2-formulas}
The Steinberg representative of $(A_k^{\mathrm M})$ is the
sequence $(a(T_k,V_k))$. Its coordinates $T_k=T_{1,k}$ and
$V_k=T_{2,k}$ are
\begin{align}
 T_k={}&\miuraq_{1,k}+\frac{\miuraq_{2,k+1}\miuraq_{2,k+2}}{\miuraq_{1,k+2}}
       +\frac{\miuraq_{2,k+1}}{\miuraq_{2,k+3}}
       +\frac{\miuraq_{1,k+1}}{\miuraq_{2,k+3}\miuraq_{2,k+2}}
       +\frac1{\miuraq_{1,k+3}},\label{tw:b2-five}\\
 V_k={}&\miuraq_{2,k}+\frac{\miuraq_{1,k}}{\miuraq_{2,k+1}}
       +\frac{\miuraq_{2,k+2}}{\miuraq_{1,k+2}}
       +\frac1{\miuraq_{2,k+3}}.\label{tw:b2-four}
\end{align}
\end{proposition}

\begin{proof}
By Proposition~\ref{tw:normalization}, a unique sequence of
upper unitriangular matrices $u_k$ sends $(A_k^{\mathrm M})$ to the
Steinberg section. In terms of the entries
\[
 u_k=\begin{pmatrix}
 1&a_k&b_k&c_k\\0&1&d_k&e_k\\0&0&1&f_k\\0&0&0&1
 \end{pmatrix},
\]
the equation $u_{k+1}A_k^{\mathrm M}\sigma(u_k)^{-1}=a(T_k,V_k)$ reads
\begin{equation}\label{tw:b2-certificate}
 u_{k+1}A_k^{\mathrm M} J_0u_k^{\mathsf T}J_0^{-1}=a(T_k,V_k).
\end{equation}
The two factors on the left are
\begin{align*}
 u_{k+1}A_k^{\mathrm M}&=\begin{pmatrix}
 \miuraq_{1,k}+a_{k+1}&a_{k+1}\miuraq_{2,k}/\miuraq_{1,k}+b_{k+1}&b_{k+1}/\miuraq_{2,k}&c_{k+1}\\
 1&\miuraq_{2,k}/\miuraq_{1,k}+d_{k+1}&d_{k+1}/\miuraq_{2,k}&e_{k+1}\\
 0&1&1/\miuraq_{2,k}&f_{k+1}\\0&0&0&1
 \end{pmatrix},\\[3pt]
 J_0u_k^{\mathsf T}J_0^{-1}
 &=\begin{pmatrix}
 1&-f_k&e_k&-c_k\\
 0&1&-d_k&b_k\\
 0&0&1&-a_k\\0&0&0&1
 \end{pmatrix}.
\end{align*}
The third row of $a(T_k,V_k)$ is $(0,1,0,0)$, and the second is
$(1,0,0,0)$. Entries $(3,3)$ and $(3,4)$ in the product determine
$d_k$ and $a_k$. With these values, entries $(2,2)$ and $(2,3)$
determine $f_k$ and $e_k$, respectively. This gives
\[
\begin{aligned}
 d_k&=\miuraq_{2,k}^{-1},
 & a_k&=\miuraq_{2,k}(b_k+f_{k+1}),\\
 f_k&=\miuraq_{2,k}/\miuraq_{1,k}+\miuraq_{2,k+1}^{-1},
 & e_k&=f_kd_k-d_{k+1}/\miuraq_{2,k}=\miuraq_{1,k}^{-1}.
\end{aligned}
\]
The $(2,4)$ entry gives
\[
 c_k=e_{k+1}-d_{k+1}f_{k+1}+\frac{\miuraq_{2,k}b_k}{\miuraq_{1,k}}
     =-\frac1{\miuraq_{2,k+1}\miuraq_{2,k+2}}+\frac{\miuraq_{2,k}b_k}{\miuraq_{1,k}} .
\]
In the first row the first two entries are $T_k$ and $-V_k$, hence
\[
\begin{aligned}
 T_k&=\miuraq_{1,k}+a_{k+1}
     =\miuraq_{1,k}+\miuraq_{2,k+1}(b_{k+1}+f_{k+2}),\\
 V_k&=T_kf_k-a_{k+1}\miuraq_{2,k}/\miuraq_{1,k}-b_{k+1}
     =\miuraq_{2,k}+\miuraq_{1,k}/\miuraq_{2,k+1}+f_{k+2}.
\end{aligned}
\]
Substitution for $f_{k+2}$ gives \eqref{tw:b2-four}. The third
entry of the first row is
\[
 T_k e_k-(a_{k+1}\miuraq_{2,k}/\miuraq_{1,k}+b_{k+1})d_k
       +b_{k+1}/\miuraq_{2,k}
   =(T_k-a_{k+1})/\miuraq_{1,k}=1 .
\]
The formula for $T_k$ still contains $b_{k+1}$. To determine this
entry of $u_{k+1}$, we use the vanishing of the $(1,4)$ entry of
$a(T_k,V_k)$. Substitution for $a_k$ and $c_k$ gives
\[
 c_{k+1}=
    -\frac{T_k}{\miuraq_{2,k+1}\miuraq_{2,k+2}}+
        b_{k+1}f_{k+1}+\miuraq_{2,k}b_k .
\]
Using the expressions for $c_{k+1},T_k,f_{k+1}$ cancels $b_{k+1}$, leaving
\[
 \miuraq_{2,k}b_k=
     \frac{\miuraq_{1,k}}{\miuraq_{2,k+1}\miuraq_{2,k+2}}
     +\frac{f_{k+2}}{\miuraq_{2,k+2}}
     -\frac{1}{\miuraq_{2,k+2}\miuraq_{2,k+3}}
   = \frac{\miuraq_{1,k}}{\miuraq_{2,k+1}\miuraq_{2,k+2}}+\frac{1}{\miuraq_{1,k+2}}.
\]
Replacing $k$ by $k+1$ in this identity and substituting into
$T_k=\miuraq_{1,k}+\miuraq_{2,k+1}(b_{k+1}+f_{k+2})$ gives \eqref{tw:b2-five}.
\end{proof}

To identify the two coordinates with characters, we fix a spectral
parameter $a\in\mathbb C^\times$. For $B_2$, the substitution
\eqref{tw:dictionary}, the Cartan matrix, and the two orbit sizes are
\[
 \miuraq_{1,k}=Y_{1,a\qaff^{2k}},\qquad \miuraq_{2,k}=Y_{2,a\qaff^{2k-1}},\qquad
 C=\begin{pmatrix}2&-1\\-2&2\end{pmatrix},\quad(o_1,o_2)=(2,1).
\]
The terms $\miuraq_{1,k}$ of $T_k$ and $\miuraq_{2,k}$ of $V_k$ thus have spectral
parameters $b=a\qaff^{2k}$ and $b=a\qaff^{2k-1}$, respectively.
Expressing \eqref{tw:b2-five} and
\eqref{tw:b2-four} in these variables gives, by
Theorem~\ref{tw:character-theorem},
\begin{align}
 \chi_{\qaff}(L(Y_{1,b}))={}&Y_{1,b}
  +\frac{Y_{2,b\qaff}Y_{2,b\qaff^3}}{Y_{1,b\qaff^4}}
  +\frac{Y_{2,b\qaff}}{Y_{2,b\qaff^5}}
  +\frac{Y_{1,b\qaff^2}}{Y_{2,b\qaff^3}Y_{2,b\qaff^5}}
  +\frac1{Y_{1,b\qaff^6}},\label{tw:b2-character-one}\\
 \chi_{\qaff}(L(Y_{2,b}))={}&Y_{2,b}
  +\frac{Y_{1,b\qaff}}{Y_{2,b\qaff^2}}
  +\frac{Y_{2,b\qaff^4}}{Y_{1,b\qaff^5}}
  +\frac1{Y_{2,b\qaff^6}}.\label{tw:b2-character-two}
\end{align}
At $b=\qaff^{-3}$ and $b=\qaff^{-4}$, respectively, these formulas
agree with \cite[Example~3.9]{HL}. Setting every Laurent variable
$Y_{j,b'}$ equal to one sums the dimensions of the $\ell$-weight
spaces and gives the two module dimensions, $5$ and $4$.

The four terms of $V_k$ also display the cancellation under both
screenings. The root monomials from \eqref{tw:screen-recurrence} are
\[
 \mathsf R_{1,k}=\frac{\miuraq_{1,k}\miuraq_{1,k+2}}{\miuraq_{2,k+1}\miuraq_{2,k+2}},\qquad
 \mathsf R_{2,k}=\frac{\miuraq_{2,k}\miuraq_{2,k+1}}{\miuraq_{1,k}}.
\]
By Lemma~\ref{tw:screen-identification}, the generators of the two
screening modules satisfy
\[
 s_{1,k+2}=\mathsf R_{1,k}s_{1,k},\qquad
 s_{2,k+1}=\mathsf R_{2,k}s_{2,k}.
\]
The first recurrence expresses $s_{1,2n}$ in terms of $s_{1,0}$ and
$s_{1,2n+1}$ in terms of $s_{1,1}$. These give the two components of
$S_1$. The second recurrence expresses every $s_{2,k}$ in terms of
$s_{2,0}$, so $S_2$ has one component. Writing $N_1,N_2,N_3,N_4$
for the monomials in \eqref{tw:b2-four} in their displayed order gives
\[
 \frac{N_1}{N_2}=\mathsf R_{2,k},\qquad
 \frac{N_2}{N_3}=\mathsf R_{1,k},\qquad
 \frac{N_3}{N_4}=\mathsf R_{2,k+2}.
\]
Since $N_1$ and $N_4$ contain only variables $\miuraq_{2,s}$,
$S_1(N_1)=S_1(N_4)=0$. The other contributions cancel in pairs by
\eqref{tw:screen-recurrence},
\begin{align*}
 S_1(N_2+N_3)
   &=N_2s_{1,k}-N_3s_{1,k+2}
     =(N_2-N_3\mathsf R_{1,k})s_{1,k}=0,\\
 S_2(N_1+N_2)
   &=N_1s_{2,k}-N_2s_{2,k+1}
     =(N_1-N_2\mathsf R_{2,k})s_{2,k}=0,\\
 S_2(N_3+N_4)
   &=N_3s_{2,k+2}-N_4s_{2,k+3}
     =(N_3-N_4\mathsf R_{2,k+2})s_{2,k+2}=0 .
\end{align*}
Hence $S_1(V_k)=S_2(V_k)=0$.

\section{Discrete universal Wronskians and the common screening kernel}
\label{sec:wronskian-descent}

\sectionoverview{In this section we prove directly that the common screening
kernel on the spectral lattice is the image of the discrete Miura pullback.
We vary the preserved flag over a fixed Steinberg sequence and prove that
every Laurent polynomial annihilated by the screenings descends to a
regular function of that sequence.}

Theorems~\ref{ade:main} and~\ref{tw:character-theorem} identify the
pulled-back Steinberg coordinates with the fundamental $q$-characters
by the Frenkel--Mukhin criterion for dominant monomials
(Lemma~\ref{lem:screened-uniqueness}). The proof of descent below uses
only the action on the initial Wronskian, the Miura coordinates
expressed as ratios of minors, and regularity on the Steinberg
section. It shows that every Laurent polynomial killed by the screenings
is a Miura pullback. Combined with the character identification, this
proves the Frenkel--Mukhin description of the $q$-character ring on
the chosen lattice \cite[Theorem~5.1]{FM}.

\subsection{Discrete universal Wronskians and their initial values}

\sectionoverview{In this subsection, we define the discrete universal
Wronskian for a sequence $(A_k)$ in Steinberg form and relate it to
a group-valued solution of $\Xi_{k+1}=A_k\sigma(\Xi_k)$. Varying the
initial value then identifies the regular and rational functions that
depend only on $(A_k)$.}

We consider the ordinary Miura transformation for simply laced $G$ and
the four twisted transformations of \eqref{tw:data}. In the first case
$\sigma=1$ and $J=(i_1,\ldots,i_\ell)$ is any ordering of the simple roots.
In the four twisted cases $J$ lists the representatives of the $\sigma$-orbits
in the order given in \eqref{tw:data}.
In both cases $c=s_{i_1}\cdots s_{i_\ell}$, and we write
$\mathcal A=\mathcal A_c$ for the Steinberg section. The domain of
$\mu_{c,\sigma}:(M_c^-)^{\mathbb Z}\to\mathcal A^{\mathbb Z}$ has
coordinate ring $\mathbb C[\miuraq_{i,k}^{\pm1}:i\in J,\ k\in\mathbb Z]$.
We identify this ring with $\mathcal Y_\Lambda$ by
\eqref{ade:spectral-dictionary} in simply laced type and by
\eqref{tw:dictionary} in twisted type.

\begin{definition}\label{def:wronskian-frame}
Let $D$ be a commutative complex algebra. For $Z_0\in H(D)$ and a
sequence $(A_k)$ in $\mathcal A(D)$, a \emph{discrete
universal Wronskian} is a sequence $(W_k)_{k\in\mathbb Z}$ in $G(D)$
such that
\begin{equation}\label{eq:discrete-universal-wronskian}
 W_{k+1}A_k=Z_0\sigma(W_k),\qquad k\in\mathbb Z.
\end{equation}
We call it \emph{normalized at $k=0$} if $W_0=1$.
\end{definition}

At lattice index $k$, left multiplication by $W_k$ changes from
coordinates in which the connection sends $g$ to $A_k\sigma(g)$
to coordinates in which it sends $g$ to $Z_0\sigma(g)$. Indeed,
$W_{k+1}A_k\sigma(g)=Z_0\sigma(W_kg)$. For any initial value
$W_0\in G(D)$, the recurrence $W_{k+1}=Z_0\sigma(W_k)A_k^{-1}$
and the inverse recurrence
$W_k=\sigma^{-1}(Z_0^{-1}W_{k+1}A_k)$ determine the whole sequence.

To express the solution in the original frame, we first remove the torus
factor from the Wronskian equation. The recurrence $R_{k+1}=Z_0\sigma(R_k)$ with $R_0=1$ determines $(R_k)$ in $H(D)$. The sequence
$\Xi_k=W_k^{-1}R_k$ satisfies
\begin{equation}\label{eq:wronskian-frame}
 \Xi_{k+1}=A_k\sigma(\Xi_k),\qquad k\in\mathbb Z,
\end{equation}
because $W_{k+1}^{-1}Z_0=A_k\sigma(W_k)^{-1}$. Conversely,
$W_k=R_k\Xi_k^{-1}$ recovers the Wronskian from any solution, and
$W_0=\Xi_0^{-1}$.

The space $\mathscr X=\mathcal A^{\mathbb Z}$ of sequences in the Steinberg section has coordinate ring $R=\mathbb C[T_{i,k}:i\in J,\ k\in\mathbb Z]$. Over $\mathscr X$ the scheme $\mathscr W_\sigma$ parametrizes pairs
$((A_k),(\Xi_k))$ satisfying \eqref{eq:wronskian-frame}. For fixed $Z_0\in H(\mathbb C)$, the map
$\Xi_k\mapsto R_k\Xi_k^{-1}$ identifies these pairs over $D$ with
sequences $(A_k)$ in $\mathcal A(D)$ equipped with a discrete universal
Wronskian for $Z_0$.

For any $(A_k)\in\mathscr X$, the solution $(F_k)$ with $F_0=1$ is
determined by $F_{k+1}=A_k\sigma(F_k)$. Every solution with
initial value $\Xi_0$ satisfies
\[
 \Xi_k=F_k\sigma^k(\Xi_0).
\]
For $k>0$ the recurrence uses only $A_0,\ldots,A_{k-1}$, and for
$k<0$ the inverse recurrence $F_k=\sigma^{-1}(A_k^{-1}F_{k+1})$
uses only $A_k,\ldots,A_{-1}$. Hence every $\Xi_k$ depends regularly
on the initial value and finitely many entries of $(A_k)$. Conversely,
$(\Xi_k)$ recovers the initial value $\Xi_0$ and each
$A_k=\Xi_{k+1}\sigma(\Xi_k)^{-1}$ by regular maps. This gives
\begin{equation}\label{eq:frame-product}
 \mathscr W_\sigma\simeq G\times\mathscr X,\qquad
 B:=\mathbb C[\mathscr W_\sigma]=\mathbb C[G]\otimes R .
\end{equation}
Under this isomorphism, right multiplication by $h\in G(D)$ on the
first factor replaces $\Xi_0$ by $\Xi_0h$. Since
$\Xi_k=F_k\sigma^k(\Xi_0)$, the induced action on solutions is
\begin{equation}\label{eq:frame-right-action}
  (\Xi_k)\cdot h=(\Xi_k\sigma^k(h)).
\end{equation}
The identity $W_k=R_k\Xi_k^{-1}$ then gives the action on Wronskians
as $W_k\mapsto R_k\sigma^k(h)^{-1}R_k^{-1}W_k$.
Both regular and rational invariants of this action are determined
by the sequence $(A_k)$. More precisely,
\begin{equation}\label{eq:frame-invariants}
 B^G=R,\qquad (\Frac B)^G=\Frac R .
\end{equation}
For $f\in B^G$, translation gives $f(g,A)=f(1,A)$, which belongs
to $R$. For rational invariants we work over the field of fractions
$K=\Frac R$. Over the
generic point of $\mathscr X$, the product \eqref{eq:frame-product}
has fiber $G_K$, the group $G$ over the field $K$. An invariant
$f\in(\Frac B)^G$ thus gives a rational function on $G_K$.
In $K(G_K\times G_K)$, invariance
and the substitution $(g,h)\mapsto(g,g^{-1}h)$ give
\[
 f(gh)=f(g)\quad\Longrightarrow\quad f(h)=f(g).
\]
The open subset $N_-HN_+$ of $G_K$ is isomorphic to the product
$N_-\times H\times N_+$ over $K$. Each $N_\pm$ is an affine space
and $H$ is a product of copies of the multiplicative group. Since $K$
is infinite, the $K$-points of this open subset are dense. Specializing $g$
to one of these points at which $f$ is regular gives
$f\in K$. Both reverse inclusions follow because the action fixes
$(A_k)$.

Over a meromorphic function field, the choice of initial value in
\eqref{eq:frame-product} must be compatible with the spectral shift.
Suppose $A(u)$ is a meromorphic connection in Steinberg form and
$Z_0\in H(\C)$ is prescribed. A meromorphic Wronskian for
this connection satisfies
\[
 W(\eta u) A(u)\sigma(W(u))^{-1}=Z_0 .
\]
This is exactly the equation for a meromorphic change of frame to the
constant connection $g\mapsto Z_0\sigma(g)$ \cite{QDE,QDEII}. Restriction
to an orbit where $A$ and $W$ are defined gives
\eqref{eq:discrete-universal-wronskian}. If $W'$
is another meromorphic Wronskian for the same $A$ and $Z_0$, then
$K=W'W^{-1}$ satisfies $K(\eta u)=Z_0\sigma(K(u))Z_0^{-1}$.
Conversely, each meromorphic $K$ satisfying this equation gives the
Wronskian $KW$. Thus, once $W$ is given, the possible values of $W'$
at $u_*$ are determined by the solutions $K$ on the chosen domain.
To write $\Xi(u)=W(u)^{-1}R(u)$ as in \eqref{eq:wronskian-frame}, one
also needs a meromorphic $H$-valued solution of
$R(\eta u)=Z_0\sigma(R(u))$ on that domain.

The Steinberg condition also characterizes $W_k$ among representatives
of the cosets $W_kN_+$ under right multiplication by $N_+$. Indeed, suppose $(g_k)$ is a sequence
in $G(D)$ such that $a_k=g_{k+1}^{-1}Z_0\sigma(g_k)$
belongs to $N_+\dot c N_+$ for every $k$. Replacing $g_k$ by $g_kn_k$, with
$n_k\in N_+(D)$, replaces $a_k$ by $n_{k+1}^{-1}a_k\sigma(n_k)$.
Theorem~\ref{thm:ade-full-cell} in simply laced type and
Proposition~\ref{tw:normalization} in the twisted cases give a unique
sequence $(n_k)$ for which the transformed $a_k$ lie in $\mathcal A(D)$.
Thus the cosets $g_kN_+$ have unique representatives $W_k=g_kn_k$
for which $W_{k+1}^{-1}Z_0\sigma(W_k)$ lies in the Steinberg section.

When the $A_k$ lie in the Steinberg section, the group identity
\eqref{eq:discrete-universal-wronskian} implies the Wronskian
equations on highest-weight vectors. Indeed, each $A_k$ has the form
$n\dot c$ with $n\in U_c\subset N_+$. For the highest-weight
vector $v_i\in V(\omega_i)$ we therefore have
$A_k^{-1}v_i=\dot c^{-1}v_i$, and
\[
 W_{k+1}v_i=Z_0\sigma(W_k)\dot c^{-1}v_i .
\]
For $\sigma=1$ this is the equation of
\cite[Definition~4.9 and (4.9)]{KZ} on the chosen spectral orbit,
with constant lift $\dot c$.

When $\sigma=1$, the relation between $W_k$ and $\Xi_k$ gives a
comparison with the sequence used in
\cite[Definition~1.1]{FLbands}. In the convention of that paper,
a $(G,c)$-band is a sequence $g(s)\in G$ whose successive quotients
$g(s)g(s+1)^{-1}$ belong to $\mathcal A_c$. With the same lift
$\dot c$, setting $g(s)=\Xi_{-s}=W_{-s}^{-1}Z_0^{-s}$ gives
$g(s)g(s+1)^{-1}=A_{-s-1}$. Conversely these quotients determine
the elements $A_k$, and $W_k=Z_0^k g(-k)^{-1}$ recovers
the Wronskian.

For $G=\SL_N$, $c=s_1\cdots s_{N-1}$ and $\sigma=1$, there is also a
Wronskian formed directly from scalar solutions. If
$\mathbf f_k=(f^{(1)}_k,\ldots,f^{(N)}_k)$ is the first row of $\Xi_k$,
the lower rows of the companion matrix in Proposition~\ref{ade:companion}
give
\begin{equation}\label{eq:casorati-frame}
 \Xi_k=\begin{pmatrix}
   \mathbf f_k\\ \mathbf f_{k-1}\\ \vdots\\ \mathbf f_{k-N+1}
 \end{pmatrix},\qquad \det \Xi_k=1 .
\end{equation}
The first row of $\Xi_{k+1}=A_k\Xi_k$ gives
\[
 \mathbf f_{k+1}=\sum_{i=1}^{N-1}(-1)^{i-1}T_{i,k}\mathbf f_{k+1-i}
             +(-1)^{N-1}\mathbf f_{k+1-N}.
\]
Conversely, arranging the values of any $N$ scalar solutions in the
rows displayed in \eqref{eq:casorati-frame} gives a matrix $M_k$
satisfying $M_{k+1}=A_kM_k$. The recurrence gives the first row of
$M_{k+1}$, and each lower row is the preceding row of $M_k$. Its
determinant, the Casorati determinant or discrete Wronskian, is
constant because $\det A_k=1$. If it is a unit $d$,
multiplying one of the scalar solutions by $d^{-1}$ gives a solution
in $\mathrm{SL}_N(D)$.
For $f^{(j)}_k=f^{(j)}(u^{(k)})$ and $u^{(k)}=q^ku_*$, the entry
in row $r$ and column $j$ is $f^{(j)}(q^{1-r}u^{(k)})$. In this sense
$\Xi_k$ is the $q$-Wronskian matrix of the scalar solutions, while the
change of trivialization in \eqref{eq:discrete-universal-wronskian}
is $W_k=R_k\Xi_k^{-1}$.

\subsection{Gauss factorization and Wronskian minors}

\sectionoverview{In this subsection, we recover the Miura variables from a solution by
factoring its values as products in $N_+B_-$. We express each variable as a ratio of principal minors and
determine which choices of solution give the same Miura variables.}

For a solution $(\Xi_k)$ and an index $k$, Proposition~\ref{prop:gauss-cell}
gives a factorization $\Xi_k=u_kb_k$, with $u_k\in N_+(D)$ and
$b_k\in B_-(D)$, precisely when $\Delta_j(\Xi_k^{-1})$ is a unit
in $D$ for every simple root $j$. For fixed $j$ and $k$, this
principal minor is a regular function on the space $\mathscr W_\sigma$
of pairs $((A_k),(\Xi_k))$. The multiplicative subset $S\subset B$ consists of all finite products of the functions $\Delta_j(\Xi_k^{-1})$, including the empty product.
We define $\mathscr W_\sigma^\circ=\operatorname{Spec}S^{-1}B$.
By the universal property of localization, a $D$-point of
$\mathscr W_\sigma$ factors through $\mathscr W_\sigma^\circ$ precisely
when $\Xi_k$ admits this factorization for every $k$. Each function
in $S^{-1}B$ has a denominator involving only finitely many minors.

On $\mathscr W_\sigma^\circ$ the same minors determine the fundamental
weight coordinates of the torus factor $[b_k]_0$. Indeed,
$\Xi_k^{-1}=b_k^{-1}u_k^{-1}$ and $u_k^{-1}$ fixes the highest-weight
vector, so $\Delta_j(\Xi_k^{-1})=\omega_j([b_k]_0)^{-1}$. Thus
\begin{equation}\label{eq:cartan-wronskians}
 Q_{j,k}=\omega_j([b_k]_0)=\Delta_j(\Xi_k^{-1})^{-1}.
\end{equation}
Substitution of $\Xi_k=u_kb_k$ and $\Xi_{k+1}=u_{k+1}b_{k+1}$ in
\eqref{eq:wronskian-frame} gives
\begin{equation}\label{eq:wronskian-miura}
 A_k^{\mathrm M}=b_{k+1}\sigma(b_k)^{-1}
    =u_{k+1}^{-1}A_k\sigma(u_k)\in B_-\cap N_+\dot c N_+=M_c^-.
\end{equation}
The Cartan projection $B_-\to H$ is a homomorphism commuting with
$\sigma$, and $\omega_j(\sigma(h))=\omega_{\sigma^{-1}(j)}(h)$
for $h\in H$. Hence, for every simple root $j$ of $G$,
\[
 \frac{Q_{j,k+1}}{Q_{\sigma^{-1}(j),k}}
 =\omega_j([A_k^{\mathrm M}]_0)
 =\omega_j\!\left(\prod_{p\in J}\alpha_p^\vee(\miuraq_{p,k})\right)
 =\prod_{p\in J}\miuraq_{p,k}^{\langle\omega_j,\alpha_p^\vee\rangle}
 =\begin{cases}\miuraq_{j,k},&j\in J,\\
                1,&j\notin J.\end{cases}
\]
In particular, for each $i\in J$,
\begin{equation}\label{eq:miura-Q-ratio}
 \miuraq_{i,k}=\frac{Q_{i,k+1}}{Q_{\sigma^{-1}(i),k}}.
\end{equation}
Since $\Xi_k^{-1}=R_k^{-1}W_k$ and $R_k\in H(D)$,
\eqref{eq:cartan-wronskians} gives $Q_{j,k}=\omega_j(R_k)/\Delta_j(W_k)$.
Using $R_{k+1}=Z_0\sigma(R_k)$ in \eqref{eq:miura-Q-ratio} yields
the formula in terms of the Wronskian itself,
\[
 \miuraq_{i,k}=\omega_i(Z_0)
          \frac{\Delta_{\sigma^{-1}(i)}(W_k)}
               {\Delta_i(W_{k+1})}.
\]
If $i$ has orbit size $o_i$, all $\sigma^{-p}(i)$ with $1\le p<o_i$
lie outside $J$. At these labels the ratio is one, so
$Q_{\sigma^{-p}(i),s}=Q_{\sigma^{-p-1}(i),s-1}$. Iterating from
$p=1$ to $o_i-1$ in the denominator of \eqref{eq:miura-Q-ratio}
gives
\begin{equation}\label{eq:orbit-cartan-quotient}
 Q_{\sigma^{-1}(i),k}=Q_{i,k+1-o_i},\qquad
 \miuraq_{i,k}=\frac{Q_{i,k+1}}{Q_{i,k+1-o_i}}.
\end{equation}
This includes $o_i=1$, when $\sigma^{-1}(i)=i$ and no iteration is needed.

The sequence $(A_k^{\mathrm M})$ in Miura form determines the solution
once the lower Borel factor $b_0$ at index zero is prescribed. Indeed, for
$(b_0,(A_k^{\mathrm M}))\in B_-(D)\times M_c^-(D)^{\mathbb Z}$,
Corollary~\ref{ade:normalization} or Proposition~\ref{tw:normalization}
gives the unique sequence $(u_k)$ in $N_+(D)$ such that
$A_k=u_{k+1}A_k^{\mathrm M}\sigma(u_k)^{-1}\in\mathcal A(D)$ for every $k$. The
recursion $b_{k+1}=A_k^{\mathrm M}\sigma(b_k)$ determines $(b_k)$ in both directions, and
\[
 \Xi_{k+1}=u_{k+1}b_{k+1}=u_{k+1}A_k^{\mathrm M}\sigma(b_k)
            =A_k\sigma(u_kb_k)=A_k\sigma(\Xi_k).
\]
Gauss factorization of the constructed solution recovers $u_k$ and
$b_k$, hence also $b_0$ and $A_k^{\mathrm M}$. In the other direction,
if $\Xi_k=u_kb_k$ is given, equation \eqref{eq:wronskian-miura} and
uniqueness of the unipotent transformation to the section recover $u_k$ from
$(A_k^{\mathrm M})$. The recurrence with initial value $b_0$ recovers
$b_k$ and hence $\Xi_k$.
Both constructions are regular at each index and yield the isomorphism
\begin{equation}\label{eq:wronskian-gauss-product}
 \mathscr W_\sigma^\circ\simeq B_-\times (M_c^-)^{\mathbb Z},\qquad
 (\Xi_k)\longmapsto (b_0,(b_{k+1}\sigma(b_k)^{-1})_k).
\end{equation}
For $h\in B_-$ the action in \eqref{eq:frame-right-action} sends
$b_k$ to $b'_k=b_k\sigma^k(h)$ and leaves $u_k$ fixed. Since
$b'_{k+1}\sigma(b'_k)^{-1}=A_k^{\mathrm M}$ and $b'_0=b_0h$, it follows that
\[
 (S^{-1}B)^{B_-}
   =\mathbb C[(M_c^-)^{\mathbb Z}]
   =\mathcal Y_\Lambda=\mathbb C[\miuraq_{i,k}^{\pm1}].
\]
Projection to the Miura sequence in \eqref{eq:wronskian-gauss-product} defines the
map $p$, while projection to $\mathscr X$ in \eqref{eq:frame-product}
defines $\pi$. Since the Miura transformation
takes $(A_k^{\mathrm M})$ to $(A_k)$ by \eqref{eq:wronskian-miura},
the square
\[
 \begin{array}{ccc}
 \mathscr W_\sigma^\circ & \xrightarrow{\ p\ } & (M_c^-)^{\mathbb Z}\vphantom{\Big|}\\
 \big\downarrow && \big\downarrow\rlap{$\scriptstyle\,\mu_{c,\sigma}$}\\
 \mathscr W_\sigma & \xrightarrow{\ \pi\ } & \mathscr X
 \end{array}
\]
commutes. The left arrow is the morphism induced by localization at
the principal minors. On coordinate rings, the equality
$\pi|_{\mathscr W_\sigma^\circ}=\mu_{c,\sigma}\circ p$ identifies
$p^*\mu_{c,\sigma}^*$ with the composite $R\hookrightarrow B\to S^{-1}B$.
For each fixed $(A_\ell)$ and $k$, the recurrence $\Xi_{\ell+1}=A_\ell\sigma(\Xi_\ell)$ has a solution in both directions with $\Xi_k=1$. At this solution each function $\Delta_j(\Xi_k^{-1})$
has value $\Delta_j(1)=1$. Thus for each $j,k$ the corresponding
minor is a nonzero element of $B$. Since $B$ is a polynomial ring
over the domain $\mathbb C[G]$ by \eqref{eq:frame-product}, every
finite product of these minors is nonzero. Thus $B\to S^{-1}B$ is
injective, and the square gives the injectivity of the Miura pullback
$\mu_{c,\sigma}^*$.

For $G=\SL_2$ the only minor needed to factor $\Xi_k$ in $N_+B_-$
can be read from its second column. If $(r_k)$ and $(s_k)$ are the
scalar solutions giving the two columns of $\Xi_k$, then
\[
 \Xi_k=\begin{pmatrix}r_k&s_k\\r_{k-1}&s_{k-1}\end{pmatrix},
 \qquad
 A_k=\begin{pmatrix}T_k&-1\\1&0\end{pmatrix}.
\]
The upper left entry of $\Xi_k^{-1}$ is $s_{k-1}$, since $\det\Xi_k=1$.
Hence the factorization exists at every index precisely when all $s_k$
are units. Writing $Q_k=Q_{1,k}$ and $T_k=T_{1,k}$ in this case,
\eqref{eq:cartan-wronskians},
\eqref{eq:miura-Q-ratio}, and $s_{k+1}=T_ks_k-s_{k-1}$ give
\[
 Q_k=s_{k-1}^{-1},\qquad \miuraq_{1,k}=s_{k-1}/s_k,\qquad
 T_k=\miuraq_{1,k}+\miuraq_{1,k+1}^{-1}.
\]
Right multiplication by $h\in\mathrm{SL}_2(D)$ preserves
$A_k=(\Xi_{k+1}h)(\Xi_kh)^{-1}$. For $h\in B_-(D)$,
its action on the second column is $s_k\mapsto h_{22}s_k$.
Hence $\miuraq_{1,k}$ is preserved as well.

\subsection{Screenings and variation of the initial value}

\sectionoverview{In this subsection, we differentiate the Miura
variables as the initial value varies in a positive simple-root
direction. We show that each such derivative is a unit multiple of
one component of a quantum screening.}

For a fixed sequence $A=(A_k)$ in the Steinberg section and a simple-root index
$j$ of $G$, the right action \eqref{eq:frame-right-action} replaces
the initial value $\Xi_0$ by $\Xi_0x_j(z)$. Differentiating at $z=0$
defines a derivation $\partial_j$ on $B=\mathbb C[\mathscr W_\sigma]$,
\[
 (\partial_j\Phi)(\Xi_0,A)=
 \left.\frac{d}{dz}\Phi(\Xi_0x_j(z),A)\right|_{z=0}.
\]
The Miura variables are regular on $\mathscr W_\sigma^\circ$, and
the quotient rule extends $\partial_j$ to its coordinate ring $S^{-1}B$.

\begin{lemma}\label{lem:frame-screenings}
For $i\in J$ and $0\le\rho<o_i$, the simple root
$j=\sigma^{-1-\rho}(i)$ corresponds to the component $\mathscr D_{i,\rho}$
of $S_i$, with $\mathscr D_{i,0}=\mathscr D_i$ when $\sigma=1$. The
derivative along $x_j(z)$ satisfies
\begin{equation}\label{eq:frame-screening-component}
 \partial_j(p^*F)=\lambda_{i,\rho}\,p^*(\mathscr D_{i,\rho}F),\qquad
 \lambda_{i,\rho}=\frac{\partial_j(p^*\miuraq_{i,\rho})}{p^*\miuraq_{i,\rho}}
 \in (S^{-1}B)^\times
\end{equation}
for every $F\in\mathcal Y_\Lambda$. Consequently, a Laurent
polynomial $F\in\mathcal Y_\Lambda$ satisfies $S_i F=0$
for every $i\in J$ if and only if $\partial_j(p^*F)=0$ for
every simple-root index $j$ of $G$.
\end{lemma}

\begin{proof}
We compute $\partial_j$ by changing the initial value to
$\Xi_0x_j(\epsilon)$ over $D=(S^{-1}B)[\epsilon]/(\epsilon^2)$.
For $j=\sigma^{-1-\rho}(i)$ and $h_k=\sigma^k(j)$,
\eqref{eq:frame-right-action} gives $\Xi'_k=\Xi_kx_{h_k}(\epsilon)$.
Each $\Delta_l((\Xi'_k)^{-1})$ is a unit in $D$ because it is a unit
modulo $\epsilon$. Thus $\Xi'_k$ has Gauss factors reducing to
$u_k,b_k$ at $\epsilon=0$. They have the form $u_k(1+\epsilon X_k)$
and $(1+\epsilon Y_k)b_k$, where $1+\epsilon X$ denotes the tangent
element at the identity with derivative $X$. Comparing their product
with $u_kb_kx_{h_k}(\epsilon)$ gives
\[
 X_k+Y_k=\Ad(b_k)e_{h_k},\qquad
 X_k\in\mathfrak n_+\otimes S^{-1}B,\quad
 Y_k\in\operatorname{Lie}(B_-)\otimes S^{-1}B.
\]
For $b_k=n_{-,k}a_k$ with $a_k=[b_k]_0$, each bracket of $e_{h_k}$
with a negative-root vector lowers its height from one. Hence
$\Ad(n_{-,k})e_{h_k}-e_{h_k}\in\operatorname{Lie}(B_-)\otimes S^{-1}B$.
Since $a_k$ acts on $e_{h_k}$ by the scalar $\alpha_{h_k}(a_k)$,
projection along $\mathfrak g=\mathfrak n_+\oplus\operatorname{Lie}(B_-)$
gives
\[
 X_k=\theta_k e_{h_k},\qquad \theta_k=\alpha_{h_k}(a_k)\in(S^{-1}B)^\times .
\]
Thus $u'_k=u_kx_{h_k}(\epsilon\theta_k)$ and
\eqref{eq:wronskian-miura} becomes
\[
 \widehat A_k^{\mathrm M}=(u'_{k+1})^{-1}A_k\sigma(u'_k)
 =x_{h_{k+1}}(-\epsilon\theta_{k+1})A_k^{\mathrm M}
   x_{h_{k+1}}(\epsilon\theta_k).
\]

We identify each Miura coordinate with its pullback under $p$.
The Miura product contains a factor $P_i$ for the sole representative
$i\in J$ of the $\sigma$-orbit of $h_{k+1}=\sigma^{k-\rho}(i)$.
Thus a factor can have simple root $h_{k+1}$ only when
$k\in\rho+m\mathbb Z$, where $m=o_i$. At such an index, moving the
left-hand root element past a factor $P_p$ preceding $P_i$ multiplies
its parameter by $\miuraq_{p,k}^{-A^G_{pi}}$ by
\eqref{ade:rank-one-distinct}. At $P_i$ the resulting parameter is
added to $\miuraq_{i,k}$ by \eqref{ade:rank-one-same}. All other Miura
coordinates are fixed. Once the left-hand element has passed through
the product, the two upper factors on the right cancel. Indeed, their
product is in $N_+$, and both $\widehat A_k^{\mathrm M}$ and the
product of the varied Miura factors lie in $B_-$. Since
$N_+\cap B_-=\{1\}$, the relative derivative is
\[
 \frac{\partial_j\miuraq_{l,k}}{\miuraq_{l,k}}=
 \begin{cases}
 \lambda_k:=-\theta_{k+1}\miuraq_{i,k}^{-1}
       \displaystyle\prod_{p\prec i}\miuraq_{p,k}^{-A^G_{pi}},
       & l=i,\quad k\in\rho+m\mathbb Z,\\
 0,&\text{otherwise}.
 \end{cases}
\]
Here $p\prec i$ is the order in $J$ and
$A^G_{ph}=\langle\alpha_h,\alpha_p^\vee\rangle$.
The Cartan projection of $b_{r+1}=A_r^{\mathrm M}\sigma(b_r)$ gives
$a_{r+1}=[A_r^{\mathrm M}]_0\sigma(a_r)$. Evaluating $\alpha_{h_{r+1}}$
and using $h_{r+1}=\sigma(h_r)$, we obtain
\[
 \frac{\theta_{r+1}}{\theta_r}
 =\alpha_{h_{r+1}}([A_r^{\mathrm M}]_0)
 =\prod_{p\in J}\miuraq_{p,r}^{A^G_{p,h_{r+1}}}.
\]
For $k\in\rho+m\mathbb Z$, multiplication from $r=k+1$ to $k+m$
and cancellation of the factors with $p\prec i$ at $k+m$ yield
\[
\begin{aligned}
 \frac{\lambda_{k+m}}{\lambda_k}
 &=\miuraq_{i,k}\miuraq_{i,k+m}
   \prod_{p\prec i}\miuraq_{p,k}^{A^G_{pi}}
   \prod_{i\prec p}\miuraq_{p,k+m}^{A^G_{pi}}
   \\[-2pt]
 &\qquad\times\prod_{r=k+1}^{k+m-1}\prod_{p\in J}
         \miuraq_{p,r}^{A^G_{p,\sigma^{r-k}(i)}}
 =\mathsf R_{i,k}.
\end{aligned}
\]
Here the exponent of $\miuraq_{i,k+m}$ is $A^G_{ii}-1=1$.
The product is $\mathsf R_{i,k}$ by \eqref{ade:geometric-root} for
$\sigma=1$ and by \eqref{tw:screen-recurrence} for $\sigma\ne1$.
Thus $\lambda_k$ and
$\mathscr D_{i,\rho}(\miuraq_{i,k})/\miuraq_{i,k}$ have the same
successive ratios on $\rho+m\mathbb Z$. At $k=\rho$ the latter equals
one, whereas $\lambda_\rho$ is a unit by its formula above. Hence
\eqref{eq:frame-screening-component} holds on every Laurent generator
with $\lambda_{i,\rho}=\lambda_\rho$, and the Leibniz rule gives it
on the whole ring. Since $p^*$ is injective by
\eqref{eq:wronskian-gauss-product}, we obtain
\[
 S_iF=0\ \text{for all }i
 \quad\Longleftrightarrow\quad \mathscr D_{i,\rho}F=0\ \text{for all }i,\rho
 \quad\Longleftrightarrow\quad \partial_j(p^*F)=0\ \text{for all }j,
\]
where the last equivalence uses the bijection
$(i,\rho)\mapsto\sigma^{-1-\rho}(i)$ with the simple roots of $G$.
\end{proof}

The factor $\lambda_{i,\rho}$ matters when interpreting the
screening as a vector field. On
$B$, $\partial_j$ is locally nilpotent because the translate of
each regular function by $x_j(z)$ is a polynomial in $z$. On the Miura
coordinate ring, $\mathscr D_{i,\rho}^{\,n}(\miuraq_{i,\rho})=\miuraq_{i,\rho}$
for every $n\ge1$, which excludes a regular additive-group action
with this infinitesimal generator.

\subsection{Regularity on the Steinberg section}

\sectionoverview{In this subsection, we prove that every Laurent polynomial annihilated
by the screenings descends to a regular function of the sequence
$(A_k)$ in the Steinberg section. We first show that it is independent of the initial value
and then exclude a pole along each irreducible divisor.}

We use the generation theorem for a split semisimple group. The
rank-one subgroups $\varphi_j(\mathrm{SL}_2)$ at its simple roots
generate $G$ as an algebraic group \cite[Proposition~21.62]{Milne}.
Since $\mathrm{SL}_2$ is generated by its upper and lower unipotent
subgroups, the subgroups $x_j(\mathbb G_a)$ together with $B_-$
also generate $G$. This applies after extension of the base field.

\begin{theorem}\label{thm:regular-descent}
For every Coxeter element in simply laced type and every pair
$(c,\sigma)$ in \eqref{tw:data},
\begin{equation}\label{eq:regular-descent}
 \mu_{c,\sigma}^*\mathbb C[\mathcal A^{\mathbb Z}]
    =\bigcap_{i\in J}\ker(S_i|_{\mathcal Y_\Lambda}).
\end{equation}
\end{theorem}

\begin{proof}
For $P\in R$, the pullback $p^*\mu_{c,\sigma}^*P=\pi^*P$ depends only
on the sequence $(A_k)$ in the Steinberg section. Right multiplication
on $\Xi_0$ fixes this sequence, so $\partial_j(p^*\mu_{c,\sigma}^*P)=0$
for every $j$. Hence
Lemma~\ref{lem:frame-screenings} gives $\mu_{c,\sigma}^*R\subseteq\bigcap_i\ker S_i$.

Conversely, for $F\in\bigcap_i\ker S_i$, the function
$\widetilde F=p^*F\in (S^{-1}B)^{B_-}$ satisfies
$\partial_j\widetilde F=0$ for every $j$ by the same lemma.
Over $K=\Frac R$ we regard $\widetilde F$ as a rational function of
the initial value $g$ on $G_K$. Differentiation of $x_j(s)x_j(t)=x_j(s+t)$
in $t$ at $t=0$ gives
\[
 \frac{d}{ds}\widetilde F(gx_j(s))=(\partial_j\widetilde F)(gx_j(s))=0
 \quad\text{in }K(G_K)(s).
\]
The constants of $d/ds$ are $K(G_K)$ in characteristic zero. Evaluation
at $s=0$ gives $\widetilde F(gx_j(s))=\widetilde F(g)$. Since $B_-$ and the
$x_j(\mathbb G_a)$ generate $G$, \eqref{eq:frame-invariants} yields
$\widetilde F\in(\Frac B)^G\cap S^{-1}B=K\cap S^{-1}B$.

It remains to show that this rational function of the Steinberg
coordinates is polynomial, that is, $K\cap S^{-1}B=R$. The ring $R$ is factorial because
the prime factors of a polynomial in finitely many Steinberg coordinates
remain prime on adjoining the other coordinates. Each element of
$K\cap S^{-1}B$ therefore has the two expressions
\[
 \frac{P}{Q}=\frac{b}{d},\qquad
 P,Q\in R,\ \gcd(P,Q)=1,\quad b\in B,\quad d\in S.
\]
Here $Q$ is a polynomial in the Steinberg coordinates, while $d$ is
a finite product of the principal minors $\Delta_j(\Xi_k^{-1})$.
Equality in $\operatorname{Frac}B$ gives an identity in $B$,
\begin{equation}\label{eq:descent-denominators}
 Pd=Qb .
\end{equation}
Suppose an irreducible $f\in R$ divides $Q$. At the generic point of
the divisor $f=0$ the sequence $(A_k)$ is defined over
$K_f=\Frac(R/(f))$ and the initial value $\Xi_0$ still varies over
$G_{K_f}$ by \eqref{eq:frame-product}. The image of $Q$ is zero and the
image of $P$ is nonzero because $f$ does not divide $P$. With bars
denoting images after extension from $R$ to $K_f$,
\eqref{eq:descent-denominators} gives
\[
 B\otimes_R K_f=K_f[G],\qquad \overline P\,\overline d=0,
 \qquad \overline P\in K_f^\times .
\]
On this fiber $\Xi_k=F_k\sigma^k(\Xi_0)$ with $F_k\in G(K_f)$
given by the solution with $F_0=1$. For each $k$ the map
$g\mapsto (F_k\sigma^k(g))^{-1}$ is an automorphism of the variety
$G_{K_f}$. The restriction of $\Delta_j(\Xi_k^{-1})$ is its pullback
of $\Delta_j$, which is nonzero since $\Delta_j(1)=1$.
Their finite product $\overline d$ is therefore nonzero
in the domain $K_f[G]$, since $G$ is geometrically integral.
This contradicts $\overline P\,\overline d=0$, so $Q$ is a unit and
$K\cap S^{-1}B=R$. Consequently
\[
 p^*F=\pi^*P_0=p^*\mu_{c,\sigma}^*P_0\qquad(P_0\in R).
\]
Injectivity of $p^*$ gives $F=\mu_{c,\sigma}^*P_0$, and injectivity
of $\mu_{c,\sigma}^*$ gives uniqueness.
\end{proof}

\section{B\"acklund transformations and the Weyl group symmetry of \texorpdfstring{$q$}{q}-characters}
\label{sec:fh-weyl}

\sectionoverview{In this section we lift the Frenkel--Hernandez
reflections to discrete universal Wronskians and give a uniform
geometric proof of the Weyl symmetry on the spectral lattice of
Theorem~\ref{thm:intro}. The pinned representatives of the simple
reflections act on the Wronskians and induce a Weyl group action on
torus-invariant rational functions.
By comparing this action with the reflection formulas in each chamber,
we deduce the Weyl group relations on the completed Laurent rings.
Differentiating along the positive simple-root subgroups, with a common
parameter in each $\sigma$-orbit, realizes the screening operators with
values in the field of rational functions on the space of Wronskians. A Laurent polynomial is fixed by the reflection
at a vertex precisely when the screening at that vertex annihilates it.
For folded characters the reflection formulas follow by restriction.}

\subsection{The group action on discrete universal Wronskians}
\label{fh:action-subsection}

\sectionoverview{Here we describe the action of $G$ on the discrete
universal Wronskians of a fixed Steinberg sequence. Gauss factorization
expresses a change of initial value as a change of Miura form and
determines how two successive changes compose.}

For a sequence $(A_k)$ in the Steinberg section and a fixed $Z_0\in H$,
the solutions and Wronskians of Section~\ref{sec:wronskian-descent}
are related by
\[
 \Xi_{k+1}=A_k\sigma(\Xi_k),\qquad
 R_0=1,\qquad R_{k+1}=Z_0\sigma(R_k),\qquad
 W_k=R_k\Xi_k^{-1}.
\]
For an ordinary connection $\sigma=1$. In the folded construction
the solutions take values in $G_\nu$ and we use $G=G_\nu$ in these
formulas. The groups $H$, $B_-$ and $N_+$ are the maximal torus,
negative Borel subgroup and positive unipotent subgroup of this $G$.
By \eqref{eq:frame-product}, the field of rational functions on the space
of discrete universal Wronskians is
$\mathcal K=\Frac B=\Frac(\mathbb C[G]\otimes R)$, where $R$ is the
polynomial ring in the Steinberg coordinates $T_{i,k}$. Thus an element
of $\mathcal K$ is a rational function of the initial value
$\Xi_0=W_0^{-1}$ and finitely many of the $T_{i,k}$.

At a complex point of $\mathscr W_\sigma$, the factorization $\Xi_k=u_kb_k$ with $u_k\in N_+$
and $b_k\in B_-$ exists precisely when every principal minor
$\Delta_j(\Xi_k^{-1})$ at that index is nonzero. Over a commutative
complex algebra the same minors must be units. The choice $W_0=1$ always
determines a Wronskian, but a minor can vanish at another index.
For example, in $\mathrm{SL}_2$ the choice $A_0=\dot s_1$ gives
$\Xi_1=\dot s_1$ and $\Delta_1(\Xi_1^{-1})=0$. On the locus where
the factorizations exist at every index, \eqref{eq:wronskian-gauss-product}
identifies a Wronskian with its Miura sequence $(A_k^{\mathrm M})$
and its initial Borel factor $b_0$. The remaining Borel factors are
determined by $b_{k+1}=A_k^{\mathrm M}\sigma(b_k)$. For $G_\nu$
the same identification follows from uniqueness of the unipotent
change of frame in the proof of Corollary~\ref{cor:folded-miura}.
We denote the projection to $(A_k^{\mathrm M})$ by $p$.

\begin{proposition}\label{fh:frame-action}
Right multiplication $\Xi_0\mapsto\Xi_0g$ defines a regular action
of $G$ on $\mathscr W_\sigma$. It fixes $(A_k)$ and acts by
\[
 \Xi_k\longmapsto\Xi_k\sigma^k(g),\qquad
 W_k\longmapsto R_k\sigma^k(g)^{-1}R_k^{-1}W_k .
\]
At an index where the original and transformed solutions both admit
Gauss factorization, the factors $v_k(g,b_k)\in N_+$ and $b_k^g\in B_-$
are uniquely specified by
\[
 b_k\sigma^k(g)=v_k(g,b_k)b_k^g .
\]
Where the factorizations exist at $k$ and $k+1$, the transformed
Miura form is
\begin{equation}\label{fh:framed-gauge}
 A_k^{\mathrm M,g}=b_{k+1}^g\sigma(b_k^g)^{-1}
     =v_{k+1}(g,b_{k+1})^{-1}A_k^{\mathrm M}\sigma(v_k(g,b_k)).
\end{equation}
Where the successive factorizations for $g,h\in G$ exist, their factors satisfy
\begin{equation}\label{fh:gauss-composition}
 b_k^{gh}=(b_k^g)^h,\qquad
 v_k(gh,b_k)=v_k(g,b_k)v_k(h,b_k^g).
\end{equation}
The invariant field $\mathcal K^{B_-}$ is the pullback under $p$ of
the field of rational functions in the Miura variables.
\end{proposition}

\begin{proof}
Under \eqref{eq:frame-product}, changing the initial value
$\Xi_0$ to $\Xi_0g$ is right translation on the $G$ factor, a regular
action. The sequence $\Xi_k^g=\Xi_k\sigma^k(g)$ has this initial
value and satisfies the same recurrence as $\Xi_k$. By uniqueness
in both directions it is the transformed sequence, with Wronskian
\[
 W_k^g=R_k(\Xi_k^g)^{-1}
      =R_k\sigma^k(g)^{-1}R_k^{-1}W_k .
\]
Wherever $b_k\sigma^k(g)$ admits Gauss factorization, its positive
factor $v_k=v_k(g,b_k)$ gives $\Xi_k^g=(u_kv_k)b_k^g$. Substitution
at $k$ and $k+1$ into $\Xi_{k+1}^g=A_k\sigma(\Xi_k^g)$ yields
\[
 b_{k+1}^g\sigma(b_k^g)^{-1}
   =(u_{k+1}v_{k+1})^{-1}A_k\sigma(u_kv_k)
   =v_{k+1}^{-1}A_k^{\mathrm M}\sigma(v_k).
\]
For a second element $h\in G$, the equality
\[
 b_k\sigma^k(gh)
 =v_k(g,b_k)b_k^g\sigma^k(h)
 =v_k(g,b_k)v_k(h,b_k^g)(b_k^g)^h
\]
and uniqueness of the factorization give \eqref{fh:gauss-composition}.
For $g\in B_-$ we have $b_k\sigma^k(g)\in B_-$, hence $v_k=1$.
In the coordinates $((A_k^{\mathrm M}),b_0)$ the action sends $b_0$ to
$b_0g$ and fixes $(A_k^{\mathrm M})$. Right translations act transitively on $B_-$, so an invariant rational
function depends only on $(A_k^{\mathrm M})$. Conversely, every
rational function pulled back under $p$ is fixed.
\end{proof}

\begin{remark}\label{fh:meromorphic-action}
On the integer lattice any initial value is allowed. Over a prescribed
field of meromorphic functions, two Wronskians $W(u)$ and $W'(u)$
for the same connection $A(u)$ and twist $Z_0$ have quotient
$L(u)=W'(u)W(u)^{-1}$ satisfying
\[
 L(\eta u)=Z_0\sigma(L(u))Z_0^{-1}.
\]
For $G=\mathrm{SL}_2$ over $\mathbb C(u)$, this equation determines $L$
explicitly when $\sigma=1$, $q=\eta$ has infinite multiplicative order, and
$Z_0=\operatorname{diag}(\zeta,\zeta^{-1})$ with $\zeta^2\notin q^{\mathbb Z}$.
Indeed, a nonzero rational solution of $f(qu)=a f(u)$ can have no zero or
pole in $\mathbb C^\times$, since its zeros and poles are unions of
$q$-orbits. Thus $f(u)=c u^m$ and $a=q^m$. For the diagonal entries
of $L$ this gives constants, and for the off-diagonal entries the
multipliers $\zeta^{\pm2}$ force both entries to vanish. Hence $L$ is
constant and diagonal.

For the constant connection $Z_0$, a preserved rational line other than
a coordinate line would have a nonzero rational slope with multiplier
$\zeta^2$ or $\zeta^{-2}$. The coordinate lines are therefore the only
two rational preserved flags. At either line an infinitesimal change
of slope $\epsilon f(u)$, with $\epsilon^2=0$, gives the same equation
for $f$. Both points of the rational Miura fiber have zero tangent
space. In this example the screening variations of the flag on the
lattice consequently have no nonzero rational realization.
\end{remark}

\subsection{Normalizer representatives and the braid relations}
\label{fh:braid-subsection}

\sectionoverview{Here we prove the braid relations for the changes of
Wronskian induced by representatives of the simple reflections. We
compute their squares, identify the functions on which they act as
involutions, and relate the root subgroups of the fixed group to the
components of each screening.}

We denote by $K$ the group whose pinned representatives act on the
Wronskians. For an ordinary connection, $K=G$ in simply laced type
and $K=G_\nu$ for folded characters.
For a connection with the diagram automorphism $\sigma$ of
\eqref{tw:data}, the solutions take values in the simply laced group
$G$ and we use $K=(G^\sigma)^\circ$. In all cases $\alpha_i^K$ and
$(\alpha_i^K)^\vee$ denote the simple root and coroot of $K$ at $i$.
In the twisted case, $i\in J$ represents an orbit
$O_i=\{\sigma^k(i):k\in\mathbb Z\}$ of simple roots of $G$.
The pinned homomorphism $\mathrm{SL}_2\to K$ at $i$ sends $z$ to
$\prod_{j\in O_i}\varphi_j(z)$. Its factors commute because distinct
vertices in $O_i$ are nonadjacent. As in \eqref{eq:folded-pinning},
$\alpha_i^K$ is the restriction of $\alpha_i$ to $(H^\sigma)^\circ$
and $(\alpha_i^K)^\vee(t)=\prod_{j\in O_i}\alpha_j^\vee(t)$. The formula
for the Cartan matrix $C$ of $\mathfrak a$ in Section~\ref{tw:section}
then gives
\[
 \langle\alpha_j^K,(\alpha_i^K)^\vee\rangle
   =\sum_{v\in O_i}\langle\alpha_j,\alpha_v^\vee\rangle=C_{ji}.
\]
Thus $\operatorname{Lie}K$ is Langlands dual to $\mathfrak a$ for the
twisted construction. For folded characters the root system governing
the group action is that of $G_\nu$ itself.

The representative $n_i\in K$ of the simple reflection at $i$ is the
image of $\bigl(\begin{smallmatrix}0&-1\\1&0\end{smallmatrix}\bigr)$
under the pinned rank-one homomorphism. For a nontrivial $\sigma$ this
gives $n_i=\prod_{j\in O_i}\dot s_j$ in $G$. In every case $n_i$
normalizes $H$ and $\sigma(n_i)=n_i$. On $\mathcal K$ the right
action of $g\in G$ induces the pullback $\rho(g)$ given by
\[
 (\rho(g)f)((A_k),\Xi_0)=f((A_k),\Xi_0g).
\]
In particular $\rho(g)\rho(h)=\rho(gh)$. We denote $\rho(n_i)$ by
$\mathcal B_i$, and the order of the product of the simple reflections
at $i,j$ in $W(K)$ by $m_{ij}$.

\begin{proposition}
\label{fh:frame-braid}
The operators $\mathcal B_i$ satisfy
\[
 \underbrace{\mathcal B_i\mathcal B_j\mathcal B_i\cdots}_{m_{ij}\ \mathrm{factors}}
 =
 \underbrace{\mathcal B_j\mathcal B_i\mathcal B_j\cdots}_{m_{ij}\ \mathrm{factors}}
 \quad(i\ne j),\qquad
 \mathcal B_i^2=\rho\bigl((\alpha_i^K)^\vee(-1)\bigr),\qquad
 \mathcal B_i^4=1 .
\]
At each index $k$, two reduced expressions of the same Weyl group element
give the same product of the successive unipotent factors in
\eqref{fh:gauss-composition}, wherever those factorizations exist.
On the torus-invariant field $\mathcal K^H$ the operators define
an action of $W(K)$. On a rational function of right $H$-weight
$\lambda$, the square $\mathcal B_i^2$ acts by
$(-1)^{\langle\lambda,(\alpha_i^K)^\vee\rangle}$.
\end{proposition}

\begin{proof}
For each reduced expression $w=s_{i_1}\cdots s_{i_p}$ the pinned
representative $n_w=n_{i_1}\cdots n_{i_p}$ depends only on $w$
\cite[\S1.4, (1.7)--(1.9)]{FZ}. In $\mathrm{SL}_2$ the matrix
representing the simple reflection has square $-I_2$. The pinned
homomorphism at $i$ sends this square to $(\alpha_i^K)^\vee(-1)$.
Applying the group action $\rho$ gives
\[
 \mathcal B_{i_1}\cdots\mathcal B_{i_p}=\rho(n_w),\qquad
 \mathcal B_i^2=\rho(n_i^2)=\rho\bigl((\alpha_i^K)^\vee(-1)\bigr),
 \qquad \mathcal B_i^4=1.
\]
The first equality gives the braid relations. By
\eqref{fh:gauss-composition}, the successive unipotent factors at $k$
multiply to $v_k(n_w,b_k)$ for every reduced expression of $w$.

Since $n_i$ normalizes $H$, $\mathcal B_i$ preserves $\mathcal K^H$.
On a function $f$ of right $H$-weight $\lambda$ its square acts by
\[
 \mathcal B_i^2 f
 =\lambda\bigl((\alpha_i^K)^\vee(-1)\bigr) f
 =(-1)^{\langle\lambda,(\alpha_i^K)^\vee\rangle} f.
\]
On $\mathcal K^H$ the weight is zero, so the braid and square relations
give the action of $W(K)$.
\end{proof}

The orders $m_{ij}$ are $2,3,4,6$ for the rank-two root systems
$A_1\times A_1,A_2,B_2,G_2$, respectively. They depend only on the
product of the two off-diagonal Cartan entries, so transposing the
Cartan matrix preserves them. The pinned representatives generate the
Tits extension of $W(K)$ in the torus normalizer of $K$. Its action on
Wronskians is right multiplication on $\Xi_0$.

The distinction between this action and its restriction to torus
invariants is visible in the Wronskian minors. Right multiplication
by $t\in H$ sends $[b_k]_0$ to $[b_k]_0\sigma^k(t)$. Hence
$Q_{j,k}=\omega_j([b_k]_0)$ has right $H$-weight
$\omega_{\sigma^{-k}(j)}$. The numerator and denominator of the
Miura ratio \eqref{eq:miura-Q-ratio} have the same weight, so the
ratio and its Weyl translates belong to $\mathcal K^H$. The Borel
subgroup fixing a function changes under conjugation. If $f$ is $B_-$-invariant, then $\rho(n_i)f$ is invariant under
$n_iB_-n_i^{-1}$, since $bn_i=n_i(n_i^{-1}bn_i)$ for
$b\in n_iB_-n_i^{-1}$.

The same group action yields Lie algebra derivations on $\mathcal K$.
For $X\in\operatorname{Lie}K$, the differential of $\rho$ at the
identity in the direction $X$ is denoted by $D_X$. The pinning gives
Chevalley generators $e_i^K,f_i^K,h_i^K$ with $h_i^K=[e_i^K,f_i^K]$.

\begin{corollary}\label{fh:framed-serre}
For $X\in\operatorname{Lie}K$, the operators on $\mathcal K$ satisfy
\[
 \begin{aligned}
 \mathcal B_iD_X\mathcal B_i^{-1}&=D_{\operatorname{Ad}_{n_i}X},&
 [D_{e_i^K},D_{f_j^K}]&=\delta_{ij}D_{h_i^K},\\
 (\operatorname{ad}D_{e_i^K})^{\,1-\langle\alpha_j^K,(\alpha_i^K)^\vee\rangle}
 D_{e_j^K}&=0 &&(i\ne j).
 \end{aligned}
\]
\end{corollary}
\begin{proof}
Differentiation of $\rho(n_i)\rho(g)\rho(n_i)^{-1}
=\rho(n_ign_i^{-1})$ at $g=1$ in the direction $X$ proves the
conjugation formula. Since $\rho$ is a group homomorphism, its
differential satisfies $[D_X,D_Y]=D_{[X,Y]}$. Applying this identity
to the Chevalley relations and then to the iterated brackets in the
Serre relations gives the other two formulas.
\end{proof}

In the twisted construction, each positive simple-root subgroup of
$G$ gives one component of an ordinary quantum screening. The
corresponding subgroup of $K$ acts at every vertex of one orbit
with the same parameter. The following proposition identifies its derivative
on the pullback of $\mathcal Y_\Lambda$ with the quantum screening operator
by embedding the screening module in $\mathcal K$.

\begin{proposition}
\label{prop:fh-orbit-group}
In the four cases of \eqref{tw:data}, the Weyl group of
$K=(G^\sigma)^\circ$ identifies with $W(G)^\sigma$. Its simple
reflection at $i$, represented in $W(G)$ by $r_i=\prod_{h\in O_i}s_h$,
corresponds to the simple reflection at $i$ in $W(\mathfrak a)$.

For every $i\in J$, the screening module on $\Lambda$ embeds in
$\mathcal K$ by a $\mathcal Y_\Lambda$-linear map $\iota_i$, where
the scalars act on $\mathcal K$ through $p^*$. For every
$F\in\mathcal Y_\Lambda$ this map satisfies
\[
 \left.\frac{d}{dz}\rho\!\left(x_i^\sigma(z)\right)(p^*F)\right|_{z=0}
       =\iota_i(S_iF),\qquad
 x_i^\sigma(z)=\prod_{h\in O_i}x_h(z).
\]
\end{proposition}

\begin{proof}
For $v$ in the $\sigma$-invariant real cocharacter space, the
pairings $\langle\alpha_h,v\rangle$ are equal for $h\in O_i$.
Since the simple roots in $O_i$ are mutually orthogonal, their
reflections commute and
\[
 r_i(v)=v-\sum_{h\in O_i}\langle\alpha_h,v\rangle\alpha_h^\vee
       =v-\langle\alpha_i^K,v\rangle(\alpha_i^K)^\vee .
\]
Thus $r_i$ restricts to the simple reflection of $K$ at $i$.
If $1\ne w\in W(G)^\sigma$, some $w\alpha_h$ is negative. Since
$w$ commutes with $\sigma$, this holds for every $h$ in that orbit
$O_i$. Right multiplication by $s_h$ then decreases length by one
and leaves the images of the other simple roots in $O_i$ unchanged.
Thus $\ell(wr_i)=\ell(w)-|O_i|$.
Induction on length gives $W(G)^\sigma=\langle r_i:i\in J\rangle$.
The sum of the positive coroots is $\sigma$-invariant and regular,
so its stabilizer in $W(G)$ is trivial. Restriction to the invariant
cocharacter space is therefore faithful and identifies $W(G)^\sigma$
with $W(K)$. The Cartan matrices of $K$ and $\mathfrak a$ are transposes,
so their Coxeter relations agree under the stated identification of
simple reflections.

For $0\le\ell<o_i$ and $h=\sigma^{-1-\ell}(i)$, differentiation
along $x_h$ gives the logarithmic derivative
$\lambda_{i,\ell}=\partial_h(p^*\miuraq_{i,\ell})/(p^*\miuraq_{i,\ell})$,
which is a unit by \eqref{eq:frame-screening-component}. Conjugating
$x_h(z)$ by $t\in H$ replaces $z$ by $\alpha_h(t)z$. Differentiation
and the $H$-invariance of $p^*\miuraq_{i,\ell}$ therefore give
\[
 \rho(t)\partial_h\rho(t)^{-1}=\alpha_h(t)\partial_h,\qquad
 \rho(t)\lambda_{i,\ell}=\alpha_h(t)\lambda_{i,\ell}.
\]
For $0\le\ell<o_i$ these characters of $H$ are distinct. Their
eigenvectors $\lambda_{i,\ell}$ are therefore linearly independent over
$\operatorname{Frac}(p^*\mathcal Y_\Lambda)$, which $H$ fixes pointwise.
By Lemma~\ref{tw:screen-identification}, $s_{i,0},\ldots,s_{i,o_i-1}$
form a basis of the screening module. Sending $s_{i,\ell}$ to
$\lambda_{i,\ell}$ and extending linearly through $p^*$ gives the
required injection $\iota_i$. Differentiating the product
$x_i^\sigma(z)=\prod_{h\in O_i}x_h(z)$ and using
\eqref{eq:frame-screening-component} gives
\[
 \left.\frac{d}{dz}\rho(x_i^\sigma(z))(p^*F)\right|_{z=0}
  =\sum_{h\in O_i}\partial_h(p^*F)
  =\sum_{\ell=0}^{o_i-1}\lambda_{i,\ell}p^*(\mathscr D_{i,\ell}F)
  =\iota_i(S_iF).
\]
\end{proof}

These Lie algebra derivations act on the coordinate ring of the
Wronskians and extend to $\mathcal K$ by the quotient rule. A related
construction of Frenkel and Hernandez, reported in \cite[\S8.5]{FLbands},
uses screening derivations to define a Lie algebra action on the
Grothendieck ring of category $\mathcal O$ for a quantum affine Borel
algebra. To compare their Weyl reflections with the action of the representatives
$n_i$, we recall the completions of the Laurent ring on which the
reflections are defined.

\subsection{The Frenkel--Hernandez action on the completed Laurent rings}
\label{fh:source-subsection}

\sectionoverview{In this subsection we recall the Weyl group action of Frenkel and Hernandez on Laurent series in the spectral variables. The choice of Weyl chamber determines which expansion of a first-order difference equation gives the action of a simple reflection.}

For the ordinary $q$-characters of Theorem~\ref{thm:intro}, the spectral
parameter of the Miura variable $\miuraq_{i,k}$ is
$b_{i,k}=a\qaff^{2k+\xi_i}$, with $\qaff$ of infinite multiplicative
order. We identify $Y_{i,b_{i,k}}$ with $\miuraq_{i,k}$ in the Laurent
ring $\mathcal Y_\Lambda$. The simple roots and fundamental weights of $\mathfrak a$ are denoted by
$\alpha_i^{\rm q}$ and $\omega_i^{\rm q}$, respectively, and its Weyl group
by $W$. The identification of $W$ with $W(K)$
matches the simple reflections indexed by the same vertex. The
spectral shift at vertex $i$ uses $\qaff_i=\qaff^{o_i}$, where $o_i$
is the size of the $\sigma$-orbit of $i$ and equals one when $\sigma=1$.

\begin{definition}\label{fh:completion}
The weight of a Laurent monomial $m=\prod_{j,b}Y_{j,b}^{n_{j,b}}$ is
$\mathrm{wt}(m)=\sum_{j,b}n_{j,b}\omega_j^{\rm q}$. For $w\in W$,
the Weyl chamber completion of \cite[Definition~2.1]{FHWeyl} on $\Lambda$
is the ring $\widehat{\mathcal Y}^{\,w}_\Lambda$ of formal sums
$\sum_m a_m m$ of Laurent monomials on $\Lambda$ with complex
coefficients satisfying the following conditions. The weights
$w\,\mathrm{wt}(m)$ for which $a_m\ne0$ lie in a finite union of
sets $\lambda-Q_+$, where $\lambda$ is a weight of $\mathfrak a$ and
$Q_+=\sum_{i\in J}\mathbb Z_{\ge0}\alpha_i^{\rm q}$.
At each weight only finitely many coefficients are nonzero.
\end{definition}

The shift $\tau(Y_{i,b})=Y_{i,b\qaff^2}$ sends $\miuraq_{i,k}$ to
$\miuraq_{i,k+1}$. It preserves weights and acts on every chamber
completion. The same Laurent polynomial gives a finite sum in
each chamber, so $\mathcal Y_\Lambda$ embeds diagonally in
$\prod_{w\in W}\widehat{\mathcal Y}^{\,w}_\Lambda$.

Frenkel and Hernandez construct the Weyl group action and describe its
invariants for all spectral parameters
\cite[Definition~3.1, Lemma~3.2, and Theorems~3.4 and~5.1]{FHWeyl}.
We state their result on $\Lambda$. For our geometric proof of the Weyl
group relations, we use the formulas for the individual reflections, the
unique solution of the defining recurrence in each chamber, and the
continuity of the resulting substitutions
\cite[Lemmas~2.8 and~3.2]{FHWeyl}. We obtain the square and braid
relations from the pinned representatives acting on universal Wronskians.

\begin{theorem}\label{fh:external}
For $i\in J$ and $w\in W$ there is a unique family
$\Sigma^w_{i,x}\in\widehat{\mathcal Y}^{\,w}_\Lambda$, indexed by
$x=b_{i,k}\qaff_i^{-1}$ for $k\in\mathbb Z$, satisfying
\begin{equation}\label{fh:external-rec}
 \Sigma^w_{i,x}=1+\mathsf R_i(x)^{-1}\Sigma^w_{i,x\qaff_i^{-2}},
 \qquad \Sigma^w_{i,x\qaff^2}=\tau\Sigma^w_{i,x}.
\end{equation}
Each $\Sigma^w_{i,x}$ is a unit. The simple reflection $s_i$ acts from
the completion for $w$ to that for $ws_i$ by a continuous isomorphism
$\Theta_i^w$. It fixes $Y_{j,b}$ for $j\ne i$, and for
$b=b_{i,k}$ its action is
\begin{equation}\label{fh:external-reflection}
 \Theta_i^w(Y_{i,b})=
 Y_{i,b}\,\mathsf R_i(b\qaff_i^{-1})^{-1}
 \frac{\Sigma^{ws_i}_{i,b\qaff_i^{-3}}}{\Sigma^{ws_i}_{i,b\qaff_i^{-1}}}.
\end{equation}
These maps commute with $\tau$ and define a $W$-action on the product
of the completions. The invariant elements of the diagonal Laurent
ring $\mathcal Y_\Lambda$ form the image under $\chi_{\qaff}$ of the
complexified Grothendieck ring of $\mathcal C_\Lambda$.
\end{theorem}

The series in \eqref{fh:external-rec} is the solution of
\cite[Lemma~2.8]{FHWeyl}. When $w\alpha_i^{\rm q}$ is positive,
iteration towards $x\qaff_i^{-2},x\qaff_i^{-4},\ldots$ gives terms
of decreasing weight in the chamber. When $w\alpha_i^{\rm q}$ is
negative, solving the recurrence towards $x\qaff_i^{2},x\qaff_i^{4},\ldots$
gives the expansion in that chamber. Both series in the quotient
\eqref{fh:external-reflection} belong to the \emph{target} completion,
indexed by $ws_i$.

On $\Lambda$, all root monomials occurring in \eqref{fh:external-rec}
belong to $\mathcal Y_\Lambda$ by \eqref{ade:spectral-dictionary}
and \eqref{tw:root-identification}. The expansions therefore use only
variables on $\Lambda$, and their dependence on $x$ gives the shift
identity in \eqref{fh:external-rec}. Formula \eqref{fh:external-reflection}
then commutes with $\tau$ on Laurent generators and, by continuity,
on each completion. The formulas are valid over $\mathbb C$, and
Theorem~\ref{thm:intro} identifies the character ring on $\Lambda$.

\subsection{The simple-root subgroup and the reflection representative}
\label{sec:fh-rank-one}

\sectionoverview{In this subsection we compute the finite change of the Miura coordinates under the simple-root subgroup and the representative of its reflection. Their actions are determined by the same Gauss factorization. Its simple-root coordinates give a rational solution of the difference equation in Theorem~\ref{fh:external}.}

For $i\in J$ we write $O_i$ for its $\sigma$-orbit. The two elements
acting on the initial value of the Wronskian are
$x_i^\sigma(s)=\prod_{h\in O_i}x_h(s)$ and
$n_i=\prod_{h\in O_i}\dot s_h$. We denote by $\Gamma_i$ the restriction
of $\rho(n_i)$ to the field $\mathcal K^H$ of functions invariant under
right translation by $H$. If $\sigma=1$ these products each have one
factor.

In the Gauss factorization $\Xi_k=u_k b_k$ the negative Borel factor
has the form $b_k=z_k t_k$, with $z_k\in N_-$ and $t_k=[b_k]_0\in H$.
Similarly, $A_k^{\mathrm M}=V_kd_k$ with
$d_k=[A_k^{\mathrm M}]_0$ and $V_k\in N_-$. For $n\in N_-$,
the image of $n$ in $N_-/[N_-,N_-]$ has a factor
$y_h(z_h(n))$ at each simple root $h$. Thus $z_h$ is an additive
coordinate on this quotient. We write $z_{h,k}=z_h(z_k)$.
The Miura factorization \eqref{eq:rank-one-factors} gives
\[
 z_h(V_k)=
 \begin{cases}
 \displaystyle\miuraq_{h,k}^{-1}
       \prod_{p\prec h}\miuraq_{p,k}^{-A^G_{ph}},&h\in J,\\[3pt]
 0,&h\notin J.
 \end{cases}
\]
Indeed, the torus factors preceding $P_h(\miuraq_{h,k})$ conjugate
$y_h(\miuraq_{h,k}^{-1})$ by
$\prod_{p\prec h}\alpha_p^\vee(\miuraq_{p,k})$, and no other factor
contributes to $z_h(V_k)$. In particular $z_i(V_k)$ is a Laurent
monomial and a unit for $i\in J$.

Two normalizations of $z_{i,k+1}$ occur in the action. Conjugation of
the root subgroup by $t_{k+1}$ introduces the factor
$\alpha_i(t_{k+1})$. For the reflection the relevant quantity is the
ratio of $z_{i,k+1}$ to the coordinate $z_i(V_k)$ of the Miura form.
We use the functions
\begin{equation}\label{eq:fh-simple-auxiliary}
 \begin{aligned}
 B_{h,k}&=\alpha_h(t_{k+1})z_{h,k+1}
      &&(h\text{ a simple root of }G),\\
 \gamma_{i,k}&=\alpha_i(t_{k+1})z_i(V_k),\qquad
 \Sigma_{i,k}=\frac{z_{i,k+1}}{z_i(V_k)}
             =\frac{B_{i,k}}{\gamma_{i,k}}
      &&(i\in J).
 \end{aligned}
\end{equation}
The rational function $\Sigma_{i,k}$ belongs to $\mathcal K^H$.
Indeed, right multiplication of the initial value by an element of
$H$ changes only the torus factor of $b_k$ and fixes $A_k^{\mathrm M}$.

\begin{proposition}
\label{prop:fh-reflection}
Right multiplication by $x_i^\sigma(s)$ and by $n_i$ fixes
$\miuraq_{j,k}$ for $j\ne i$. At vertex $i$ their actions on the rational functions are
\begin{equation}\label{eq:fh-normalizer-miura}
 \begin{aligned}
 \frac{\rho(x_i^\sigma(s))(\miuraq_{i,k})}{\miuraq_{i,k}}
    &=\frac{1+sB_{i,k-o_i}}{1+sB_{i,k}},\\
 \frac{\Gamma_i(\miuraq_{i,k})}{\miuraq_{i,k}}
    &=1-\frac{1}{\Sigma_{i,k}},\qquad
 \Gamma_i(\Sigma_{i,k})=1-\Sigma_{i,k}.
 \end{aligned}
\end{equation}
The functions in \eqref{eq:fh-simple-auxiliary} satisfy
\begin{equation}\label{eq:fh-rankone-recurrence}
 B_{i,k}-B_{i,k-o_i}=\gamma_{i,k},\qquad
 \gamma_{i,k+o_i}=\mathsf R_{i,k}\gamma_{i,k},\qquad
 \Sigma_{i,k+o_i}=1+\mathsf R_{i,k}^{-1}\Sigma_{i,k},
\end{equation}
where $\mathsf R_{i,k}$ is the root monomial in
\eqref{ade:geometric-root} or \eqref{tw:screen-recurrence}.
\end{proposition}

\begin{proof}
Substitution of $b_k=z_kt_k$ and $A_k^{\mathrm M}=V_kd_k$ into
$b_{k+1}=A_k^{\mathrm M}\sigma(b_k)$ separates its torus and
unipotent parts. This gives
\begin{equation}\label{eq:fh-unipotent-recursion}
 t_{k+1}=d_k\sigma(t_k),\qquad
 z_{k+1}=V_k\operatorname{Ad}_{d_k}(\sigma(z_k)).
\end{equation}
On the simple negative-root space at $h$, $\operatorname{Ad}_{d_k}$ acts
by $\alpha_h(d_k)^{-1}$. Applying the additive coordinate $z_h$ to
the second identity and multiplying by $\alpha_h(t_{k+1})$ gives
$B_{h,k}-B_{\sigma^{-1}(h),k-1}=\alpha_h(t_{k+1})z_h(V_k)$.
At the indices $k-r$ and roots $\sigma^{-r}(i)$, the left sides
telescope. On the right only $r=0$ contributes, since $i$ is the unique
root of $O_i$ in the Miura product. Thus
\[
 B_{i,k}-B_{i,k-o_i}
  =\sum_{r=0}^{o_i-1}\alpha_{\sigma^{-r}(i)}(t_{k-r+1})
                  z_{\sigma^{-r}(i)}(V_{k-r})=\gamma_{i,k}.
\]
For $h=\sigma^{-k-1}(i)$, the calculation in
Lemma~\ref{lem:frame-screenings} gives
$\partial_h(\miuraq_{i,k})/\miuraq_{i,k}
=-\alpha_i(t_{k+1})z_i(V_k)=-\gamma_{i,k}$.
The same $h$ acts at $k+o_i$. The recurrence for these logarithmic
derivatives in that lemma gives
$\gamma_{i,k+o_i}=\mathsf R_{i,k}\gamma_{i,k}$. Dividing
$B_{i,k+o_i}=B_{i,k}+\gamma_{i,k+o_i}$ by $\gamma_{i,k+o_i}$ now yields
\[
 \Sigma_{i,k+o_i}
 =\frac{B_{i,k}+\gamma_{i,k+o_i}}{\gamma_{i,k+o_i}}
 =1+\mathsf R_{i,k}^{-1}\Sigma_{i,k}.
\]

For $g=x_i^\sigma(s)$ or $g=n_i$, $\sigma(g)=g$. In
\eqref{fh:framed-gauge} the positive unipotent factor of $b_kg$ can
be computed in the Levi subgroup generated by $H$ and the root
subgroups at $\pm\alpha_h$ for $h\in O_i$. Indeed, the other negative
root subgroups belong to the unipotent radical of the corresponding
negative parabolic. This radical is normalized by the Levi subgroup and
can therefore be moved to the right into $B_-$. The rank-one groups
at the vertices of $O_i$ commute. For $r=z_{h,k}$ and $u=s\alpha_h(t_k)$ their
two matrix products, up to a right torus factor, are
\[
 \begin{pmatrix}1&0\\r&1\end{pmatrix}
 \begin{pmatrix}1&u\\0&1\end{pmatrix}
 =\begin{pmatrix}1&u\\r&1+ru\end{pmatrix},\qquad
 \begin{pmatrix}1&0\\r&1\end{pmatrix}
 \begin{pmatrix}0&-1\\1&0\end{pmatrix}
 =\begin{pmatrix}0&-1\\1&-r\end{pmatrix}.
\]
Factoring either matrix $M$ as an upper unipotent matrix times a lower
triangular one gives $M_{12}/M_{22}$ for its upper unipotent parameter
where $M_{22}\ne0$.
At root $i$ and index $k+1$ this gives
\[
 v=\frac{s\alpha_i(t_{k+1})}{1+sB_{i,k}}
 \quad(g=x_i^\sigma(s)),\qquad
 v=\frac{1}{z_{i,k+1}}\quad(g=n_i).
\]
In \eqref{fh:framed-gauge}, the inverse factor $x_i(-v)$ multiplies
$A_k^{\mathrm M}$ on the left. Moving it past the factors preceding
$P_i(\miuraq_{i,k})$ multiplies $-v$ by
$\miuraq_{p,k}^{-A^G_{pi}}$ at each $p\prec i$ by
\eqref{ade:rank-one-distinct}. At $P_i(\miuraq_{i,k})$ the resulting
parameter is added to $\miuraq_{i,k}$ by \eqref{ade:rank-one-same}. Hence
\[
 \rho(g)\miuraq_{i,k}
 =\miuraq_{i,k}-v\prod_{p\prec i}\miuraq_{p,k}^{-A^G_{pi}}
 =\miuraq_{i,k}(1-vz_i(V_k)).
\]
The other upper factors have roots in $O_i\setminus\{i\}$, none of
which occur in the Miura product. By \eqref{ade:rank-one-distinct},
passing an upper factor through a factor at a different root preserves
that factor's Miura coordinate. After all the factors have been
reordered, the remaining upper factor, including
$\sigma(v_k(g,b_k))$ from \eqref{fh:framed-gauge}, is trivial. Indeed,
both the transformed element and its Miura product lie in $B_-$ and
$B_-\cap N_+=\{1\}$. Thus all $\miuraq_{j,k}$ with $j\ne i$ are fixed.
The two values of
$1-vz_i(V_k)$ are
\[
 1-\frac{s\gamma_{i,k}}{1+sB_{i,k}}
  =\frac{1+sB_{i,k-o_i}}{1+sB_{i,k}},\qquad
 1-\frac{z_i(V_k)}{z_{i,k+1}}=1-\Sigma_{i,k}^{-1}.
\]
Since $B_{i,k}-B_{i,k-o_i}=\gamma_{i,k}$ and
$\Sigma_{i,k}=B_{i,k}/\gamma_{i,k}$, the multiplier for $x_i^\sigma(s)$
has value $1-\Sigma_{i,k}^{-1}$ at $s=\infty$. The formula for
the reflection gives $\Sigma_{i,k}=\miuraq_{i,k}/(\miuraq_{i,k}-\Gamma_i\miuraq_{i,k})$.
Applying $\Gamma_i$ and using $\Gamma_i^2=1$ on $\mathcal K^H$
(Proposition~\ref{fh:frame-braid}) gives
$\Gamma_i\Sigma_{i,k}=1-\Sigma_{i,k}$.
\end{proof}

To compare with \eqref{fh:external-reflection}, the index $k$ corresponds
to $b=b_{i,k}$, the index $k-o_i$ to $b\qaff_i^{-2}$, and the root
monomial $\mathsf R_{i,k-o_i}$ to $\mathsf R_i(b\qaff_i^{-1})$ by
\eqref{ade:spectral-dictionary} and \eqref{tw:root-identification}.
Thus the recurrence for $\Sigma_{i,k}$ is precisely
\eqref{fh:external-rec} at $x=b\qaff_i^{-1}$. For the simple reflection
from $w$ to $ws_i$ its solution in the target completion satisfies
\[
 1-\frac{1}{\Sigma^{ws_i}_{i,b\qaff_i^{-1}}}
   =\mathsf R_i(b\qaff_i^{-1})^{-1}
      \frac{\Sigma^{ws_i}_{i,b\qaff_i^{-3}}}
           {\Sigma^{ws_i}_{i,b\qaff_i^{-1}}}.
\]
The right-hand side is $\Theta_i^w(Y_{i,b})/Y_{i,b}$.
Thus the same recurrence determines the multiplier of the rational
reflection and, upon expansion in the target chamber, the multiplier
of the Frenkel--Hernandez reflection.

\begin{remark}
\label{rem:fh-first-order}
The parameter $s$ of $x_i^\sigma(s)$ also realizes the formal
deformation in \cite[\S5.2, (5.63)]{FHWeyl}. Indeed,
$D_{i,k}(s)=\gamma_{i,k}^{-1}+s\Sigma_{i,k}
=(1+sB_{i,k})/\gamma_{i,k}$ satisfies
\[
 D_{i,k+o_i}(s)=s+\mathsf R_{i,k}^{-1}D_{i,k}(s),\qquad
 \rho(x_i^\sigma(s))(\miuraq_{i,k})
    =\miuraq_{i,k}\bigl(1-s/D_{i,k}(s)\bigr).
\]
Differentiation at $s=0$ sends $\miuraq_{i,k}$ to
$-\miuraq_{i,k}\gamma_{i,k}$. The identity
$-\gamma_{i,k+o_i}=\mathsf R_{i,k}(-\gamma_{i,k})$
is the defining recurrence of the screening symbols. Here the
differentiation concerns only the subgroup parameter $s$, with the
quantum parameter $\qaff$ fixed.
\end{remark}

\begin{remark}
\label{rem:fh-bax-rational}
For a $Z$-twisted Miura--Pl\"ucker $(G,q)$-oper with simply laced $G$,
the rational simple-root change of frame in \cite[Proposition~7.1]{FKSZ}
replaces $Q_i^+$ by $Q_i^-$ and $Z$ by $s_i(Z)$ under the hypotheses
of that proposition. It preserves $Q_j^+$ for $j\ne i$.
After exchanging the two Borels to use the convention of this paper,
we use the same torus changes of frame for the two Miura forms at
each index on the orbit. Their respective
Miura coordinates at $i$ are $y_k$ and $y'_k$. The parameter added
to $y_k$ at its simple-root factor is $y'_k-y_k$. On passing through
that factor its value becomes $(y'_k-y_k)/(y_k y'_k)$ by
\eqref{ade:rank-one-same}. The torus multipliers from the other
factors in \eqref{ade:rank-one-distinct} then give
\[
 \frac{y'_k-y_k}{y_k}=
 \mathsf R_{i,k}^{-1}\frac{y'_{k+1}-y_{k+1}}{y'_{k+1}} .
\]
Where $y'_k\ne y_k$, the function $S_k=y_k/(y_k-y'_k)$ therefore
satisfies
\[
 y'_k=y_k(1-S_k^{-1}),\qquad S_{k+1}=1+\mathsf R_{i,k}^{-1}S_k .
\]
These are the two equations for the reflected coordinate and
$\Sigma_{i,k}$ in Proposition~\ref{prop:fh-reflection} with $o_i=1$.
The rational B\"acklund transformation thus supplies a solution of this recurrence
on the chosen multiplicative orbit. Formula~\eqref{fh:external-reflection}
uses the solution determined by the target Weyl chamber.
\end{remark}

\subsection{Simultaneous evaluation in the Weyl chambers}
\label{fh:simultaneous}

\sectionoverview{In this subsection we extend the comparison for one
simple reflection to successive reflections. For each Weyl chamber we
construct the factor $z_0\in N_-$ in the decomposition
$\Xi_0=u_0z_0t_0$ of the initial value over the completed Laurent ring. Evaluation at that factor
identifies the action of each normalizer representative with the
corresponding Frenkel--Hernandez reflection. This comparison allows us
to deduce the relations on the completed Laurent rings from the
relations among the representatives acting on Wronskians.}

In the coordinates of \eqref{eq:wronskian-gauss-product}, the initial
Borel factor is $b_0=z_0t_0$ with $z_0\in N_-$ and $t_0\in H$.
Right multiplication of the initial solution by $H$ fixes both $z_0$
and the Miura variables. The ring
\[
 D_0=\mathcal Y_\Lambda\otimes_{\mathbb C}\mathbb C[N_-]
 \ \subset\ \mathcal K^H
\]
therefore consists of Laurent polynomials in the Miura variables
with coefficients that are regular functions of $z_0$. Equation
\eqref{eq:fh-unipotent-recursion} expresses every $z_k$ in these
coordinates. The recurrence and its inverse are regular in the $N_-$ coordinates,
with Laurent polynomial coefficients in the Miura variables. Consequently the shift
$\tau(\miuraq_{j,k})=\miuraq_{j,k+1}$ extends to an automorphism
of $D_0$ sending the coordinates of $z_k$ to those of $z_{k+1}$.

\begin{lemma}\label{fh:chamber-sol}
For each $w\in W$ there is a unique sequence
$z_k^w\in N_-(\widehat{\mathcal Y}^{\,w}_\Lambda)$ satisfying
\[
 z_{k+1}^w=V_k\operatorname{Ad}_{d_k}\bigl(\sigma(z_k^w)\bigr),
 \qquad z_{k+1}^w=\tau z_k^w \qquad(k\in\mathbb Z),
\]
where $d_k=[A_k^{\rm M}]_0$ and $V_k=A_k^{\rm M}d_k^{-1}$.
\end{lemma}

\begin{proof}
Since $z_k^w=\tau^k z_0^w$, the equation to solve is
$\tau z_0^w=V_0\operatorname{Ad}_{d_0}(\sigma(z_0^w))$.
An ordered product of the negative-root subgroups gives coordinates
on $N_-$. We use increasing height of the corresponding positive
roots. The negative-root subgroups indexed by positive roots of height at
least $p$ generate a subgroup $N_{-,\ge p}$. The quotient
$N_{-,\ge p}/N_{-,\ge p+1}$ is the additive group of
$\bigoplus_{\operatorname{ht}(\beta)=p}\mathfrak g_{-\beta}$.
Thus, once the coordinates of smaller height are known, the coordinate
$x_\beta$ of $z_0^w$ for $\operatorname{ht}(\beta)=p$ satisfies
\[
 \tau x_\beta=a_\beta x_{\sigma^{-1}\beta}+P_\beta,\qquad
 a_\beta=c_\beta\beta(d_0)^{-1}.
\]
Here $c_\beta\in\mathbb C^\times$ records the action of $\sigma$ on
root vectors, and $P_\beta$ is a polynomial in the smaller-height
coordinates with Laurent polynomial coefficients.

If the $\sigma$-orbit $O$ of $\beta$ has size $m$, iterating this
equation $m$ times gives
\[
 \tau^m x_\beta=Cx_\beta+P,\qquad
 C=\prod_{j=0}^{m-1}\tau^{m-1-j}(a_{\sigma^{-j}\beta}).
\]
Here $P$ depends only on translates under $\tau$ of the known
terms $P_\gamma$ and the coefficients $a_\gamma$. To determine
which series solves this equation in chamber $w$, we compute the
weight of the monomial factor of $C$. On the $\sigma$-invariant
root lattice of $G$ we use the map
$\psi(\lambda)=\sum_{j\in J}\langle\lambda,\alpha_j^\vee\rangle
\omega_j^{\rm q}$. Since
$d_0=\prod_{j\in J}\alpha_j^\vee(\miuraq_{j,0})$ and $\tau$ preserves
weights, $-\mathrm{wt}(C)=\psi(\sum_{\gamma\in O}\gamma)$.
The Cartan formula of Section~\ref{tw:section} and orthogonality within
$O_i$ give
\[
 \psi\!\left(\sum_{h\in O_i}\alpha_h\right)=\alpha_i^{\rm q},\qquad
 \psi\!\left(\Bigl(\prod_{h\in O_i}s_h\Bigr)\lambda\right)
     =s_i\psi(\lambda).
\]
The second identity uses the equality of the pairings of $\lambda$
with $\alpha_h^\vee$ for $h\in O_i$. For nontrivial $\sigma$,
Proposition~\ref{prop:fh-orbit-group} identifies $w$ with an element
of $W(G)$ commuting with $\sigma$. For $\sigma=1$ we already have
$w\in W(G)$.
Consequently $wO$ is a single $\sigma$-orbit, and its roots have the
same sign because $\sigma$ permutes the positive simple roots.
In their sum, simple roots from the same orbit have equal
coefficients. Applying $\psi$ therefore gives a nonzero sum of
simple roots of $\mathfrak a$ with coefficients of one sign. It follows
that $w\,\mathrm{wt}(C)$ lies in $Q_+\setminus\{0\}$ or its negative.

Applying $\tau^{-m}$ to the equation for $x_\beta$ and writing
$T(f)=\tau^{-m}(C)\tau^{-m}f$ gives
$(1-T)x_\beta=\tau^{-m}P$. On $\widehat{\mathcal Y}^{\,w}_\Lambda$,
\[
 (1-T)^{-1}=
 \begin{cases}
   \displaystyle\sum_{a\ge0}T^a,&w\,\mathrm{wt}(C)\in -Q_+ ,\\[4pt]
   \displaystyle-\sum_{a\ge1}T^{-a},&w\,\mathrm{wt}(C)\in Q_+ .
 \end{cases}
\]
In each series, successive powers lower chamber weights by a fixed
nonzero element of $Q_+$. The supports therefore remain in the
original finite union of translated negative cones, and only finitely
many powers contribute at each weight. Telescoping gives a two-sided
inverse. This determines $x_\beta$ for one root of $O$. Since each $a_\gamma$
is a unit, the equation for $\tau x_\gamma$ recovers
$x_{\sigma^{-1}\gamma}$ successively around $O$. The $m$-shift equation
ensures that the recovered value at $\beta$ agrees with the initial one. Induction on $p$ proves existence and uniqueness of $z_0^w$.
\end{proof}

Evaluation of regular functions on $N_-$ at $z_0^w$ defines a
homomorphism $\Phi_w:D_0\to\widehat{\mathcal Y}^{\,w}_\Lambda$ that
fixes the Miura variables. It commutes with $\tau$ by
Lemma~\ref{fh:chamber-sol}.

\Needspace{7\baselineskip}
\begin{theorem}\label{fh:intertwining}
There is a $W$-stable localization $D$ of $D_0$ inside $\mathcal K^H$
to which every $\Phi_w$ extends. For each $w\in W$ and $i\in J$, chamber
evaluation intertwines the reflections as in the commutative diagram
\begin{equation}\label{fh:intertwiners}
 \begin{array}{ccc}
 D & \xrightarrow{\ \Gamma_i\ } & D\\
 \big\downarrow\rlap{$\scriptstyle\,\Phi_w$} & &
 \big\downarrow\rlap{$\scriptstyle\,\Phi_{ws_i}$}\\
 \widehat{\mathcal Y}^{\,w}_\Lambda & \xrightarrow{\ \Theta_i^w\ } &
 \widehat{\mathcal Y}^{\,ws_i}_\Lambda .
 \end{array}
\end{equation}
\end{theorem}

\begin{proof}
A denominator can be inverted before chamber evaluation if its image
is a unit in every target ring. This gives the multiplicative set $S$
and the localization $D$
\[
 S=\bigcap_{v\in W}\Phi_v^{-1}\bigl((\widehat{\mathcal Y}^{\,v}_\Lambda)^\times\bigr),
 \qquad D=S^{-1}D_0\subset\mathcal K^H .
\]
Every $\Phi_v$ extends to $D$, and $\tau S=S$ since $\Phi_v\tau=\tau\Phi_v$.
For $i\in J$, the unit $z_i(V_k)\in\mathcal Y_\Lambda^\times$ gives
$\Sigma_{i,k}=z_{i,k+1}/z_i(V_k)\in D_0$. By
\eqref{eq:fh-rankone-recurrence} and the uniqueness in
\eqref{fh:external-rec} (\cite[Lemma~2.8]{FHWeyl}),
\[
 \Phi_v(\Sigma_{i,k})=\Sigma^v_{i,b_{i,k}\qaff_i^{-1}}
 \in(\widehat{\mathcal Y}^{\,v}_\Lambda)^\times .
\]
Thus $\Sigma_{i,k},z_{i,k+1}\in S$. For $h\notin J$ we have $z_h(V_k)=0$,
so \eqref{eq:fh-unipotent-recursion} gives
$z_{h,k+1}=\alpha_h(d_k)^{-1}z_{\sigma^{-1}(h),k}$.
As $\alpha_h(d_k)$ is a Laurent monomial, this equality preserves
membership in $S$. Iterating until $\sigma^{-1}$ reaches the unique
representative in $J$ therefore gives $z_{h,k}\in S$ for every
simple root $h$ and every $k$.

We next check $\Gamma_i(D_0)\subset D$ on the unipotent coordinates
and the Laurent generators separately. Since $n_i$ normalizes $H$,
$b_0n_i=z_0n_i(n_i^{-1}t_0n_i)$. Thus the transformed $N_-$ factor
is obtained from the Gauss factorization of $z_0n_i$. In the standard
Levi subgroup for $O_i$, the rank-one factors at $h\in O_i$ commute.
Their Gauss factorizations require only $z_{h,0}^{-1}$, by the
rank-one calculation in Proposition~\ref{prop:fh-reflection}. The
remaining negative-root factors belong to the unipotent radical of
the opposite parabolic. The Levi subgroup normalizes this radical, so
these factors can be moved to the right without further division. Hence
\[
 \Gamma_i\bigl(\mathbb C[N_-]\bigr)
 \subset \mathbb C[N_-][z_{h,0}^{-1}\mid h\in O_i]\subset D .
\]
On the Laurent generators, \eqref{eq:fh-normalizer-miura}
and \eqref{eq:fh-rankone-recurrence} give
\[
 \frac{\Gamma_i(\miuraq_{j,k})}{\miuraq_{j,k}}=
 \begin{cases}
 \mathsf R_{i,k-o_i}^{-1}\Sigma_{i,k-o_i}/\Sigma_{i,k},&j=i,\\
 1,&j\ne i .
 \end{cases}
\]
Both multipliers lie in $D^\times$, proving the claimed inclusion.

For $v=ws_i$, the homomorphisms $E=\Phi_v\Gamma_i$ and
$E'=\Theta_i^w\Phi_w$ are now defined on $D_0$. Evaluating the
displayed multipliers by $\Phi_v$ gives \eqref{fh:external-reflection},
so $E=E'$ on $\mathcal Y_\Lambda$. It remains to compare their images of
the coordinates of $z_0$. Both maps commute with $\tau$.
For $\Gamma_i$ this follows from $\sigma(n_i)=n_i$, and for $\Theta_i^w$
from the shift of the spectral parameter in \eqref{fh:external-rec}
and \eqref{fh:external-reflection}. Since $V_k$ and $d_k$ depend only on the Miura variables, the two
sequences $E(z_k)$ and $E'(z_k)$ satisfy the same transformed
recurrence \eqref{eq:fh-unipotent-recursion}. In both sequences the
value at $k$ is the $k$th shift of the initial value.

To compare the initial values, we use an integral linear form $\ell$
strictly positive on the simple roots. The maximal height in chamber $v$ is
\[
 \deg_v X=\max_{m\in\operatorname{supp}X}\ell(v\,\mathrm{wt}(m))
 \qquad(0\ne X\in\widehat{\mathcal Y}^{\,v}_\Lambda).
\]
The support condition in Definition~\ref{fh:completion} bounds these
heights above, and the shift preserves them. In
\eqref{fh:external-reflection} the two $\Sigma$-series have unique
leading monomials of the same weight, since their arguments differ
by a spectral shift. Thus the leading term of $\Theta_i^w(Y_{i,b})$
has weight $\omega_i^{\rm q}-\alpha_i^{\rm q}=s_i\omega_i^{\rm q}$.
For $j\ne i$ the generator $Y_{j,b}$ is fixed and
$s_i\omega_j^{\rm q}=\omega_j^{\rm q}$. Multiplication and inversion
therefore give a unique leading monomial of weight $s_i\,\mathrm{wt}(M)$
in the image of each Laurent monomial $M$. Its product with the
maximal-height part of $X$ is nonzero. Since $v s_i=w$,
\[
 \deg_v\!\bigl(\Theta_i^w(M)X\bigr)
   =\ell(w\,\mathrm{wt}(M))+\deg_v X .
\]
If $E(z_0)\ne E'(z_0)$, there is a positive root $\beta$ of least
root height for which the two coordinates at $-\beta$ differ. Their
difference $x$ is nonzero. With $m$ the size of the $\sigma$-orbit
of $\beta$, the inhomogeneous part of the $m$-fold recurrence
in Lemma~\ref{fh:chamber-sol} depends only on
coordinates at smaller root heights and on the Miura variables.
These agree for $E$ and $E'$. Subtraction therefore gives
\[
 \tau^m x=\Theta_i^w(C)x,\qquad
 w\,\mathrm{wt}(C)\in\pm(Q_+\setminus\{0\}),
\]
where $C$ is the coefficient of the $\beta$-coordinate after the $m$
shifts, as in that lemma. The shift preserves $\deg_v x$, whereas
multiplication by $\Theta_i^w(C)$ changes it by
$\ell(w\,\mathrm{wt}(C))\ne0$. Consequently $E(z_0)=E'(z_0)$, so the
two homomorphisms agree on $D_0$.

For $s\in S$ the inclusion $\Gamma_i(D_0)\subset D$ gives
$\Gamma_i(s)=a/b$ with $a\in D_0$ and $b\in S$. Equality on $D_0$ yields, for every $w$,
\[
 \Phi_{ws_i}(a)=\Phi_{ws_i}(b)\,\Theta_i^w(\Phi_w(s))
 \in(\widehat{\mathcal Y}^{\,ws_i}_\Lambda)^\times .
\]
As $ws_i$ ranges over all chambers, this shows $a\in S$ and hence
$\Gamma_i(S)\subset D^\times$.
Localization extends \eqref{fh:intertwiners} to $D$ and gives
$\Gamma_i(D)\subset D$. Since $\Gamma_i^2=1$ on $\mathcal K^H$,
this inclusion is an equality for every $i$.
\end{proof}

Theorem~\ref{fh:intertwining} now gives the Weyl group relations.
For a Laurent polynomial in the Miura variables, each composition
of chamber maps evaluates the corresponding composition of the
$\Gamma_i$ in the final chamber. In the square relation the final
chamber is the initial one, and $\Gamma_i^2=1$ on $D$. In a braid
relation both compositions end in the same chamber and the two
products agree on $D$ by Proposition~\ref{fh:frame-braid}.
Continuity and the density of the Laurent polynomials extend these
identities to each completed Laurent ring.

\begin{remark}\label{fh:primitive-operators}
The extension of the character variables by the highest
$\ell$-weight variables $\Psi_{i,b}$ of the prefundamental
representations admits two actions. The normalized operators of
\cite[Theorem~7.6]{GHLII} define a Weyl group action. The operators
$\Theta_i'$ of \cite[Proposition~3.17]{FHExtended} satisfy the
braid relations, and each has infinite order
\cite[Remark~7.7]{GHLII}. By Proposition~\ref{fh:frame-braid}, the
pinned representatives on Wronskians have fourth power one. On the completed rings of $Y$-variables, Theorem~\ref{fh:intertwining}
identifies their action with the Frenkel--Hernandez maps from $w$
to $ws_i$ and back, which are mutual inverses.
\end{remark}

\subsection{Fixed functions and the screening equations}
\label{fh:invariants-subsection}

\sectionoverview{Here we characterize the rational functions of the Miura variables
fixed by a simple reflection. The criterion is that the corresponding
screening annihilates the function. For Laurent polynomials it is also equivalent to
invariance under the corresponding Frenkel--Hernandez reflection in
any Weyl chamber.}

By Proposition~\ref{fh:frame-action}, pullback identifies the field
$F_\Lambda=\operatorname{Frac}\mathcal Y_\Lambda$ with $\mathcal K^{B_-}$.
Each screening extends to this field by the quotient rule, with values
in its screening module tensored with $F_\Lambda$. A Laurent polynomial
$f$ also defines an element of each chamber completion by the same finite
sum. Thus in $\Theta_i^w(f)=f$ the left side is the image of $f$ from the
chamber $w$ and the right side is its finite sum in the chamber $ws_i$.

To recover invariance under $\Gamma_i$ from its chamber evaluation, we
need algebraic independence of the series that enter the formula for
$\Gamma_i$ on the Miura variables.

\begin{lemma}\label{lem:fh-chamber-independence}
For each $v\in W$, the chamber completion $\widehat{\mathcal Y}^{\,v}_\Lambda$
is a domain. For $i\in J$, the $o_i$ elements
$\Phi_v(\Sigma_{i,0}),\ldots,\Phi_v(\Sigma_{i,o_i-1})$ of its fraction
field are algebraically independent over $F_\Lambda$.
\end{lemma}

\begin{proof}
An integral linear form $h$ positive on the simple roots assigns to each
monomial $m$ the height $h(v\,\mathrm{wt}(m))$. The weights $v\,\mathrm{wt}(m)$ of a series in chamber $v$ lie in a
finite union of translates $\lambda-Q_+$. Their heights are bounded
above, and only finitely many weights occur at each height. Since the
series has finitely many monomials per weight
(Definition~\ref{fh:completion}), its component of each height is a
Laurent polynomial. If
$F$ and $G$ are nonzero with highest nonzero components $F_a$ and
$G_b$, then
\[
 (FG)_{a+b}=F_aG_b\ne0
 \qquad\text{in }\mathcal Y_\Lambda .
\]
This proves the domain assertion and embeds $F_\Lambda$ in the
fraction field of the completion.

For $0\le\rho<o_i$, the two chamber expansions of
\eqref{eq:fh-rankone-recurrence} give
\[
 H_\rho:=\Phi_v(\Sigma_{i,\rho})=
 \begin{cases}
 \displaystyle\sum_{n\ge0}\prod_{p=1}^n\mathsf R_{i,\rho-po_i}^{-1},
     &v\alpha_i^{\rm q}>0,\\[5pt]
 \displaystyle-\sum_{n\ge1}\prod_{p=0}^{n-1}\mathsf R_{i,\rho+po_i},
     &v\alpha_i^{\rm q}<0 .
 \end{cases}
\]
The recurrence is expanded toward smaller indices in the first case
and toward larger indices in the second. To distinguish the $H_\rho$, we use the Euler
derivation $E_N=\miuraq_{i,N}\partial/\partial\miuraq_{i,N}$. It
multiplies each monomial by its exponent at $(i,N)$ and therefore acts
termwise in the completion. The only
variables at vertex $i$ in $\mathsf R_{i,k}$ are $\miuraq_{i,k}$
and $\miuraq_{i,k+o_i}$, both with exponent one. For $n\ge1$, the
$n$th product in the first expansion thus contains
$\miuraq_{i,\rho-no_i}$ with exponent $-1$, and the $n$th product in
the second contains $\miuraq_{i,\rho+no_i}$ with exponent $1$.
For $N=\rho-no_i$ in the first case and $N=\rho+no_i$ in the second,
the derivative of this product is nonzero. It cannot cancel against
another summand, since the products have distinct weights
$-n\alpha_i^{\rm q}$ or $n\alpha_i^{\rm q}$, respectively.
At vertex $i$, the series $H_{\rho'}$ uses only the indices
$\rho'+po_i$ with $p\in\mathbb Z$. Thus
\[
 E_N H_\rho\ne0,\qquad
 E_N H_{\rho'}=0\quad(0\le\rho'<o_i,\ \rho'\ne\rho).
\]

A relation over $F_\Lambda$ would, after clearing denominators, give a
nonzero polynomial $P\in\mathcal Y_\Lambda[X_0,\ldots,X_{o_i-1}]$
of least total degree with $P(H_0,\ldots,H_{o_i-1})=0$. Some $X_\rho$
occurs in $P$. For $n$ sufficiently large the corresponding variable
$\miuraq_{i,N}$ occurs in none of its coefficients, so
\[
 0=E_N\bigl(P(H_0,\ldots,H_{o_i-1})\bigr)
   =(\partial_{X_\rho}P)(H_0,\ldots,H_{o_i-1})\,E_N H_\rho .
\]
The domain property and $E_N H_\rho\ne0$ force the first factor to
vanish. Since $\partial_{X_\rho}P$ is nonzero in characteristic zero
and has smaller total degree, this contradicts minimality.
\end{proof}

The evaluation $\Phi_v$ fixes the Miura variables. By the lemma it is
injective on $\mathcal Y_\Lambda[\Sigma_{i,0},\ldots,\Sigma_{i,o_i-1}]$
and thus extends to a field embedding
\[
 F_\Lambda(\Sigma_{i,0},\ldots,\Sigma_{i,o_i-1})
   \lhook\joinrel\longrightarrow
 \operatorname{Frac}\widehat{\mathcal Y}^{\,v}_\Lambda .
\]
For $f\in\mathcal Y_\Lambda$, equations \eqref{eq:fh-normalizer-miura}
and \eqref{eq:fh-rankone-recurrence} place $\Gamma_i f$ in the field
on the left.

\begin{theorem}\label{fh:fixed-functions}
For $f\in\mathcal Y_\Lambda$, $i\in J$, and $w\in W$,
\[
 S_i f=0
 \quad\Longleftrightarrow\quad \Gamma_i f=f
 \quad\Longleftrightarrow\quad
 \Theta_i^w(f)=f \ \text{in }\widehat{\mathcal Y}^{\,ws_i}_\Lambda .
\]
For rational functions, the common kernel is
\[
 \bigcap_{i\in J}\ker(S_i|_{F_\Lambda})
   =\mathcal K^G=F_\Lambda\cap\mathcal K^K
   =F_\Lambda\cap(\mathcal K^H)^W .
\]
Its intersection with $\mathcal Y_\Lambda$ is
$\mu_{c,\sigma}^*\mathbb C[\mathcal A^{\mathbb Z}]$.
\end{theorem}

\begin{proof}
For $f\in\mathcal Y_\Lambda$, both $f$ and $\Gamma_i f$ belong to
$F_\Lambda(\Sigma_{i,0},\ldots,\Sigma_{i,o_i-1})$. Evaluation in
chamber $v=ws_i$ is injective on this field by
Lemma~\ref{lem:fh-chamber-independence} and fixes $f$. Hence
\eqref{fh:intertwiners} gives
\[
 \Theta_i^w(f)-f=\Phi_v(\Gamma_i f-f),\qquad
 \Theta_i^w(f)=f\ \Longleftrightarrow\ \Gamma_i f=f .
\]

To relate the rational reflection to the screening, we identify each
$f\in F_\Lambda$ with its pullback $p^*f\in\mathcal K^{B_-}$.
Its stabilizer $G_f$ under right translation therefore contains $B_-$. For $h\in O_i$, the factors of
$n_i$ away from $h$ commute with the rank-one subgroup at $h$. Thus
$n_i y_h(\mathbb G_a)n_i^{-1}=x_h(\mathbb G_a)$, and
\[
 n_i\in G_f
 \quad\Longleftrightarrow\quad
 x_h(\mathbb G_a)\subset G_f\ \text{for every }h\in O_i .
\]
Indeed, when $n_i$ fixes $f$, so does the conjugate of each
$y_h(\mathbb G_a)\subset B_-$ by $n_i$. Conversely, the two root
subgroups at $h$ generate its rank-one subgroup, which contains
$\dot s_h$. Their representatives have product $n_i$.

The group law gives
$\frac{d}{dt}\rho(x_h(t))f=\rho(x_h(t))\partial_h f$.
In characteristic zero a rational function of $t$ with zero derivative
is constant. Here evaluation at $t=0$ gives $f$, so
$x_h(\mathbb G_a)\subset G_f$ if and only if $\partial_h f=0$.
For $h=\sigma^{-1-\ell}(i)$, the quotient rule extends
\eqref{eq:frame-screening-component} to $F_\Lambda$ and gives
\[
 \partial_h(p^*f)
   =\lambda_{i,\ell}\,p^*(\mathscr D_{i,\ell}f),
 \qquad \lambda_{i,\ell}\ne0 .
\]
As $\ell$ runs from $0$ to $o_i-1$, $h$ runs through $O_i$, and
$\mathscr D_{i,\ell}f$ runs through the coefficients of $S_i f$ in a
basis of the screening module. Consequently
\[
 \Gamma_i f=f
 \quad\Longleftrightarrow\quad
 \mathscr D_{i,\ell}f=0\quad(0\le\ell<o_i)
 \quad\Longleftrightarrow\quad S_i f=0 .
\]
As $i$ ranges over $J$, these positive simple-root subgroups together
with $B_-$ generate $G$. Therefore $\bigcap_i\ker(S_i|_{F_\Lambda})
=\mathcal K^G$. Since $n_i\in K\subset G$ and the $\Gamma_i$ generate
the $W$-action on $\mathcal K^H$, we have
\[
 \mathcal K^G\subset F_\Lambda\cap\mathcal K^K
 \subset F_\Lambda\cap(\mathcal K^H)^W
 =\bigcap_i\ker(S_i|_{F_\Lambda}) .
\]
On $\mathcal Y_\Lambda$, Theorem~\ref{thm:regular-descent} gives the
claimed pullback of the Steinberg coordinate ring. Theorem~\ref{thm:intro}
identifies its generators with the fundamental $q$-characters on $\Lambda$.
\end{proof}

\subsection{Weyl symmetry for folded characters}
\label{fh:folded-action}

\sectionoverview{Here we construct the Weyl action for the ordinary
Miura transformation of the fixed group $G_\nu=(\widetilde G^\nu)^\circ$.
We identify the folded $q$-character ring as the ring of invariant
Laurent polynomials and compare the spectral shifts with those in the
ordinary $q$-characters of non-simply-laced type.}

For each $\nu$-orbit $O$ of simple roots, the product
$r_O=\prod_{j\in O}\widetilde s_j$ consists of commuting reflections.
The argument of Proposition~\ref{prop:fh-orbit-group}, applied to $\nu$,
identifies these products with the simple reflections of
$W(G_\nu)=W(\widetilde G)^\nu$. We use the folded spectral lattice
$b_{O,k}=a\qaff^{2k+\xi_O}$ of Theorem~\ref{thm:intro-folded}.
For each $w\in W(G_\nu)$, Definition~\ref{fh:completion} applied to
the weights and positive roots of $G_\nu$ gives the folded chamber
completion. The folded Laurent ring embeds diagonally in the product
of these completions.

\begin{proposition}\label{fh:folded-weyl}
For $w\in W(G_\nu)$ the substitution $\pi(\widetilde Y_{j,b})=Y_{O,b}$
for $j\in O$ intertwines the product of the simply laced
Frenkel--Hernandez reflections at the vertices of $O$ with a map
$\overline\Theta_O^w$ from the folded completion at $w$ to the one at
$wr_O$. It fixes the variables at every other orbit and satisfies
\begin{equation}\label{fh:folded-reflection}
 \overline\Theta_O^w(Y_{O,b})
   =Y_{O,b}\left(1-
          \bigl(\overline\Sigma^{\,wr_O}_{O,b\qaff^{-1}}\bigr)^{-1}\right),
 \qquad
 \overline\Sigma^{\,wr_O}_{O,c}
   =1+\mathsf R_O(c)^{-1}\overline\Sigma^{\,wr_O}_{O,c\qaff^{-2}} .
\end{equation}
Here $\mathsf R_O$ is the folded root monomial of
\eqref{eq:folded-root-monomial}, and the series are the
shift-covariant solution in the target chamber $wr_O$.
These maps define an action of $W(G_\nu)$ on the product of the
folded chamber completions. The Laurent polynomials in the diagonal
subring that are fixed by this action form the
folded $q$-character ring of Theorem~\ref{thm:intro-folded}.
\end{proposition}

\begin{proof}
For $w\in W(\widetilde G)^\nu$, restriction of weights to $H_\nu$ commutes
with $w$. In each cone $\lambda-\widetilde Q_+$ it is given by
\[
 \left(\lambda-\sum_j n_j\widetilde\alpha_j\right)\bigg|_{H_\nu}
 =\lambda|_{H_\nu}-\sum_{O'}\left(\sum_{j\in O'}n_j\right)\alpha_{O'},
 \qquad n_j\ge0,
\]
where $\widetilde Q_+$ is the positive root cone of $\widetilde G$.
A prescribed folded weight fixes each sum $\sum_{j\in O'}n_j$, hence
has finitely many preimages in the cone. A series in chamber $w$
has support in a finite union of such cones after applying $w$, and
has only finitely many monomials at each weight. Thus at every folded
weight $\pi$ sums only finitely many coefficients, giving a continuous
map between the completions at $w$.

We write $v=wr_O$ and $\widetilde\Theta_{r_O}^w$ for the composition
of the simply laced chamber maps for the vertices of $O$, starting at $w$. If $j,h\in O$ are distinct, they are
nonadjacent. Thus $\widetilde s_h\widetilde\alpha_j=\widetilde\alpha_j$,
and the reflection at $h$ fixes every variable of
$\widetilde{\mathsf R}_j(c)$. Its action preserves the chamber expansion
of the series for $j$. Since $v$ commutes with $\nu$, the signs of
$v\widetilde\alpha_j$ agree for all $j\in O$. By
\eqref{eq:folded-root-monomial}, their series consequently have the
same image $\overline\Sigma^v_{O,c}$. For $j\in O'$ and $c=b\qaff^{-1}$,
\eqref{fh:external-reflection} gives
\[
 \pi\widetilde\Theta_{r_O}^w(\widetilde Y_{j,b})=
 \begin{cases}
  Y_{O',b},&O'\ne O,\\[3pt]
  \displaystyle Y_{O,b}\mathsf R_O(c)^{-1}
       \frac{\overline\Sigma^v_{O,c\qaff^{-2}}}{\overline\Sigma^v_{O,c}}
   =Y_{O,b}\bigl(1-(\overline\Sigma^v_{O,c})^{-1}\bigr),&O'=O .
 \end{cases}
\]
The right side depends only on the orbit of $j$, so it defines
$\overline\Theta_O^w$ on Laurent generators. The quotient of the two
shifted series has a single leading monomial of weight zero. Thus the
image of $Y_{O,b}$ has leading weight
$\operatorname{wt}(Y_{O,b})-\alpha_O=r_O\operatorname{wt}(Y_{O,b})$.
Every further term has strictly lower weight in chamber $v$. The
leading monomial is a unit, so inversion gives an expansion with
opposite leading weight and the same decrease in chamber $v$.
The generators at $O'\ne O$ are fixed and their weights are fixed by
$r_O$. Multiplying these expansions gives, for any Laurent monomial $m$,
\[
 v\,\mathrm{wt}(m')\in w\,\mathrm{wt}(m)-Q_+
 \qquad\bigl(m'\in\operatorname{supp}\overline\Theta_O^w(m)\bigr),
\]
where $Q_+$ is the positive root cone of $G_\nu$. For a fixed output
weight $\mu=v\,\mathrm{wt}(m')$, a contributing source weight in
$\lambda-Q_+$ lies between $\mu$ and $\lambda$ in the root order.
This interval contains only finitely many weights, and a source
series contains only finitely many monomials at each of them.
Consequently $\overline\Theta_O^w$ extends continuously, and
$\overline\Theta_O^w\pi=\pi\widetilde\Theta_{r_O}^w$ on the completions.
Every folded Laurent polynomial lifts through $\pi$. Applying the
intertwining identity to such a lift transfers the square and braid
relations for the elements $r_O\in W(\widetilde G)^\nu$ to folded
Laurent polynomials. Continuity and density extend the relations to
the folded completions.

In particular, $\pi(\chi_{\qaff}(V))$ is invariant for every ambient
character on the lattice, by Theorem~\ref{fh:fixed-functions}.
Conversely, suppose a folded Laurent polynomial $f$ is invariant.
For each $O$ the identity $\overline\Theta_O^{r_O}(f)=f$ holds in the
identity chamber. In this chamber
\[
 \overline\Sigma_k=\overline\Sigma^1_{O,b_{O,k}\qaff^{-1}},\qquad
 \mathsf R_{O,k}=\mathsf R_O(b_{O,k}\qaff),\qquad
 \overline\Sigma_0=\sum_{n\ge0}\prod_{\ell=1}^n\mathsf R_{O,-\ell}^{-1}.
\]
We next show that invariance forces $S_Of=0$. For this purpose
we first prove that $\overline\Sigma_0$ is transcendental over the
fraction field $F$ of the folded Laurent ring. The product indexed by
$n\ge1$ has weight $-n\alpha_O$ and exponent $-1$ at
$Y_{O,b_{O,-n}}$. As in Lemma~\ref{lem:fh-chamber-independence},
multiplication of the components of maximal height shows that the folded
completion is a domain. Suppose $P(\overline\Sigma_0)=0$ is a nonzero
polynomial relation of least degree over $F$. For large $n$, none of
the coefficients of $P$ involves $Y_{O,b_{O,-n}}$. The Euler derivative
$E_n$ in this variable annihilates the coefficients, while the distinct
weights of the products ensure that $E_n\overline\Sigma_0\ne0$. Hence
\[
 0=P'(\overline\Sigma_0)E_n\overline\Sigma_0
 \quad\Longrightarrow\quad P'(\overline\Sigma_0)=0 .
\]
In characteristic zero $P'$ is nonzero of smaller degree, a
contradiction.

We can therefore replace $\overline\Sigma_0$ by an indeterminate $T$
in the reflection identity. The affine polynomials $P_k(T)$ determined by
\[
 P_0(T)=T,\qquad P_{k+1}(T)=1+\mathsf R_{O,k}^{-1}P_k(T)
 \quad(k\in\mathbb Z)
\]
satisfy $P_k(\overline\Sigma_0)=\overline\Sigma_k$. Replacing each
$\overline\Sigma_k$ by $P_k(T)$ in the substitution
\eqref{fh:folded-reflection} sends $f$ to a rational function $f_T\in F(T)$.
At $T=\overline\Sigma_0$, its value is
$\overline\Theta_O^{r_O}(f)=f$. Transcendence makes evaluation on
$F(T)$ injective, so $f_T=f$. The coefficient $a_k$ of $T$ in $P_k$
satisfies $a_0=1$ and $a_{k+1}=\mathsf R_{O,k}^{-1}a_k$. Thus each
$a_k$ is nonzero and $P_k(T)^{-1}=a_k^{-1}T^{-1}+O(T^{-2})$ at
$T=\infty$. Substitution in each monomial of $f$ gives
\[
 f_T=f-T^{-1}\sum_{k\in\mathbb Z}a_k^{-1}Y_{O,b_{O,k}}
           \frac{\partial f}{\partial Y_{O,b_{O,k}}}+O(T^{-2}).
\]
The sum is finite because $f$ is a Laurent polynomial. It vanishes by
$f_T=f$. Moreover, \eqref{eq:folded-screening} and the recurrence
for $a_k$ give $s_{O,b_{O,k}}=a_k^{-1}s_{O,b_{O,0}}$. The sum is
therefore the coefficient of $s_{O,b_{O,0}}$ in $S_Of$. Since this
holds for every $O$, all folded screenings annihilate $f$.
Corollary~\ref{cor:folded-miura} identifies this common kernel with
the folded character ring.
\end{proof}

The same restriction describes the action on Wronskians before chamber
evaluation. For a Steinberg sequence in $G_\nu$ and an initial solution
value in $G_\nu$, the recurrence $\Xi_{k+1}=A_k\Xi_k$ keeps every $\Xi_k$
in $G_\nu$. By \eqref{eq:folded-pinning}, a simple-root subgroup
of $G_\nu$ is the product $\prod_{j\in O}\widetilde x_j(z)$ with the
same parameter at every vertex of $O$, and its reflection representative is
$\prod_{j\in O}\dot{\widetilde s}_j$. Uniqueness of Gauss
factorization identifies the changes of Miura form by these elements
in the ambient and fixed groups. Differentiation of the root subgroup
gives the folded screening by \eqref{eq:folded-orbit-screening}.

The difference between the folded and ordinary spectral shifts can be seen
for the involution of $A_3$ and for triality of $D_4$. In each diagram the
outer vertices form the orbit indexed by $1$, of size $m=2$ and $m=3$,
respectively. The central vertex has index $2$. The root monomials
containing $m$ neighboring factors are
\[
 \mathsf R^{\rm ord}_1(b)
  =\frac{Y_{1,b\qaff^{-m}}Y_{1,b\qaff^m}}
         {\prod_{s=0}^{m-1}Y_{2,b\qaff^{m-1-2s}}},
 \qquad
 \mathsf R^{\rm fold}_2(b)
  =\frac{Y_{2,b\qaff^{-1}}Y_{2,b\qaff}}{Y_{1,b}^{\,m}} .
\]
In the ordinary root monomial of \eqref{tw:root-identification}, the
$m$ neighboring factors occur at distinct spectral parameters. In the
folded monomial of \eqref{eq:folded-root-monomial}, identification of
the outer variables gives an $m$th power at a single parameter. The
Cartan matrices are transposes, so the multiple factors occur at opposite
vertices, while the braid length is four for the $A_3$ involution and
six for $D_4$ triality in both constructions.

\section{Difference opers and transfer matrices}
\label{sec:difference-opers}

\sectionoverview{In this section we evaluate the Laurent polynomials of
the Miura transformation on functions and on commuting operators. We
compare the flag and Steinberg descriptions of meromorphic
$(G,\eta,\sigma)$-opers, including ordinary $(G,q)$-opers when
$\sigma=1$ and $\eta=q$. We derive a local regularity criterion from the $QQ$-relations and compute the Steinberg
coordinates of twisted differential monodromy. In an XXZ-type model,
evaluation at ratios of Baxter operators gives the fundamental transfer
matrices.}

\subsection{Specialization to a difference algebra}
\label{sec:tau-specialization}

\sectionoverview{In this subsection, we specialize a lattice sequence to the successive
translates of a single group element under an algebra automorphism.
We prove that a unique change of frame in $N_+$ takes the connection
to Steinberg form and compute its coordinates by evaluating fundamental
$q$-characters.}

Let $G$ be simply laced with $\sigma=1$, or let $(G,\sigma)$ be one
of the four pairs in \eqref{tw:data}. The ordered list
$J=(i_1,\ldots,i_\ell)$ contains one representative from each orbit
of $\sigma$ on the simple roots and determines the section $\mathcal A=\mathcal A_c$ for
$c=s_{i_1}\cdots s_{i_\ell}$. A difference algebra $(D,\tau)$ is a
commutative complex algebra with a $\mathbb C$-algebra automorphism
$\tau$. For a disk centered at $0$, examples are the ring of holomorphic
germs and the field of meromorphic germs at $0$, with
$(\tau f)(u)=f(\eta u)$ for $\eta\in\C^\times$. Composition with
$u\mapsto\eta u$ and with $u\mapsto\eta^{-1}u$ is defined on germs.
An identity involving finitely many such compositions holds on a disk
small enough for all the indicated arguments.
On $G(D)$ the automorphism acts on a point
$g:\mathbb C[G]\to D$ by postcomposition, $\tau g=\tau\circ g$.
Since $\sigma$ is defined over $\mathbb C$, its action on $G(D)$
commutes with $\tau$.

In a trivialization, an element $A\in G(D)$ represents a
$(G,\tau,\sigma)$-connection. A change of frame by $g\in G(D)$ gives
\begin{equation}\label{eq:difference-gauge}
 A\longmapsto A^g=(\tau g)A\sigma(g)^{-1}.
\end{equation}
For a multiplicative shift $(\tau f)(u)=f(\eta u)$, the connection
sends $h$ in the fiber at $u$ to $A(u)\sigma(h)$ in the fiber at $\eta u$. A group-valued solution over $D$ is an element
$\Xi\in G(D)$ satisfying $\tau \Xi=A\sigma(\Xi)$. For the multiplicative shift this equation reads
$\Xi(\eta u)=A(u)\sigma(\Xi(u))$. Equation~\eqref{eq:difference-gauge}
sends $\Xi$ to a solution in the new frame, since
\[
 \tau(g\Xi)=(\tau g)A\sigma(\Xi)=A^g\sigma(g\Xi).
\]
If $\sigma$ has order $r$, iterating the solution equation $r$ times gives
\[
 \tau^r \Xi=A^{[r]}\Xi,\qquad
 A^{[r]}=\prod_{j=0}^{r-1}\sigma^j(\tau^{r-1-j}A),
\]
where the factors are ordered from left to right by increasing $j$.
The factor with index $j$ in $(A^g)^{[r]}$ is
\[
 (\sigma^j\tau^{r-j}g)\,(\sigma^j\tau^{r-1-j}A)\,
 (\sigma^{j+1}\tau^{r-1-j}g)^{-1}\qquad(0\le j<r).
\]
Multiplication in increasing order of $j$ cancels the intermediate
changes of frame and gives $(A^g)^{[r]}=(\tau^r g)A^{[r]}g^{-1}$.
For $(\tau f)(u)=f(\eta u)$ the composite is therefore an ordinary
$q$-connection with $q=\eta^r$.

The uniqueness of the Steinberg representative in
Theorem~\ref{thm:ade-full-cell} and Proposition~\ref{tw:normalization}
also gives a representative when the lift of $c$ depends on the base
point. If $L=\dot c h$ for $h\in H(D)$, the associated cell and
section are $N_+(D)LN_+(D)$ and $U_c(D)L$. For the ordinary shift
$u\mapsto qu$ over $\mathbb C(u)$ and a lift $L(u)$ prescribed by
the singularity polynomials, \cite[Theorem~8.2]{FKSZ} proves the
following result for $G\ne E_6$. In our Borel convention, for every
$g$ in this cell there is a unique $v\in N_+(\mathbb C(u))$ such that
$v(qu)g(u)v(u)^{-1}$ belongs to $U_c(\mathbb C(u))L(u)$.

\begin{proposition}
\label{prop:difference-section}
For every $g\in N_+(D)LN_+(D)$ there is a unique $u\in N_+(D)$
such that $(\tau u)g\sigma(u)^{-1}\in U_c(D)L$.
\end{proposition}

\begin{proof}
For $k\in\mathbb Z$, the elements $g_k=\tau^k g$ lie in
$N_+(D)L_kN_+(D)$, where $L_k=\dot c h_k$ and $h_k=\tau^k h$.
We first reduce this sequence with a change of frame $u_k$ at each
index. To apply the reduction for the pinned lift $\dot c$, we use the
torus sequence
\[
 t_0=1,\qquad t_{k+1}=c\bigl(\sigma(t_k)h_k^{-1}\bigr),
 \qquad c(t)=\dot c t\dot c^{-1}.
\]
Since $c$ and $\sigma$ act by automorphisms on $H$, the recurrence
determines $t_k$ in both directions. It gives
$t_{k+1}L_k\sigma(t_k)^{-1}=\dot c$. Conjugation by $H$ preserves
$N_+$ and $U_c$, so $x\mapsto t_{k+1}x\sigma(t_k)^{-1}$ identifies
the cell $N_+(D)L_kN_+(D)$ and its section $U_c(D)L_k$ with
$M_c(D)$ and $\mathcal A_c(D)$, respectively. In these coordinates a
change of frame $u_k$ is represented by $t_ku_kt_k^{-1}$.
Theorem~\ref{thm:ade-full-cell} for $\sigma=1$ and
Proposition~\ref{tw:normalization} for the other cases therefore give
a unique sequence $(u_k)$ satisfying
\[
 u_{k+1}g_k\sigma(u_k)^{-1}\in U_c(D)L_k\qquad(k\in\mathbb Z).
\]

It remains to show that $u_k=\tau^ku_0$. Applying $\tau$ to the
displayed condition gives
\[
 (\tau u_{k+1})g_{k+1}\sigma(\tau u_k)^{-1}\in U_c(D)L_{k+1}.
\]
Thus $(\tau u_k)_k$ and $(u_{k+1})_k$ both carry $(g_{k+1})_k$
to the sections $U_c(D)L_{k+1}$. Uniqueness gives
$u_{k+1}=\tau u_k$, so $u_0$ has the required property. For any other
$v$ with this property, the sequence $(\tau^kv)_k$ satisfies the same
reduction condition as $(u_k)_k$, hence $v=u_0$.
\end{proof}

Thus $U_c(D)L$ meets each $N_+(D)$-orbit in $N_+(D)LN_+(D)$
in exactly one point. For $L=\dot c$ this gives the section
$\mathcal A_c(D)$ used in the character formula.
For units $\boldsymbol{\miuraq}=(\miuraq_i)_{i\in J}$ in $D$, the translates of the
element $A^{\mathrm M}=\prod_{p=1}^{\ell}P_{i_p}(\miuraq_{i_p})$ in Miura form
have Miura variables $\miuraq_{i,k}=\tau^k(\miuraq_i)$. Their invertibility gives
a homomorphism from the Laurent ring of
\eqref{ade:spectral-dictionary} or \eqref{tw:dictionary},
\begin{equation}\label{eq:difference-evaluation}
 \operatorname{ev}_{\boldsymbol{\miuraq}}:\mathcal Y_\Lambda\longrightarrow D,\qquad
 Y_{i,a\qaff^{2k+\xi_i}}\longmapsto\tau^k(\miuraq_i).
\end{equation}

\begin{corollary}\label{cor:difference-specialization}
There is a unique $u\in N_+(D)$ for which
$A=(\tau u)A^{\mathrm M}\sigma(u)^{-1}\in\mathcal A(D)$. With the spectral
parameters of Theorems~\ref{ade:main} and~\ref{tw:character-theorem},
its coordinates are
\begin{equation}\label{eq:difference-character}
 T_i(A)=\operatorname{ev}_{\boldsymbol{\miuraq}}\!\left(\chi_{\qaff}(L(Y_{i,a\qaff^{\xi_i}}))\right).
\end{equation}
\end{corollary}

\begin{proof}
Proposition~\ref{prop:difference-section} with $L=\dot c$ gives the
unique $u$. For the sequence $A_k^{\mathrm M}=\tau^k A^{\mathrm M}$,
the changes of frame are $u_k=\tau^k u$. Its Steinberg representative
at $k=0$ is therefore $A=(\tau u)A^{\mathrm M}\sigma(u)^{-1}$.

Equations~\eqref{ade:raw-character} and \eqref{tw:raw-character}
express the coordinate at $k=0$ as a Laurent polynomial in the
independent variables $\miuraq_{j,k}$. For the sequence above,
$\miuraq_{j,k}=\tau^k(\miuraq_j)$. Evaluating the polynomial by
\eqref{eq:difference-evaluation} gives \eqref{eq:difference-character}.
This specialization also applies to meromorphic Miura variables,
since it requires only their shifted values.
\end{proof}

The module $L(Y_{i,a\qaff^{\xi_i}})$ in \eqref{eq:difference-character}
specifies a Laurent polynomial. In $\mathcal Y_\Lambda$ this is the
unique polynomial annihilated by every screening whose only dominant
monomial is $Y_{i,a\qaff^{\xi_i}}$, with coefficient one
(Lemma~\ref{lem:screened-uniqueness}). The screening equations thus characterize
the polynomial before any values are assigned to the spectral
variables. Its evaluation in $D$ uses only the units $\miuraq_j$
and their translates under $\tau$.

For the ordinary transformation on the connected fixed subgroup $G_\nu$
in Corollary~\ref{cor:folded-miura}, the same uniqueness argument
gives $u_{k+1}=\tau u_k$. Given units $(\miuraq_O)_O$ indexed by the
$\nu$-orbits on the vertices of $\widetilde G$, the Steinberg coordinate
for the folded simple root at $O$ and index zero
is therefore the folded character in \eqref{eq:folded-character}
evaluated at $Y_{O',b_{O',k}}=\tau^k(\miuraq_{O'})$ for all $O'$ and $k$.

\subsection{Equivalence of the flag and Steinberg descriptions of an oper}
\label{sec:flag-steinberg}

\sectionoverview{In this subsection we prove that the flag and Steinberg
descriptions of a meromorphic $(G,\eta,\sigma)$-oper are equivalent for
a fixed torus trivialization and lift. We determine how the Steinberg form
changes with the trivialization. When the connection also preserves a
flag defining a Miura reduction, we express the condition for a holomorphic
Steinberg representative in terms of fundamental $q$-characters.}

\subsubsection{The fixed lift and the torus trivialization}

We compare the flag and Steinberg descriptions over a field $D$,
either the field of meromorphic functions on a
connected domain in $\C$ stable under $u\mapsto\eta^{\pm1}u$ or the
field of meromorphic germs at $0$. In both cases $(\tau f)(u)=f(\eta u)$.
We work in a meromorphic trivialization of
the principal $G$-bundle. If $A\in G(D)$ represents the connection
and $F_+=gB_+$ is the flag defining the $(G,\eta,\sigma)$-oper, its
image in the fiber at $\eta u$ is $A\sigma(g)B_+$. The relative-position
condition of Definition~\ref{def:eta-sigma-oper} therefore reads
\[
  (\tau g)^{-1} A\sigma(g)\in B_+(D)\dot c B_+(D).
\]

A meromorphic frame adapted to $F_+$ exists over $D$. Indeed, the
Bruhat cell containing the flag is parametrized by a product of root
subgroups acting on a Weyl group representative, which gives a lift
of the flag to $G(D)$. In adapted frames the isomorphisms of
$(G,\eta,\sigma)$-opers are precisely the changes of frame in
$B_+(D)$ acting by \eqref{eq:difference-gauge}. Fixing a trivialization of the associated
torus bundle, the quotient of the $B_+$-reduction by $N_+$, restricts
these changes of frame to $N_+(D)$. With
$U_c=N_+\cap\dot c N_-\dot c^{-1}$, the element $A$ representing
the connection in an adapted frame has a unique factorization
\[
 A=v\dot c h n,\qquad v\in U_c(D),\quad h\in H(D),\quad n\in N_+(D).
\]
Indeed, the decomposition of $N_+$ in the proof of
Proposition~\ref{prop:cell-parametrizations} gives the factorization.
The Gauss factors of $\dot c^{-1}A$ are
$\dot c^{-1}v\dot c$, $h$, and $n$, which proves uniqueness.
In this torus trivialization we call $L=\dot c h$ the lift of $c$
determined by $A$. Changes of frame in $N_+(D)$ preserve
the lift because they preserve the cell $N_+(D)LN_+(D)$. The
Steinberg description uses this lift to write the connection directly
in terms of meromorphic functions, one for each factor of $c$.

\begin{definition}
\label{def:steinberg-meromorphic-oper}
Let $c=s_{i_1}\cdots s_{i_\ell}$ have pinned lift
$\dot c=\dot s_{i_1}\cdots\dot s_{i_\ell}$. For a given torus
trivialization and lift $L=\dot c h$, $h\in H(D)$, a meromorphic
$(G,\eta,\sigma)$-oper in \emph{Steinberg form} is specified by
functions $T_{i_1},\ldots,T_{i_\ell}\in D$. In the
standard trivialization its flag $F_+$ is the constant coset $B_+\in G/B_+$,
and its connection sends $(u,g)$ to $(\eta u,S(u)\sigma(g))$, where
\[
 S(u)=
 \bigl(x_{i_1}(T_{i_1}(u))\dot s_{i_1}\bigr)\cdots
 \bigl(x_{i_\ell}(T_{i_\ell}(u))\dot s_{i_\ell}\bigr)\,h(u).
\]
If $h$ is holomorphic with values in $H$, this Steinberg form is
\emph{holomorphic} when each $T_{i_p}$ is holomorphic on the base domain.
\end{definition}

By Proposition~\ref{prop:cell-parametrizations}, the product
$S h^{-1}$ ranges over $U_c(D)\dot c$ as the functions $T_{i_p}$
vary. Thus $S$ ranges over
$U_c(D)L$, and its parameters are recovered from the principal minors
$T_{i_p}=\Delta_{i_p}(S h^{-1})$. In particular, for $h=1$ the map
$u\mapsto S(u)$ takes values in the fixed Steinberg section $\mathcal A_c$.

\begin{theorem}
\label{thm:flag-steinberg}
For fixed torus trivialization and lift $L$, every isomorphism class
of meromorphic $(G,\eta,\sigma)$-opers of Definition~\ref{def:eta-sigma-oper}
contains a unique Steinberg form $S$ of
Definition~\ref{def:steinberg-meromorphic-oper}, with standard flag $B_+$.
If $A$ represents the connection in a frame adapted to $F_+$ and inducing this
trivialization, there is a unique $n\in N_+(D)$ with
$S=(\tau n)A\sigma(n)^{-1}$.

When the torus trivialization varies, every meromorphic
$(G,\eta,\sigma)$-oper has a representative
$v\dot c h$, where $v\in U_c(D)$ and $h\in H(D)$. All isomorphisms
between these representatives are given by elements $t\in H(D)$ acting by
\begin{equation}\label{eq:flag-torus-action}
 (v,h)\longmapsto
 \bigl((\tau t)v(\tau t)^{-1},\ \dot c^{-1}(\tau t)\dot c\,h\,\sigma(t)^{-1}\bigr).
\end{equation}
\end{theorem}

\begin{proof}
In a frame adapted to $F_+$ and the chosen torus
trivialization, $A$ lies in $N_+(D)LN_+(D)$. Two such frames differ
by an element of $N_+(D)$. Proposition~\ref{prop:difference-section}
therefore gives a unique representative $A^n\in U_c(D)L$ of the
isomorphism class. Conversely, every element of $U_c(D)L$ lies in
$B_+(D)\dot c B_+(D)$, so the standard flag satisfies the required
relative-position condition.

For $S=v\dot c h$ and $t\in H(D)$,
\[
 S^t=\bigl((\tau t)v(\tau t)^{-1}\bigr)\dot c\,
       \bigl(\dot c^{-1}(\tau t)\dot c h\sigma(t)^{-1}\bigr).
\]
The first factor lies in $U_c(D)$ because $H$ normalizes $U_c$.
The last factor is the new torus component, giving
\eqref{eq:flag-torus-action}. Every isomorphism preserving the standard
flag has the form $b=nt$, with $n\in N_+(D)$ and $t\in H(D)$. The
formula shows that $S^t$ is already in the section for its new lift.
Since a change of frame by $n$ preserves that lift, if $S^b=(S^t)^n$
is in Steinberg form it belongs to the same section. The identity
already takes $S^t$ to that section, so uniqueness in
Proposition~\ref{prop:difference-section} gives $n=1$.
\end{proof}

The relative-position condition alone allows an arbitrary $L$.
Even over $\mathbb C(z)$ it may be impossible to replace that
lift by $\dot c$ through a change of frame preserving $F_+$.
For example, for the shift $(\tau f)(z)=f(qz)$, the connection of an
$(\mathrm{SL}_2,q)$-oper with the standard flag can be represented by
\[
 \begin{pmatrix}0&-z^{-1}\\z&0\end{pmatrix}.
\]
Indeed, this matrix sends the standard flag to a flag in relative position
$s_1$ with it.
A change of frame in $B_+(\mathbb C(z))$ with torus factor
$\operatorname{diag}(s,s^{-1})$ changes the lower-left entry to
$z/((\tau s)s)$. To make it one would require $(\tau s)s=z$.
The order of the left side at zero is even, whereas the right
side has order one. The choice of lift is therefore a genuine
condition for using the fixed section $\mathcal A_c$ over this field.

For an ordinary $(G,q)$-oper with $h=1$, the ordered product in
Definition~\ref{def:steinberg-meromorphic-oper} is the Steinberg map of \cite[Definition~9.1]{EP}, with the
same pinning, ordering and lifts. Composing a meromorphic function
on $\mathbb P^1$ valued in $H/W$ with this map gives the multiplicative
Hitchin section on the prescribed Hitchin base \cite[Definition~9.2]{EP}.
The same argument works for either a multiplicative or an
additive shift, with $\tau$ the pullback by the chosen shift.
In the version with prescribed singularities
\cite[Definition~9.6]{EP}, each simple lift has an additional torus
factor. Wherever these factors are defined over $D$, the product of
the modified lifts gives $L$. Conjugation by the preceding lifts sends
the $p$th simple root subgroup to the root subgroup for
$s_{i_1}\cdots s_{i_{p-1}}\alpha_{i_p}$, rescaling its parameter by a
unit. These roots exhaust $R_c$, so the products with arbitrary
parameters in $D$ give $U_c(D)L$. Elliott and Pestun restrict the
parameters to polynomials of prescribed degrees and leading
coefficients. They also observe that the conjugation giving this
presentation changes the framing at infinity
\cite[after Definition~9.6]{EP}. The equivalence above identifies
the meromorphic descriptions of the connection. To identify the global
locus of \cite[Remarks~2.19 and~5.7]{EP}, one would also have to
identify the image of the framed Hitchin section under twistor
rotation. In that moduli problem the connection must be an
isomorphism off the prescribed divisor, have the prescribed local
coweight at each point of the divisor, and take the specified value
in the framing at infinity \cite[Remark~2.4 and Definitions~2.6 and~2.9]{EP}.

\subsubsection{Miura reductions and holomorphic representatives}

For a Miura $(G,\eta,\sigma)$-oper with pinned lift $L=\dot c$, we
compute the Steinberg coordinates using the preserved flag $F_-$ of
Definition~\ref{def:eta-sigma-oper}. We assume that $F_-$ is opposite
to $F_+$ at the generic point and work in a frame in which $F_+=B_+$.
The flags in
$G/B_-$ opposite to the standard $B_+$ form the open orbit
$N_+ B_-/B_-\simeq N_+$. Indeed $B_+=N_+ H$ and $N_+\cap B_-=\{1\}$.
Consequently $F_-=pB_-$ for a unique $p\in N_+(D)$ in this trivialization.
Since $\sigma$ preserves $B_-$, the preservation equation gives
\[
 (\tau p) B_-=A\sigma(p)B_-
 \quad\Longleftrightarrow\quad A^{\mathrm M}:=(\tau p)^{-1}A\sigma(p)\in B_-(D).
\]
Both $\tau p$ and $\sigma(p)$ belong to $N_+(D)$, hence
$A^{\mathrm M}\in B_-\cap N_+\dot cN_+=M_c^-$ in the frame adapted to the two flags.
Conversely every $A^{\mathrm M}\in M_c^-(D)$ with $F_+=B_+$ and $F_-=B_-$ gives
such a Miura reduction. For an ordinary $(G,q)$-oper in Steinberg
form, the definition of triangularization at infinity in
\cite[Definition~9.4]{EP} is the same equation over
$\mathbb C[[z^{-1}]]$, with the shift denoted there by
$z\mapsto q^{-1}z$. In this case $p$ is a power-series point of $N_+$,
so $F_-$ is opposite to $F_+$ also at infinity.
Their upper unipotent element is $p^{-1}$ and their generalized
eigenvalues are the torus factor of $A^{\mathrm M}$. After forgetting
$F_-$, the Miura transformation
assigns to $A^{\mathrm M}$ the representative $(\tau u)A^{\mathrm M}\sigma(u)^{-1}\in\mathcal A(D)$
of Corollary~\ref{cor:difference-specialization}.

The fixed-lift equivalence also gives a criterion for the existence
of a holomorphic representative. Let $R$ be the ring of holomorphic
functions on the shift-invariant domain defining $D$, so that
$\tau(R)=R$.

\begin{corollary}
\label{cor:holomorphic-representatives}
Let $L=\dot c h$ with $h\in H(R)$ and let $A\in N_+(D)LN_+(D)$.
There is a change of frame $g\in N_+(D)$ for which
$(\tau g)A\sigma(g)^{-1}\in N_+(R)LN_+(R)$ if and only if the
Steinberg representative of $A$ lies in $U_c(R)L$. In particular,
for $L=\dot c$ and meromorphic Miura variables, such a change of
frame exists precisely when every evaluated fundamental
$q$-character in \eqref{eq:difference-character} is holomorphic.
\end{corollary}

\begin{proof}
If a meromorphic change of frame sends $A$ to $A'\in N_+(R)LN_+(R)$,
Proposition~\ref{prop:difference-section} over $R$ sends $A'$ to
$S\in U_c(R)L$. Composing the two changes of frame and using uniqueness
over $D$ identifies $S$ with the Steinberg representative of $A$.
Conversely, if that representative belongs to $U_c(R)L$, the
unipotent change of frame taking $A$ to it has the required property,
since $U_c(R)L\subset N_+(R)LN_+(R)$. For $L=\dot c$, the
$T_i$ are affine coordinates on $\mathcal A_c$ by
Proposition~\ref{prop:cell-parametrizations}. Thus a point of
$\mathcal A_c(D)$ lies in $\mathcal A_c(R)$ precisely when all its
$T_i$ belong to $R$. Corollary~\ref{cor:difference-specialization}
identifies these functions with the evaluated characters.
\end{proof}

For example, triangularizing a holomorphic connection in Steinberg
form may introduce poles in individual character monomials. Their sum
recovers the original coordinate by \eqref{eq:difference-character},
which explains the cancellation of apparent poles in this situation
\cite[Remark~9.10]{EP}. For arbitrary meromorphic Miura variables,
holomorphicity imposes a condition. In type $A_1$, the choice
$\miuraq_1=Q/(\tau Q)$ gives $T_1=(Q+\tau^2 Q)/(\tau Q)$. At a simple
zero of $\tau Q$ for a polynomial $Q$, the pole cancels precisely
when $Q+\tau^2 Q$ vanishes there.

\subsubsection{Meromorphic solutions and preserved flags}

Holomorphicity of the Steinberg coordinates leaves open whether the
connection admits a meromorphic change of frame to a constant. An
ordinary $(G,q)$-oper is $Z$-twisted if a meromorphic change of frame
makes its connection equal to the prescribed constant $Z\in H$. For a
Miura $(G,q)$-oper, $Z$-twistedness requires such a frame adapted
to the specified preserved flag \cite[Definitions~3.1--3.2]{FKSZ}.

Even a holomorphic Miura $(\SL_2,q)$-oper with connection in
$M_c^-=B_-\cap N_+\dot cN_+$ need not admit a meromorphic change of
frame to a constant torus element near $0$. For $q$ of infinite order,
the constant Miura variable $\miuraq_1=1$
gives $T_1=\miuraq_1+(\tau \miuraq_1)^{-1}=2$. Its Steinberg representative is
\[
 A=\begin{pmatrix}2&-1\\1&0\end{pmatrix}.
\]
Suppose $W(qu)A=ZW(u)$ held with $W$ invertible over the field
of meromorphic germs at $0$ and with $Z$ diagonal. If $v(u)=\sum_n v_nu^n$ is a row of $W$ and
$\lambda$ the corresponding diagonal entry of $Z$, then
$q^n v_n A=\lambda v_n$ for each $n$. The only left eigenline of $A$
is $\C(1,-1)$. Every nonzero $v_n$ lies on that line, so both rows
of $W$ annihilate $(1,1)^{\mathsf T}$, contradicting invertibility.
Nevertheless $W_k=Z^k W_0 A^{-k}$ satisfies
\eqref{eq:discrete-universal-wronskian} for every $W_0\in\SL_2(\C)$.
Here the Steinberg coordinate is holomorphic at $0$ although its
lattice Wronskian has no realization meromorphic at $0$.

A difference field extension allows us to ask which Miura coordinates
of a preserved flag can be recovered over $D$ when a solution of the
connection equation is only available over the extension. For
$A\in M_c(D)$ in the chosen trivialization, a solution $\Xi\in G(D')$
over a difference field extension $D'/D$ supplies a preserved flag
$F_-=\Xi B_-$ because
$\tau\Xi=A\sigma(\Xi)$. If the two flags are opposite over $D'$,
then $\Xi=pb$ with $p\in N_+(D')$, $b\in B_-(D')$, and $F_-=pB_-$.
The solution equation gives
\[
 A_{\Xi}^{\mathrm M}=(\tau p)^{-1}A\sigma(p)=(\tau b)\sigma(b)^{-1}.
\]
For $b_k=\tau^k b$ one has $b_{k+1}\sigma(b_k)^{-1}=\tau^k A_{\Xi}^{\mathrm M}$.
The Cartan identities \eqref{eq:miura-Q-ratio} and \eqref{eq:orbit-cartan-quotient}
at $k=0$ therefore give, for $i\in J$,
\begin{equation}\label{eq:meromorphic-orbit-ratio}
 \Delta_i(A_{\Xi}^{\mathrm M})=\frac{\tau Q_i}{Q_{\sigma^{-1}(i)}}
          =\frac{\tau Q_i}{\tau^{1-o_i}Q_i},\qquad
 Q_j=\omega_j([b]_0).
\end{equation}
For any $h\in G(D')$,
\[
 \tau(\Xi h)=A\sigma(\Xi h)\ \Longleftrightarrow\ \tau h=\sigma(h),
 \qquad \Xi hB_-=\Xi B_-\ \Longleftrightarrow\ h\in B_-(D').
\]
Thus solutions differing by $h$ with $\tau h=\sigma(h)$ preserve the
same flag precisely when $h\in B_-(D')$. Formula
\eqref{eq:meromorphic-orbit-ratio} computes the Miura variables determined by this flag. If its variables belong to
$D$, Corollary~\ref{cor:difference-specialization} constructs the
Steinberg representative over $D$ even when $\Xi$ is defined over $D'$.

\subsection{Cancellation of poles at the zeros of the \texorpdfstring{$Q$}{Q}-functions}
\label{sec:qq-pole-cancellation}

\sectionoverview{In this subsection we evaluate the Steinberg coordinates at
ratios of holomorphic $Q$-functions. We use the screening equations to
group the terms with a common denominator and the $QQ$-relations to
prove cancellation at a shifted zero of $Q_i^+$, including a zero of
higher multiplicity.}

The connected spectral domain $U\subset\C$ is invariant under
$u\mapsto\eta^{\pm1}u$. We write $\mathcal O(U)$ and $\mathcal M(U)$
for its ring of holomorphic functions and field of meromorphic functions. At vertex $i$ the screening
relation couples the indices $k$ and $k+o_i$, so the corresponding
arguments differ by $h_i=\eta^{o_i}$. In \eqref{tw:data}, $o_i$ is
the size of the $\sigma$-orbit. For the ordinary connections in the
simply laced and folded formulas $o_i=1$, including at a vertex of
$G_\nu$ coming from a $\nu$-orbit of size greater than one.

At each vertex the data are $Q_i^\pm\in\mathcal O(U)$, with
$Q_i^+\not\equiv0$, and a prescribed scalar $s_i\in\mathcal M(U)^\times$
that includes the twist factor. The substitution
\begin{equation}\label{eq:qq-miura-substitution}
 \miuraq_{i,k}(u)=s_i(\eta^k u)
       \frac{Q_i^+(\eta^k u)}{Q_i^+(h_i\eta^k u)}
\end{equation}
takes values in $\mathcal M(U)^\times$. Thus
Corollary~\ref{cor:difference-specialization} and its folded analogue
give meromorphic Steinberg coordinates $T_{j,k}(u)$. If $0\in U$ and
$Q_i^+$ has order $m_i\ge0$ at $0$, then the quotient of the
two shifted $Q_i^+$ in \eqref{eq:qq-miura-substitution} extends to $0$
with value $h_i^{-m_i}$. Indeed, writing $Q_i^+(u)=u^{m_i}p_i(u)$
with $p_i(0)\ne0$ cancels the common factor $u^{m_i}$. Thus, if the
prescribed $s_i$ are holomorphic and nonzero there, all Miura
variables are units at $0$ and every Steinberg coordinate is
holomorphic there by Laurent polynomial evaluation. To study a
possible pole at a nonzero $u$ with $Q_i^+(h_i\eta^t u)=0$, we use the
$i$th screening. In \eqref{ade:geometric-root},
\eqref{eq:folded-root-monomial} and \eqref{tw:screen-recurrence},
the root monomial is
$\mathsf R_{i,t}=\miuraq_{i,t}\miuraq_{i,t+o_i}/M_{i,t}$, where
$M_{i,t}$ is a product of Miura variables at neighboring vertices.
The corresponding generator of the screening kernel is
\begin{equation}\label{eq:qq-two-term-generator}
 E_{i,t}=\miuraq_{i,t}(1+\mathsf R_{i,t}^{-1})
        =\miuraq_{i,t}+\frac{M_{i,t}}{\miuraq_{i,t+o_i}} .
\end{equation}
Its two summands acquire the same denominator $Q_i^+(h_i\eta^t u)$
under \eqref{eq:qq-miura-substitution}.

For fixed $i$, the fundamental character computing $T_{j,k}$ can be
written as a finite sum of the form
\[
 \sum_{\mathbf m} c_{\mathbf m}
                \prod_{t\in \mathcal I}E_{i,t}^{m_t},
 \qquad \mathbf m\in\mathbb Z_{\ge0}^{\mathcal I},\qquad
 c_{\mathbf m}\in\C[\miuraq_{\ell,n}^{\pm1}:\ell\ne i,\ n\in\Z].
\]
Here $\mathcal I\subset\Z$ is finite. Indeed, the character belongs to the
$i$th screening kernel.
Theorem~\ref{thm:FM} gives such an expression for an ordinary
character, and Theorem~\ref{thm:fhr-folded-screening} gives one for a
folded character on $a\qaff^{\mathbb Z}$. By the three root monomial formulas, if the
initial variable $Y_{i,b}$ of a generator belongs to the lattice
$b_{\ell,n}=a\qaff^{2n+\xi_\ell}$, all its variables belong to that
lattice. The same holds after multiplying all spectral parameters
by a common scalar. Two such translates are equal or disjoint as
sets of pairs $(\ell,b)$. Setting all variables off the chosen lattice
equal to one therefore sends generators on the other translates to
constants and preserves $T_{j,k}$.

The evaluation of the neighboring product $M_{i,0}$ includes a scalar
$M_i[\mathbf s](u)$, obtained on replacing each $\miuraq_{\ell,n}$
by $s_\ell(\eta^n u)$. The remaining factors are shifted
$Q^+$-functions and give the quotient
\begin{equation}\label{eq:qq-neighbor-quotient}
 M_{i,0}(u)=M_i[\mathbf s](u)\frac{H_i(u)}{H_i(h_i u)} .
\end{equation}
In the simply laced and folded formulas,
\[
 H_i(u)=\prod_{\ell\prec i}Q_\ell^+(u)^{m_{\ell i}}
        \prod_{i\prec\ell}Q_\ell^+(\eta u)^{m_{\ell i}} .
\]
Here $m_{\ell i}$ is the adjacency multiplicity in the neighboring
product $M_{i,0}$. In the folded formula it counts the neighbors of
one vertex in the $\nu$-orbit indexed by $i$ that lie in the orbit
indexed by $\ell$. For the four sets of data in \eqref{tw:data},
the products following \eqref{tw:screen-recurrence} give the
following factor of $H_i$ for each neighboring vertex $\ell$.
The integer $n$ is its first index in $M_{i,0}$.
\[
\begin{array}{c|c|c}
 & \text{factors of }M_{i,0}\text{ at }\ell
 & \text{factor of }H_i(u)\\ \hline
 o_i=o_\ell & \miuraq_{\ell,n} & Q_\ell^+(\eta^n u)\\
 o_i>o_\ell=1 & \displaystyle\prod_{d=0}^{o_i-1}\miuraq_{\ell,n+d}
              & Q_\ell^+(\eta^n u)\\
 o_i=1<o_\ell & \miuraq_{\ell,n}
              & \displaystyle\prod_{d=0}^{o_\ell-1}Q_\ell^+(\eta^{n+d}u)
\end{array}
\]
In each row, substitution of \eqref{eq:qq-miura-substitution} in
the middle column and division by its scalar factors gives the last
column at $u$ divided by its value at $h_i u$.
For instance, in type $G_2$ one has
$H_1(u)=Q_2^+(\eta u)$ and
$H_2(u)=Q_1^+(u)Q_1^+(\eta u)Q_1^+(\eta^2u)$.

The prescribed singularity factor $\Lambda_i\in\mathcal M(U)^\times$
and the neighboring product $H_i$ determine $F_i=\Lambda_i H_i$.
For constants $\alpha_i,\beta_i\in\C^\times$, the $QQ$-relation
and the compatibility of the scalar factors are, respectively,
\begin{align}
 \alpha_i Q_i^-(u)Q_i^+(h_i u)
 -\beta_i Q_i^-(h_i u)Q_i^+(u)&=F_i(u),
 \label{eq:qq-local-equation}\\
 \frac{M_i[\mathbf s](u)}{s_i(u)s_i(h_i u)}
 &=\frac{\alpha_i\Lambda_i(u)}{\beta_i\Lambda_i(h_i u)} .
 \label{eq:qq-local-scalars}
\end{align}
In particular $F_i$ is holomorphic by \eqref{eq:qq-local-equation}.
Equation~\eqref{eq:qq-local-scalars} fixes the relative scalar in the
two summands of \eqref{eq:qq-two-term-generator}. In
\eqref{eq:ordinary-folded-qq} one has $\eta=q$,
$\alpha_i=\widetilde\gamma_i$, $\beta_i=\gamma_i$. When
$\Lambda_i=1$ at every vertex, $s_\ell=\zeta_\ell^{-1}$ gives
$\zeta_i^2\prod_{\ell\ne i}\zeta_\ell^{-m_{\ell i}}$ on both sides of
\eqref{eq:qq-local-scalars}.

For $u_0\in U$, the local ring $\mathcal O_{U,u_0}$ consists of
holomorphic germs at $u_0$. Its units are the germs with nonzero value
there.

\begin{proposition}
\label{prop:qq-pole-cancellation}
Suppose the holomorphic functions $Q_\ell^\pm$, with
$Q_\ell^+\not\equiv0$, and the prescribed meromorphic scalar factors
satisfy \eqref{eq:qq-local-equation} and \eqref{eq:qq-local-scalars}
at vertex $i$. At $u_0\in U$, assume the following conditions hold for
one expression of the fundamental character in the generators
$E_{i,t}$, $t\in\mathcal I$.
\begin{enumerate}
\item The functions $s_j(\eta^k u)$ and, for $t\in \mathcal I$,
$s_i(\eta^t u)$ and $s_i(h_i\eta^t u)$ are units. Each substituted
variable $\miuraq_{\ell,n}$ occurring in a coefficient $c_{\mathbf m}$
or in a product $M_{i,t}$ is also a unit.
\item For each $t\in \mathcal I$, at $x_t=h_i\eta^t u_0$ one has
$Q_i^+(x_t)=0\ \Longrightarrow\ F_i(x_t)\ne0$.
\end{enumerate}
Then $T_{j,k}(u)/s_j(\eta^k u)$ is holomorphic at $u_0$. The zeros
of $Q_i^+$ at the points $x_t$ may have any multiplicity.
\end{proposition}

\begin{proof}
For each $t\in\mathcal I$ we check holomorphicity of the generator
$E_{i,t}$ at $u_0$. With $w=\eta^t u$ and $x=h_iw$, substitution of
\eqref{eq:qq-miura-substitution} in \eqref{eq:qq-two-term-generator}
gives
\[
 E_{i,t}(u)=
 \frac{s_i(w)Q_i^+(w)+M_{i,t}(u)Q_i^+(h_i x)/s_i(x)}
      {Q_i^+(x)}.
\]
If $Q_i^+(x_t)\ne0$, this quotient lies in $\mathcal O_{U,u_0}$
by the first hypothesis. If $Q_i^+(x_t)=0$, equations
\eqref{eq:qq-neighbor-quotient} and \eqref{eq:qq-local-scalars} give
$M_{i,t}(u)/(s_i(w)s_i(x))=\alpha_i F_i(w)/(\beta_i F_i(x))$.
Substituting this ratio yields the meromorphic identity
\begin{equation}\label{eq:qq-generator-numerator}
 \frac{E_{i,t}(u)}{s_i(w)}
 =\frac{\beta_i F_i(x)Q_i^+(x/h_i)
              +\alpha_i F_i(x/h_i)Q_i^+(h_i x)}
             {\beta_i F_i(x)Q_i^+(x)} .
\end{equation}
Multiplying \eqref{eq:qq-local-equation} at $x$ and $x/h_i$ by
$\beta_i Q_i^+(x/h_i)$ and $\alpha_i Q_i^+(h_i x)$, respectively,
and adding gives
\begin{align}
 &\beta_i F_i(x)Q_i^+(x/h_i)+\alpha_i F_i(x/h_i)Q_i^+(h_i x)
 \notag\\
 &\qquad=Q_i^+(x)\bigl(
    \alpha_i^2 Q_i^-(x/h_i)Q_i^+(h_i x)
   -\beta_i^2 Q_i^-(h_i x)Q_i^+(x/h_i)\bigr).
 \label{eq:qq-exact-divisibility}
\end{align}
The factor in parentheses is holomorphic at $x_t$ and
$\beta_iF_i(x_t)\ne0$. Substitution in \eqref{eq:qq-generator-numerator} cancels the
entire factor $Q_i^+(x)$, regardless of its order of vanishing.
Since $s_i(w)$ is a unit, $E_{i,t}$ is holomorphic at $u_0$.
The coefficients $c_{\mathbf m}$ are holomorphic there by the first
hypothesis. Thus the polynomial expression for $T_{j,k}$ in the
generators $E_{i,t}$ is holomorphic at $u_0$, and division by the
unit $s_j(\eta^k u)$ preserves holomorphicity.
\end{proof}

For a variable $\miuraq_{\ell,n}$ the first hypothesis
therefore holds, in particular, when its scalar factor and both
shifted $Q_\ell^+$ in \eqref{eq:qq-miura-substitution} are
holomorphic and nonzero there. The second hypothesis concerns the
remaining denominator $F_i(x)$ in \eqref{eq:qq-generator-numerator}
after cancellation of $Q_i^+(x)$. If shifted zeros at adjacent
vertices coincide, the neighboring product $H_i(x)$ can vanish,
and regularity then requires control of the quotient at this
additional zero. For $Q$-functions supplied by a differential
equation, the two local hypotheses give cancellation independently
of any identification with Baxter eigenvalues. If the hypotheses hold
at every shifted $Q^+$-zero occurring in the character outside the zeros
and poles of the shifted scalar factors, the normalized coordinate is
holomorphic away from these scalar zeros and poles. At all other
points of that complement, holomorphicity follows directly from
\eqref{eq:qq-miura-substitution}.

At a simple zero $x_0$ of $Q_i^+$ with
$F_i(x_0)F_i(x_0/h_i)Q_i^+(x_0/h_i)\ne0$, vanishing of the
numerator of \eqref{eq:qq-generator-numerator} at $x_0$ is
equivalent to the Bethe equation
\[
 \frac{\alpha_i Q_i^+(h_i x_0)}{\beta_i Q_i^+(x_0/h_i)}
       =-\frac{F_i(x_0)}{F_i(x_0/h_i)} .
\]
This is the equation obtained by eliminating $Q_i^-$ in
\cite[(6.16) and Theorem~6.4]{FKSZ}, also used for the folded model
in \cite[Proposition~4.4]{FHR}. At a zero of multiplicity $m$,
\eqref{eq:qq-exact-divisibility} gives divisibility of the numerator
by $(x-x_0)^m$ in the ring of holomorphic germs. Thus cancellation
with multiplicity follows directly from the functional identity.

\subsection{Steinberg coordinates of twisted monodromy}
\label{sec:central-monodromy}

\sectionoverview{In this subsection we apply the Miura transformation to the
difference connection obtained by rotating the independent variable of a
differential equation. We express its Steinberg coordinates in terms of
solutions normalized at zero and at infinity.}

In the $q$DE/ODE/IM correspondence, the loop realization of an affine
oper gives an ordinary differential equation in $z$ for each fixed
value of a spectral parameter $u$. For example, the affine oper of
Feigin--Frenkel--Hernandez in untwisted simply laced type, with
singularities only at zero and infinity, has differential operator
\[
  \partial_z+\frac{f+\ell+(z+z^\kappa u)v_\theta}{z},
  \qquad 0<\kappa<1,
\]
where $f$ is the sum of the negative simple-root generators,
$\ell\in\operatorname{Lie}(H)$, and $v_\theta$ is a highest-root
vector. The Lie algebra elements act in the chosen representation of
$G$. On the logarithmic cover of the punctured $z$-plane, a positive
turn about zero multiplies $z^\kappa$ by $e^{2\pi i\kappa}$. Changing
$u$ to $qu$ with $q=e^{-2\pi i\kappa}$ therefore preserves the
product $z^\kappa u$ and hence the differential equation
\cite{QDE}.

For the general construction, the triple $(G,\sigma,c)$ and the ordered
set $J$ of simple-root representatives are as in Theorem~\ref{thm:intro}.
The differential operator has the form $\partial_z+\mathscr A(z;u)$,
with $\mathscr A(z;u)\in\operatorname{Lie}(G)$. In a faithful
representation $\rho$ of $G$, a fundamental solution at fixed $u$ is
an invertible matrix whose columns satisfy
$\partial_z y=-d\rho(\mathscr A(z;u))y$, where $d\rho$ is the
induced Lie algebra representation.
We use the two $G$-valued fundamental solutions of the analytic
construction in \cite{QDE,QDEII}.
The Frobenius frame $\mathcal Y_{\rm F}(z;u)$ is specified by the local
expansions of its columns at $z=0$. The sectorial frame
$\mathcal Y_{\rm as}(z;u)$ is specified by their asymptotic behavior as
$z$ tends to infinity within a chosen sector. We assume that both
continue to a common connected domain in the logarithmic cover and
depend meromorphically on $u$ in a connected spectral domain $U$.

At each $u$ where the frames are defined, $\mathcal S_u$ denotes the
set of $G$-valued fundamental solutions on this common domain. If
$Y_1,Y_2\in\mathcal S_u$, their differential equations give
$\partial_z(Y_1^{-1}Y_2)=0$. Thus $Y_2=Y_1g$ for a unique $g\in G$
independent of $z$, and each normalization gives coordinates on
$\mathcal S_u$. The element $C(u)$ defined by
\[
 \mathcal Y_{\rm as}(z;u)=\mathcal Y_{\rm F}(z;u)C(u)
\]
takes the sectorial coordinate $g$ of a solution to its Frobenius
coordinate $C(u)g$. In a representation of $G$, its matrix is the
central connection matrix.

For the twisted differential families, \cite[\S7.2]{QDE} describes how
fractional rotation identifies $C(u)$ with the universal
$(G,\eta,\sigma)$-Wronskian under the analytic hypotheses of \cite{QDEII}.
To express this rotation in the two frames, we use the symmetry
of the differential equation. With $r=\operatorname{ord}(\sigma)$ and
$\varepsilon=e^{2\pi i/r}$, the notation $\varepsilon z$ means that
$\arg z$ increases by $2\pi/r$ on the logarithmic cover. We assume
that $\eta U=U$ for some $\eta\in\C^\times$ and that the differential
operator satisfies
\[
 \varepsilon\mathscr A(\varepsilon z;\eta u)
       =\sigma\bigl(\mathscr A(z;u)\bigr),
\]
where $\sigma$ also denotes the induced automorphism of the Lie algebra.
For $Y\in\mathcal S_{\eta u}$, continuation to the rotated domain
and the chain rule show that $z\mapsto\sigma^{-1}(Y(\varepsilon z))$
solves the equation at $u$. This defines the pullback from
$\mathcal S_{\eta u}$ to $\mathcal S_u$.
Its inverse $F_u:\mathcal S_u\to\mathcal S_{\eta u}$ is the map of the
difference connection associated with one rotation. The normalized
Frobenius data are assumed to give $F_u$ the form $g\mapsto Z_0\sigma(g)$
for a fixed $Z_0\in H(\C)$. In sectorial coordinates we write the same
map as $g\mapsto\mathscr M(u)\sigma(g)$.
We assume $\mathscr M(u)\in B_+\dot cB_+$ for generic $u$. With this assumption,
$F_u$ and the flag $F_+$ represented by $B_+$ in sectorial coordinates
define a $(G,\eta,\sigma)$-oper. For
$\sigma=1$ this gives an ordinary $(G,q)$-oper with $q=\eta$.
A full turn in $z$ consists of $r$ such rotations. Its map on
solutions is $F_{\eta^{r-1}u}\circ\cdots\circ F_u$ from $\mathcal S_u$
to $\mathcal S_{qu}$, where $q=\eta^r$. In sectorial coordinates
we denote its action by $g\mapsto\mathscr M^{[r]}(u)g$.
As before, $(\tau f)(u)=f(\eta u)$.

\begin{proposition}
\label{prop:monodromy-central-minors}
The central connection matrix satisfies
$C(\eta u)\mathscr M(u)=Z_0\sigma(C(u))$ for one rotation.
For a full turn the monodromy in sectorial coordinates is
\[
 \mathscr M^{[r]}(u)=C(qu)^{-1}
     \bigl(Z_0\sigma(Z_0)\cdots\sigma^{r-1}(Z_0)\bigr)C(u).
\]
Where $C(u)$ and $C(\eta u)$ belong to $B_-N_+$, the factorization
$C=v_-n$ has $v_-\in B_-$ and $n\in N_+$. In the frame
$v_-$ relative to Frobenius coordinates, the flag $F_+$ of this
$(G,\eta,\sigma)$-oper is $B_+$, the preserved flag is $B_-$, and the
connection is given by
\[
 A_-=(\tau v_-)^{-1}Z_0\sigma(v_-)
     =(\tau n)\mathscr M\sigma(n)^{-1}\in B_-.
\]
For each simple root $i$, its torus component satisfies
\begin{equation}\label{eq:alternative-Q-ratio}
 \omega_i([A_-(u)]_0)=\omega_i(Z_0)\,
       \frac{\Delta_{\sigma^{-1}(i)}(C(u))}{\Delta_i(C(\eta u))}.
\end{equation}
\end{proposition}

\begin{proof}
A solution with sectorial coordinate $g$ at $u$ has Frobenius
coordinate $C(u)g$. The image under $F_u$ has sectorial coordinate
$\mathscr M(u)\sigma(g)$ at $\eta u$. Computing its Frobenius
coordinate in the two frames gives
\[
 C(\eta u)\mathscr M(u)\sigma(g)=Z_0\sigma(C(u)g).
\]
At $g=1$ this is $(\tau C)\mathscr M=Z_0\sigma(C)$. After $r$ iterations, with factors ordered from left to right by
increasing $j$, we obtain
\[
 (\tau^r C)\prod_{j=0}^{r-1}\sigma^j(\tau^{r-1-j}\mathscr M)
   =\left(\prod_{j=0}^{r-1}\sigma^j(Z_0)\right)\sigma^r(C).
\]
The product of the $r$ factors involving $\mathscr M$ represents the successive rotations
$F_{\eta^{r-1}u}\circ\cdots\circ F_u$ in sectorial coordinates and
hence equals $\mathscr M^{[r]}(u)$ when evaluated at $u$.
Since $(\tau^r C)(u)=C(qu)$ and $\sigma^r(C)=C$, left multiplication
by $C(qu)^{-1}$ gives the asserted expression for $\mathscr M^{[r]}(u)$.
Its inverse represents the full-turn pullback from $\mathcal S_{qu}$ to
$\mathcal S_u$ \cite{QDE}.

In Frobenius coordinates $F_+=CB_+$ and the connection
preserves $B_-$ because $Z_0\in H$ and $\sigma(B_-)=B_-$. Where
$C=v_-n$, the relations $CB_+=v_-B_+$ and $v_-B_-=B_-$ make both
flags standard in the frame $v_-$. Substituting this factorization into $(\tau C)\mathscr M=Z_0\sigma(C)$ gives
\[
 (\tau n)\mathscr M\sigma(n)^{-1}
    =(\tau v_-)^{-1}Z_0\sigma(v_-)=A_-\in B_- .
\]
The torus projection is multiplicative on $B_-$ and commutes with
$\tau$ and $\sigma$. Thus
\[
 \omega_i([A_-]_0)=\omega_i(Z_0)
 \frac{\omega_{\sigma^{-1}(i)}([v_-]_0)}
      {\omega_i([\tau v_-]_0)} .
\]
Proposition~\ref{prop:gauss-cell} gives
$\omega_j([v_-]_0)=\Delta_j(C)$, proving \eqref{eq:alternative-Q-ratio}.
\end{proof}

We now compute the Steinberg coordinates of $F_u$ for a single
rotation. For this purpose we assume that each principal minor
$\Delta_j(C)$ is generically nonzero, so that the two flags are opposite
for generic $u$.
The points where $C$ fails to be holomorphic in $G$, a principal minor
vanishes, or $\mathscr M$ fails to belong to $B_+\dot cB_+$ form a discrete
subset of $U$. The union of its translates by all powers of $\eta$
is countable. For a nonzero $u_*$ outside this union, we write
$u^{(k)}=\eta^ku_*$, $C_k=C(u^{(k)})=v_{-,k}n_k$, and
$\mathscr M_k=\mathscr M(u^{(k)})$.

The change of frame $n_k$ obtained from $C_k$ puts $F_+$ and the
preserved flag in standard position. The resulting connection maps lie in $B_-$.
A torus change of frame $s_k$ preserves both flags and puts the maps
in $M_c^-$. Once in Miura form, the sequence has a unique unipotent
change of frame $u_k$ taking it to the Steinberg section.
To determine $s_k$, we use the Bruhat factorization of
Section~\ref{sec:flag-steinberg}. At $u^{(k)}$ this factorization
places $\mathscr M_k$ in $N_+\dot c hN_+$ for a unique $h\in H$.
The identity
\[
 \dot c h=c(h)\dot c,\qquad c(h)=\dot c h\dot c^{-1},
\]
places the torus factor on the left, so that $\mathscr M_k\in
N_+h_k\dot cN_+$ with $h_k=c(h)$. The recurrence for $s_k$ then
involves the automorphism
$\gamma(s)=c(\sigma(s))$.

\begin{proposition}
\label{prop:monodromy-torus-trivialization}
Every $s_0\in H$ determines a unique sequence $(s_k)_{k\in\Z}$ in $H$
such that
\[
 A_k^{\rm M}=s_{k+1}A_-(u^{(k)})\sigma(s_k)^{-1}\in M_c^- .
\]
The sequence and its Miura variables are given by
\[
 s_{k+1}=\gamma(s_k)h_k^{-1},\qquad
 \miuraq_{i,k}=
 \frac{\omega_i(s_{k+1})\,\omega_i(Z_0)}
      {\omega_{\sigma^{-1}(i)}(s_k)}
 \frac{\Delta_{\sigma^{-1}(i)}(C_k)}{\Delta_i(C_{k+1})}
 \quad (i\in J).
\]
If the orbit of $u_*$ has finite length $L$, this change of frame
descends from $\Z$ to the orbit precisely when
\[
 s_0^{-1}\gamma^L(s_0)=\prod_{j=0}^{L-1}\gamma^{L-1-j}(h_j).
\]
Over the field of meromorphic functions on $U$, a torus change of
frame gives the pinned lift $\dot c$ precisely when it is represented
by an $H$-valued meromorphic function $s$ satisfying
$s(\eta u)=\gamma(s(u))h(u)^{-1}$, where
$\mathscr M(u)\in N_+h(u)\dot cN_+$ generically.
\end{proposition}

\begin{proof}
By Proposition~\ref{prop:monodromy-central-minors},
$A_-(u^{(k)})=n_{k+1}\mathscr M_k\sigma(n_k)^{-1}$ belongs to
$B_-\cap N_+h_k\dot cN_+$. Conjugation by $H$ preserves $N_+$,
and $\dot c\sigma(s_k)=\gamma(s_k)\dot c$. Hence
\[
 A_k^{\rm M}\in B_-\cap N_+
      \bigl(s_{k+1}h_k\gamma(s_k)^{-1}\bigr)\dot cN_+,
 \qquad
 A_k^{\rm M}\in M_c^-\ \Longleftrightarrow\
 s_{k+1}=\gamma(s_k)h_k^{-1}.
\]
The equivalence uses uniqueness of the torus factor in the Bruhat
cell. Since $\gamma$ is an automorphism of $H$, the recurrence
determines every $s_k$ from $s_0$ in both directions. In $B_-$ the
projection to $H$ gives
$[A_k^{\rm M}]_0=s_{k+1}[A_-(u^{(k)})]_0\sigma(s_k)^{-1}$.
Applying $\omega_i$ yields
\[
 \miuraq_{i,k}=\frac{\omega_i(s_{k+1})}
                    {\omega_{\sigma^{-1}(i)}(s_k)}
                \,\omega_i([A_-(u^{(k)})]_0),
\]
which gives the minor formula by \eqref{eq:alternative-Q-ratio}.

On an orbit of length $L$ we have $h_{k+L}=h_k$. Thus $s_{k+L}$
and $s_k$ satisfy the same recurrence, and they agree for every $k$
precisely when $s_L=s_0$. Iteration of the recurrence gives
\[
 s_L=\gamma^L(s_0)\prod_{j=0}^{L-1}\gamma^{L-1-j}(h_j)^{-1},
\]
which proves the periodicity criterion. For a meromorphic change of
frame $s(u)$, the same double-coset calculation gives torus factor
$(\tau s)h\gamma(s)^{-1}$. The transformed lift is $\dot c$ precisely when
$\tau s=\gamma(s)h^{-1}$.
\end{proof}

\begin{proposition}
\label{prop:monodromy-character}
With the torus sequence and Miura variables $\miuraq_{j,n}$ of
Proposition~\ref{prop:monodromy-torus-trivialization}, the
integer-indexed sequence representing twisted monodromy in Steinberg
form has coordinates
\[
 T_i(A_k)=
 \left.\chi_{\qaff}\!\left(L(Y_{i,b_{i,k}})\right)
 \right|_{Y_{j,b_{j,n}}=\miuraq_{j,n}},
 \qquad i\in J,\quad k\in\Z,
\]
where $a\in\C^\times$ and $b_{j,n}=a\qaff^{2n+\xi_j}$. On an
orbit of length $L$, if $s_L=s_0$, both the unipotent change of frame
and the Steinberg representative $(A_k)$ are $L$-periodic.
\end{proposition}

\begin{proof}
Proposition~\ref{prop:monodromy-torus-trivialization} gives a
representative $(A_k^{\rm M})$ of the monodromy in Miura form.
Its principal minor $\Delta_j(A_n^{\rm M})$ equals $\miuraq_{j,n}$.
The substitution in Theorem~\ref{thm:intro} is therefore exactly
$Y_{j,b_{j,n}}\mapsto\miuraq_{j,n}$, giving the asserted formula. If $s_L=s_0$, then
$(s_k)$ and $(A_k^{\rm M})$ are $L$-periodic. Translation by $L$
therefore gives another unipotent change of frame reducing the same
Miura sequence to Steinberg form. Uniqueness in
Theorem~\ref{thm:intro} gives $u_{k+L}=u_k$ and hence $A_{k+L}=A_k$.
\end{proof}

\begin{remark}
\label{rem:central-connection-discrete-wronskian}
When $\mathscr M_k\in\mathcal A_c$ for every $k$,
Proposition~\ref{prop:monodromy-central-minors} gives
$C_{k+1}\mathscr M_k=Z_0\sigma(C_k)$. Thus $(C_k)$ itself is the
discrete universal Wronskian for $(\mathscr M_k)$ with twist $Z_0$
and initial value $C_0$.

In general, the three changes of frame have product $g_k=u_ks_kn_k$,
which takes sectorial coordinates to Steinberg coordinates.
Since $A_k=g_{k+1}\mathscr M_k\sigma(g_k)^{-1}$,
\[
 (C_{k+1}g_{k+1}^{-1})A_k
 = C_{k+1}\mathscr M_k\sigma(g_k)^{-1}
 = Z_0\sigma(C_kg_k^{-1}).
\]
By Definition~\ref{def:wronskian-frame}, $W_k=C_kg_k^{-1}$ is
the discrete universal Wronskian for $(A_k)$ with twist $Z_0$ and
initial value $C_0g_0^{-1}$.
\end{remark}

The substitution of Proposition~\ref{prop:monodromy-character}
uses the minors of the central connection matrix directly. Expressing
them through the normalized $Q$-functions of \cite{QDEII} introduces
the scalar factors and the shifts of the spectral parameter prescribed
there at each vertex. Equations~\eqref{eq:ade-forward-fundamental}
and~\eqref{tw:forward-formula} express the same coordinates through
the lowest-weight minors of the group elements $A_k^{\rm M}$.

\subsection{Transfer matrices of XXZ-type models}
\label{sec:transfer-realization}

\sectionoverview{In this subsection we evaluate the $q$-character at
ratios of commuting Baxter operators. Restriction to a common
eigenline gives the scalar Baxter formula, which identifies the
Steinberg coordinates with the transfer eigenvalues for fundamental
auxiliary modules. We also give the corresponding formula for folded
characters.}

The quantum space $W$ of the model is a tensor product of
finite-dimensional simple type-one $U_{\qaff}(L\mathfrak a)$-modules.
We write $W_\lambda$ for its weight space of weight $\lambda$.
In Section~\ref{sec:central-monodromy} the differential equation had
independent variable $z$ and spectral parameter $u$. Here $z$ denotes
the spectral parameter of the transfer matrix, and $\mathbf u=(u_j)_{j\in K}$
denotes the formal Cartan twist of \cite[\S5]{FH}.
The Cartan matrix of $\mathfrak a$ is $C=(C_{ij})$. For independent formal parameters $\rho_i$, $i\in K$, we use an algebraic closure $\mathbb F$ of the fraction field of $\mathbb C[[\rho_i:i\in K]]$.
We fix twist parameters $u_j\in\mathbb F^\times$ satisfying
$\rho_i=\prod_j u_j^{C_{ji}}$. For a finite-dimensional auxiliary
module $V$, its character
$\chi_{\qaff}(V)$ is a Laurent polynomial determined by $V$. The
transfer matrix $t_V(z,\mathbf u)$ is an operator on $W$ obtained by
taking the twisted trace of the universal $R$-matrix over $V$. In this trace,
the diagonal block for each weight $\mu$ of $V$ contributes its trace
multiplied by the monomial $\prod_j u_j^{\langle\mu,\alpha_j^\vee\rangle}$
\cite[Definition~5.1]{FH}. We use the finite-dimensional auxiliary modules
and the positive prefundamental modules $R^+_{i,b}$ in the dual convention
of \cite[\S3.4]{FH}. The coefficients of their transfer matrices preserve
$W_\lambda$ and generate a commutative unital $\mathbb F$-algebra
$\mathcal T_\lambda\subset\operatorname{End}_{\mathbb F}
(W_\lambda\otimes\mathbb F)$ \cite[Definition~5.1 and Theorem~5.3]{FH}.

The positive prefundamental transfer matrix on $W_\lambda$ has the factorization
\[
 \boldsymbol{\mathcal Q}_i^+(x)
 :=t_{R^+_{i,1}}(x,\mathbf u)\big|_{W_\lambda}
 = f_i(x)\mathbf Q_i(x),
\]
where $\mathbf Q_i$ is a polynomial with coefficients in
$\mathcal T_\lambda$ \cite[Theorem~5.9]{FH}.
The scalar series $f_i$ is determined by the product
$\prod_{j,b}Y_{j,b}^{n_{j,b}}$ of the highest monomials of the tensor
factors of $W$. In terms of the quantum Cartan matrix $\mathsf C(t)$,
formula~(5.19) of \cite{FH} reads
\[
 f_i(x)=\exp\!\left(\sum_{j,b}n_{j,b}
          \sum_{r\ge1}\frac{(x/b)^r}{r}
                       \bigl(\mathsf C(\qaff^r)^{-1}\bigr)_{ij}\right).
\]
Here $\mathsf C_{ii}(t)=t^{d_i^{\mathrm{sym}}}+t^{-d_i^{\mathrm{sym}}}$
and $\mathsf C_{ij}(t)=(t^{C_{ij}}-t^{-C_{ij}})/(t-t^{-1})$ for
$i\ne j$. On $W_\lambda$ the one-dimensional Borel module of weight
$\varpi_i$ has transfer matrix $c_i(\lambda)\operatorname{id}$ with
$c_i(\lambda)=u_i\prod_j\qaff_j^{-(C^{-1})_{ji}\lambda_j}$, where
$\lambda_j=\langle\lambda,(\alpha_j^{\rm q})^\vee\rangle$ and the rational
powers use a fixed logarithm of $\qaff$ \cite[Example~5.2]{FH}.

The generalized Baxter relation of Frenkel and Hernandez specifies the
substitution we use. For every finite-dimensional type-one auxiliary
module $V$, replacing $Y_{j,b}$ in $\chi_{\qaff}(V)$ by
\[
 [\varpi_j]\,\frac{[R^+_{j,b\qaff_j^{-1}}]}{[R^+_{j,b\qaff_j}]}
\]
gives $[V]$ after clearing denominators in the Grothendieck ring
of the quantum Borel modules in the dual convention used here.
Here $[\varpi_j]$ denotes the class of the one-dimensional Borel
module of weight $\varpi_j$. The transfer matrix is additive on short
exact sequences and multiplicative on tensor products, and spectral
translation sends $[R^+_{j,b}]$ to $\boldsymbol{\mathcal Q}_j^+(bz)$ on $W_\lambda$
\cite[Theorem~4.8, (5.17), and Theorem~5.3]{FH}. For the formal twist,
$\boldsymbol{\mathcal Q}_j^+(0)$ has entries in $\mathbb C[[\rho_i:i\in K]]$ and
specializing all $\rho_i$ to zero gives the identity on $W_\lambda$
\cite[Proposition~5.5]{FH}. This last assertion permits division at
the level of transfer matrices.

\begin{corollary}\label{cor:transfer-operator-realization}
For the groups and orderings of Theorem~\ref{thm:intro}, the spectral lattice has parameters $b_{i,k}=a\qaff^{2k+\xi_i}$ with $a\in\C^\times$ fixed. Over
$D=\mathcal T_\lambda[[z]]$ with $(\tau F)(z)=F(\qaff^2z)$, the ratios
\begin{equation}\label{eq:FH-baxter-functions}
 \mathbf B_i(x)=c_i(\lambda)\,
    \frac{f_i(x\qaff_i^{-1})\mathbf Q_i(x\qaff_i^{-1})}
         {f_i(x\qaff_i)\mathbf Q_i(x\qaff_i)}
\end{equation}
are units and give, for every finite-dimensional type-one auxiliary module $V$,
\[
 t_V(z,\mathbf u)|_{W_\lambda}
    =\left.\chi_{\qaff}(V)\right|_{Y_{j,b}=\mathbf B_j(bz)}
    \qquad\text{in }D .
\]
For the Miura variables $\miuraq_i(z)=\mathbf B_i(b_{i,0}z)$, the Steinberg
representative $A$ of $A^{\mathrm M}=\prod_{p=1}^{\ell}P_{i_p}(\miuraq_{i_p})$,
with factors in the order defining $c$, has coordinates
\[
 T_i(\tau^k A)=t_{L(Y_{i,b_{i,k}})}(z,\mathbf u)|_{W_\lambda}\qquad(i\in J,\ k\in\mathbb Z).
\]
\end{corollary}

\begin{proof}
The constant term $M=\boldsymbol{\mathcal Q}_i^+(0)$ specializes
to the identity when all $\rho_j=0$. Its determinant therefore
specializes to $1$, and the characteristic polynomial
$p(t)=t p_1(t)+p(0)$ of $M$ has nonzero constant term in $\mathbb F$.
Cayley--Hamilton gives
\[
 M^{-1}=-p(0)^{-1}p_1(M)\in\mathcal T_\lambda,\qquad
 \boldsymbol{\mathcal Q}_i^+(x)\in\mathcal T_\lambda[[x]]^\times .
\]
Thus the two shifted series in \eqref{eq:FH-baxter-functions} are
units, and their ratio has constant term $\operatorname{id}$. Under
the transfer homomorphism the one-dimensional class $[\varpi_i]$
maps to $c_i(\lambda)\operatorname{id}$. Consequently, after applying
the homomorphism to the generalized Baxter relation with denominators
cleared, division by the shifted prefundamental series replaces
each $Y_{i,b}$ by $\mathbf B_i(bz)$. This proves the identity for $t_V$.

For $\miuraq_j(z)=\mathbf B_j(b_{j,0}z)$ one has
$(\tau^s\miuraq_j)(z)=\mathbf B_j(b_{j,s}z)$. Corollary~\ref{cor:difference-specialization}
therefore identifies $T_i(A)$ with the transfer matrix for
$L(Y_{i,b_{i,0}})$. Translation by $\tau^k$ replaces each $b_{j,s}$
by $b_{j,k+s}$ and yields the identity for $T_i(\tau^k A)$.
\end{proof}

To pass from the operator identity to its eigenvalues, we restrict to a
common eigenline $\mathbb Fv$ of the transfer-matrix coefficients in
$W_\lambda\otimes\mathbb F$. Every joint generalized eigenspace contains
such a line because $\mathbb F$ is algebraically closed and $W_\lambda$
is finite-dimensional. On this line we write
\[
 \mathbf Q_j(x)v=Q_{j,v}(x)v,\qquad
 t_V(z,\mathbf u)v=\Lambda_{V,v}(z)v.
\]
Here $Q_{j,v}$ is the Baxter polynomial of the chosen joint eigenvalue.
By Corollary~\ref{cor:transfer-operator-realization}, $\boldsymbol{\mathcal Q}_j^+(0)$ is
invertible, so its eigenvalue $f_j(0)Q_{j,v}(0)$ is nonzero. Hence
$\mathbf B_j(x)$ acts on this line by the scalar series
\begin{equation}\label{eq:scalar-baxter-ratios}
 B_{j,v}(x)=c_j(\lambda)
 \frac{f_j(x\qaff_j^{-1})Q_{j,v}(x\qaff_j^{-1})}
      {f_j(x\qaff_j)Q_{j,v}(x\qaff_j)} ,
 \qquad B_{j,v}(0)=c_j(\lambda)\ne0 .
\end{equation}

Restriction of the operator identity to $\mathbb Fv$ thus gives the
scalar Baxter formula of \cite[Theorem~5.11]{FH}.

\begin{theorem}
\label{thm:fh-scalar-baxter}
Let $W$ be a tensor product of finite-dimensional simple type-one
$U_{\qaff}(L\mathfrak a)$-modules with formal twist $\mathbf u$.
Suppose $0\ne v\in W_\lambda\otimes\mathbb F$ is a common eigenvector
of the transfer-matrix coefficients, with scalar Baxter ratios $B_{j,v}$
given by \eqref{eq:scalar-baxter-ratios}. For every finite-dimensional
type-one auxiliary $U_{\qaff}(L\mathfrak a)$-module $V$, the eigenvalue
of $t_V(z,\mathbf u)$ on $v$ is
\begin{equation}\label{eq:scalar-transfer-character}
 \Lambda_{V,v}(z)=
  \left.\chi_{\qaff}(V)\right|_{Y_{j,b}=B_{j,v}(bz)}
  \qquad\text{in }\mathbb F[[z]] .
\end{equation}
\end{theorem}

\begin{proposition}
\label{prop:scalar-transfer-steinberg}
For $b_{j,n}=a\qaff^{2n+\xi_j}$, let $(A_{k,v})$ be the Steinberg
representative of the Miura sequence with coordinates
$\miuraq_{j,n}=B_{j,v}(b_{j,n}z)$. Then
\[
 T_i(A_{k,v})=\Lambda_{L(Y_{i,b_{i,k}}),v}(z)
 \qquad\text{for every }i,k .
\]
\end{proposition}

\begin{proof}
For $V=L(Y_{i,b_{i,k}})$, Theorem~\ref{thm:intro} expresses
$T_i(A_{k,v})$ as the evaluation of $\chi_{\qaff}(V)$ at the Miura
variables. The prescribed values $B_{j,v}(b_{j,n}z)$ are units by
\eqref{eq:scalar-baxter-ratios}, and \eqref{eq:scalar-transfer-character}
identifies the same evaluation with $\Lambda_{V,v}(z)$.
\end{proof}

\begin{remark}\label{rem:monodromy-transfer-comparison}
For $\mathrm{SL}_2$ and $\mathrm{SL}_3$, the fundamental transfer
eigenvalues appear explicitly in the coefficients of the scalar
$q$-difference equations in \cite[\S\S3.3 and~5.2]{KSZ}.
To compare the transfer eigenvalues with twisted monodromy for general
$G$, we use Proposition~\ref{prop:monodromy-character}, which expresses
the Steinberg coordinate $T_i(A_k)$ of twisted monodromy in its Miura
variables $\miuraq_{j,n}$ for the chosen Coxeter lift. Suppose these variables and the
Baxter ratios of a joint eigenvector $v$ are units in a common
commutative $\mathbb C$-algebra after identification of the spectral
parameters. In the highest-weight convention,
$\miuraq_{j,n}=B_{j,v}(b_{j,n}z)$ gives
$T_i(A_k)=\Lambda_{L(Y_{i,b_{i,k}}),v}(z)$
by \eqref{eq:scalar-transfer-character}. In the lowest-weight
convention, \eqref{eq:ade-forward-fundamental} and
\eqref{tw:forward-formula} give
$T_i(A_k)=\Lambda_{L(Y_{i,a\qaff^{-2k-\xi_i}}),v}(z)$ when
$\miuraq_{j,n}^{-1}=B_{\iota(j),v}(a\qaff^{H_{\mathfrak a}-2n-\xi_j}z)$.
In either convention the two sides evaluate the same fundamental character,
so for fixed $i,k$ it suffices to match the finitely many variables
occurring in that character.
\end{remark}

For meromorphic $QQ$-data, Proposition~\ref{prop:qq-pole-cancellation}
proves local cancellation under its nonvanishing hypotheses. When the
Miura variables are the Baxter ratios of an XXZ model,
Proposition~\ref{prop:scalar-transfer-steinberg} identifies the coordinate
with a transfer eigenvalue. The denominators of the normalized $R$-matrices
then locate its possible poles after scalar normalization for all joint
eigenvectors at once. For $V=L(Y_{i,b_{i,k}})$, these denominators come from its $R$-matrices with the individual factors of $W=W_1\otimes\cdots\otimes W_m$. With spectral parameter $z$ on $V$, the
$R$-matrix for $(V,W_\ell)$ is $r_\ell(z)R_\ell^{\rm norm}(z)$,
where $r_\ell$ is a scalar series and $R_\ell^{\rm norm}$ is a rational
matrix. Its denominator $d_\ell(z)$ is the monic polynomial of least
degree for which $d_\ell R_\ell^{\rm norm}$ has Laurent polynomial
entries \cite[p.~2693]{KKOP}.

\begin{corollary}
\label{prop:transfer-common-denominator}
For $V=L(Y_{i,b_{i,k}})$ and $W=W_1\otimes\cdots\otimes W_m$, the products of the scalar factors and denominators are $S_{V,W}=\prod_{\ell=1}^m r_\ell$ and $d_{V,W}=\prod_{\ell=1}^m d_\ell$, respectively.
For every joint eigenvector $v$, the Steinberg coordinate obtained
from its Baxter ratios satisfies
\[
 \frac{d_{V,W}(z)}{S_{V,W}(z)}\,T_i(A_{k,v})
       \in\mathbb F[z,z^{-1}] .
\]
In particular, $S_{V,W}^{-1}T_i(A_{k,v})$ is regular at every
$z_0\ne0$ with $d_{V,W}(z_0)\ne0$, even if several shifted Baxter
polynomials vanish there.
\end{corollary}

\begin{proof}
Multiplication of the ordered product of $R$-matrices by
$d_{V,W}/S_{V,W}$ replaces each factor by $d_\ell R_\ell^{\rm norm}$,
whose entries are Laurent polynomials. The twist is independent
of $z$, so the twisted trace over $V$ still has Laurent polynomial
entries. In any $\mathbb F$-basis of $W_\lambda\otimes\mathbb F$
beginning with the joint eigenvector $v$, its first diagonal entry is
$d_{V,W}S_{V,W}^{-1}\Lambda_{V,v}$. By
Proposition~\ref{prop:scalar-transfer-steinberg}, this is
$d_{V,W}S_{V,W}^{-1}T_i(A_{k,v})$. Division by $d_{V,W}$ preserves
regularity at every nonzero point where $d_{V,W}$ does not vanish.
\end{proof}

For the folded characters of Theorem~\ref{thm:intro-folded}, the
transfer matrices come from the quantum loop algebra associated with
the simply laced Lie algebra $\operatorname{Lie}\widetilde G$. The
quantum space $W$ of \cite[\S5]{FHR} is a finite-dimensional simple
type-one module. The diagram automorphism $\nu$ fixes its highest
monomial and the twist $\mathbf u$.
We work at twists where the transfer matrices in the Baxter relation
are defined and the positive prefundamental eigenvalues used as
denominators are nonzero series.

The positive and negative prefundamental modules of \cite[\S5.2]{FHR}
give the transfer matrices on $W_\lambda$
\[
 \boldsymbol{\mathcal Q}_j^\pm(x)
 =t_{R^\pm_{j,1}}(x,\mathbf u)\big|_{W_\lambda}.
\]
On a joint generalized eigenspace $E$ of weight $\lambda$, we write
$p^\pm_{j,E}(x)$ for the scalar eigenvalue series of
$\boldsymbol{\mathcal Q}_j^\pm(x)$. The subspace $W(\mathbf u)$
of \cite[Definition~5.6 and Lemma~5.8]{FHR} is the sum of all such
$E$ with
\[
 \nu(\lambda)=\lambda,\qquad
 p^\pm_{\nu(j),E}(x)=p^\pm_{j,E}(x)\quad\text{for every }j.
\]
By \cite[Theorem~5.13]{FHR}, every finite-dimensional auxiliary
transfer matrix preserves $W(\mathbf u)$, and its generalized
eigenvalues there are evaluations of folded characters, up to the
overall scalar factor in that theorem. In the normalization of
\eqref{eq:scalar-baxter-ratios} the ratio at an orbit $O$ is
\[
 B_{O,E}(x)=c_j(\lambda)
       \frac{p^+_{j,E}(x\qaff^{-1})}{p^+_{j,E}(x\qaff)},
       \qquad j\in O .
\]
Indeed, the shift is $\qaff$ at every simply laced vertex, and the
invariance of $\lambda$ and $\mathbf u$ gives
$c_{\nu(j)}(\lambda)=c_j(\lambda)$. Thus the displayed expression is
independent of $j\in O$. To evaluate the folded character on the
spectral lattice of Theorem~\ref{thm:intro-folded}, we use these ratios
at $b_{O,n}z$, where $b_{O,n}=a\qaff^{2n+\xi_O}$.

\begin{proposition}
\label{prop:folded-transfer-steinberg}
For each joint generalized eigenspace $E$ occurring in $W(\mathbf u)$,
the Miura sequence in $G_\nu=(\widetilde G^\nu)^\circ$ with coordinates
$\miuraq_{O,n}=B_{O,E}(b_{O,n}z)$ has a Steinberg representative
$(A_{k,E})$. For every $k\in\mathbb Z$ and $i\in O$, the coordinate
$T_O(A_{k,E})$ equals the generalized eigenvalue on $E$ of the
transfer matrix with auxiliary module $L(\widetilde Y_{i,b_{O,k}})$
in the normalization of \eqref{eq:scalar-transfer-character}.
\end{proposition}

\begin{proof}
For $j\in O$, the nonzero series $p^+_{j,E}(x)$ has the form
$x^m v_j(x)$ with $v_j(0)\ne0$. Its two shifts have the same order at zero.
Cancelling the factor $x^m$ in their ratio gives
\[
 B_{O,E}(x)=c_j(\lambda)\qaff^{-2m}
     \frac{v_j(x\qaff^{-1})}{v_j(x\qaff)}
\]
with nonzero constant term. The Miura variables are therefore units,
and Theorem~\ref{thm:intro-folded} gives the unique Steinberg
representative and the identity
\[
 T_O(A_{k,E})=
 \left.\pi\!\left(\chi_{\qaff}(L(\widetilde Y_{i,b_{O,k}}))\right)
 \right|_{Y_{O',b}=B_{O',E}(bz)} .
\]

The coefficients of the auxiliary transfer matrix commute with those
of the prefundamental transfer matrices and therefore preserve their
joint generalized eigenspace $E$. In a basis of $E$ simultaneously
triangularizing these commuting operators, every diagonal entry of
$\boldsymbol{\mathcal Q}_j^+(x)$ is $p^+_{j,E}(x)$. Taking a diagonal
entry is multiplicative on upper triangular matrices, so the
generalized Baxter relation with denominators cleared holds for
each diagonal entry. Division in the field of Laurent series by
the nonzero shifted series $p^+_{j,E}$ gives the substitution
\[
 \widetilde Y_{j,b}\ \longmapsto\
 c_j(\lambda)\frac{p^+_{j,E}(bz\qaff^{-1})}{p^+_{j,E}(bz\qaff)}
   =B_{O',E}(bz)\qquad(j\in O').
\]
Since $\pi(\widetilde Y_{j,b})=Y_{O',b}$, this is precisely the
evaluation of the folded variable at $B_{O',E}(bz)$. For
$V=L(\widetilde Y_{i,b_{O,k}})$ every diagonal entry of the normalized
$t_V$ equals $T_O(A_{k,E})$.
\end{proof}

\appendix

\section*{Appendix. The Drinfeld current presentation of the quantum loop algebra}
\label{app:quantum-presentation}

\sectionoverview{In this appendix, we give the defining relations of the
quantum loop algebra in the normalization of \eqref{eq:cartan-currents}.
We relate the two customary choices of Cartan generators and specify
the coproduct used for tensor products of modules.}

The level-zero current presentation of $U_{\qaff}(L\mathfrak a)$
has generators $k_i^{\pm1}$, $h_{i,m}$ ($m\in\mathbb Z\setminus\{0\}$)
and $x_{i,r}^{\pm}$ ($r\in\mathbb Z$), indexed by the finite vertices
$i\in K$ \cite[Theorem~1.2]{FM}. The relations use the quantum integers
$[n]_w=(w^n-w^{-n})/(w-w^{-1})$ and the quantum binomial coefficients
\[
 [n]_w!=\prod_{j=1}^{n}[j]_w,\qquad
 \begin{bmatrix}n\\a\end{bmatrix}_w
 =\frac{[n]_w!}{[a]_w![n-a]_w!}\qquad(0\le a\le n).
\]
The choices $w=\qaff$ and $w=\qaff_i=\qaff^{d_i^{\mathrm{sym}}}$ have
infinite order, so these denominators are nonzero. The $k_i$ are
invertible, and all the Cartan generators $k_i^{\pm1}$ and $h_{j,m}$
commute with one another. With $B_{ij}=d_i^{\mathrm{sym}}C_{ij}$,
their relations with the root generators are
\[
 k_i x_{j,r}^{\pm}k_i^{-1}=\qaff^{\pm B_{ij}}x_{j,r}^{\pm},\qquad
 [h_{i,m},x_{j,r}^{\pm}]
   =\pm\frac{[mB_{ij}]_{\qaff}}{m}x_{j,r+m}^{\pm}.
\]
The conjugation factor is $\qaff^{B_{ij}}=\qaff_i^{C_{ij}}$, so the
ordinary weight of $x_{j,r}^\pm$ is $\pm\alpha_j^{\rm q}$. To express
the Cartan currents with coefficient $\qaff_i-\qaff_i^{-1}$ in the
exponent of \eqref{eq:cartan-currents}, one rescales the generators by
$\bar h_{i,m}=(\qaff-\qaff^{-1})h_{i,m}/(\qaff_i-\qaff_i^{-1})$.
The rescaled commutator is
\[
 [\bar h_{i,m},x_{j,r}^{\pm}]
 =\pm\frac{\qaff^{m B_{ij}}-\qaff^{-m B_{ij}}}
              {m(\qaff_i-\qaff_i^{-1})}x_{j,r+m}^{\pm}
 =\pm\frac{[m C_{ij}]_{\qaff_i}}{m}x_{j,r+m}^{\pm}.
\]

The symbol $\phi_{i,n}^\pm$ denotes the coefficient of $u^n$ in the corresponding Cartan current $\phi_i^\pm(u)$ of \eqref{eq:cartan-currents}.
The expansion at zero gives
$\phi_{i,n}^+=0$ for $n<0$, and the expansion at infinity gives
$\phi_{i,n}^-=0$ for $n>0$. The root generators satisfy, for all
$r,s\in\mathbb Z$,
\begin{align*}
 [x_{i,r}^+,x_{j,s}^-]
   &=\delta_{ij}\frac{\phi_{i,r+s}^+-\phi_{i,r+s}^-}{\qaff_i-\qaff_i^{-1}},\\
 x_{i,r+1}^{\pm}x_{j,s}^{\pm}
   -\qaff^{\pm B_{ij}}x_{j,s}^{\pm}x_{i,r+1}^{\pm}
   &=\qaff^{\pm B_{ij}}x_{i,r}^{\pm}x_{j,s+1}^{\pm}
     -x_{j,s+1}^{\pm}x_{i,r}^{\pm}.
\end{align*}
At $r=s=0$ the mixed relation gives the usual
finite-type commutator
$[x_{i,0}^+,x_{j,0}^-]=\delta_{ij}(k_i-k_i^{-1})/(\qaff_i-\qaff_i^{-1})$,
since $\phi_{i,0}^+=k_i$ and $\phi_{i,0}^-=k_i^{-1}$.

For $i\ne j$ and $p=1-C_{ij}$ the quantum Serre relations use
the group $\mathfrak S_p$ of permutations of $1,\ldots,p$ and read
\[
 \sum_{\pi\in\mathfrak S_p}\sum_{a=0}^{p}
 (-1)^a\begin{bmatrix}p\\a\end{bmatrix}_{\qaff_i}
 x_{i,r_{\pi(1)}}^{\pm}\cdots x_{i,r_{\pi(a)}}^{\pm}
 x_{j,s}^{\pm}
 x_{i,r_{\pi(a+1)}}^{\pm}\cdots x_{i,r_{\pi(p)}}^{\pm}=0
\]
for all $r_1,\ldots,r_p,s\in\mathbb Z$, with empty products equal to one.

The Drinfeld realization theorem identifies this algebra with the
level-zero Drinfeld--Jimbo algebra of the untwisted affine Cartan
matrix of $\mathfrak a$ \cite[\S1.2]{FM}. In the Chevalley presentation
the generators $E_j,F_j,K_j^{\pm1}$ are indexed by the finite nodes
and the affine node. The tensor product in $K_0(\mathcal C)$ uses the
coproduct
\[
 \Delta(K_j)=K_j\otimes K_j,\qquad
 \Delta(E_j)=E_j\otimes1+K_j\otimes E_j,\qquad
 \Delta(F_j)=F_j\otimes K_j^{-1}+1\otimes F_j .
\]

The current relations are homogeneous for the loop grading, so rescaling
by loop degree defines the spectral automorphism $\tau_z$ used in
\eqref{eq:spectral-translation} \cite[(1.4)]{FM}.

\end{document}